%% file: GGP2.tex
\documentclass[10pt,leqno]{article}

\usepackage[francais,english]{babel}
\usepackage{amsmath,amssymb,amscd,amsxtra}
\usepackage{a4wide}
\usepackage[latin1]{inputenc}
\usepackage{entete-hanoi}
\usepackage[all]{xy} 
\usepackage[T1]{fontenc}
\usepackage{hyperref}
\usepackage{xcolor}
\usepackage{mathdots}

\title{The global Gan-Gross-Prasad conjecture for unitary groups. II. \\
  From Eisenstein series to Bessel periods.} 
\author{Raphaël Beuzart-Plessis and Pierre-Henri Chaudouard}
\date{}

\begin{document}
\selectlanguage{english}
\maketitle

\begin{abstract}
  We state and prove an extension of the global Gan-Gross-Prasad conjecture and the Ichino-Ikeda conjecture to the  case of some Eisenstein series on unitary groups $U_n\times U_{n+1}$.  Our theorems are  based on a comparison of the Jacquet-Rallis trace formulas. A new point is the expression of some interesting spectral contributions in these formulas in terms  of integrals of relative characters. As an application of our mains theorems,  we prove the global Gan-Gross-Prasad and the  Ichino-Ikeda conjecture for Bessel periods of unitary groups.
\end{abstract}

\tableofcontents

\section{Introduction}

\subsection{Arthur parameters and weak base change}

\begin{paragr}
  In some sense this paper is a sequel of \cite{BCZ} where we proved the global Gan-Gross-Prasad, see \cite[section 24]{GGP}, and the Ichino-Ikeda conjectures for a product of unitary groups $U(n)\times U(n+1)$, see \cite{IIk} and \cite{NHar}. The goal of the present paper is two-fold: first we state and prove an extension of these two conjectures to the case of some Eisenstein series.   Second we show that this extension, when applied to some specific Eisenstein series, implies  the global Gan-Gross-Prasad conjecture and its refinement à la  Ichino-Ikeda  for general Bessel periods of unitary groups. To state our results we first review the notion of Arthur parameter.
\end{paragr}

\begin{paragr}[Hermitian Arthur parameter.] --- \label{S:Aparam}Let  $E/F$ be a quadratic extension of number fields and $c$ be the non-trivial element of the Galois group $\Gal(E/F)$. Let $\AAA$ be the ring of adèles of $F$. Let $n\geq 1$ be an integer.  Let $G_n$ be the group of automorphims of the $E$-vector space $E^n$. We view $G_n$ as an $F$-group by Weil restriction. For an automorphic representation $\Pi$ of $G_n(\AAA)$ we denote by $\Pi^*$ its conjugate-dual. Let's introduce some definitions.  A {\em discrete Hermitian Arthur parameter} of $G_n$ is an irreducible automorphic representation $\Pi$ of $G_n(\AAA)$ such  that

  \begin{itemize}
  \item  $\Pi$ is  isomorphic  to the full induced representation $\Ind_Q^{G_n}(\Pi_1\boxtimes \ldots \boxtimes \Pi_r)$ where $Q$ is a parabolic subgroup of $G_n$ with Levi factor $G_{n_1}\times\ldots \times G_{n_r}$ where $n_1+\ldots+n_r=n$;
  \item  $\Pi_i$ is a  conjugate self-dual cuspidal automorphic representation of $G_{n_i}(\AAA)$  and the  Asai $L$-function $L(s,\Pi_i,\As^{(-1)^{n+1}})$ has a pole at $s=1$ for $1\leq i\leq r$ .
  \item the representations $\Pi_i$  are mutually non-isomorphic for $1\leq i\leq r$;
   \end{itemize}
  The integer $r$ and the representations $(\Pi_i)_{1\leq i \leq r}$ are unique (up to a permutation).  We set $S_\Pi=(\mathbb{Z}/2\mathbb{Z})^{r}$.

  For our purpose, we need more general Arthur parameters of $G_n$, which we call {\em {regular} Hermitian Arthur parameters} and which are by definition the automorphic representations $\Pi$ of $G_n$ such that
    \begin{itemize}
  \item  $\Pi$ is  isomorphic  to the full induced representation $\Ind_Q^{G_n}(\Pi_1\boxtimes \ldots \boxtimes \Pi_r \boxtimes \Pi_0\boxtimes \Pi_r^* \boxtimes \ldots \boxtimes \Pi_1^*)$ where $Q$ is a parabolic subgroup of $G_n$ with Levi factor $M_Q=G_{n_1}\times\ldots \times G_{n_r}\times G_{n_0}\times G_{n_r}\times \ldots \times G_{n_1}$ where $n_0+2(n_1+\ldots+n_r)=n$;
  \item $\Pi_0$ is a discrete Hermitian Arthur parameter of $G_{n_0}$;
  \item  $\Pi_i$ is a  cuspidal automorphic representation of $G_{n_i}(\AAA)$ (with character central trivial on $A_{G_{n_i}}^\infty$) for $1\leq i\leq r$;
  \item the representations $\Pi_1,\ldots,\Pi_r,\Pi_1^*,\ldots,\Pi_r^*$ are  mutually non-isomorphic .
  \end{itemize}

  The representation $\Pi_0$ is uniquely determined by $\Pi$ and is called the discrete component of $\Pi$. We set:
  \begin{align*}
      S_\Pi=S_{\Pi_0}.
  \end{align*}

  The parabolic subgroup $Q$ depends on the ordering on the  representations  $\Pi_1,\ldots,\Pi_r,\Pi_1^*,\ldots,\Pi_r^*$: we fix one.

  Let $\ago_{Q,\CC}^*$ be the complex vector space of unramified characters of $Q(\AAA)$. We have the real subspaces $\ago_Q^*$ and $i\ago_Q^*$ respectively of real and unitary characters.  Let $w$ be the permutation that exchanges the two blocks $G_{n_i}$  corresponding to $\Pi_i$ and $\Pi_i^*$ for all $1\leq i\leq r$. We set
$$\ago_{\Pi,\CC}^*=\{\la\in \ago_{Q,\CC}^*\mid w\la=-\la\}.$$
For any $\la\in \ago_{\Pi,\CC}^*$, we define $\Pi_\la$ as the full induced representation
\begin{align*}
    \Ind_Q^{G_n}\left((\Pi_1\boxtimes \ldots \boxtimes \Pi_r \boxtimes \Pi_0\boxtimes \Pi_r^* \boxtimes \ldots \boxtimes \Pi_1^*)\otimes\la\right).
\end{align*}
 If $\la\in i\ago_\Pi^*=\ago_{\Pi,\CC}^*\cap i \ago_Q^*$  then $\Pi_\la$ is irreducible.
\end{paragr}

 \begin{paragr}[Unitary groups and (weak) base change.] --- \label{S:weakBC}For any integer $n\geq 1$, let $\hc_n$ be the set of isomorphism classes of non-degenerate $c$-Hermitian spaces $h$ over $E$ of rank $n$. We identify any    $h \in \hc_n$ with a representative and we denote by $U(h)$ its automorphism group. Let $h\in \hc$ and  $P\subset U(h)$ be  a parabolic subgroup  with Levi factor $M_P$. There exist a decomposition $n_0+2(n_1+\ldots+n_r)=n$ and $h_{n_0} \in \hc_{n_0}$ such that $M_P$ is identified with $G_{n_1}\times \ldots\times G_{n_r} \times U(h_{n_0})$.
 Let  $\sigma$ be a cuspidal automorphic subrepresentation of $M_P(\AAA)$ (with central character trivial on the central subgroup $A_P^\infty$ defined in §\ref{S:goodcpt}).
   Accordingly we have $\sigma=\Pi_1\boxtimes \ldots \boxtimes \Pi_r\boxtimes \sigma_0$ with $\Pi_i$ a  cuspidal automorphic representation of $G_{n_i}$ (with central character trivial on $A_{G_{n_i}}^\infty$).

   We shall say that a {regular} Hermitian Arthur parameter $\Pi$ of $G_n$ is a {\em weak base-change} of $(P,\sigma)$ if there exist  a parabolic subgroup $Q$ of $G_n$ with Levi factor $M_Q=G_{n_1}\times\ldots \times G_{n_r}\times G_{n_0}\times G_{n_r}\times \ldots \times G_{n_1}$ and a discrete Hermitian Arthur parameter  $\Pi_0$ of $G_{n_0}$ such that
   \begin{enumerate}
   \item  $\Pi$ is  isomorphic  to the full induced representation $\Ind_Q^{G_n}(\Pi_1\boxtimes \ldots \boxtimes \Pi_r \boxtimes \Pi_0\boxtimes \Pi_r^* \boxtimes \ldots \boxtimes \Pi_1^*)$;
   \item for almost all places of $F$ that split in $E$, the local component $\Pi_{0,v}$ is the split local base change of $\sigma_{0,v}$.
   \end{enumerate}
   Note that this implies that the representations $\Pi_1,\ldots,\Pi_r,\Pi_1^*,\ldots,\Pi_r^*$ are  mutually non-isomorphic and that $\Pi_0$ is in fact the discrete component of $\Pi$.
  If condition 2 above is satisfied, we shall also say that $\Pi_0$ is the weak base change of $\sigma_0$.

If $\Pi$ is a  weak base-change of $(P,\sigma)$,  we can naturally identify the space $\ago_{P,\CC}^*$ of unramified characters of $P(\AAA)$ with $\ago_{\Pi,\CC}^*$ and so  we will not distinguish between the two spaces. Thus for $\la\in \ago_{\Pi,\CC}^*$ we can consider the full induced representation $\Sigma_\la=\Ind_P^{U_h}(\sigma\otimes\la)$.

\end{paragr}

\begin{paragr} \label{S-intro:product} We can extend the notions above  to the case of a product. Let  $n,n'\geq 1$ be integers.  A \emph{regular} Hermitian Arthur parameter of $G_{n}\times G_{n'}$ is then an  automorphic representation of the form $\Pi=\Pi_n\boxtimes \Pi_{n'}$ where $\Pi_k$ is a \emph{regular} Hermitian Arthur parameter of $G_k$ for $k=n,n'$. Then we set $S_\Pi=S_{\Pi_n}\times S_{\Pi_{n'}}$ and $\ago_{\Pi,\CC}^*=\ago_{\Pi_n,\CC}^*\times \ago_{\Pi_{n'},\CC}^*  $ etc. For $\la=(\la_n,\la_{n'})\in \ago_{\Pi,\CC}^*$, we set $\Pi_\la=\Pi_{n,\la_n}\boxtimes \Pi_{n',\la_{n'}}$. A  parameter $\Pi$ is discrete if both $\Pi_n$ and $\Pi_{n'}$ are discrete.

Let $h\in\hc_n$ and $h'\in \hc_{n'}$. Let $P=P_n\times P_{n'}$ be a parabolic subgroup of $U(h)\times U(h')$.   We say that a {regular} Hermitian Arthur parameter $\Pi=\Pi_n\boxtimes \Pi_{n'}$ of $G_n\times G_{n'}$ is a weak base-change of $(P,\sigma)$ if $\Pi_n$ and $\Pi_{n'}$ are respectively  weak base-changes of $(P_n,\sigma_n)$ and  $(P_{n'},\sigma_{n'})$ where $\sigma=\sigma_n\boxtimes \sigma_{n'}$.
 \end{paragr}

\subsection{An extension of the Gan-Gross-Prasad conjecture to some  Eisenstein series}\label{ssec:ext-GGP}

\begin{paragr}[Corank $1$ and regular Hermitian Arthur parameter.] --- \label{S:AparamG} Let $n\geq 1$. Consider the ``corank $1$'' case $G=G_n\times G_{n+1}$.  Let  $\Pi=\Pi_n\boxtimes \Pi_{n+1}$ be a {regular} Hermitian Arthur parameter of $G$. We can write  $\Pi_k=\Ind_{Q_k}^{G_k}(\Pi_{1,k}\boxtimes \ldots \boxtimes \Pi_{r_k,k})$ for some parabolic subgroup $Q_k\subset G_k$ for $k=n,n+1$ and some cuspidal automorphic representations of $G_{n_{i,k}}(\AAA)$ with $n_{1,k}+\ldots+n_{r_k,k}=k$. We shall say that the parameter $\Pi$ is $H$-\emph{regular} if for  all $1\leq i\leq r_n$ and $1\leq j\leq r_{n+1}$  the representation $\Pi_{i,n}$ is not isomorphic to the contragredient of $\Pi_{j,n+1}$. 

  \begin{remarque}
    \label{rq:regularity}
     In the core of the paper, $H$ will stand for the diagonal subgroup $G_n$ of $G$ and the term $H$-regular refers to the fact that $H$-regular Hermitian Arthur parameter features some particularly nice properties with respect to the (regularized)  Rankin-Selberg period over that subgroup (roughly stemming from the fact that the Rankin-Selberg $L$-function of $\Pi$ has no poles).
     A discrete Hermitian Arthur parameter is necessarily {$H$-regular}. Otherwise we would get a self-conjugate cuspidal representation $\Pi_i$ of $G_{n_i}$ for some $n_i\geq 1$ such that both Asai $L$-functions $L(s,\Pi_i,\As^{+})$ and   $L(s,\Pi_i,\As^{-})$ have a pole at $s=1$: this is not possible.
  \end{remarque}
 
On the unitary side, let $h_0\in \hc_1$ be the element of rank $1$ given by the norm $N_{E/F}$. Then we attach to any $h\in \hc_n$ the following algebraic groups over $F$:
\begin{itemize}
  \item the product of unitary groups $U_h=U(h)\times U(h\oplus h_0)$ where $h\oplus h_0$ denotes the orthogonal sum;
\item the unitary group $U'_h$ of automorphisms of $h$ viewed as a subgroup of $U_h$ by the obvious diagonal embedding.
\end{itemize}

\end{paragr}

\begin{paragr} \label{S-intro:period} Let $P=M_PN_P\subset U_h$ be  a parabolic subgroup  with Levi factor $M_P$ and unipotent radical $N_P$. Let  $\sigma$ be a cuspidal automorphic subrepresentation of $M_P(\AAA)$ with central character trivial on $A_P^\infty$.  Let $\Ac_{P,\sigma}(U_h)$ be the space of automorphic forms on the quotient $A_P^\infty M_P(F) N_P(\AAA)\back U_h(\AAA)$ such that for all $g\in U_h(\AAA)$
  $$m\in M_P(\AAA)\mapsto \delta_P(m)^{-\frac12} \varphi(mg)$$
  belongs to the space of $\sigma$.  Here  $N_P$ is the unipotent radical of $P$  and $\delta_P$ is the modular character of $P(\AAA)$.  The  representation of $U_h(\AAA)$ on  $\Ac_{P,\sigma}(U_h)$ is isomorphic to  the  induced representation $\Sigma=\Ind_P^{U_h}(\sigma)$.  Let $\varphi\in \Ac_{P,\sigma}(U_h)$. For $\la\in \ago_{\Pi,\CC}^*$, we introduce the Eisenstein series $E(\varphi,\la)$ and the Ichino-Yamana regularized period
 \begin{align}\label{eq-intro:IY}
   \pc_{U'_h}(\varphi,\la)=\int_{[U_h']} \La_u^TE(x,\varphi,\la)\,dx
 \end{align}
where $[U_h']=U_h'(F)\back U_h'(\AAA)$ is equipped with the Tamagawa measure, $\La_u^T$ is the truncation operator introduced by Ichino-Yamana in \cite{IYunit} depending on an auxiliary parameter $T$ whose definition is recalled in §\ref{S:IY-truncU}. The integral is absolutely convergent. Moreover if the base change of $\Sigma$ is a {$H$-regular} Arthur parameter (which will be our assumption) then the integral does not depend on $T$ (see proposition \ref{prop:trunc-perU} below). In this case  $\pc_{U'_h}(\varphi,\la)$ is a meromorphic function, which is regular outside the singularities of  the Eisenstein series. In particular, it is holomorphic on $i\ago_{\Pi}^*$.
\end{paragr}

\begin{paragr}[The Gan-Gross-Prasad conjecture for some Eisenstein series.] --- \label{S-intro:GGP}
  
    \begin{theoreme}\label{thm:GGP}
      Let $\Pi$ be a {$H$-regular} Hermitian Arthur parameter of $G$ and let $\la\in i\ago_{\Pi}^*$. The following two statements are equivalent:
      \begin{enumerate}
      \item The complete Rankin-Selberg $L$-function of $\Pi_\la$ (including Archimedean places) satisfies
      $$L(\frac12,\Pi_\la)\not=0;$$
      \item There exist $h\in \hc_n$,  a parabolic subgroup  $P\subset U_h$ with Levi factor $M_P$ and $\sigma$ an irreducible  cuspidal automorphic subrepresentation of $M_P(\AAA)$ such that $\Pi$ is a weak base change of $(P,\sigma)$ and the period integral  $\varphi\mapsto \pc_{U'_h}(\varphi,\la)$ induces a non-zero linear form on $\Ac_{P,\sigma}(U_h)$.
      \end{enumerate}
    \end{theoreme}

    \begin{remarque}\label{rq:GGP}
     The Levi subgroup $M_P$ is determined up to conjugation by the parameter $\Pi$. Moreover we have $P=U_h$ if $\Pi$ is discrete. In this case, the theorem is proved in \cite[Theorem 1.8]{BPLZZ}  if $\Pi$ is cuspidal and in  \cite[Theorem 1.1.5.1]{BCZ} for a general discrete Hermitian parameter. The novelty of the theorem is to consider \emph{non-discrete}  Arthur parameters and thus periods of \emph{proper} Eisenstein series on unitary groups. 
          \end{remarque}
  \end{paragr}

\begin{paragr}[Factorization of periods  of  some Eisenstein series à la Ichino-Ikeda.] --- \label{S-intro:IIEis}Let $h\in \hc_n$. Let $P$ be   a parabolic subgroup of  $ U_h$ with Levi factor $M_P$ and let $\sigma$ be an irreducible  cuspidal automorphic subrepresentation of $M_P(\AAA)$ such that the weak base change of $(P,\sigma)$  is a regular Hermitian Arthur parameter $\Pi$.  We have a restricted tensor product decomposition  $\sigma=\bigotimes_{v\in V_F}'\sigma_v$ over the set $V_F$ of places of $F$.  We assume that $\sigma$ is tempered that is, for every place $v$, the local representation $\sigma_v$ is tempered. Let $\la\in i\ago_\Pi^*$. We define  $\Pi_\la$ and $\Sigma_\la$ as above. Let $\Sigma_{\la,v}=\Ind_P^{U_h}(\sigma_v\otimes\la)$ and $\Pi_{\la,v}$ be their local components.

We set
\begin{align*}
     \mathcal{L}(s,\Sigma_\la)= (s-\frac12)^{-\dim(\ago_{\Pi}^*)} \prod_{i=1}^{n+1}L(s+i-1/2,\eta^i)\frac{L(s,\Pi_\la)}{L(s+1/2,\Pi_\la,\As')}
\end{align*}
where $\eta$ denotes the quadratic idele class character associated to the extension $E/F$, $L(s,\eta^i)$ is the completed Hecke $L$-function associated to $\eta^i$ and $L(s,\Pi_\la,\As')$ is the $L$-function associated to $\As^{(-1)^n}\boxtimes\As^{(-1)^{n+1}}$. Note that with our hypothesis the function $L(s,\Pi_\la,\As')$ has a pole of order $\dim(\ago_{\Pi}^*)$ at $s=1$. Thus the function 
 $(s-1)^{-\dim(\ago_{\Pi}^*)} L(s,\Pi_\la,\As')$ is holomorphic and non-vanishing at $s=1$. In particular, the function  $\mathcal{L}(s,\Sigma_\la)$ is holomorphic at $s=\frac12$.
  
 We denote by $\mathcal{L}(s,\Sigma_{\la,v})$ the corresponding quotient of local $L$-factors, namely for $s$ in some half-space we have:
 \begin{align*}
    \mathcal{L}(s,\Sigma_\la)= (s-\frac12)^{-\dim(\ago_{\Pi}^*)} \prod_{v\in V_F} \mathcal{L}(s,\Sigma_{\la,v}).
 \end{align*}
For each place $v$ of $F$, we define a {\em local normalized period} $\pc^\natural_{h,\sigma_v}: \Sigma_v\times \Sigma_v\to \CC$ as follows:
$$\displaystyle \pc^\natural_{h,\Sigma_{\la,v}}(\varphi_v,\varphi'_v)=\mathcal{L}(\frac12,\Sigma_{\la,v})^{-1}\int_{U'_h(F_v)} (\Sigma_{\la,v}(h_v)\varphi_v,\varphi'_v)_v dh_v,\;\; \varphi_v,\varphi'_v\in \Sigma_v.$$
It depends on the choice of a Haar measure $dh_v$ on $U'_h(F_v)$ as well as an invariant inner product on $\sigma_v$ which gives in the usual way an invariant product on $\Sigma_v$ denoted by $(\cdot,\cdot)_v$.  By  the temperedness assumption, the integral is absolutely convergent \cite[Proposition 2.1]{NHar} and the local factor $\mathcal{L}(s,\Sigma_{\la,v})$ has neither zero nor pole at $s=\frac12$. 
  
We introduce on $\Ac_{P,\sigma}(U_h)$ the Petersson inner product  given by
  $$\displaystyle (\varphi,\varphi)_{\Pet}=\int_{A_P^{\infty}M_P(F)N_P(\AAA)\back U_h(\AAA)} \lvert \varphi(g)\rvert^2 dg,\;\; \varphi\in \sigma.$$
  Recall that we have normalized the period integral $\pc_{U_h'}(\la)$ by choosing the invariant Tamagawa measures on $[U_h']$. We also normalize the Petersson  product  by using the quotient of Tamagawa measures. We assume that the local Haar measures $dh_v$ on $U'_h(F_v)$ are such that the product $\prod_v dh_v$ gives the Tamagawa measure on $U_h'(\AAA)$.
  
  \begin{theoreme}\label{thm:II}
 Let $\Pi$ and  $(P,\sigma)$ as above.  For $\la\in i\ago_\Pi^*$ and every nonzero factorizable vector $\varphi=\otimes_v' \varphi_v\in \Ac_{P,\sigma}(U_h)\simeq \otimes'_{v\in V_F} \Sigma_v$, we have
  \begin{align}\label{intro-eq:II-E}
  \frac{\lvert \pc_{U_h'}(\varphi,\la)\rvert^2}{(\varphi, \varphi)_{\Pet}}=\lvert S_\Pi\rvert^{-1} \mathcal{L}(\frac12,\Sigma_\la)\prod_v \frac{\pc^\natural_{h,\Sigma_{\la,v}}(\varphi_v,\varphi_v)}{(\varphi_v,\varphi_v)_v}.
  \end{align}
  \end{theoreme}

  \begin{remarque}By \cite[Theorem 2.12]{NHar} and our choice of measures, almost all factors in the right-hand side are equal to $1$. As in remark \ref{rq:GGP}, the statement reduces to \cite[Theorem 1.1.6.1]{BCZ} for a discrete Hermitian Arthur parameter $\Pi$  and even to  \cite[Theorem 1.10]{BPLZZ}  if $\Pi$ is moreover simple.
  \end{remarque}
  
  \end{paragr}

  \subsection{The case of Bessel periods}\label{ssect: case of Bessel}

  \begin{paragr}
    Let $n\geq m\geq 0$ be two  integers of the same parity. We have $n=m+2r$ for some $r\geq 0$. Recall that we denote by $h_0$ the $1$-dimensional Hermitian space given by the norm $N_{E/F}$. Let $h_s\in \hc_2$ be the orthogonal sum of $h_0$ and $-h_0$.  For any $h\in \hc_m$,  we define $\tilde h\in \hc_n$ to be the orthogonal sum of $h$ and $r$ copies of $h_s$ denoted by $h_s^1,\ldots, h_s^r$ . For each $1\leq i\leq r$, let  $(x_i,y_i) $ be a hyperbolic basis of $h_s^i$ that is we have $h_s^i(x_i,x_i)=h_s^i(y_i,y_i)=0$ and $h_s^i(x_i,y_i)=1$. We consider also the orthogonal sum $h_{n+1}= \tilde h \oplus h_0 \in \hc_{n+1}$. We denote by $v_0$ the vector of $h_{0}$ corresponding to $1\in E$. We have a diagonal embedding
    \begin{align*}
    U(h)  \hookrightarrow \gc_h=U(h)\times U(h_{n+1})
    \end{align*}
    for which the image of $U(h)$  in $U(h_{n+1})$ is the subgroup which acts by the identity on $h_0\oplus h_s^1\oplus \ldots\oplus  h_s^r$.
    
    Let $B\subset U(h_{n+1})$ be the stabilizer of the isotropic flag
    \begin{align}\label{intro-eq:iso-flag}
      (0)\subsetneq \vect(x_1)\subsetneq \vect(x_1,x_2)\subsetneq \ldots \subsetneq \vect(x_1,\ldots ,x_r).
    \end{align}
    Let $N$ be the unipotent radical of $B$. Then the group $U(h)$ normalizes $N$. Let $\bc_h=U(h)\ltimes (\{1\}\times N)$: this is the so-called Bessel  subgroup of $\gc_h$.
  \end{paragr}

  \begin{paragr}[Bessel periods.] ---
    Let $\psi:\AAA/F\to \CC^\times$ a non-trivial continuous character. We define a character $\psi_N:[N]\to \CC^\times$ by
    \begin{align*}
       \psi_N(u)=\psi\left(\sum_{i=1}^{r-1} h_{n+1}(ux_{i+1},y_i)+h_{n+1}(uv_0,y_r) \right),\;\;\; u\in [N].
    \end{align*}
    This character extends uniquely to a character $\psi_{\mathcal{B}_h}: [\mathcal{B}_h]\to \CC^\times$ that coincides with $\psi_N$ on $[N]$ and is trivial on $[U(h)]$.
    Let $\sigma$ be a cuspidal automorphic subrepresentation of $\gc_h(\AAA)$.  We define the \emph{global Bessel period} for $\varphi\in \sigma$ by the absolute convergent integral
    \begin{align*}
      \pc_{\bc_h,\psi}(\varphi)=\int_{[\bc_h]} \varphi(g) \psi_{\bc_h}(g) dg.
    \end{align*}
  \end{paragr}

      \begin{paragr}[The Gan-Gross-Prasad conjecture for Bessel periods.] --- \label{S-intro:GGP-Bessel} Let $G^\flat=G_{m}\times G_{2n+1}$.  We can now state our first theorem about Bessel periods.

         \begin{theoreme} \label{thm:GGP-Bessel}
Let $\Pi$ be a discrete Hermitian Arthur parameter of $G^\flat$.
    The following assertions are equivalent:
    \begin{enumerate}
    \item The complete Rankin-Selberg $L$-function of $\Pi$  satisfies 
      \begin{align*}
        L(\frac12,\Pi)\not=0;
      \end{align*}
    \item There exist a Hermitian form $h\in \hc_m$ and  an automorphic  cuspidal subrepresentation $\sigma$ of $\gc_h(\AAA)$  such that  its  weak base to $G^\flat$ is $\Pi$ and the Bessel period 
\begin{align*}
        \varphi \mapsto \pc_{\bc_h,\psi}(\varphi)
      \end{align*}
      does not vanish identically on $\sigma$.
    \end{enumerate}
  \end{theoreme}

\begin{remarques}
    \begin{itemize}
    \item The case $r=0$ is just a particular case of Theorem \ref{thm:GGP-Bessel}.
    \item Assume $m=0$. Then the  $L$-function is the constant function of value $1$. So the assertion $1$ is automatically satisfied. On the other hand, the group $\gc_h$ is the quasi-split unitary group $U_{2r+1}$ of rank $2r+1$. Moreover the Bessel subgroup is a maximal unipotent subgroup of $U_{2r+1}$. Then the Bessel period is the so-called Fourier-Whittaker coefficient. The theorem is proved in the work of Ginzburg-Rallis-Soudry, see \cite{GRS}. 
\item The direction $2\Rightarrow 1$ is also proved by D. Jiang-L. Zhang, see \cite[Theorem 5.7]{JZ}.
    \end{itemize}
  \end{remarques}
\end{paragr}

\begin{paragr}   In our approach, Theorem \ref{thm:GGP-Bessel} is a consequence of Theorem \ref{thm:GGP}.  To explain this, we may and shall assume $r>0$. We start with  a discrete Hermitian Arthur parameter $\Pi$  of $G^\flat$. It can be written $\Pi=\Pi_m\boxtimes\Pi_{n+1}$ where  $\Pi_m$ and $\Pi_{n+1}$ are  respective discrete parameters  of $G_m$ and $G_{n+1}$. 
Let $\al_1,\ldots,\al_r$ be $r$ characters of $E^\times \back \AAA_E^1$ such that the characters $\al_1,\ldots,\al_r,\al_1^*,\ldots,\al_r^*$ are two by two distinct (we recall that $\al_i^*$ denotes the conjugate-dual of $\al_i$). Let $Q_n\subset G_n$ be a parabolic subgroup of Levi factor $G_1^{r}\times G_{m}\times G_1^r$. Then

  \begin{align*}
    \tilde\Pi=\Ind _{Q_n}^{G_n}(\al_1  \boxtimes \ldots \boxtimes\al_r   \boxtimes \Pi_m\boxtimes  \al_1^*  \boxtimes \ldots \boxtimes\al_r^* )  \boxtimes  \Pi_{n+1}
  \end{align*}
  is a {regular} Hermitian Arthur parameter of $G=G_n\times G_{n+1}$. Even if  $\tilde\Pi$ is not discrete, it is at least {$H$-regular} in the sense of § \ref{S:AparamG}: this is an obvious consequence of remark \ref{rq:regularity}  and the assumption on the characters $\al_i$. We have an identification $\CC^r \simeq \ago_{\tilde\Pi,\CC}^*$ such that if $\la_s\in  \ago_{\tilde\Pi,\CC}^*$ is the image of $(s,\ldots,s)$ with $s\in \CC$ we have
   \begin{align*}
    \tilde\Pi_{\la_s}=\Ind _{Q_n}^{G_n}(\al_1 |\cdot|_E^s  \boxtimes \ldots \boxtimes\al_r |\cdot|_E^s   \boxtimes \Pi_m \boxtimes \al_1^* |\cdot|_E^{-s}  \boxtimes \ldots \boxtimes\al_r^* |\cdot|_E^{-s})  \boxtimes  \Pi_{n+1}.
  \end{align*}
For simplicity we set $\tilde\Pi_{s}=\tilde\Pi_{\la_s}$.   By elementary properties of Rankin-Selberg $L$-function, it is clear that assertion 1 of  \ref{thm:GGP-Bessel} is equivalent to $1'$: 
  \begin{enumerate}
      \item[$1'$.] There exists $s\in i\RR$ such that  $L(\frac12,\tilde\Pi_{s})\not=0$.
  \end{enumerate}

  Let $h\in \hc_m$ and $\sigma$ be an automorphic cuspidal subrepresentation of $\gc_h$ whose weak base to $G^\flat$ is $\Pi$. Let $P_n\subset U(\tilde h)$ be the parabolic subgroup stabilizing the isotropic flag
  $$\displaystyle 0\subsetneq \vect(x_r)\subsetneq \vect(x_r,x_{r-1})\subsetneq \ldots \subsetneq \vect(x_r,\ldots,x_1).$$
  (Note that this flag is opposite position to \eqref{intro-eq:iso-flag}.) and set $P=P_n\times U(h_{n+1})$; a parabolic subgroup of $U_{\tilde h}=U(\tilde h)\times U(h_{n+1})$.  Then $G_1^r \times \gc_h$ is a Levi factor $M_P$ of $P$. Set $\tilde \sigma= \al_1\boxtimes \ldots \boxtimes\al_r \boxtimes \sigma$. This is an automorphic  cuspidal representation of $M_P(\AAA)$ and $\tilde\Pi$ is the weak base change of $(P,\sigma)$. Let $\varphi\in \Ac_{P,\tilde \sigma}(U_{\tilde h})$. As in subsection \ref{ssec:ext-GGP}, we denote by $U'_{\tilde h}$ the ``diagonal'' subgroup of $U_{\tilde h}$. In the case at hand, the restriction of the Eisenstein series $E(\varphi,\la)$ to $[U'_{\tilde h}]$ is rapidly decreasing for any $\la \in \ago_{\tilde \Pi,\CC}$ where the Eisenstein series is regular and for any such $\la$, we have
  \begin{align*}
    \pc_{U'_{\tilde h}}(\varphi,\la)=\int_{[U_h']} E(x,\varphi,\la)\,dx
  \end{align*}
where the left-hand side is defined according to \eqref{eq-intro:IY} and the right-hand side is absolutely convergent (see Proposition \ref{prop:trunc-perU} assertion 3). Moreover, the map $s\mapsto \pc_{U'_{\tilde h}}(\varphi,\la_s)$ is meromorphic and holomorphic on $i\RR$. We prove in Proposition \ref{prop:reduc-GGP-cork1} that the map  $\varphi\in  \Ac_{\sigma}(\gc_h)\mapsto \pc_{\bc_h,\psi}(\varphi)$ does not vanish identically if and only if there is $s\in \CC$ such that the map $\varphi \mapsto \pc_{U'_{\tilde h}}(\varphi,\la_s)$ does not vanish identically on $\Ac_{P,\tilde \sigma}(U_{\tilde h})$. This last fact is eventually a consequence of some unfolding identity that roughly takes the following form
$$\displaystyle \pc_{U'_{\tilde h}}(\varphi,\la_s)=\int_{\bc'(\AAA)\backslash U'_{\tilde h}(\AAA)} \pc_{\bc_h,\psi}(\varphi_s(h))dh$$
for $\varphi\in \Ac_{P,\tilde \sigma}(U_{\tilde h})$ where $\varphi_s$ stands for the corresponding element of $\Ind_P^{U_{\tilde h}}(\tilde \sigma\otimes \lambda_s)$ (given through the choice of a suitable Iwasawa decomposition $U_{\tilde h}(\AAA)=P(\AAA)K$ that is implicit in the definition of the Eisenstein series $E(\varphi,\lambda_s)$) and $\bc'=U(h)\rtimes V$ with $V$ the unipotent radical of the parabolic subgroup of $U(\tilde h)$ stabilizing the isotropic subspace $\vect(x_1,\ldots,x_r)$. It should be emphasized however that this identity doesn't make sense per se as the Eulerian integral on the right hand side is not absolutely convergent in general. More precisely, it has to be ``interpreted in the sense of $L$-functions'' which requires some nontrivial unramified computations of local integrals involving Bessel functions. We refer the reader to Section \ref{sect: application} and more specifically \ref{ssect: unr computation} and \ref{ssect: reduction} for details.

It follows that condition 2 of Theorem  \ref{thm:GGP-Bessel} holds for $h\in \hc_m$ and $\sigma$ if and only if the following assertion holds:
\begin{itemize}
\item[2'.]  There exists $s\in i\RR$ such that  $\varphi \mapsto \pc_{U'_{\tilde h}}(\varphi,\la_s)$ does not vanish identically on $\Ac_{P,\sigma}(U_{\tilde h})$. 
\end{itemize}
It is then straightforward to deduce Theorem \ref{thm:GGP-Bessel} from Theorem  \ref{thm:GGP}.
  
\end{paragr}

\begin{paragr}[Local Bessel periods.] --- From now on, we fix $h\in \hc_m$ and  a decomposition of the character $\psi=\otimes_{v\in V_F} \psi_v$ from which we get a decomposition  $\psi_{\mathcal{B}_h}= \otimes_{v\in V_F}\psi_{\mathcal{B}_h,v}$ where $\psi_{\mathcal{B}_h,v}$ is a character of $\mathcal{B}_h(F_v)$.
Let $v$ be a place of $F$.   The integral 
  \begin{align*}
    \int_{\bc_h(F_v)} f_v(g_v) \psi_{\bc_h,v}(g_v) dg_v
  \end{align*}
  is well-defined for a smooth and compactly supported function $f_v$ on $\gc_h(F_v)$ and extends to a continuous linear form $ f_v\mapsto \pc_{\bc_h,\psi_v}(f_v)$ on the space of tempered functions, see subsection \ref{Sect local periods}. It depends on the choice of a Haar measure on $\bc_h(F_v)$.

  Let  $\sigma_v$ be a  tempered irreducible  representation of $\gc_h(F_v)$ equipped with an invariant inner product $(\cdot,\cdot)_v$. Let $\varphi_v$ and $\varphi_v'$ be vectors of $\sigma_v$. The associated matrix coefficient defined by $f_{ \varphi_v,\varphi_v' }(g)=   (\sigma_v(g) \varphi_v,\varphi'_v)_v$ for all $g  \in \gc_h(F_v) $ belongs to this space and we set:
  \begin{align*}
     \pc_{\bc_h,\psi_v}( \varphi_v,\varphi'_v)= \pc_{\bc_h,\psi_v}( f_{ \varphi_v,\varphi_v' }).
  \end{align*}
\end{paragr}

\begin{paragr}[The Ichino-Ikeda conjecture for Bessel periods.] ---  Let $\sigma $ be a tempered automorphic  cuspidal subrepresentation $\sigma$ of $\gc_h(\AAA)$. Tempered means  that  we have a decomposition $\sigma=\otimes_{v\in V_F}' \sigma_v$ with $\sigma_v$ tempered for all $v$. We also assume that the weak base change of $\sigma$ to $G^\flat$ is a discrete Hermitian parameter $\Pi$. As in §\ref{S-intro:IIEis}, we define the ratio of $L$-functions  $\mathcal{L}(s,\sigma)$ and its local counterparts  $\mathcal{L}(s,\sigma_v)$ for $s\in \CC$.  Explicitly we have:

\begin{align*}
  \lc(s,\sigma_v)=\prod_{i=1}^{n+1}L(s+i-1/2,\eta^i_v)\frac{L(s,\Pi_v)}{L(s+1/2,\Pi_v,\As')}
\end{align*}
where $\As'=\As^{(-1)^m}\otimes\As^{(-1)^{n+1}}$ and  $\mathcal{L}(s,\sigma)$ is the product of the local factors in some half-plane. We use the local factor to define the normalized local Bessel period
  \begin{align*}
     \pc_{\bc_h,\psi_v}^\natural( \varphi_v,\varphi'_v)= \lc(\frac12 ,\sigma_v)^{-1}\pc_{\bc_h,\psi,v}( f_{ \varphi_v,\varphi_v' }).
  \end{align*}

We assume that the product of local measures on  $\bc_h(F_v)$ gives the Tamagawa measure on $\bc_h(\AAA)$. On $\sigma$, we use the Petersson inner product $(\cdot, \cdot)_{\Pet}$ normalized by the Tamagawa measure on $\gc_h(\AAA)$.

\begin{theoreme}\label{thm:II-Bessel}
 Let $\sigma$ and $\Pi$ as above.  For every nonzero factorizable vector $\varphi=\otimes_v' \varphi_v\in \sigma$, we have
  \begin{align}\label{intro-eq:II-Bessel}
  \frac{\lvert \pc_{\bc_h,\psi}(\varphi)\rvert^2}{(\varphi, \varphi)_{\Pet}}=\lvert S_\Pi\rvert^{-1} \mathcal{L}(\frac12,\sigma)\prod_v \frac{\pc_{\bc_h,\psi_v}^\natural(\varphi_v,\varphi_v)}{(\varphi_v,\varphi_v)_v}.
  \end{align}
  \end{theoreme}

  \begin{remarques}
    \begin{enumerate}
    \item In the right-hand side, almost all factors are equal to $1$, see \cite[Theorem 2.2]{Liu}.
    \item The statement has been conjectured by Y. Liu in a more general context,  see \cite[conjecture 2.5]{Liu}.
    \item  For $m=0$, the group $\gc_h$ is the quasi-split unitary group $U_{2r+1}$ of rank $2r+1$. The theorem has been conjectured by Lapid and Mao, \cite[conjecture 1.1]{LaMao}.
    \item The proof we give is along the same lines as for Theorem \ref{thm:GGP-Bessel}, namely it is eventually deduce it from Theorem \ref{thm:II} in a similar fashion.
    \end{enumerate}
  \end{remarques}

\end{paragr}

\subsection{On some spectral contributions of the Jacquet-Rallis trace formulas}

\begin{paragr} In this subsection, we explain some new  ingredients that play a role in the proof of  Theorems \ref{thm:GGP} and \ref{thm:II}. As many other contributions on the subject (among them, see \cite{Z1},  \cite{Zhang2}, \cite{Xue}, \cite{RBP},  \cite{BPPlanch}, \cite{BPLZZ}, \cite{BCZ}), we follow the strategy of the seminal paper \cite{JR} of Jacquet and Rallis. More precisely, besides the local harmonic analysis performed in the mentioned papers, our work is based on the  geometric comparison, fully established in \cite{CZ}, of the {\em relative trace formulas} constructed in  \cite{Z3} of the unitary groups $U_h$ for $h\in \hc_n$ and the corresponding group $G$. However to be able to exploit this comparison, we need to obtain more  tractable expressions for the spectral contributions we are interested in.
\end{paragr}  

\begin{paragr}
  Let's first explain our result in the unitary case namely for the group $U=U_h$ and its subgroup $U'=U_h'$. Let $\Xgo(U)$ be the set of cuspidal data of $U$. According to the work of Zydor, see  \cite[section 4]{Z3}, the contribution of $\chi\in \Xgo(U)$ to the relative trace formula for the group $U$ is built upon the absolutely convergent integral
  \begin{align*}
    \int_{[U']\times [U']} K^T_{f,\chi}(x,y) \, dxdy.
  \end{align*}
  Here $K^T_{f,\chi}$ is a suitably modified version à la Arthur of the $\chi$-part $K_{f,\chi}$ of the automorphic kernel $K_f(x,y)=\sum_{\gamma\in U(F)}f(x^{-1}\gamma y)$ associated to a Schwartz function $f$ on $U(\AAA)$,  see \eqref{eq:KchiTU} below for the precise definition. It depends on a truncation paramerer $T$. It turns out that the integral above is an exponential-polynomial function in $T$ whose purely polynomial part is constant and gives by definition the $\chi$-contribution denoted by $J^U_\chi(f)$ of the relative trace formula, see Theorem \ref{thm:jfDefU} for this slight extension of Zydor's work to Schwartz test functions. The problem however is to get an expression for $J_\chi(f)$ that reflects the Langlands spectral decomposition of $K_{f,\chi}$ and that is related to the periods \eqref{eq-intro:IY} we are interested in.  The starting point is the following new independent characterization of  $J^U_\chi(f)$:  the integral
  \begin{align}\label{intro-eq:KLaTu}
    \int_{[U']\times [U']}( K_{f,\chi} \Lambda^T_u)(x,y) \, dxdy
  \end{align}
  is absolutely convergent and is asymptotic to an  polynomial exponential in the variable  $T$ whose purely polynomial term is constant and equal to $J_\chi^U(f)$, see Corollary \ref{cor:LaTu-cstterm}. Here $K_{f,\chi} \Lambda^T_u$ means that we have applied the Ichino-Yamana truncation operator $\Lambda^T_u$ already mentioned in §\ref{S-intro:period} to the right variable of the kernel  $K_{f,\chi}$. Let's now assume that the cuspidal datum $\chi$ is  $(U,U')$-regular in the sense of §\ref{S:UU'reg}.  Then the expression \eqref{intro-eq:KLaTu} does not depend on $T$ and thus is equal to $J_\chi^U(f)$. To state our result we fix a representative $(M_P,\sigma)$ where $M_P$ is a Levi factor of a parabolic subgroup $P=M_PN_P$ of $U$ and $\sigma$ is a cuspidal automorphic representation of $M_P(\AAA)$. Let $\Ac_{P,\sigma,\cusp}(U_h)$ be the space of automorphic forms on the quotient $A_P^\infty M_P(F) N_P(\AAA)\back U(\AAA)$ such that for all $g\in U(\AAA)$
  $$m\in M_P(\AAA)\mapsto \delta_P(m)^{-\frac12} \varphi(mg)$$
  belongs to the $\sigma$-isotypic component of the space of cuspidal automorphic forms on the quotient  $A_P^\infty M_P(F) \back M_P(\AAA)$. Working throughout Langlands spectral decomposition of $K_{f,\chi}$, we get (see Theorem \ref{thm:JchiU}):
\begin{align}\label{eq-intro:JchiU}
  J^U_\chi(f)= \int_{i\ago_P^*}   J_{P,\sigma}^{U}(\la, f)\,d\la.
\end{align}
where  the right-hand side is the absolutely convergent integral of the relative character defined by:
 \begin{align*}
    J_{P,\sigma}^{U}(\la, f)=\sum_{\varphi\in \bc_{P,\sigma}}   \pc_{U'}(I_P(\la,f)\varphi,\la)\overline{   \pc_{U'}(\varphi,\la)  }.
  \end{align*}
Here the periods  $\pc_{U'}(\cdot,\la)$ are those defined in \eqref{eq-intro:IY} and  $I_P(\la,f)$ denotes the induced action of $f$ twisted by $\la$. The sum is over some orthonormal basis $\bc_{P,\sigma}$ of  $\Ac_{P,\sigma,\cusp}(U)$, see §\ref{S:rel-charU}, for the Petersson inner product.
\end{paragr}

\begin{paragr} Let's turn to the linear case namely $G=G_n\times G_{n+1}$. In this case, we have to consider two subgroups namely $H=G_n$  diagonally embedded in $G$ and $G'=G_n'\times G_{n+1}'$ where $G'_n=\GL(n,F)$ is naturally embedded in $G_n=\GL(n,E)$. Let $\chi$ be a cuspidal datum of $G$ and let $f$ be a Schwartz function on $G(\AAA)$. As before we denote by $K_{f,\chi}$ the $\chi$-part of the automorphic kernel. According to \cite[Theorem 1.2.4.1]{BCZ}, the contribution $I_\chi(f)$, as defined by Zydor in \cite{Z3}, is also the constant term of the polynomial exponential (in the variable  $T$) which is asymptotic to  the absolutely convergent integral
  \begin{align*}
      \int_{[H]} \int_{[G'] }   \Lambda_r^TK_\chi (h,g)\,\eta_{G'}(g)dgdh.
  \end{align*}
  Here $ \eta_{G'}$ is the quadratic character defined in §\ref{S:Gn} and $\La_r^T$ is a truncation operator (in the parameter $T$) introduced by Ichino-Yamana and well-suited for the study of Rankin-Selberg period. Assume that $\chi$ is represented by a pair $(M,\pi)$ where $M$ is the Levi factor of a parabolic subgroup $P$ of $G$. We assume also that $\chi$ is $(G,H)$-regular and Hermitian in the sense of §\ref{S:GHgeneric}. Then we have (see Theorem \ref{thm:Ichi})
   \begin{align*}
    I_\chi(f)=2^{-\dim(\ago_L)}\int_{i\ago_M^{L,*}}I_{P,\pi}(\la,f)\,d\la.
   \end{align*}
   Here $L$ is a Levi subgroup of $G$ containing $M$ and determined by $\pi$ (see §\ref{S:choiceL}). For $\la\in i\ago_M^{L,*}$,  the relative character $I_{P,\pi}(\la,f)$ is given by one of the two expressions:
   \begin{align*}
     I_{P,\pi}(\la,f)&=\sum_{\varphi\in \bc_{P,\pi}}  \PP(E(I_P(\la,f)\varphi,\la))\cdot \overline{J(\xi, \varphi,\la)}\\
     &=\sum_{\varphi\in \bc_{P,\pi}} \frac{Z^{RS}(0, W(I_P(\la,f)\varphi,\la))  \overline{\beta_\eta(W(\varphi,\la))}       }{ \bg W(\varphi,\la),W(\varphi,\la)\bd_{\Pet}}.
   \end{align*}
   The sums are over some orthonormal bases $\bc_{P,\pi}$ for the Petersson inner product of  the space $\Ac_{P,\pi,\cusp}(G)$ (defined as above). The first expression is built upon $\PP(E(\varphi,\la))$ and $J(\xi, \varphi,\la)$. The former is the regularized Rankin-Selberg period à la Ichino-Yamana, see \cite{IY}, of the Eisenstein series associated to the pair $(P,\pi)$. The latter is the intertwining (Flicker-Rallis) period of Jacquet-Lapid-Rogawski \cite{JLR}. The second expression uses Whittaker functionals $W(\varphi,\la)$ associated to Eisenstein series and linear forms on the Whittaker models of $\Ind_{P(\AAA)}^{G(\AAA)}(\pi\otimes\la)$ equipped with the Petersson product $\bg \cdot, \cdot\bd_{\Pet}$. The linear forms $Z^{RS}$ and $\beta_\eta$ are counterparts of the period and the intertwining period. In the Rankin-Selberg case, the link is recalled in proposition \ref{prop:ZRS}: it is based on \cite{IY} which  generalizes the classical Rankin-Selberg theory. In the Flicker-Rallis case, the precise relation is given in proposition \ref{prop:Flicker-Whitt} which  generalizes the work of Flicker \cite{Flicker}. Besides some reductions based on \cite[section 9]{BCZ}, the bulk of the proof of proposition is the object of section \ref{sec:Inter}. Note also that we prove in section \ref{sec:Inter} a result that is of independent interest: we express a basic intertwining period of Jacquet-Lapid-Rogawski in terms of an integral of a Whittaker functional, see theorem \ref{thm:J-W}.  The second expression of the relative character is better suited for the proof of Theorem \ref{thm:II}.   In section \ref{chap: Alt proof}, we give an alternative proof of this spectral expansion of $I_\chi(f)$, for $(G,H)$-regular $\chi$ , that is based on the theory of Zeta integrals.

   Finally let us remark that if $\chi$ is a $(U,U')$-regular cuspidal datum of $U$ attached to a cuspidal representation of  $M_{P_n}(\AAA)\times U(h\oplus h_0)(\AAA)$ for  some parabolic subgroup $P_n$ of $U(h)$ then the modified kernel $K^T_{f,\chi}(x,y) $ coincides with the usual $\chi$-part of the kernel. Then one can directly get the expression \eqref{eq-intro:JchiU} in which the $U'$-periods  are absolutely convergent. This includes in particular the Eisenstein series needed for the deduction of the Gan-Gross-Prasad and Ichino-Ikeda conjectures for general Bessel periods, as outlined in subsection \ref{ssect: case of Bessel}.

 \end{paragr}

 \subsection{Organization of the paper}

 \begin{paragr}
   The reader will find the main notations and some prelimary results in Section \ref{sec:prelim}. The Section \ref{sec:spectral-unitary} is devoted to the study of some spectral contributions of the Jacquet-Rallis trace formula for unitary groups. The main general results is Theorem \ref{thm:jfDefU} which gives a more tractable expression to compute spectral contributions. The proof of Theorem \ref{thm:jfDefU} is the bulk of Subsections \ref{ssec:Truncation operator} and \ref{ssec:proofmainthm}. Then in Subsection \ref{ssec:UU'reg-contrib}, the spectral contribution for some cuspidal data are explicitly given in terms of relative characters  (see Theorem \ref{thm:JchiU}). In Section \ref{sec:GH-contrib-JR}, we turn to the Jacquet-Rallis trace formula for general linear groups. The main result is Theorem \ref{thm:Ichi} which expresses  the spectral contribution for some cuspidal data in terms of relative characters. In Subsection \ref{ssec:char-Whittaker}, we show that these relative characters can expressed in terms of Whittaker functionals, see Theorem \ref{thm:IPi=Ipi}. The main result of Section \ref{sec:Inter} is Theorem \ref{thm:J-W} which relates some basic  intertwining period of Jacquet-Lapid-Rogawski to some integral of a Whittaker functional: it is used in the proof of Theorem \ref{thm:IPi=Ipi} but it is also of independent interest. Its proof occupies the whole part of Section \ref{sec:Inter}. As explained above, we provide in Section \ref{chap: Alt proof} an alternative proof for the description of the spectral contributions to the Jacquet-Rallis trace formula for general linear groups that can be obtained by combining Theorem \ref{thm:Ichi} with Theorem \ref{thm:IPi=Ipi}. In Section \ref{sec:Preuve Eis}, we explain the proof of Theorems \ref{thm:GGP} and \ref{thm:II}  based on the comparison of the Jacquet-Rallis trace formulas and the results obtained before. The aim of the final Section \ref{sect: application} is to establish the reduction of Theorems \ref{thm:GGP-Bessel} and \ref{thm:II-Bessel} to special cases of Theorems \ref{thm:GGP} and \ref{thm:II} respectively. Its most technical part is in Subsection \ref{ssect: unr computation} where necessary unramified computation is performed. Finally, Appendix \ref{appendix Weyl character formula} presents a, probably well-known, extension of Weyl's character formula to non-connected groups that is necessary for the unramified computation.
 \end{paragr}

\section{Preliminaries}\label{sec:prelim}

\subsection{Algebraic and adelic groups}

\begin{paragr}
  We shall try to follow the usual notations of Arthur and the main notations of  the previous article \cite{BCZ}. For the reader's convenience, we briefly recall our choices.
\end{paragr}

\begin{paragr}
  We denote by $F$ a number field, $V_F$ (resp. $V_{F,\infty}$) the set of its places (resp. Archimedean places)  and $\AAA$ its ring of adèles. For $v\in V_F$, let $F_v$ be its completion at $v$. For any finite subset $S\subset V_F$, we set $F_S=\otimes_{v\in S}F_v$ and $F_\infty=F_{V_{F,\infty}}$. We denote by  $ | \cdot | $  the morphism $ \AAA^\times \to \RR_+^\times $ given by the product of normalized absolute values $|\cdot|_v$ on each $ F_v $. 
\end{paragr}

\begin{paragr}  Let   $G$ be a reductive group defined over $F$.   All the subgroups of $G$ we consider are assumed to be defined over $F$. We fix  $P_0\subset G$  a minimal parabolic subgroup and $M_0$ a Levi factor of $P_0$. A parabolic subgroup of $G$ which contains $P_0$, resp. $M_0$, is said to be standard, resp. semi-standard. Let $P$ be a  semi-standard parabolic subgroup of $G$. It has a Levi decomposition $P=M_PN_P$ such that $M_P$ is a semi-standard Levi factor (that is $M_0\subset M_P$) and $N_P$ is the unipotent radical of $P$. Such a group $M_P$ is called a semi-standard Levi subgroup of $G$. It is said to be standard if moreover $P$ is standard. Let $X^*(P)$ be the group of rational characters of $P$ defined over $F$. Attached to $P$ are real vector spaces $\ago_P^*=X^*(P)\otimes_\ZZ\RR$ and $\ago_P=\Hom_\ZZ(X^*(P),\RR)$ in canonical duality:
\begin{align}\label{eq:pairing}
\bg \cdot,\cdot\bd : \ago_P^*\times \ago_P\to \RR.
\end{align}
  If $P\subset Q\subset G$, we have natural  maps  $\ago_Q^*\to \ago_P^*$ and $\ago_P\to \ago_Q$. The kernel of the second one is denoted by $\ago_P^Q$. We have natural decomposition $\ago_P=\ago_P^Q\oplus \ago_Q$ and dually $\ago_P^*=\ago_P^{Q,*}\oplus \ago_Q^*$. We put a subscript $\CC$ to denote the extension of scalars to $\CC$. Then one has a decomposition
$$\ago_{P,\CC}^{Q,*}=\ago_{P}^{Q,*}\oplus i\ago_{P}^{Q,*}$$
where $i^2=-1$.  We shall denote by $\Re$ and $\Im$ the associated projections. The complex conjugate is then defined by $\bar\la=\Re(\la)-i\Im(\la)$. Note that the spaces $\ago_P^Q$, $\ago_P^{Q,*}$ depend only on the Levi factors $M_Q$ and $M_P$ and thus are also denoted by $\ago_{M_P}^{M_Q}$, $\ago_{M_P}^{M_Q,*}$ etc.

Let  $\Ad_{P}^Q$ be the adjoint action of $M_P$ on the Lie algebra of $M_Q\cap N_P$. Let  $\rho_P^Q$ be the unique element of $\ago_{P}^{Q,*}$ such that for every $m\in M_P(\AAA)$ we have
\begin{align*}
  |\det(\Ad_P^Q(m))|=\exp(\bg 2\rho_P^Q,H_P(m)\bd).
\end{align*}
We set $\rho_P=\rho_P^G$. For every $g\in G(\AAA)$, we set
\begin{align*}
   \delta^Q_P(g)=\exp(\bg 2\rho^Q_P,H_P(g)\bd).
\end{align*}
 Usually we replace the subscript $P_0$ simply by $0$ e.g. $\ago_0=\ago_{P_0}$, $\rho_0=\rho_{P_0}$ etc.
\end{paragr}

\begin{paragr}\label{S:root}Let $P $ be a  standard parabolic subgroups of $G$.   Let $A_P=A_{M_P}$ be   the maximal central $F$-split torus of $M_P$. Let $\Delta_{0}^P$ be the set of simple roots of  $A_{0}$ in  $M_P\cap P_0$.  Let $\Delta_P^Q$ be the image of $\Delta_{0}^Q\setminus \Delta_{0}^P$  by the projection $\ago_{0}^*\to \ago_P^*$. It is a basis of $\ago_P^{Q,*}$ whose dual basis is the subset of coweights $\hat\Delta_P^{Q,\vee}\subset \ago_P^Q$. We have also the set of coroots $\Delta_P^{Q,\vee}$ which is a basis of $\ago_P^{Q}$. We have the dual basis given by the  set of simple weights $\hat{\Delta}_P^Q \subset \ago_P^{Q,*}$.  The sets $\Delta_P^Q$ and  $\hat{\Delta}_P^Q$ determine open  cones  in  $\ago_{0}$ whose characteristic functions are denoted respectively by $\tau_P^Q$ and  $\hat{\tau}_P^Q$. We set

\begin{align*}
     \ago_P^{*,Q+}=\left\{\lambda\in \ago_P^*\mid \langle \lambda,\alpha^\vee\rangle \geqslant 0,\; \forall \alpha^\vee\in \Delta_P^{Q,\vee} \right\}.
\end{align*}
We define similarly $\ago_P^{Q+}$ using roots instead of coroots. If $Q=G$, the exponent $G$ is omitted.
\end{paragr}

\begin{paragr} \label{S:Weyl}Let $W$ be the Weyl group of $(G,A_{0})$ namely  the  quotient of the normalizer of  $A_0$ in $G(F)$by $M_0$. For  $P=M_P N_P$ and  $Q=M_QN_Q$ two standard parabolic subgroups of $ G $, we denote by $W(P,Q)$ or $W(M_P,M_Q)$ the set of $w\in W$ such that  $w\Delta^P_0=\Delta_0^Q$. For $w\in W(P,Q)$, we have  $wM_Pw^{-1}=M_Q$. When $P=Q$, the group $W(P,P)$ is simply denoted by $W(P)$ or $W(M_P)$.  We will also write $W^{M_P}$ for the Weyl group of $(M_P,A_0)$.
\end{paragr}

\begin{paragr}\label{S:goodcpt}  Let $ K=\prod_{v\in V_F} K_v\subset G (\AAA) $ be a ``good'' maximal compact subgroup in good position relative to $M_0$ (called ``admissible'' in \cite[p. 9]{arthur2}). We write
\begin{align*}
K=K_\infty K^\infty
\end{align*}
where  $K_\infty=\prod_{v\in V_{F,\infty}} K_v$ and $K^\infty=\prod_{v\in V_F\setminus V_{F,\infty}} K_v$. We have a homomorphism $H_P:P(\AAA)\to \ago_P$ such that  $\bg \chi, H_P(p) \bd = \log | \chi (p) | $ for any $p\in P(\AAA)$ and $\chi\in X^*(P)$. By Iwasawa decomposition $G(\AAA)=P(\AAA)K$, it extends to a map $H_P:G(\AAA)\to \ago_P$ which is left-invariant by $M_P(F)N_P(\AAA)$ and right-invariant by $K$. We denote $A_P^\infty=A_{M_P}^\infty$ the neutral component of the group of real points of the maximal $\QQ$-split torus of the Weil restriction $\Res_{F/\QQ}(A_P)$. Then the restriction of $H_P$ to $A_P^\infty$ is an isomorphism.

We set  $[G]_P=M_P(F)N_P(\AAA)\back G(\AAA)$ and $[G]_{P,0}=A_P^\infty M_P(F) N_P(\AAA)\back G(\AAA)$. Let $[G]^1_P$ be the subset of $[G]_P$ where the map $H_P$ vanishes. If $P=G$ we shall omit the subscript $P$.
\end{paragr}

\begin{paragr} \label{S:hauteurs}We fix a height  $\|\cdot\|$ on $G(\AAA)$ as in \cite[section 2.4]{BCZ}.  Let $P\subset G$ be a standard parabolic subgroup. We set for $x\in [G]_P$
	\begin{equation*}
\| x\|_P=\inf_{\gamma\in M_P(F)N_P(\AAA)} \| \gamma x\|.
	\end{equation*}
      \end{paragr}

      \begin{paragr}[Haar measures.] --- \label{S:Haar}Let's explain briefly our choice and notations  of measures, see \cite[section 2.3]{BCZ} for more details.          We fix  a non-trivial continuous additive character $\psi':\AAA/F\to \CC^\times$. For each place $v\in V_F$, the local component $\psi'_v$ of $\psi'$ determines an autodual Haar measure on $F_v$. The choice of an invariant rational volume form on $G$ then determines a Haar measure $dg_v$ on $G(F_v)$. We have an Artin-Tate $L$-function $L_G(s)=  \prod_{v\in V_F} L_{G,v}(s)$, see \cite{Gros} and more generally  $L_G^S(s)=  \prod_{v\in V_F\setminus S} L_{G,v}(s)$ for $S\subset V_F$ finite. Then $\Delta^{S,*}_G$ (simply  $\Delta^*_G$  if $S$ is empty) is defined to be  the leading coefficient in the Laurent expansion of $L_G(s)$ at $s=0$. We also set  $\Delta_{G,v}=L_{G,v}(0)$. The  Tamagawa measure $dg$ on $G(\AAA)$ is defined by $dg=dg_S\times dg^S$ where $dg_S=\prod_{v\in S} dg_v$ and $dg^S=(\Delta_G^{S,*})^{-1}\prod_{v\notin S} \Delta_{G,v}dg_v$ for  $S\subset V_F$ finite.

We equip $\ago_P$ with the Haar measure that gives a covolume $1$ to the lattice $\Hom(X^*(P),\ZZ)$ and  $i\ago_{P}^{*}$ with  the dual Haar measure. The group $ A_P ^ \infty $  is equipped with the Haar measure compatible with the isomorphism $ A_P ^ \infty \simeq  \ago_ {P} $ induced by the map $H_P$. The groups $\ago_P^G\simeq \ago_P/\ago_G$ and  $i\ago_P^{G,*}\simeq i\ago_P^*/i\ago_G^*$ are provided with the quotient Haar measures.

 The homogeneous space $[G]$ (resp. $[G]^1\simeq [G]_0$) is equipped with the quotient of the Tamagawa measure on $G(\AAA)$ by the counting measure on $G(F)$ (resp. by the product of the counting measure on $G(F)$ with the Haar measure we fixed on $A_G^\infty$). For $P$ a standard parabolic subgroup, we equip similarly $[G]_{P}$ with the quotient of the Tamagawa measure on $G(\AAA)$ by the product of the counting measure on $M_P(F)$ with the Tamagawa measure on $N_P(\AAA)$. Since the action by left translation of $a\in A_P^\infty$ on $[G]_P$ multiplies the measure by $\delta_P(a)^{-1}$, taking the quotient by the Haar measure on $A_P^\infty$ induces a ``semi-invariant'' measure on $[G]_{P,0}=A_P^\infty\backslash [G]_P$.
\end{paragr}

\subsection{Space of functions}

\begin{paragr}
For two positive functions  $f$ and $g$ on a set $X$, we write $ f(x)\ll g(x), \; x\in X$ if there exists a constant $C>0$ such that $f(x)\leqslant Cg(x)$ for every $x\in X$. We write $ f(x)\sim g(x), \; x\in X$ if $f(x)\ll g(x)$ and $g(x)\ll f(x)$.
\end{paragr}

\begin{paragr}
For every $C\in \RR\cup \{-\infty\}$ with $D>C$, we set $\cH_{>C}=\{z\in \CC\mid \Re(z)>C \}$ 
\end{paragr}

\begin{paragr}[Schwartz space.] ---  \label{S:Schwartz} As before  $G$ is a reductive group defined over $F$.   We let $\ggo_\infty$ be the Lie algebra of $G(F_\infty)$ and $\uc(\ggo_\infty)$ be the enveloping algebra of its complexification and $\zc(\ggo_\infty)\subset \uc(\ggo_\infty)$ be its center.
  
  We shall briefly review three useful  locally convex topological spaces of functions, see \cite[section 2.5]{BCZ} for more details. First  $\Sc(G(\AAA))$ is  the Schwartz space of $G(\AAA)$: it  contains  the dense subspace $\Cc(G(\AAA))$ of smooth and compactly supported functions.

Second the  Schwartz space $\Sc([G]_P)$ of $[G]_P$ is the space of smooth functions $\varphi:[G]_P\to \CC$ such that for every $N>0$ and $X\in \uc(\ggo_\infty)$ we have
$$\displaystyle \lVert \varphi\rVert_{N,X}=\sup_{x\in [G]_P}\lVert x\rVert_{P}^{N}\lvert (R(X)\varphi)(x)\rvert<\infty.$$

Third the space of functions of uniform moderate growth  on $[G]_P$ is defined as
\begin{equation*}
\tc([G]_P)=\bigcup_{N>0} \tc_N([G]_P).
\end{equation*}
where  $\tc_N([G]_P)$ is the space of smooth functions $\varphi: [G]_P\to \CC$ such that for every $X\in \uc(\ggo_\infty)$ we have
\begin{align*}
 \lVert \varphi\rVert_{-N,X}=\sup_{x\in [G]_P}\lVert x\rVert_{P}^{-N}\lvert (R(X)\varphi)(x)\rvert <\infty.
\end{align*}
\end{paragr}

\begin{paragr}[Automorphic forms.] ---  \label{S:automor}For this §,  we refer the reader to \cite[section 2.7]{BCZ} for more details.  The space $\Ac_P(G)$ of  automorphic forms on $[G]_P$ is  the subspace of $\zc(\ggo_\infty)$-finite functions in $\tc([G]_P)$. Let $\Ac_{P,\cusp}(G)\subset \Ac_P(G)$ be the subspace of cuspidal automorphic forms: these are the functions   $\varphi\in \Ac_P(G)$ such that $\varphi_Q=0$  for every proper parabolic subgroup $Q\subsetneq P$. The constant term  $\varphi_Q$ is defined by
  \begin{align*}
    \varphi_Q(x)=\int_{[N_Q]} \varphi(ux) du
  \end{align*}
   for all $x\in G(\AAA)$.
   A cuspidal  automorphic representation  $\pi$ of $M_P(\AAA)$ is  a topologically irreducible subrepresentation of $\Ac_{\cusp}(M_P)$. For every $\lambda\in \ago_{P,\CC}^*$, we define  $\pi_\lambda=\pi\otimes \lambda$ as the space of functions of the form
   \begin{align*}
     m\in [M_P]\mapsto \exp(\langle \lambda, H_P(m)\rangle) \varphi(m)
   \end{align*}
   for $ \varphi\in \pi$.

   We denote by $\Ac_{\pi,\cusp}(M_P)$  the $\pi$-isotypic component of $\Ac_{\cusp}(M_P)$.  The normalized smooth induction  $\Ind_{P(\AAA)}^{G(\AAA)}(\Ac_{\pi,\cusp}(M_P))$ is denoted by $\Ac_{P,\pi,\cusp}(G) $ and is identified with the space of automorphic forms $\varphi\in \Ac_{P}(G)$ such that
$$m\in [M_P]\mapsto \exp(-\langle \rho_P, H_P(m)\rangle) \varphi(mg)$$
belongs to $\Ac_{\pi,\cusp}(M_P)$  for every $g\in G(\AAA)$.  The algebra $\Sc(G(\AAA))$ acts on $\Ac_{P,\pi,\cusp}(G)$  by right convolution. For every $\lambda\in \ago_{P,\CC}^*$, we  denote by $I(\la)$ the action on $\Ac_{P,\pi,\cusp}(G)$ we get by transport from the action of $\Sc(G(\AAA))$ on $\Ac_{P,\pi_\la, \cusp}$ and the identification $\Ac_{P,\pi,\cusp}\to \Ac_{P,\pi_\la,\cusp}$ given by $\varphi\mapsto \exp(\langle \lambda, H_P(.)\rangle)\varphi$.
Assume that  the central character of $\pi$ is unitary. Then we equip $\Ac_{P,\pi,\cusp}(G)$  with the Petersson inner product given by 
\begin{align*}
  \| \varphi\|^2_{\Pet}= \displaystyle \langle \varphi,\varphi\rangle_{\Pet}=\int_{[G]_{P,0}} \lvert \varphi(g)\rvert^2 dg
\end{align*}

Note that for every cuspidal  automorphic representation  $\pi$, there exists a unique  $\lambda\in \ago_{P,\CC}^*$ such that the central character is trivial on $A_M^\infty$. Most of the time,  the cuspidal  automorphic representations we consider  are implicitly assumed to have a central character which is trivial on $A_M^\infty$.

Let $\hat K$ be the set of isomorphism classes of irreducible unitary representations of $K$. A $K$-basis $\bc_{P,\pi}$ of  $\Ac_{P,\pi,\cusp}(G)$ is by definition the union over of $\tau\in \hat K$ of orthonormal bases $\bc_{P,\pi,\tau}$  for the Petersson inner product of  the finite dimensional  subspaces $\Ac_{P,\pi,\cusp}(G,\tau)$ of functions in $\Ac_{P,\pi,\cusp}(G)$  which transform under $K$ according to $\tau$.

For any $\varphi\in \Ac_{P,\cusp}(G)$, $g\in G(\AAA)$ and $\la\in \ago_{P,\CC}^{*}$, we denote by
  \begin{align*}
    E(g,\varphi,\la)=\sum_{\delta\in P(F)\back G(F)} \exp(\bg \la,H_P(\delta g)\bd \varphi(\delta g)
  \end{align*}
  the Eisenstein series which is analytic. The right-hand side is convergent for $\Re(\la)$ in a suitable cone. Let $P$ and $Q$ be  standard parabolic subgroups of $G$. For any $w\in W(P,Q)$ and $\la\in  \ago_{P,\CC}^{*}$, we have the intertwining operator
  \begin{align*}
    M(w,\la):\Ac_P(G)\to \Ac_{Q}(G).
  \end{align*}
For more details and continuity properties of these constructions we refer the reader to \cite[§§ 2.7.3, 2.7.4]{BCZ}.
\end{paragr}

\begin{paragr}[Cuspidal datum.] --- \label{S:cusp-datum} Let $P\subset G$ be a standard parabolic subgroup of $G$. Let $L^2([G]_P)$ be the space of square integrable functions on $[G]_P$.
  Let  $\Xgo(G)$ be the set of cuspidal data of $G$ Recall that $\Xgo(G)$ is the quotient of  the set of pairs $(M_P,\pi)$ such that 
\begin{itemize}
\item $P$ is a standard parabolic subgroup of $G$;
\item $\pi$ is the isomorphism class of a cuspidal automorphic representation of $M_P(\AAA)$ with central character trivial on $A_P^\infty$.
\end{itemize}
by the following equivalence relation: $(M_P,\pi)\sim (M_Q,\tau)$ if there exists $w\in W(P,Q)$ such that $w\pi w^{-1}\simeq \tau$.  Then for any standard parabolic subgroup $P$ of $G$, we have  the Langlands decomposition (see \cite[section 2.9]{BCZ} for more details):
\begin{align*}
 L^2([G]_P)=\widehat{\bigoplus_{\chi\in \Xgo(G)}} L^2_{\chi}([G]_P).
\end{align*}
The Schwartz algebra acts on   $L^2([G]_P)$  and for each $\chi\in \Xgo(G)$ on  $L^2_{\chi}([G]_P)$ by  right convolution. For each $f\in \Sc(G(\AAA))$ we get integral operators whose kernela are  denoted respectively by $K_{f,P}$, resp.   $K_{f,P,\chi}$. We shall use also the decomposition
\begin{align*}
 L^2([G]_{P,0})=\widehat{\bigoplus_{\chi\in \Xgo(G)}} L^2_{\chi}([G]_{P,0}).
\end{align*}
\end{paragr}

\subsection{Around Arthur's partition}

\begin{paragr}
  Let $T_1,T\in \ago_0$ and $P$ be  a standard parabolic subgroup of $G$. We define
  \begin{align*}
    A_{P_0}^{P,\infty}(T_1)=\{ a\in A_0^\infty\mid \langle \alpha, H_0(a)\rangle\geqslant \langle \alpha, T_1\rangle , \; \forall \alpha\in \Delta^P_0 \}
  \end{align*}
and
\begin{align*}
    A_{P_0}^{P,\infty}(T_1,T)=\{ a\in  A_0^{P,\infty}(T_1)\mid \langle \varpi, H_0(a)\rangle\leqslant \langle \alpha, T\rangle , \; \forall \varpi\in \hat\Delta^P_0 \}.
  \end{align*}
\end{paragr}

\begin{paragr}[Siegel domains.] --- \label{S:Siegel} Let $T_-\in \ago_0^{G,+}$  and $\omega_0\subset P_0(\AAA)^1$ be a compact subset such that $P_0(\AAA)^1=P_0(F)\omega_0$. Let $P$ be a standard parabolic subgroup of $G$. We define
\begin{equation*}
 \sgo_P(\om_0,T_-,K)=\omega_0 A_{P_0}^{P,\infty}(T_-)K.
\end{equation*}
There exists $T_-\in -\ago_0^{G,+}$ such that for all standard parabolic subgroup of $G$ we have  $G(\AAA)=P(F)\sgo_P(\om_0,T_-,K)$.We will fix $T_-$ and $\om_0$ and  we set $\sgo_P= \sgo_P(\om_0,T,K)$: this is a {\em Siegel domain} of $[G]_P$.
\end{paragr}

\begin{paragr}\label{S:partition}
  For any standard parabolic subgroup $P$ of $G$ and $T\in \ago_0^+$, we consider the characteristic function $F^P(\cdot,T)$ of $\om_0 A_0^{P,\infty}(T_-,T)K\subset \sgo_P$. This function descends on $[G]_P $. There is  a point $T_0\in \ago_0^{G,+}$ such that for all $T\in T_0+\ago_0^{G,+}$ and all standard parabolic subgroup $Q$ of $G$ we have the following formula which gives a partition of $[G]_Q$ (see \cite[Lemma 6.4]{ar1} and also \cite[Proposition 3.6.3]{labWal}):
\begin{align}\label{eq:partition}
  \sum_{P_0\subset P\subset Q} \sum_{\delta \in P(F)\back Q(F)}F^P(\delta g,T) \tau_P^Q(H_P(\delta g)-T)=1.
\end{align}
Note that the relation implies a simpler definition of $F^P(\cdot,T)$ for $T\in T_0+\ago_0^{G,+}$: the function  $F^P(\cdot,T)$ is the characteristic function of the following set:
\begin{align}\label{eq:FQ-explicit}
  \{g\in [G]_P \mid \forall \varpi\in \hat\Delta^P_0, \delta\in P(F)  \, \bg \varpi,H_0(\delta g)-T\bd \leq 0\}. 
\end{align}
Indeed this is a  straightforward consequence of \cite[Lemma 2.1] {ar-unip} and the definition of the operator $\Lambda^T$ there.

We shall use several times the following simple lemma.

\begin{lemme}
  \label{lem:Ftau-Siegel}
  Let $P\subset Q$ be parabolic subgroup of $G$. Let $g\in G(\AAA)$ be such that $F^P( g,T) \tau_P^Q(H_P( g)-T)=1$. Then we have
  \begin{align*}
    \forall \al\in \Delta_0^Q\setminus \Delta_0^P, \ \bg \al, H_0(g)\bd  > \bg \al, T\bd.
  \end{align*}
In particular, if $g\in \sgo_P$ then $g\in \sgo_Q$.
\end{lemme}

\begin{preuve}
  Let $\al\in \Delta_0^Q\setminus \Delta_0^P$. We  write $\al=\al_P+\al^P$ according to the decomposition $\ago_0^*= \ago_P^*\oplus \ago_0^{P,*}$.  Since  we have $\bg \al^P,\be^\vee\bd =\bg \al,\be^\vee\bd\leq 0$ for all $\be\in \Delta_0^P$,  we know that $\al^P$ is a non positive linear combination of elements of $\hat\Delta_0^P$. Thus the condition $F^{P}(g,T) =1$ implies  that 
\begin{align*}
  \bg \al^P, H_0(g)\bd  \geq \bg \al^P, T\bd.
\end{align*}
The condition $\tau_P^{Q}(H_P(g)-T)=1$ implies that we have  $\bg \al_P, H_0(g)\bd  > \bg \al_P, T\bd$. We conclude that $\bg \al, H_0(g)\bd  > \bg \al, T\bd$ for all $\al\in \Delta_0^Q\setminus \Delta_0^P$. The last assertion is then obvious since $T$ is positive.
\end{preuve}
\end{paragr}

\begin{paragr}\label{S:fct-d}
Let $\phi:\ago_0\to \CC$ be a function and let $T\in T_0+\ago_0^{G,+}$. 
  For any standard parabolic subgroup  $Q$  of $G$, we define  for all $g\in [G]_Q$
  \begin{align*}
    d_Q(\phi,g,T)= \sum_{P_0\subset P\subset Q} \sum_{\delta \in P(F)\back Q(F)}F^P(\delta g,T) \tau_P^Q(H_P(\delta g)-T)\exp(\phi(H_P(\delta g))).
  \end{align*}
We will mainly use $d_Q(\la,g,T)$ for $\la\in \ago_0^*$. 
  
  \begin{proposition}\label{prop:equiv}Let $\la\in \ago_0^*$ and $Q\subset G$ a parabolic subgroup.
    \begin{enumerate}
    \item Let $\tc \subset T_0+\ago_0^{G,+}$ be  a compact subset 
 We have
$$\exp(\bg \la, H_0(g)\bd) \sim d_Q(\la,g,T)$$
for all $g \in \sgo_Q$ and $T\in\tc$. In particular, all functions $d_Q(\la,T)$ are equivalent for $T\in \tc$.
\item Let $\kc$ be a compact subset of $G(\AAA)$ and  let   $T\in T_0+\ago_0^{G,+}$, we have
$$d_Q(\la,g,T)\sim d_Q(\la ,gc,T)$$
for all $g \in [G]_Q$, $c\in \kc$.
\end{enumerate}
  \end{proposition}

  \begin{preuve}
 1. If we write $\la=\la_Q+\la^Q$ according to the decomposition $\ago_0^*=\ago_Q^{*}+\ago_0^{Q,*}$, we see that  $d_Q(\la,g,T)=\exp(\bg \la_Q,H_Q(g)\bd) d_Q(\la^Q,g,T)$ for all $T\in T_0+\ago_0^+$ and all $g\in [G]_Q$ . Thus we may and shall assume that $\la\in \ago_0^{Q,*}$.  
Let $g\in \sgo_Q$ and $T\in \tc$. There exists  $P \subset Q$  is such that $F^P(g,T)\tau_P^Q(H_P(g)-T)=1$. Then, by definition, we get $d_Q(\la,g,T)=\exp(\bg \la, H_P(g)\bd) $. It suffices to prove that $H_0(g)-H_P(g)$ stays in a compact subset which depends only on $\tc$ for any $g\in \sgo_P$ such that  $F^P(g,T)=1$. In fact the latter  condition  implies $\bg \varpi, H_0(g)\bd \leq \bg \varpi, T\bd$ for all $\varpi\in \hat\Delta^P_0$ and $g\in \sgo_P$ implies that $\bg \al, H_0(g)\bd \geq \bg \al, T_-\bd$ for all $\al\in \hat\Delta^P_0$. Hence the projection of $H_0(g)$ on $\ago_0^{P}$ stays in a fixed compact subset but this projection is nothing else but $H_0(g)-H_P(g)$.

2. First we observe that $H_0(kc)$ stays in a fixed compact for $k\in K$ and $c\in\kc$. By assertion 1, we may replace $T$ by any element in $T_0+\ago_0^{G,+}$. In particular, we may and shall  assume that $T$ is such that $T-H_0(kc)\in  T_0+\ago_0^{G,+}$ for all $k\in K$ and $c\in \kc$. For any $g\in G(\AAA)$, we shall denote $k(g)$ an element of $K$  such that $gk(g)^{-1}\in P_0(\AAA)$.

Let $g\in G(\AAA)$ and  $c\in \kc$. Then there exist  a unique parabolic sugroup  $P \subset Q$   and $\delta\in Q(F)$ such that $F^P(\delta gc,T)\tau_P^Q(H_P(\delta gc)-T)=1$. Observe that $H_0(\delta gc)=H_0(\delta g)+H_0(k(\delta g) c)$. We deduce that on the one hand  that $\tau_P^Q(H_P(\delta g)-(T-H_0(k(\delta g)c))=1$ is equivalent to  $\tau_P^Q(H_P(\delta gc)-T)=1$ and on the other hand that we also have $F^P(\delta g,T-H_0(k(\delta g)c))=1$: indeed for all $\varpi\in \hat\Delta^P$ we have
\begin{align*}
  \bg \varpi,H_0(\delta g)\bd&=  \bg \varpi,H_0(\delta gc)\bd -  \bg \varpi,H_0(k(\delta g) c)\bd \\
&\leq   \bg \varpi,T  -  H_0(k(\delta g) c)\bd 
\end{align*}
since $F^{P}(\delta gc,T)=1$. 
From this, we get
\begin{align*}
  d_Q(\la, g,  T-H_0(k(\delta g)c))=\exp(\bg \la, H_P(\delta g)\bd)=d_Q(\la,gc,T)\exp(-\bg \la, H_0(k(\delta g)c)\bd).
\end{align*}
Using assertion 1 and the fact that  $H_0(k(\delta g) c)$ stays in a fixed compact set, we can easily conclude.

\end{preuve}

  \begin{lemme}\label{lem:egalite-d}Let $Q\subset R$ be standard parabolic subgroups of $G$. 
There exists $T_1\in T_0+\ago_0^{G,+}$ such that for all $T\in T_1+\ago_0^{G,+}$, all $\la\in \ago_0^*$ and all   $g\in G(\AAA)$ such that $F^{Q}(g,T)\tau_Q^{R}(H_Q(g)-T)=1$, we have 
$$d_Q(\la,g,T_0)=d_R(\la,g,T_0).$$
    \end{lemme}

    \begin{preuve} 
Both the hypothesis and the conclusion are invariant under left translation by $Q(F)$. It suffices to prove the result for $g\in \sgo_Q$ such that $F^{Q}(g,T)\tau_Q^{R}(H_Q(g)-T)=1$. First, there exists a  unique standard parabolic subgroup $P\subset Q$ such that  $F^{P}(g,T_0)\tau_P^{Q}(H_Q(g)-T_0)=1$. In the proof of proposition \ref{prop:equiv} we have shown that $H_0(g)-H_P(g)$ stays in a fixed compact subset for $g\in \sgo_Q$ such that $F^{P}(g,T_0)=1$. In particular we can assume that  for such elements $g$ we have
\begin{align}\label{eq:T_1}
  \bg \al, T_1 \bd > \bg \al_P,T_0\bd +\bg \al, H_0(g)-H_P(g)\bd.
\end{align}

By lemma \ref{lem:Ftau-Siegel}, we have  $\bg \al, H_0(g)\bd  > \bg \al, T\bd$ for all $\al\in \Delta_0^R\setminus \Delta_0^Q$. So we have for all   $\al\in \Delta_0^R\setminus \Delta_0^Q$
\begin{align*}
  \bg \al_P, H_0(g)\bd =  \bg \al, H_P(g)\bd&=    \bg \al, H_0(g)\bd+ \bg \al, H_P(g)-H_0(g)\bd\\
&>  \bg \al, T\bd+ \bg \al, H_P(g)-H_0(g)\bd\\
&\geq \bg \al, T_1 \bd + \bg \al, H_P(g)-H_0(g)\bd\\
&>\bg \al_P,T_0\bd
\end{align*}
by \eqref{eq:T_1}. In particular, we see that we have $\tau_P^{R}(H_Q(g)-T_0)=1$. Thus by definition $d_R(\la,g,T_0)=\exp(\bg \la,H_P(g)\bd)=d_Q(\la,g,T_0)$.
    \end{preuve}

\begin{proposition}\label{prop:unicity} Let $Q\subset R$ be standard parabolic subgroups of $G$ such that $\Delta_0^R\setminus  \Delta_0^Q =\{\al\}$ for some simple root $\al$. Let $c$ such that $c > \bg \al, T_{0,P}\bd$ for all  parabolic subgroups $P\subset Q$. 

For any $g\in G(\AAA)$ there is at most one element $\delta \in Q(F)\back R(F)$ such that $d_Q(\al,\delta g,T_0)>\exp(c)$.
 \end{proposition}

\begin{preuve}
  Let $g\in G(\AAA)$ such that $d_Q^R(g,T_0)>\exp(c)$. Using left translations by $Q(F)$, we may and shall assume that there exists a parabolic subgroup    $P\subset Q$  such that $F^P(g,T_0) \tau_P^Q(H_P( g)-T_0)=1$. Then the condition $d_Q^R(\al, g,T_0)>\exp(c)$ is then equivalent to  $\bg \al, H_P(g)\bd >c$ and so $\bg \al, H_P(g)\bd >\bg \al, T_{0,P}\bd$. So we have also $F^P(g,T_0) \tau_P^R(H_P( g)-T_0)=1$. Now the uniqueness follows from the partition \eqref{eq:partition} applies to $[G]_R$.

\end{preuve}
\end{paragr}


\input{Contrib-unit}


\input{GHregular}


\input{Inter-Whittaker}


\input{FRspectral}


\input{Preuve}


\input{Besselcomputation}

\bibliography{biblio}
\bibliographystyle{alpha}

\end{document}

%% file: Contrib-unit.tex
\section{On the spectral expansion of the Jacquet-Rallis trace formula for unitary groups}\label{sec:spectral-unitary}

\subsection{Notations}\label{S:nota-SpecU}

\begin{paragr}
  Let  $E/F$ be a quadratic extension of number fields and $c$ be the non-trivial element of the Galois group $\Gal(E/F)$. Let $\AAA$ be the ring of adèles of $F$. Let $n\geq 0$ be an integer. Let $\hc_n$ be the set of isomorphism classes of pairs $(V,h)$  where $V$ is a $E$-vector space of dimension $n$ and $h$ a non-degenerate $c$-Hermitian form on $V$. For any  $(V,h)\in \hc_n$, we  identify $(V,h)$ with a representative and we shall denote by $U(V,h)$ or simply $U(h)$  its automorphisms group.  We will fix $(V_0,h_0)\in \hc_1$.
\end{paragr}

\begin{paragr}
  We attach to any $h\in \hc_n$ the following algebraic groups over $F$:
\begin{itemize}
\item the unitary group $U'_h=U(h)$ of automorphisms of $h$;
\item  the unitary group $U''_h= U(h\oplus h_0) $ where $h\oplus h_0$ denoted the orthogonal sum;
\item the product of unitary groups $U_h=U_h'\times U_h''$.
\end{itemize}
We have an embedding 
\begin{align}
  \label{eq:U(h)}
U_h'\hookrightarrow U_h'' \ \ \ g\mapsto \begin{pmatrix}
 g & 0 \\  0 & 1
\end{pmatrix}
\end{align}
and  a diagonal embedding $U'_h\hookrightarrow U_h$.  
\end{paragr}

\begin{paragr}\label{S:unitary-choices} Let $n\geq 1$ an integer and $(V,h)\in \hc_n$. We fix a parabolic subgroup $P_0'$ of $U_h'$ and a Levi factor $M_0'$ of $P_0'$ both defined over $F$ and minimal for these properties. As a parabolic subgroup, $P_0'$ is the stabilizer of a flag of totally isotropic subspaces of $V$. Let $P_0''$ be the parabolic subgroup  of $U''_h$ stabilizing the same flag. Let $M_0''$ be the unique Levi factor of $P_0''$ that contains $M_0'$. Then $P_0=P_0'\times P_0''$  is a parabolic subgroup of $U_h$ with Levi factor  $M_0=M_0'\times M_0''$ .  We fix also a pair $(B_0'',T_0'')$ consisting of a minimal parabolic $F$-subgroup of $U_h''$ and a Levi factor defined over $F$. We may and shall assume $B_0''\subset P_0''$ and $T_0''\subset M_0''$. We fix maximal good compact subgroups $K_h'\subset U_h'(\AAA)$ and $K_h''\subset U_h''(\AAA)$ respectively in good position relative to $M_0'$ and $T_0''$ in the sense of §\ref{S:goodcpt}. We set $K_h=K_h'\times K_h''\subset U_h(\AAA)$. Any $g\in U_h=U_h'\times U_h''$ will be written $(g',g'')$ without any further comment. 

From now on $(V,h)$ is fixed and  we shall omit it in the notations: thus we have $U=U_h$, $U'=U'_h $ and so on. 
\end{paragr}

\begin{paragr} Let $\fc_0'$  be the the set of parabolic subgroups of  $U'$  that contain $P_0'$.  For any $P'\in \fc_0'$, let $P''$ be the parabolic subgroups of  $U'$ which stabilizes the flag of totally isotropic subspaces of $V$ that defines $P'$. Note that $P''\cap U'=P'$. We   get a one-to-one map 
\begin{align}\label{eq:P'-P}
  P'\mapsto P=P'\times P''
\end{align}
from $\fc_0'$ onto a subset, denoted by $\fc_0$ of parabolic subgroups of  $U$. For any standard parabolic subgroup $P'$ of $U'$, resp. $P''$ of $U''$, resp. $P$ of $U$, we will denote by $\Delta_0^{P'}$, resp. $\Delta_0^{P''}$, resp. $\Delta_0^{P}$ the set of simple roots of $A_{M_0'}$, resp. $A_{T_0'''}$, resp. $A_{M_0'\times T_0''}$ in $P_0'\cap M_{P'}$, resp. $B_0''\cap M_{P''}$, resp. $(P_0'\times B_0'')\cap M_P$. We set $\ago_0'=\ago_{M_0'}$.   The inclusion $P_0'\subset P_0''$ gives an identification  $\ago_{P_0''}=\ago_{P_0'}=\ago_0'$. Using the  map $a_{B_0''}^*\to \ago_{P_0''}^*=\ago_{P_0'}^*$ dual to the inclusion $A_{P_0''}\subset A_{T_0''}$, we see that any $\la\in \ago_{B_0''}^*$  defined  a linear map $\ago_0'\to \RR$ still denoted by $\la$. Let $P'\in \fc_0'$ and $P''\subset U''$ be the associated parabolic subgroup. Observe that if $\la$ is a simple root  $\tilde\al\in \Delta_0^{U''}\setminus  \Delta_0^{P''}$, the  linear map $\ago_0'\to \RR$ we get is  equal, up to a positive constant, to a unique root  $\al\in \Delta_0^{U'}\setminus  \Delta_0^{P'}$.

Following §\ref{S:Siegel},  we have the notion of Siegel domains. They depend on auxiliary choices $\om_0',T_-'$ for $U'$ and  $\om_0'',T_-''$ for $U''$. The Siegel domains for $U=U'\times U''$ will be the product of Siegel domains for $U'$ and $U''$. We may and shall assume that $T_-'$ and $T_-''$ are chosen so that
\begin{align}
  \label{eq:inclusionAT}
A_{P_0'}^{P',\infty}(T_-')\subset A_{B_0''}^{P'',\infty}(T_-'').
\end{align}
\end{paragr}

\begin{paragr} We fix a point $T_0\in \ago_{0}'$  as in §\ref{S:partition} and we set 
$$d_{P'}(\la)=d_{P'}(\la,T_0)$$
for any $P'\in \fc'_0$ and $\la\in \ago_0^*$. The precise choice of $T_0$ is irrelevant in the sequel and we won't use this notation anymore.  We proceed in the same way to define $d_{P''}(\la)$.

\begin{proposition}\label{prop:equiv-unit}
  \begin{enumerate}
  \item Let $\la\in \ago_{B_0''}^*$. We have:
    \begin{align}\label{eq:dP-dP'}
      d_{P''}(\la,x)\sim d_{P'}(\la,x) \ \ x\in [U']_{P'}.
    \end{align}
  \item For $\tilde\al\in \Delta_0^{U''}\setminus  \Delta_0^{P''}$, there exists $r>0$ such that
    \begin{align*}
       d_{P''}(\tilde\al,x)\sim d_{P'}(\al,x)^r \ \ x\in [U']_{P'}.
    \end{align*}
where $\al \in \Delta_0^{U'}\setminus  \Delta_0^{P'}$ is deduced from $\tilde\al$ as above.
  \end{enumerate}
\end{proposition}

\begin{preuve}
  1. Both sides of \eqref{eq:dP-dP'} are invariant by left $P'(F)$-translation. Thus it suffices to prove the equivalence on $\sgo_{P'}= \om_0'  A_{P_0'}^{P',\infty}(T_-') K'$. Note that if $a\in  A_{P_0'}^{P',\infty}(T_-')$ then $a^{-1}\om_0'a $ is included in a fixed compact subset $\kc\subset U'(\AAA)$ that does not depend on $a$.  Thus by proposition \ref{prop:equiv} assertion 2, we have $d_{P''}(xa k)\sim d_{P''}(a)$  for $x\in \om_0'$, $a\in   A_{P_0'}^{P',\infty}(T_-')$ and $k\in K'$. By the inclusion \eqref{eq:inclusionAT} and proposition \ref{prop:equiv} assertion 1 we have  $d_{P''}(a)\sim\exp(\bg \la, H_{B_0''}(a)\bd)$ for $a\in   A_{P_0'}^{P',\infty}(T_-')$. Then we have $H_{B_0''}(a)=H_{P_0'}(a)=H_{P_0'}(xak)$ and $\exp(\bg \la, H_{B_0''}(a)\bd) \sim d_{P'}(\la,xak)$  by proposition \ref{prop:equiv} assertion 1  applied to $P'$.

2. is a consequence of 1 and the fact that the restriction of $\tilde\al$ to $\ago_0'$ is $r\al$ for some $r>0$.
\end{preuve}

\end{paragr}

\begin{paragr}[Sufficiently positive $T$.] ---   \label{S:suffic-pos}We will fix an euclidean norm $\lVert .\rVert$ on $\ago_0'$,  invariant by the Weyl group. For  $T\in \ago_{0}'$, let
$$\displaystyle d(T)=\inf_{\alpha\in \Delta_{0}^{U'}} \bg \alpha, T\bd .$$
We will fix $C>0$ and $\eps>0$ respectively   large and small enough constants. We shall throughout assume that $T$ is ``sufficiently positive'' that is we assume $T$ satisfies the inequality $d(T)\geqslant \max(\eps \lVert T\rVert,C)$. In particular, we shall assume that lemma \ref{lem:egalite-d} holds for any sufficiently positive $T$ and our functions $d_{P'}(\la)$ and $d_{P''}(\la)$.
\end{paragr}

\subsection{Truncated kernel}\label{ssec:truncated ker}

\begin{paragr}
  Let $f\in \Sc(U(\AAA))$. For any cuspidal datum $\chi\in \Xgo(U)$ and any  standard parabolic subgroup $P$ of $U$ we get kernels  $K_{f,P}$, resp.   $K_{f,P,\chi}$, see §\ref{S:cusp-datum}. If $P=U$ we shall omit the subscript $U$.
\end{paragr}

\begin{paragr}
  For any $T\in \ago_0'$, $x,y\in U'(\AAA)$ and $\chi\in \Xgo(U)$ we set:
\begin{align}
\label{eq:KchiTU}
 K_{f,\chi}^{T}(x,y) = \sum_{ P'\in \fc_0'} \epsilon_{P'} \sum_{\gamma \in P'(F)\bsl U'(F)} \sum_{\delta  \in P' (F) \bsl U'(F)} \htau_{P'}(H_{P'}(\delta y ) - T)K_{f,P, \chi}(\gamma x, \delta  y)
\end{align}
where  we set
\begin{align}
  \label{eq:epsU}
\epsilon_{P'}=(-1)^{\dim(\ago_{P'})}
\end{align}
 and $K_{f,P, \chi}$ is the kernel attached to the parabolic subgroup $P$ of $U$ image of $P'$ by \eqref{eq:P'-P}. Recall that we set $\htau_{P'}=\htau_{P'}^{U'}$.

\begin{remarque}
  \label{rq:KTchiU}
This is the kernel used in \cite{Z3} for compactly supported functions. Since we consider more general test functions we need a comment here. First the sum over $\delta$ may be taken in a \emph{finite} set which depends on $y$ (see \cite{ar1} Lemma 5.1). Second  by \cite[Lemma 2.10.1.1]{BCZ} for any $N>0$ there exist $N'>0$ and a continuous semi-norm $\|\cdot\|_{\Sc}$ on $\Sc(U(\AAA))$ such that we have
\begin{align}\label{eq:sumKPchi}
  \sum_{\chi\in \Xgo(U)} |K_{f,P,\chi}(x,y)| \leq \|f\|_{\Sc} \|x\|_P^{-N} \|y\|^{N'}_{P} 
\end{align}
for any standard parabolic subgroup $P$ of $U$ and all $x,y\in [U]_P$ and $f\in \Sc(U(\AAA))$. It follows that the sum over $\gamma$ is absolutely convergent.
\end{remarque}

\end{paragr}

\begin{paragr} The next theorem is an extension to Schwartz test functions of the work of Zydor in  \cite[section 4]{Z3}. The proof will be given in subsection \ref{ssec:proofmainthm} below.

\begin{theoreme}\label{thm:jfDefU} Let $T \in \ago_{0}'$ sufficiently positive.
  \begin{enumerate}
  \item We have
\[
\sum_{\chi \in \Xgo(G)} \int_{[U']\times [U']}| K^T_{f,\chi}(x,y) | \, dxdy < \infty
\]
\item As a function of $T$, the integral
\begin{align}\label{eq:polynexp}
J_\chi^{U,T}(f)=\int_{[U']\times [U']} K^T_{f,\chi}(x,y) \, dxdy 
\end{align}
coincides with  an  exponential-polynomial function in $T$ whose purely polynomial part is constant and denoted by $J^U_\chi(f)$.
\item The distributions $J_\chi^U$ are continuous, left and right $U'(\AAA)$-invariant.
\item The sum
  \begin{align}\label{eq:expansion-spec}
    J^U(f)= \sum_{\chi} J^U_{\chi}(f)
  \end{align}
is absolutely convergent and defines a continuous distribution $J^U$.
  \end{enumerate}
\end{theoreme}
\end{paragr}

\subsection{Truncation operator and distributions $J_\chi^U$}\label{ssec:Truncation operator}

\begin{paragr}
  The goal of this section is to state theorem \ref{thm:asym-LaTU} which gives the asymptotics of the distributions $J_\chi^{U,T}(f)$ defined in theorem \ref{thm:jfDefU} when the parameter goes to infinity.  The theorem will be useful  for subsequent computations in section \ref{ssec:UU'reg-contrib}.
\end{paragr}

\begin{paragr}[The Ichino-Yamana truncation operator.] --- \label{S:IY-truncU} Let $T\in \ago_{0}'$ sufficiently positive. In \cite{IYunit}, Ichino-Yamana defined a truncation operator which transforms functions of uniform moderate growth on $[U'']$ into rapidly decreasing  functions on $[U']$, see \cite[lemma 2.2]{IYunit}. By applying it to the \emph{right} component of $[U]=[U']\times [U'']$, we get a truncation operator which we denote by $\Lambda^T_u$. It associates to any function $\varphi$ on $[U]$ the function on $[U']$ defined by the following formula: for any  $x\in [U']$:
\begin{align}  \label{eq:LaTdU}
    (\Lambda^T_u\varphi)(x)= \sum_{ P' \in \fc_0' }  \epsilon_{P'} \sum_{\delta \in P'(F)\back U'(F) }   \hat\tau_{P'}(H_{P'}(\delta x)-T) \varphi_{U'\times P''}(  \delta  x)
\end{align}
where we follow  notations involved in \eqref{eq:KchiTU}. Moreover $\varphi_{U'\times P''}$ is the constant term of $\varphi$ along the parabolic subgroup $U'\times P''$ of $U$ where $P'\times P''$ is the image of $P'$ by \eqref{eq:P'-P}.
We shall recall and precise the main properties of $\Lambda_u^T$. 

\begin{remarque}
  To avoid confusions, we emphasize that in \eqref{eq:LaTdU} the map $\varphi_{U'\times P''}$ is evaluated at $ \delta  x\in U'(\AAA)$ where $U'$ is viewed  as a diagonal subgroup of $U$.
\end{remarque}
\end{paragr}

\begin{paragr}[Properties of $\La_u^T$.]

  \begin{proposition}\label{prop:LaTu-continu} 
    \begin{enumerate}
    \item The map $\varphi\mapsto \Lambda^T_u\varphi$ induces a linear continuous map from $\tc([U])$ to $\Sc([U'])$.
    \item For every $r>0$ there exists a continuous semi-norm $\|.\|$ on $\tc([U])$   such that
  \begin{align*}
   | (\Lambda^T_u\varphi )(x)- F^{U'}(x,T)\varphi(x) |  \leq e^{-r \|T\|} \|x\|_{U'}^{-N} \|\varphi\|
  \end{align*}
for all $x\in [U']$, $\varphi\in \tc([U])$ and $T\in \ago_0'$ sufficiently positive.
    \end{enumerate}
  \end{proposition}

  \begin{preuve} The assertion 1  can be easily extracted from \cite[(proof of) lemma 2.2]{IYunit}. For the convenience of the reader we give some details. First, it suffices to show that for any $N>0$ there exists a continuous semi-norm $\|.\|$ on $\tc([U])$ such that for all $\varphi\in \tc([U])$ we have:
    \begin{align*}
      \| \Lambda^T_u\varphi \|_{\infty,N} \leq   \|\varphi\|
    \end{align*}
where $\|\psi\|_{\infty,N}=\sup_{x\in [U']} \|x\|_{U'}^{N} |\psi(x)|$ for any function $\psi$ on $[U']$. We can write for $x\in [U']$:

\begin{align}\label{eq:dec-LaTu}
  (\Lambda^T_u\varphi)(x)=\sum_{P_1'\subset P_2' }   \sum_{\delta \in P_1'(F)\back U'(F)} F^{P_1'}(\delta x,T)  \sigma_1^2(H_0(\delta  x)-T) \varphi_{1,2}(\delta x)
\end{align}
where $\sigma_1^2=\sigma_{P_1'}^{P_2'}$ is the eponymous function with values in $\{0,1\}$ introduced by Arthur in \cite[section 6]{ar1} for the group $U'$ and
\begin{align}\label{eq:phi12}
   \varphi_{1,2}=\sum_{P_1'\subset P'\subset P_2'} \epsilon_{P'} \varphi_{U'\times P''}    
\end{align}
with $\varphi_{U'\times P''}  $ as above. Note that if $P_1'=P_2'$ then $\sigma_1^2=0$ unless $P_1'=P_2'=U'$. In this case, the corresponding term is $F^{U'}(x,T) \varphi(x)$ for which the result is obvious. The other cases are deduced from the next lemma which also gives assertion 2.
\end{preuve}

\begin{lemme}\label{lem:tranche-inst} Assume $P_1'\subsetneq P_2'$. For every $N>0$ and $r>0$ there exists $C>0$ and a continuous semi-norm $\|.\|$ on $\tc([U])$   such that
  \begin{align*}
    \sum_{\delta \in P_1'(F)\back U'(F)} F^{P_1'}(\delta x,T)  \sigma_1^2(H_0(\delta  x)-T) |\varphi_{1,2}(\delta x)| \leq C e^{-r \|T\|} \|x\|_{U'}^{-N} \|\varphi\|
  \end{align*}
for all $x\in [U']$, $\varphi\in \tc([U])$ and $T\in \ago_0'$ sufficiently positive.
\end{lemme}
\end{paragr}

\begin{paragr}[Proof of Lemma \ref{lem:tranche-inst}.] ---  It suffices to prove that for every $N>0$ and $r>0$ there exists $C>0$ and a continuous semi-norm $\|.\|$ on $\tc([U])$   such that  
 \begin{align*}
     F^{P_1'}(x,T)  \sigma_1^2(H_0(  x)-T) |\varphi_{1,2}(x)| \leq C e^{-r \|T\|} \|x\|_{P'_1}^{-N} \|\varphi\|
  \end{align*}
for all $x\in P_1'(F)\back U'(\AAA)$, $\varphi\in \tc([U])$ and $T\in \ago_0'$ sufficiently positive. It suffices to prove this majorization for $x$ in  the Siegel set $\sgo_{P'_1}$. Then  the condition  $F^{P_1'}( x,T)\sigma_1^2(H_0(  x)-T)=1$ implies that $x$ belongs also to the Siegel set $\sgo_{P'_2}$, see lemma \ref{lem:Ftau-Siegel}.  Then we can deduce from  \cite[proofs of lemme I.2.10 and corollaire I.2.11]{MWlivre} that for any $N>0$ and any $\la$ which is a linear combination  of elements of  $\Delta_0^{P''_{2}}\setminus \Delta_0^{P''_{1}}$ with positive coefficients  there exists a continuous semi-norm $\|\cdot \|_{N,\la}$ on $\tc([U])$ such that 
\begin{align*}
     F^{P_1'}(x,T)  \sigma_1^2(H_0(  x)-T) |\varphi_{1,2}(x)| \leq      \exp(-\bg \la, H_{P_1''}(x)\bd)         \|x\|_{P'_1}^{N} \|\varphi\|
  \end{align*}
for all $x\in \sgo_{P'_2}$, $\varphi\in \tc([U])$ and $T\in \ago_0'$ sufficiently positive. Taking into account that $H_{P_1'}(x)-H_{P_1''}(x)$ stays in a fixed compact subset for $x\in U'(\AAA)$ and the fact that the restriction to $\ago_{P_1'}$ of elements of $\Delta_0^{P''_{2}}\setminus \Delta_0^{P''_{1}}$ are (up to some positive contants) the restriction of elements of $\Delta_0^{P'_{2}}\setminus \Delta_0^{P'_{1}}$, we see that next lemma \ref{lem:Fsigma-d} enables us to conclude.

  \begin{lemme}
    \label{lem:Fsigma-d}
 For every $N\geq 0$ and $r>0$, there exists $t>0$ and $C>0$ such that for any 
 \begin{align}
   \label{eq:la=}
\la=\sum_{\al\in \Delta_0^{P'_{2}}\setminus \Delta_0^{P'_{1}}}x_\al \al       \text{   with  }  x_\al>t
 \end{align}
 we have 
\begin{align*}
  F^{P_1'}(x,T)  \sigma_1^2(H_0(x)-T) \exp(-\bg \la,H_{P_1'}(x)\bd ) \leq C e^{-r \|T\|} \|x\|_{P_1'}^{-N}
  \end{align*}
for all $x\in \sgo_{P_1'}$  and $T\in \ago_0'$ sufficiently positive.
  \end{lemme}

  \begin{preuve} Any $H_1\in\ago_{P_1'}$ can  be written $H_1=H_1^2+H_2$ according to the decomposition $\ago_{P_1'}= \ago_{P_1'}^{P_2'}\oplus \ago_{P_2'}$. By \cite[Corollary 6.2]{ar2}, there is $c_2>0$ such that 
\begin{align*}
  \|H_1\| \leq c_2 (1+\|H_1^2\|).
\end{align*}
for all $H_1\in \ago_{P_1'}$ such that $\sigma_1^2(H_1)\not=1$.  In particular, we see that there exists $\mu$ in the open convex cone generated\footnote{By this, we mean the interior of the convex cone  generated by $\Delta_0^{P'_{2}}\setminus \Delta_0^{P'_{1}}$.} by $\Delta_0^{P'_{2}}\setminus \Delta_0^{P'_{1}}$  and $m>0$ such that:
\begin{align}
  \label{eq:maj-normP1bis}
\|x\|_{P_1'}\leq \exp(\bg \mu,H_{P_1'}(x)\bd)  e^{m\|T\|}
\end{align}
for $x\in \sgo_{P_1'}$ such that
\begin{align}
  \label{eq:non-Fsigma-0}
F^{P_1'}( x,T)  \sigma_1^2(H_0(  x)-T) \not=0.
\end{align}

Thus we are reduced to prove the statement for $N=0$.

Let $\la$ be defined as in  \eqref{eq:la=} with $t>0$. The projection $\la^1$  of $\la$ on $\ago_0^{P_1',*}$ belongs to the (closed) convex cone generated by $-\hat\Delta_0^{P_1'}$. Since for all $x\in \sgo_{P_1'}$ that satisfy \eqref{eq:non-Fsigma-0}, we have on the one hand $\bg \varpi, H_0(x)-T\bd \leq 0$ for all $\varpi\in \hat\Delta_0^{P_1'}$ and  on the other hand
\begin{align*}
  \bg \la, H_{P_1'}(x)\bd =\bg \la-\la^1, H_{P_1'}(x)\bd > \bg \la -\la^1,T\bd .
\end{align*}
We get 
  \begin{align*}
  \bg \la, H_0(x)\bd &=\bg \la, H_{P_1'}(x) \bd+ \bg \la^1, H_{P_0}(x) \bd\\
&   > \bg \la-\la^1, T \bd+ \bg \la^1, T \bd=  \bg \la, T \bd \\
&\geq  t \eps|\Delta_{P_1'}^{P_2'}| \|T\|
\end{align*}
where the last lower bound holds for sufficiently positive $T$ (see § \ref{S:suffic-pos}). The conclusion is then clear.
  \end{preuve}
\end{paragr}

\begin{paragr} We can apply the operator $\Lambda_u^T$ to the right variable of the kernel $K_{f,\chi}(x,y)$: we get a function on $[U]\times [U']$ denoted by $K_{f,\chi}\Lambda^T_u$.

   \begin{theoreme}\label{thm:asym-LaTU} 
\begin{enumerate}
     \item  There exists a continuous semi-norm $\|.\|$ on $\Sc(U(\AAA))$ such that  for any $T\in \ago_0'$ sufficiently positive and  $f\in \Sc(U(\AAA))$ we have 
       \begin{align}
         \label{eq:LaTuKchi}
\sum_{\chi\in \Xgo(U) }\int_{[U']\times [U']} \left| (K_{f,\chi}\Lambda^T_u)(x,y)  \right|\, dxdy  \leq \|f\|
       \end{align}
\item For all $r>0$, there exists a continuous semi-norm $\|.\|$ on $\Sc(U(\AAA))$ such that  for any $T\in \ago_0'$ sufficiently positive and $f\in \Sc(U(\AAA))$ we have 
  \begin{align}
    \label{eq:asym-LaTu}
\sum_{\chi\in \Xgo(U) }   \left| J_\chi^{U,T}(f) - \int_{[U']\times [U']} (K_{f,\chi}\Lambda^T_u )(x,y) \, dxdy  \right|  \leq  e^{-r\|T\|}\|f\|.
  \end{align}
\end{enumerate}
\end{theoreme}

 From theorem \ref{thm:jfDefU} assertion 2 and from theorem \ref{thm:asym-LaTU}, we get: In particular, we see from assertion 2 and 

\begin{corollaire}\label{cor:LaTu-cstterm}
The absolutely convergent integral
\begin{align*}
  \int_{[U']\times [U']}( K_{f,\chi} \Lambda^T_u)(x,y) \, dxdy
\end{align*}
is  asymptotic to an  exponential-polynomial in the variable  $T$ whose purely polynomial term is constant and equal to $J_\chi^U(f)$.
\end{corollaire}

\end{paragr}

\subsection{Proof of main theorems}\label{ssec:proofmainthm}

\begin{paragr}[Proof of theorem \ref{thm:asym-LaTU} assertion 1.] --- This is a simple combination of the basic properties of the truncation operator $\La^T_u$ given in proposition \ref{prop:LaTu-continu} and the majorization \eqref{eq:sumKPchi} (for $P=U$)  of the kernel.  
\end{paragr}

\begin{paragr}[Asymptotics of several truncated kernels.] We shall use  the two following results.

  \begin{theoreme}\label{thm:conv-FRTF}
    For every $N_1,N_2, r >0$, there exists a continuous semi-norm $\|.\|$ on $\Sc(U(\AAA))$ such that
\begin{equation}\label{eq:asym KT}
\displaystyle \sum_{\chi\in \Xgo(G)} \left\lvert K^T_{\chi,f}(x,y)- K_{f,\chi}(x,y)F^{U'}(y,T)\right\rvert \leq e^{-r \| T\|} \| x\|^{-N_1}_{[U']} \| y\|^{-N_2} _{[U']}  \, \lVert f\rVert
\end{equation}
for $f\in \Sc(U(\AAA))$, $(x,y)\in [U']\times [U']$ and $T\in \ago_{0}$ sufficiently positive. 
\end{theoreme}

The proof of theorem \ref{thm:conv-FRTF} will be given in §\ref{S:preuve-conv-FRTF} below.

\begin{corollaire}
  \label{cor:conv-LaTu-RTF}
    For every $N_1,N_2, r >0$, there exists a continuous semi-norm $\|.\|$ on $\Sc(U(\AAA))$ such that
\begin{equation}\label{eq:asym KT2}
\displaystyle \sum_{\chi\in \Xgo(G)} \left\lvert K^T_{\chi,f}(x,y)- (K_{f,\chi}\La^T_u)(x,y)\right\rvert \leq e^{-r \| T\|} \| x\|^{-N_1}_{[U']} \| y\|^{-N_2} _{[U']}  \, \lVert f\rVert
\end{equation}
for $f\in \Sc(U(\AAA))$, $(x,y)\in [U']\times [U']$ and $T\in \ago_{0}$ sufficiently positive. 
\end{corollaire}

\begin{preuve}
  The corollary is a straightforward consequence of  theorem \ref{thm:conv-FRTF} and the   majorization of
  \begin{align*}
    (K_{f,\chi}\La^T_u)(x,y)-K_{f,\chi}(x,y)  F^{U'}(y,T)
  \end{align*}
  we get by using assertion 2 of proposition \ref{prop:LaTu-continu} and the majorization \eqref{eq:sumKPchi} for $P=U$ (recall that we have $\|x\|_U\sim \|x\|_{U'}$ for $x\in [U']$).

\end{preuve}
\end{paragr}

\begin{paragr}[Proof of theorem \ref{thm:jfDefU}.] --- First we mention that all the statements but the continuity are stated and proved for compactly supported functions, see \cite[theorems 4.1,  4.5 and  4.7]{Z3}. We just need the extension to Schwartz functions. The assertion 1 of theorem \ref{thm:jfDefU}  is a direct consequence of corollary \ref{cor:conv-LaTu-RTF} and assertion 1 of theorem \ref{thm:conv-FRTF}.
  
The assertion 2 of theorem \ref{thm:jfDefU} can be proved as in \cite[proof of theorems  4.5]{Z3}. We take for granted  the  obvious extension of corollary \ref{cor:conv-LaTu-RTF} and  theorem \ref{thm:asym-LaTU} assertion 1 to  the auxiliary kernels of \cite[section 4.1]{Z3} relative to some parabolic subgroups of $U$. Indeed it can be proved with the same technics and is left to the reader. In particular, the formula of  \cite[proposition 4.4]{Z3} holds for Schwartz functions and each term in right-hand side of the formula (among them $J_\chi^U(f)$) is a continuous distribution. From this we deduce assertion 3 of theorem \ref{thm:jfDefU} as in \cite[theorem 4.7]{Z3}.
\end{paragr}

\begin{paragr}[Proof of theorem \ref{thm:asym-LaTU} assertion 2.] ---  It is an obvious application of corollary \ref{cor:conv-LaTu-RTF}. 
\end{paragr}

\begin{paragr}[Proof of theorem \ref{thm:conv-FRTF}.] ---\label{S:preuve-conv-FRTF}
  We shall use the notation of the proof of proposition \ref{prop:LaTu-continu}. We start from the expression
  \begin{align*}
    K_{f,\chi}^{T}(x,y) -K_{f,\chi}(x,y)F^{U'}(y,T)\\
= \sum_{ P_1'\subsetneq  P_2'}      \sum_{\delta  \in P'_1 (F) \back U'(F)}  \sum_{\gamma\in P_2'(F)\back U'(F)}F^{P_1'}(\delta y,T) \sigma_{1}^{2}(H_{P_1'}(\delta y)-T)   K_{1,2, f, \chi}(\gamma x,\delta y)
  \end{align*}
where $\sigma_{1}^{2}=\sigma_{P_1'}^{P_2'}$ is as in \eqref{eq:dec-LaTu}, the sum is over $P_1',P_2'\in \fc_0'$ and we set
  \begin{align*}
    K_{1,2, f,\chi}(x,y)=\sum_{P_1'\subset P\subset P_2'}   \epsilon_{P'}  \sum_{\gamma \in P'(F)\back P_2'(F)}   K_{f,P, \chi}(\gamma x,y).
  \end{align*}

From now on, we fix $P_1'\subsetneq  P_2'$. Let $\al\in \Delta^{P_2'}_0\setminus  \Delta^{P_1'}_0$ and let $P_1'\subsetneq \, ^\al P_1'\subset P_2'$ be defined by $ \Delta^{^\al P_1'}_0=  \Delta^{P_1'}_0\cup\{\al\}$. For any $\, ^\al P_1'\subset P'\subset P_2'$ we denote by $P_\al'$ the parabolic subgroup $P_1'\subset P_\al'\subsetneq P'$ defined by  $\Delta_0^{P_\al'}=\Delta^{P'}_0\setminus  \{\al\}$.

We denote by $P$ and $P_\al$ the parabolic subgroups associated to $P'$ and $P'_\al$ (see \eqref{eq:P'-P}). For any $\, ^\al P_1'\subset P'\subset P_2'$, we set

\begin{align*}
  K^\al_{f,P,\chi}(x,y)= K_{f,P,\chi}(x,y)- \sum_{\gamma \in P_\al'(F)\back P'(F)} K_{f,P_\al,\chi}(\gamma x,y).
\end{align*}
Note that we have
 \begin{align}\label{eq:K12=K12al}
     K_{1,2, f,\chi}(x,y)=\sum_{\,^\al P_1'\subset P'\subset P_2'}   \epsilon_{P'}  \sum_{\delta \in P'(F)\back P_2'(F)}    K^\al_{f,P,\chi}(\delta x,y).
  \end{align}

We fix $P'$ as in the sum above and we start by majorizing each  $K^\al_{P,\chi,f}( x,y)$.

\begin{lemme}\label{lem:KPal}
  There  exists $n_0>0$ such that for any $n_1,n_2>0$ there is a continuous semi-norm  $\|\cdot\|$ on $\Sc(U(\AAA))$ such that
\begin{align*}
   \sum_{\chi\in \Xgo(U)}  |K^{\al}_{f,P,\chi,}(x,y)| \leq    \|f\| d_{P_1'}(\al,y)^{-n_1}\|y\|^{n_2+n_0}_{P'_1}  \|x\|_{P'}^{-n_2} .
\end{align*}
for $x\in P'(F)\back U'(\AAA)$ and $y \in P_1'(F)\back U'(\AAA)$ such that
\begin{align}\label{eq:Fsigma=1}
  F^{P_1'}(y,T)\sigma_1^2(H_{P_1'}(y)-T)=1.
\end{align}
\end{lemme}

\begin{preuve}
  Let $y \in  U'(\AAA)$. By using left translation by $P_1'(F)$, we may and shall assume that $y\in \sgo_{P'_1}$. Assume that $y$ satisfies also the condition \eqref{eq:Fsigma=1}. By lemma \ref{lem:Ftau-Siegel}, we have   $y\in \sgo_{P_\al'}$. By  lemma \ref{lem:egalite-d}, we also have  $d_{P_\al'}(\al,y)=d_{P_1'}(\al,y)$.  Recall that $\|y\|_{P'_1} \sim \|y\|$   for $y\in \sgo_{P'_1}$ and $d_{P_\al'}(\al,y)\sim \exp(\bg \al, H_{P_1'}(y)\bd)$  for $y\in \sgo_{P'_\al}$. Thus we may and shall freely replace $ d_{P_1'}(\al,y)$ by either $d_{P_\al'}(\al,y)$ or  $\exp(\bg \al, H_{P_1'}(y)\bd)$ and $\|y\|_{P'_1} $ by $\|y\|$ in the inequality we have to prove.  
  
  Let's observe that we have
  \begin{align*}
    \sum_{\gamma\in P_\al(F)\back P(F)}K_{f,P_\al,\chi}(\gamma x,y)=\int_{[N_\al]} K_{f,P,\chi}(x,n y)\, dn
  \end{align*}
where $N_\al=N_{P_\al}$. We have $P_\al=P_\al'\times P_\al''$ where $P_\al'\subset P'$ and $P_\al''\subset P''$ are maximal parabolic subgroups.  We set $N_\al''=N_{P_\al''}$. We denote by $\tilde\al$ the unique element of $\Delta_0^{P''}\setminus \Delta_0^{P_\al''}$.

Thus we see that $K^\al_{P,\chi,f}(x,y)$ is the sum of three terms:
\begin{align}
\label{eq:lem1}  K^{1,\al}_{f,P,\chi}(x,y)=K_{f,P,\chi}(x,y)-\int_{[N_\al'']} K_{f,P,\chi}(x,n y)\,dn\\
\label{eq:lem2}   K^{2,\al}_{f,P,\chi}(x,y)=\int_{[N_\al'']} K_{f,P,\chi}(x,n y)\,dn-\int_{[N_\al]} K_{f,P,\chi}(x,n y)\,dn\\
\label{eq:lem3}K^{3,\al}_{f,P,\chi}(x,y)=\sum_{\gamma\in \Om_{P,\al}} K_{f,P_\al, \chi}(\gamma x,y)
\end{align}
where $\Om_{P,\al}$ is the complement in $P_\al(F)\back P(F)$ of the diagonal image of $P'(F)$.

Let $J\subset U(\AAA_f)$ be a compact open subgroup and $S(U(\AAA))^J\subset S(U(\AAA))$ be the subspace of right-$J$-invariant functions. By a variant of  \cite[preuve du lemme I.2.10]{MWlivre}, for all $n_1,n_2>0$ there exists a finite family $(X_i)_{i\in I}$ of elements in the complexified universal envelopping algebra of the complexified Lie algebra of $U$ such that 
\begin{align*}
  \sum_{\chi\in \Xgo(U)}   |K^{1,\al}_{f,P,\chi}(x,y)|\leq \exp(- n_1\bg \al, H_{P_\al'}(y)\bd)\|y\|_{}^{n_2} \sup_{ z''\in [U'']_{P''} } \left(\|z''\|_{P''}^{-n_2} \sum_{i\in I} \sum_{\chi\in \Xgo(U)}   | K_{R(X_i)f,P,\chi}(x,(y,z''))|\right)
\end{align*}
for $x\in [U]_P$ and $y\in \sgo_{P'}$.

By \cite[Lemma 2.10.1.1]{BCZ}, there  exists $n_0>0$ such that for any $n_2>0$ there is a continuous semi-norm  $\|\cdot\|$ on $\Sc(U(\AAA))$ such that
\begin{align*}
   \sum_{\chi\in \Xgo(U)} |K_{f,P,\chi}(x,z)| \leq \|f\| \|x\|_P^{-n_2} \|z\|^{n_2+n_0}_{P} 
\end{align*}
for all $x,z\in [U]_P$. We conclude  that there exists $n_0$ such for all $n_1,n_2>0$  there is a semi-norm $\|\cdot\|$ on $\Sc(U(\AAA))$ such that 
\begin{align}\label{eq:K1al}
  \sum_{\chi\in \Xgo(U)}|K^{1,\al}_{f,P,\chi}(x,y)|\leq   \exp(- n_1\bg \al, H_{P_\al'}(y)\bd) \|y\|^{n_2+n_0} \|x\|_{P'}^{-n_2}
\end{align}
for all $x\in [U']_{P'}$ and $y\in \sgo_{P'}$. In the same way, one proves that the bounds  \eqref{eq:K1al} holds for $K^{2,\al}_{f,P,\chi}(x,y)$.

We introduce the weight function 
$$w(x)=\min( d_{P_\al''}(\tilde\al,x'),d_{P_\al''}(\tilde\al,x''))$$
for $x=(x',x'')\in [U]_{P_\al}=[U]_{P_\al'}\times [U]_{P_\al''}.$.
In the following we shall view $\al+\tilde\al$ as an element of $\ago_{P_0'}^*\oplus \ago_{B_0''}^*$. By \cite[Lemma 2.10.1.1]{BCZ}, there  exists $n_0>0$ such that for any $n_1,n_2>0$ there is a continuous semi-norm  $\|\cdot\|$ on $\Sc(U(\AAA))$ such that
\begin{align*}
   \sum_{\chi\in \Xgo(U)} |K_{f,P_\al,\chi}(x,z)| \leq    \|f\| w(z)^{-n_1}\|z\|^{n_2+n_0}_{P_\al}  w(x)^{n_1}\|x\|_{P_\al}^{-n_2} 
\end{align*}
for all $x,z\in [U]_{P_\al}$. Thus for all  $x\in P'(F)\back U'(\AAA)$ and $z\in [U']_{P_\al'}$  we have:
 \begin{align*}
   \sum_{\chi\in \Xgo(U)} \sum_{\gamma \in \Om_{P,\al}}|K_{f,P_\al,\chi}(\gamma x, z)| \leq    \|f\|  w(z)^{-n_1} \|z\|^{n_2+n_0}_{P_\al}  \sum_{\gamma \in \Om_{P,\al}}w(\gamma x)^{n_1}\|\gamma x\|_{P_\al}^{-n_2} .
\end{align*}

By proposition \ref{prop:unicity}, there is $c>0$ such that for all $\gamma=(\gamma',\gamma'')\in P(F)=P'(F)\times P''(F)$  and any $x\in U'(\AAA)$ (viewed as a diagonal element of $U(\AAA)$) such that $w(\gamma x)>c$ we have $\gamma''\in P_\al''(F)\gamma'$. Thus the class of $\gamma$ does not belong to the subset $\Om_{P,\al}$. In this way for $x\in U'(\AAA)$ we have:

$$\sum_{\gamma \in \Om_{P,\al}}w(\gamma x)^{n_1}\|\gamma x\|_{P_\al}^{-n_2} \leq c^{n_1} \sum_{\gamma \in P_\al(F)\back P(F)} \|\gamma x\|_{P_\al}^{-n_2} .
$$
By proposition \ref{prop:equiv-unit}, we have  $w(z)\sim d_{P_\al'}(\al,z)^r$ for some $r>0$ and all $z\in [U']_{P'_\al}$. It is then easy to see  that there  exists $n_0>0$ such that for any $n_1,n_2>0$ there is a continuous semi-norm  $\|\cdot\|$ on $\Sc(U(\AAA))$ such that
\begin{align*}
   \sum_{\chi\in \Xgo(U)}  |K^{3,\al}_{f,P,\chi}(x,y)| \leq    \|f\| d_{P_\al'}(\al,y)^{-n_1}\|y\|^{n_2+n_0}  \|x\|_{P'}^{-n_2} 
\end{align*}
for $y\in \sgo_{P_\al'}$ and $x\in P'(F)\back U'(\AAA)$. The conclusion is clear.
\end{preuve}

\begin{lemme}\label{lem:sumK12chi}
There  exists $n_0>0$ such that for any $n_2>0$ and any $\la$ in the open convex cone generated by $\Delta_0^{P_2'}\setminus \Delta_0^{P_1'}$ there is a continuous semi-norm  $\|\cdot\|$ on $\Sc(U(\AAA))$ such that
  \begin{align*}
     \sum_{\chi\in \Xgo(U)}|K_{1,2, f,\chi}(x,y)|\leq \|f\| d_{P_1'}(-\la,y)\|y\|^{n_2+n_0}_{P'_1}  \|x\|_{P'_2}^{-n_2} .
  \end{align*}
 for all $T$ sufficiently positive and all $x\in P_2'(F)\back U'(\AAA)$  and $ y\in P_1'(F)\back U'(\AAA)$ and $y$ such that $F^{P_1'}(y,T) \sigma_{1}^{2}(H_{P_1'}(y)-T) =1$.
\end{lemme}

\begin{preuve}
Let $P_1'\subset P'\subset P_2'$ and $\al\in \Delta_0^{P_2'}\setminus \Delta_0^{P_1'}$. 
By lemma \ref{lem:KPal} we have 
\begin{align*}
   \sum_{\chi\in \Xgo(U)}  |K^{\al}_{f,P,\chi}(x,y)| \leq    \|f\| d_{P_1'}(\al,y)^{-n_1}\|y\|^{n_2+n_0}_{P'_1}  \|x\|_{P'}^{-n_2} .
\end{align*}
 Moreover $\|y\|_{P'}\leq \|y\|_{P'_1}$. By \eqref{eq:K12=K12al}, we see that there  exists $n_0>0$ such that for any $n_1,n_2>0$ there is a continuous semi-norm  $\|\cdot\|$ on $\Sc(U(\AAA))$ such that

\begin{align*}
  \sum_{\chi\in \Xgo(U)}|K_{1,2, f,\chi}(x,y)|\leq \|f\| d_{P_1'}(\al,y)^{-n_1}\|y\|^{n_2+n_0}_{P'_1}  \|x\|_{P_2'}^{-n_2} 
  \end{align*}

The result follows easily.
\end{preuve}

Let $n_0$ be as in  lemma \ref{lem:sumK12chi}. For any $n_1,n_2,r>0$  there are $\la$ and $C>0$  as in lemma \ref{lem:Fsigma-d}  such that
\begin{align*}
  F^{P_1'}(y,T) \sigma_{1}^{2}(H_{P_1'}(y)-T) d_{P_1'}(\al,y)  \|y\|_{P_1'}^{n_0+n_2}  \leq C\exp(-r\|T\|) \|y\|^{-n_1}.
\end{align*}
Indeed, we recall that we may take $y\in \sgo_{P_1'}$  and  that $d_{P_1'}(\al,y)$ is equivalent to $\exp(\bg \al, H_{P_1'}(y)\bd)$ for such $y$ (see proposition \ref{prop:equiv}).

As a consequence we get by   lemma \ref{lem:sumK12chi} that there exists  a continuous semi-norm  $\|\cdot\|$ on $\Sc(U(\AAA))$ such that 
\begin{align*}
   F^{P_1'}(y,T) \sigma_{1}^{2}(H_{P_1'}(y)-T)   \sum_{\chi\in \Xgo(U)}|K_{1,2, f,\chi}(x,y)|\leq \|f\| e^{-r\|T\|} \|x\|_{P'_2}^{-n_2} \|y\|^{-n_2 }_{P'_1}  
  \end{align*}
  for all $T$ sufficiently positive, all $x\in P_2'(F)\back U'(\AAA)$  and $ y\in P_1'(F)\back U'(\AAA)$ and all $f\in\Sc(U(\AAA))$ .  It is then straightforward to get theorem \ref{thm:conv-FRTF}.

\end{paragr}

\subsection{The $(U,U')$-regular contribution in the Jacquet-Rallis trace formula}
\label{ssec:UU'reg-contrib}

\begin{paragr}
  The goal of the section  is to get theorem \ref{thm:JchiU} below which gives a computation of  the distributions $J_\chi^{U}$ of theorem \ref{thm:jfDefU} in terms of relative characters for some specific cuspidal data which we are going to define.
\end{paragr}

\begin{paragr} \label{S:UU'reg}  Recall that we have  $U=U'\times U''$. Let $\chi=(\chi',\chi'')\in\Xgo(U)=\Xgo(U')\times \Xgo(U'')$ be  a cuspidal datum.  Let $(M=M'\times M'',\pi=\pi'\boxtimes\pi'')$ be a representative where $M=M_P$ is standard Levi subgroup of  a standard parabolic subgroup $P$ of $U$.  For any integer $r$, we set $G_r=\Res_{E/F} \GL(r,E)$. We can find hermitian forms $h'$ and $h''$ respectively  of rank $n'$ and $n''$, integers $n_1',\ldots, n_{r'}'$ and $n_1'',\ldots, n_{r''}''$ and  for  $1\leq i\leq r'$ cuspidal representations $\pi_i'$ of $G_{n_i'}(\AAA)$ (with central character trivial on $A_{G_{n_i'}}^\infty $)  and   for  $1\leq i\leq r''$ cuspidal representations $\pi_i''$ of $G_{n_i''}(\AAA)$ (with central character trivial on $A_{G_{n_i''}}^\infty $), cuspidal representations  $\sigma'$ and $\sigma''$ respectively of $U(h')(\AAA)$ and  $U(h'')(\AAA)$  such that
such that
\begin{itemize}
\item $n'+2(n_1'+\ldots+n_{r'}')=n$ and  $n''+2(n_1''+\ldots+n_{r''}'')=n+1$;
\item $M'\simeq G_{n_1'}\times \ldots \times G_{n_{r'}'} \times U(h')$  and $M''\simeq G_{n_1''}\times \ldots \times G_{n''_{r''}} \times U(h'')$;
\item $\pi'=\pi_1'\boxtimes \ldots\boxtimes \pi_{r'}'\boxtimes \sigma'$ and $\pi''=\pi_1''\boxtimes \ldots\boxtimes \pi_{r''}''\boxtimes \sigma''$ accordingly.
\end{itemize}

 We shall say that $\chi$  is 
 \begin{itemize}
 \item \emph{$U$-regular} (or simply regular) if both $\chi'$ and $\chi''$ are regular. We say  that $\chi'$ is regular if the representations $\pi_1',\ldots,\pi_{r'}', (\pi_1')^{*},\ldots, (\pi_{r'}')^{*}$ 
   are two by two distinct (the same definition applies to $\chi''$). Here $\pi^*$ means the conjugate dual of the representation $\pi$.
 \item \emph{$U'$-regular} if for each $1\leq i\leq r'$ and $1\leq j\leq r''$  such that $n_i'=n_j''$ the representations $(\pi_i')^\vee$ (the contragredient of $\pi_i'$)  and $\pi_j''$ are neither isomorphic nor conjugate dual;
 \item \emph{$(U,U')$-regular} if it is both $U$-regular  and $U'$-regular.
 \end{itemize}

\end{paragr}

\begin{paragr} Let $\varphi\in \Ac_{P,\pi,\cusp}(U)$. Let  $\la\in \ago_{P,\CC}^*$ such that the Eisenstein series $E(\varphi,\la)$  on $[U]$ is regular at $\la$. For $T\in \ago_0'$, we denote by $\Lambda^T_u E(\varphi,\la)$ the function on $[U']$ obtained from $E(\varphi,\la)$ by truncation by  the operator defined in \eqref{eq:LaTdU}. The following proposition is basic to our calculation. 

  \begin{proposition} \label{prop:trunc-perU}Let $T\in \ago_0'$ sufficiently positive.
\begin{enumerate}
\item  The integral
\begin{align}\label{eq:trunc-perU}
     \int_{[U']} \Lambda^T_u E(x,\varphi,\la)\, dx
   \end{align}
is absolutely convergent.
\item If the spectral datum $\chi$ defined by $(M,\pi)$ is $U'$-regular then the integral \eqref{eq:trunc-perU} does not depend on $T$.
\item Write $P=P'\times P''$. If $P'=U'$ or $P''=U''$ then  we have:
  \begin{align}\label{eq:trunc=int-U}
     \int_{[U']} \Lambda^T_u E(x,\varphi,\la)\, dx=\int_{[U']}  E(x,\varphi,\la)\, dx
  \end{align}
  where the left-hand side is absolutely convergent.
\end{enumerate}
  \end{proposition}
  
  \begin{preuve}
    1.    The absolute convergence follows from the uniform moderate growth  of Eisenstein series and the basic properties of the truncation operator $\La^T_u$ recalled in proposition \ref{prop:LaTu-continu}.

    2. To analyse the dependence on $T$, we shall use the following formula (adapted from \cite[eq. (2.2)]{IYunit}): for any $T'\in \ago_0'$ and any smooth function $\varphi$ on $[U]$ we have for any $x\in U'(\AAA)$
    \begin{align}\label{eq:LaT+T'}
      (\Lambda^{T+T'}_u\varphi)(x)=\sum_{R'\in \fc_0'} \sum_{\delta \in R'(F)\back U'(F)} \Gamma_{R'}(H_{R'}(\delta x)-T_{R'},T') (\Lambda^{T,R'}_u\varphi)(\delta x)
    \end{align}
where the function $\Gamma_{R'}(\cdot,T')$ is compactly supported on $\ago_R$ (for the precise definition see p.9 of \cite{IYunit}); moreover we set  for any $x\in U'(\AAA)$:
\begin{align*}
  (\Lambda^{T,R'}_u\varphi)(x)=\sum_{Q'\subset R'} (-1)^{\dim(\ago_{Q'}^{R'})}\sum_{\delta \in Q'(F)\back R'(F)}\hat\tau_{Q'}^{R'}(H_{Q'}(\delta x)-T) \varphi_{U'\times Q''}(\delta x) 
\end{align*}
where the sum is over the set of standard parabolic subgroups $Q'$ of $R'$  and $\varphi_{U'\times Q''}$ is the constant term along $U'\times Q''$ such that  $Q'\times Q''$ is the image of $Q'$ by the map \eqref{eq:P'-P}. For $R'=U'$, we recover the operator $\Lambda^{T}_u$.

Using \eqref{eq:LaT+T'} and some Iwasawa decomposition, we see that assertion 2 follows from the following vanishing statement: for any \emph{proper} standard parabolic subgroups $R'$ of $U'$ we have for some Haar measure on $K$

\begin{align*}
  \int_{[M_{R'}]} \int_K   \exp(-\bg 2\rho_{R'},H_{R'}(x)\bd) \Gamma_{R'}(H_{R'}(x)-T,T')  \Lambda^{T,R'}_uE_{R}( x,\varphi,\la)\, dx=0
\end{align*}
where $E_{R}( \varphi,\la)$ is the constant term of $E(\varphi,\la)$ along the parabolic subgroup $R$, the image of $R'$ by the map \eqref{eq:P'-P}.

Using the usual computation of the constant of (cuspidal) Eisenstein series, we are reduced to prove that 
\begin{align}\label{eq:reducRunit}
  \int_{[M_{R'}]^1}   \Lambda^{T,M_R'}_u E^{R}( x,\varphi,\la)\, dx=0
\end{align}
for all parabolic subgroup $P\subset R$, all $\varphi\in \Ac_{P,\pi,\cusp}(U)$ such that the class of  $(M_P,\pi)$ is $U'$-regular. Here $E^{R}(\varphi,\la)$ denotes the Eisenstein series relative to $R$.

Let's prove this last claim. The reasonning here is very similar to that of \cite[proof of proposition 5.1.4.1]{BCZ} so we will be quite brief. There exist a hermitian form $h'$ of rank $m$ and  integers $n_1,\ldots,n_r$ such that $m+2(n_1+\ldots+n_r)=n$ 
\begin{align*}
  M_{R}\simeq  (G_{n_1}\times\ldots\times G_{n_r})\times  U(h') \times (G_{n_1}\times\ldots\times G_{n_r})\times U(h'\oplus h_0)\\
M_{R'}\simeq  (G_{n_1}\times\ldots\times G_{n_r})\times  U(h').
\end{align*}
The  embedding $M_{R'}\subset M_{R}$ is given by the product of the diagonal embeddings $G_{n_1}\times\ldots\times G_{n_r}\subset (G_{n_1}\times\ldots\times G_{n_r})^2$ and $ U(h')\subset  U(h')\times U(h'\oplus h_0)$. Then the operator  $ \Lambda^{T,M_R'}_u$ is the product of
\begin{itemize}
\item the usual Arthur operator attached to $G_{n_1}\times\ldots\times G_{n_r}$ viewed as an operator on functions on $[G_{n_1}\times\ldots\times G_{n_r}]\times [G_{n_1}\times\ldots\times G_{n_r}]$ acting on the second factor;
\item the operator  $ \Lambda^{T,U(h')}_u$ defined relatively to the embedding $ U(h')\subset  U(h')\times U(h'\oplus h_0)$.
\end{itemize}
We see that in \eqref{eq:reducRunit} each integral over $G_{n_i}$ can be interpreted as a scalar product of an Eisenstein series and a truncated Eisenstein series. By  Langlands' computation of this scalar product and our $U'$-regularity assumption, we get the expected vanishing.

3. If $P'=U'$ or $P''=U''$ then on the one hand the cuspidal datum defined by $(M,\pi)$ is automatically $U'$-regular. Thus the left-hand side of  \eqref{eq:trunc=int-U} does not depend on $T$. On the other hand, the restriction of $E(\varphi,\la)$ to $[U']$ is rapidly decreasing. So the right-hand side of \eqref{eq:trunc=int-U} is absolutely convergent. Now the equality \eqref{eq:trunc=int-U}, is a straightforward consequence of assertion 2 of proposition \ref{prop:LaTu-continu} and the dominated convergence theorem.
  \end{preuve}
\end{paragr}

\begin{paragr}[Ichino-Yamana regularized period.] --- 
  Assume that the cuspidal datum defined by $(M,\pi)$ is $U'$-regular. Following \cite{IYunit}, we define the regularized period of $E(\varphi,\la)$ for $\varphi\in \Ac_{P,\pi,\cusp}(U)$ and $\la\in \ago_{P,\CC}^*$ denoted by
  \begin{align*}
    \pc_{U'}(\varphi,\la)
  \end{align*}
to be the integral \eqref{eq:trunc-perU} attached to the Eisenstein series $E(\varphi,\la)$ for any $T\in \ago_0'$ sufficiently positive. It is meromorphic in $\la$ and holomorphic when  $E(\varphi,\la)$ is regular, in particular on $i\ago_P^*$. The map $\varphi\mapsto \pc_{U'}(\varphi,\la)$ is continuous.
\end{paragr}

\begin{paragr}[Relative character.] --- \label{S:rel-charU} We keep our assumption on $(M,\pi)$. Let $\la\in \ago_{P,\CC}^*$. We define the relative character $J_{P,\pi}^U(\la, f)$ by
  \begin{align*}
    J_{P,\pi}^{U}(\la, f)=\sum_{\varphi\in \bc_{P,\pi}}   \pc_{U'}(I_P(\la,f)\varphi,\la)\overline{   \pc_{U'}(\varphi,\la)  }
  \end{align*}
  where $f\in \Sc(U(\AAA))$ and $ \bc_{P,\pi}$ is  a $K$-basis of $\Ac_{P,\pi,\cusp}(U)$ in the sense of §\ref{S:automor}. Outside the singularities of the involved Eisenstein series, the sum is absolutely convergent. It is holomorphic on $i\ago_{P}^*$. It  does not depend on the choice of $\bc_{P,\pi}$  and it defines a continuous linear form on $\Sc(U(\AAA))$ (see \cite[proposition 2.8.4.1]{BCZ})). We have the following simple functional equation.

  \begin{proposition}
    \label{prop:JU-eqF} Let $P$ and $P_1$ be standard parabolic subgroup of $U$ of respective standard Levi factors $M$ and $M_1$.   Assume that the pairs $(M,\pi)$ and $(M_1,\pi_1)$ define the same $U'$-regular cuspidal datum. Then for $w\in W(M,M_1)$ such that $\pi_1=w\pi$ and  $\la\in i\ago_{P}^*$ we have:
      \begin{align*}
    J_{P,\pi}^{U}(\la, f)=   J_{P_1,\pi_1}^{U}(w\cdot \la, f)
      \end{align*}
    \end{proposition}

    \begin{preuve}
      This is an immediate consequence of the functional equation $E(\varphi,\la)=E(M(w,\la)\varphi,w\la)$ of Eisenstein series and the fact that for  $\la\in i\ago_{P}^*$  the intertwining operator sends a $K$-basis of $\Ac_{P,\pi,\cusp}(U)$ to a $K$-basis of $\Ac_{P,\pi_1,\cusp}(U)$.
    \end{preuve}
  \end{paragr}

  \begin{paragr}[Two auxiliary results.] --- For further use, we state and prove two useful results. 

    \begin{lemme}
   \label{lem:expansion-unit-Eis}   
Let $\chi\in \Xgo(U)$ be a $U$-regular cuspidal datum. For any standard parabolic subgroup $P$ of $U$ and any representative $(M_P,\pi)$ of $\chi$ we have for all $f\in \Sc(U(\AAA))$ and $x,y\in  [U]$ 
\begin{align}\label{eq::expansion-unit-Eis}
  K_{f,\chi}(x,y)=\int_{i\ago_P^{*}}  \sum_{\varphi \in \bc_{P,\pi}} E(x,I_P(\la,f)\varphi,\la)\overline{E(y,\varphi,\la)}\, d\la
\end{align}
where  $\bc_{P,\pi}$ is a $K$-basis as above. 
    \end{lemme}

    \begin{preuve} We start from  a representative  $(M_P,\pi)$ of $\chi$. Let $Q$ be a standard parabolic subgroup. Since $\chi$ is $U$-regular, the space  $\Ac_{Q}(U)\cap L^2_\chi([U]_{Q,0})$ (see §§\ref{S:automor},\ref{S:cusp-datum}) is equal to
\begin{align}\label{eq:AcQ}
   \hat \oplus_{w\in W(P,Q)}\Ac_{Q,w\pi,\cusp}(U).
\end{align}
In particular, this space is zero  unless $Q$ and $P$ are associated.      Using the functional equation of Eisenstein series, we see that the general expansion of the kernel $K_{f,\chi}$ (extended to the Schwartz space in \cite[lemma  2.10.2.1]{BCZ}) can be simplified and gives the expression \eqref{eq::expansion-unit-Eis}.
    \end{preuve}

    \begin{proposition}\label{prop:abs-cv-unit}
       There exists a semi-norm on $\Sc(U(\AAA))$ such that for all $f\in \Sc(U(\AAA))$
  \begin{align} \label{eq:sum-abs-val}
\sum_{(M,\pi) }  \int_{i\ago_{P}^*}   |J^U_{P,\pi}(\la,f)|\, d\la < \|f\|
  \end{align}
       where the sum is over a set of representatives $(M,\pi)$ of  $(U,U')$-regular cuspidal data of $U$ and $P$ is the  standard parabolic subgroup of which $M$ is the standard Levi factor.
     \end{proposition}
   
     \begin{preuve}
       By uniform boundedness principle, it suffices to prove that the expression \eqref{eq:sum-abs-val} is finite for each $f$. Let $(M,\pi)$ be as in the sum \eqref{eq:sum-abs-val} and $M=M_P$. By definition, see  §\ref{S:automor},  a $K$-basis $\bc_{P,\pi}$ is the union over the set $\hat K$  of finite families $\bc_{P,\pi,\tau}$. Let $T$ be a sufficiently positive point in  $\ago_0'$ such that for any $\varphi\in\Ac_{P,\pi,\cusp}(U)$
\begin{align}
  \label{eq:period=truncated}
\pc_{U'}(\varphi,\la)=\int_{[U']} \Lambda^T_u E(x,\varphi,\la)\, dx.
\end{align}
By proposition \ref{prop:LaTu-continu} and the fact that $\int_{[U']} \|x\|_{[U]}^{-N}dx<\infty$ for large enough $N$,    there exists $X,Y\in \uc(\ugo_\CC)$ and $N>0$ such that we have for all $\la\in i\ago_{P}^*$
\begin{align}\label{eq:bound-relative-char}
 | \sum_{\varphi\in \bc_{P,\pi,\tau}}  \pc_{U'}(I_P(\la,f)\varphi,\la)\overline{   \pc_{U'}(\varphi,\la)  }|\\\nonumber\leq \sup_{x,y\in [U]}  \left(  (\|x\|_{U}\|y\|_{U})^{-N}  | \sum_{\varphi \in \bc_{P,\chi,\tau}} (R(X)E)(x,I_P(\la,f)\varphi,\la)\overline{R(Y)E(y,\varphi,\la)}|\right).
\end{align}
However, by \cite[lemma  2.10.2.1]{BCZ}, there exists $C>0$ such that we have for all $x,y\in [U]$
\begin{align}\label{eq:bound-Eis}
  (\|x\|_{U}\|y\|_{U})^{-N} \sum_{(M,\pi) }  \int_{i\ago_{P}^*} \sum_{\tau\in \hat K}   | \sum_{\varphi \in \bc_{P,\pi,\tau}} (R(X)E)(x,I_P(\la,f)\varphi,\la)\overline{R(Y)E(y,\varphi,\la)}| < C
\end{align}
where the sum is as in \eqref{eq:sum-abs-val}. The basis used in  \cite[lemma  2.10.2.1]{BCZ} is slightly  different but since we are using regular cuspidal data, we can recover our formulation  (see the decomposition \eqref{eq:AcQ} in the proof of lemma \ref{lem:expansion-unit-Eis}). The conclusion is then clear.

     \end{preuve}

\end{paragr}

\begin{paragr}

We can now state and prove the main result of this section.

  \begin{theoreme}\label{thm:JchiU}
Let $\chi\in \Xgo(U)$ be a $(U,U')$-regular cuspidal datum. For any standard parabolic subgroup $P$ of $U$ and any representative $(M_P,\pi)$ of $\chi$ we have
\begin{align*}
  J^U_\chi(f)= \int_{i\ago_P^*}   J_{P,\pi}^{U}(\la, f)\,d\la
\end{align*}
where the integral in the right-hand side is absolutely convergent.
\end{theoreme}

\begin{preuve} Let $T_1,T_2$ be two sufficiently positive points in  $\ago_0'$ such that \eqref{eq:period=truncated} holds for $T=T_1,T_2$ . Let $\Lambda^{T_1}_uK_{f,\chi}\Lambda^{T_2}_u$ be the function of two variables we get when  we apply on $K_{f,\chi}$ the truncation operators $\Lambda^{T_1}_u$ and $\Lambda^{T_2}_u$ respectively on the left and right variables. Using the bound \eqref{eq:bound-Eis},  the properties of truncation operators given in proposition \ref{prop:LaTu-continu} and lemma \ref{lem:expansion-unit-Eis},  we check that we have for all $x,y\in [U']$
  \begin{align}\label{eq:LaT1T2-debut}
    (\Lambda^{T_1}_uK_{f,\chi}\Lambda^{T_2}_u) (x,y)=\int_{i\ago_P^{*}}  \sum_{\varphi \in \bc_{P,\pi}} (\Lambda^{T_1}_uE)(x,I_P(\la,f)\varphi,\la)\overline{\Lambda^{T_2}_uE(y,\varphi,\la)}\, d\la.
  \end{align}
  Using a variant of the majorization \eqref{eq:bound-relative-char} and the bound \eqref{eq:bound-Eis}, we see that
  \begin{align*}
      \int_{i\ago_{P}^*} \sum_{\tau\in \hat K}     \int_{[U']\times [U']} |\sum_{\varphi \in \bc_{P,\chi,\tau}}\Lambda^{T_1}_u E(x,I_P(\la,f) \varphi,\la)  \overline{   \Lambda^{T_1}_u E(y,\varphi,\la)}|\, dxdy<\infty.
  \end{align*}

From this, one checks that
\begin{align}\label{eq:LaT1T2}
\int_{[U']\times [U']   }  (\Lambda^{T_1}_uK_{f,\chi}\Lambda^{T_2}_u) (x,y)\, dxdy=\int_{i\ago_P^{*}}  J^U_{P,\pi}(\la, f)\,d\la
\end{align}
where both sides are absolutely convergent. Note that the right-hand side of \eqref{eq:LaT1T2} hence the left-hand side depend neither on $T_1$ nor on $T_2$.

By proposition \ref{prop:LaTu-continu} assertion 2, $\Lambda^{T_1}_uK_{f,\chi}\Lambda^{T_2}_u$ converges to $K_{f,\chi}\Lambda^{T_2}_u$ pointwise when $d(T_1)\to+\infty$ (see §\ref{S:suffic-pos}). By proposition \ref{prop:LaTu-continu} assertion 1, we see  that $(\Lambda^{T_1}_uK_{f,\chi}\Lambda^{T_2}_u) $ is bounded by  a finite sum of contributions $|K_{f_i,\chi}\Lambda^{T_2}_u|$ for certain functions $f_i\in \Sc(U(\AAA))$.  These contributions are absolutely integrable convergent by theorem \ref{thm:asym-LaTU} assertion 1. So we deduce from the  dominated convergence theorem that
\begin{align*}
\lim_{ d(T_1)\to+\infty}  \int_{[U']\times [U']   }  (\Lambda^{T_1}_uK_{f,\chi}\Lambda^{T_2}_u) (x,y)\, dxdy=\int_{[U']\times [U']   }  (K_{f,\chi}\Lambda^{T_2}_u) (x,y)\, dxdy.
\end{align*}
We saw that the left-hand side  depends  neither on $T_1$ nor on $T_2$. So the right-hand side is also  independent of  $T_2$. We can conclude from corollary \ref{cor:LaTu-cstterm} that the right-hand side is in fact equal to $J_\chi^U(f)$.
\end{preuve}
\end{paragr}

%% file: GHregular.tex
\section{The $(G,H)$-regular contribution in the Jacquet-Rallis trace formula}\label{sec:GH-contrib-JR}

\subsection{Statement and proof}
\label{ssec:statement-GH}

\begin{paragr}\label{S:Gn}
     Let  $E/F$ be a quadratic extension of number fields. Let $\AAA$ be the ring of adèles of $F$ and $\eta=\eta_{E/F}$ be the quadratic character of the group $\AAA^\times$ attached to $E/F$.  
     Let $n\geq 1$ and $G'_n=\GL_{n,F}$ be the algebraic group of $F$-linear automorphisms of $F^n$. We view as an $F$-subgroup of  $G_n=\Res_{E/F} (G_n' \times_FE)$.  We denote by $c$  the Galois involution of $G_n$. Let $\eta_{G'_n}$ be the character of $G'_n(\AAA)$ given by 
$$\eta_{G'_n}(h)=\eta(\det(h))^{n+1}$$
for all $h\in G'_n(\AAA)$. Let $(B_n',T_n')$ be a pair where $B_n'$ is the Borel subgroup $G_n'$ of  upper triangular matrices and $T_n'$ is the maximal torus of $ G' _n$ of diagonal matrices. Let $(B_n,T_n)$ be the pair deduced from $(B_n',T_n')$ by extension of scalars to $E$ and restriction to $F$: it is a pair of a minimal parabolic subgroup of $ G_n $ and its Levi factor. Let $ K_n \subset G_n (\AAA) $ and $ K '_n= K_n \cap G'_n (\AAA) \subset G '_n(\AAA) $ be the ``standard'' maximal compact subgroups.  We set:
\begin{align*}
  \ago_{n+1}^+=\ago_{B_{n+1}}^{G_{n+1}+}
\end{align*}
where the right-hand side is defined in § \ref{S:root}.

 We set  $G=G_n\times G_{n+1}$ and $G'=G_n'\times G_{n+1}'$ (see §\ref{S:Gn}).  Let $c$ be  the Galois involution of $G$ whose fixed points set if $G'$. The reductive groups  $G$ and $G'$ are equipped with the pairs $(B_n\times B_{n+1},T_n\times  T_{n+1})$ and $(B_n'\times B_{n+1}',T_n'\times  T_{n+1}')$ . Let $K=K_n\times K_{n+1}\subset G(\AAA)$ and $K'= K \cap G' (\AAA)$.  We denote by $\eta_{G'}$ the character $\eta_{G'_n}\boxtimes \eta_{G'_{n+1}}$ of $G'(\AAA)$. Let $H$ be the image of the diagonal embedding
  \begin{align*}
    G_n\hookrightarrow G_n\times G_{n+1}.
  \end{align*}

  Let $\pi$ be a cuspidal automorphic representation of $G(\AAA)$ with central character trivial on $A_G^\infty$.  As in §\ref{S:UU'reg}, we  denote by $\pi^*$ the conjugate-dual representation of $G(\AAA)$. We shall say that $\pi$ is self conjugate-dual if  $\pi\simeq\pi^*$ and that $\pi$ is $G'$-distinguished, resp. $(G',\eta)$-distinguished, if the linear form (called the Flicker-Rallis period)
  \begin{align}\label{eq:lesperiodes}
    \varphi \mapsto \int_{[G']_0} \varphi(h)\,dh, \text{ resp.}  \int_{[G']_0} \varphi(h)\eta(\det(h))\,dh
  \end{align}
  does not vanish identically on $\Ac_{\pi,\cusp}(G)$. Then $\pi$ is self conjugate-dual if and only if $\pi$ is either $G'$-distinguished or $(G',\eta)$-distinguished and it cannot be both (see  \cite{Flicker}). Note that since we are working with general linear groups we shall omit the subscript $\cusp$ and we shall write simply $\Ac_{\pi}(G)$ for the space $\Ac_{\pi,\cusp}(G)$. The same rule is applied to the space $\Ac_{P,\pi}(G)$ if $\pi$ is a  cuspidal automorphic representation of $M_P(\AAA)$. 
\end{paragr}

\begin{paragr}[The $(G,H)$-regular and Hermitian cuspidal data.] ---\label{S:GHgeneric}
  Let $\chi=(\chi_n,\chi_{n+1})\in \Xgo(G)=\Xgo(G_n) \times \Xgo(G_{n+1})$  be a cuspidal datum. Let $(M,\pi)$ be a representative in the class of $\chi$ with $M=M_P$ for a standard parabolic subgroup $P$ of $G$. We write $M=M_n\times M_{n+1}$ and $\pi=\pi_n\boxtimes \pi_{n+1}$ accordingly. Let  $j\in \{n,n+1\}$. We write $M_j=G_{n_{1,j}}\times\ldots\times G_{n_{r_j,j}}$  for some integer $r_j\geq 1$ and $\pi_j=\sigma_{1,j} \boxtimes\ldots \boxtimes\sigma_{r_j,j}$ accordingly. 

  We shall say that $\chi$ is

  \begin{itemize}
  \item   $G$-regular (or simply regular) if for all $j\in \{n,n+1\}$ and $1\leq i,i'\leq r_j$ such that  $n_{i,j}=n_{i',j}$  and $\sigma_{i,j}=\sigma_{i',j}$  we have  $i'=i$.  
  \item $H$-regular if  for $1\leq i \leq r_n$ and all $1\leq j \leq r_{n+1}$, if $n_{i,n}=n_{j,n+1}$ the representation $\sigma_{i,n}$ is not isomorphic to the contragredient of  $\sigma_{j,n+1}$;
  \item $(G,H)$-regular if it is both $G$-regular and $H$-regular;
  \item Hermitian if $\pi=\pi^*$ and if the representation  $\sigma_{i,j}$ is $\eta_{G'_{j}}$-distinguished for all $1\leq i\leq r_j$ and $j\in  \{n,n+1\}$ such that $\sigma_{i,j}=\sigma_{i,j}^*$.
  \end{itemize}
  
\end{paragr}

\begin{paragr}\label{S:choiceL}
  We fix a $(G,H)$-regular and Hermitian cuspidal datum $\chi$. With the notations as above, we may and shall choose the representative $(M,\pi)$ such that for all $j\in \{n,n+1\}$ there exists an integer $s_j\geq 0$ such that the following conditions are satisfied:
  \begin{enumerate}
  \item $2 s_j\leq r_j$  and for all odd $i$ such that $1\leq i <2s_j$ we have $\sigma_{i+1,j}=\sigma_{i,j}^*$ (in particular $\sigma_{i,j}\not=\sigma_{i,j}^*$);
  \item for all $i>2s_j$, we have $\sigma_{i,j}=\sigma_{i,j}^*$.
  \end{enumerate}
Let $L=L_n\times L_{n+1}$ be the standard Levi subgroup of $G$ such that for all $j\in \{n,n+1\}$  we have
\begin{align}\label{eq:leLevi}
L_j=G_{n_{1,j}+n_{2,j}}\times \ldots \times G_{n_{2s-1,j}+n_{2s,j}} \times G_{n_{2s+1,j}} \times \ldots \times G_{n_{r,j}}.
\end{align}
Let  $\xi\in W(M)$ such that $\xi^2=1$ and $\ago_M^{L}$ is the kernel of $\xi+\Id$ for the natural action of $W(M)$ on  $\ago_M$. We denote by $Q$ the standard parabolic subgroup of Levi $L$.
\end{paragr}

\begin{paragr}[Intertwining period.] --- \label{S:intertwining period} Let $\tilde{\xi}\in G(F)$ such that $\tilde\xi c(\tilde\xi)^{-1}=\xi$ where $c$ is the Galois involution of $G$. We define the $F$-subgroups $P_{\tilde \xi}$,   $M_{\tilde \xi}$ and $N_{\tilde \xi}$ of $G'$   respectively by  :
  \begin{align}
\label{eq:Pxi}    P_{\tilde \xi}=G'\cap \tilde\xi^{-1}P \tilde\xi.\\
\label{eq:Mxi}  M_{\tilde \xi}=G'\cap \tilde\xi^{-1}M \tilde\xi.\\
\label{eq:Nxi}N_{\tilde \xi}=G'\cap \tilde\xi^{-1}N \tilde\xi.
  \end{align}
   We have the Levi decomposition $P_{\tilde \xi}=M_{\tilde \xi} N_{\tilde \xi}$ where $M_{\tilde \xi}$ is reductive and  $N_{\tilde \xi}$ is  unipotent. The map $a\mapsto H_P({\tilde \xi} a {\tilde \xi}^{-1})$ identifies $A_{M_{\tilde \xi}}^\infty  $ with  the subspace $\ago_{M}^{L}$. In particular, we have $\bg \la, H_P({\tilde \xi} a {\tilde \xi}^{-1})\bd=0$ for any $a\in A_{M_{\tilde \xi}}^\infty  $  and $\la\in \ago_{M,\CC}^{L,*}$. We define the intertwining period for all $\varphi\in \Ac_{P,\pi}(G)$ and  $\la\in \ago_{M,\CC}^{L,*}$ by 

\begin{align}\label{eq:Jxila}
    J(\xi,\varphi,\la)=\int_{ A_{M_{\tilde \xi}}^\infty  M_{\tilde\xi}(F) N_{\tilde\xi} (\AAA)\back G'(\AAA)} \exp(\bg \la,H_P(\tilde\xi h)\bd) \varphi(\tilde\xi h)\eta_{G'}(h) \,dh.
\end{align}

  Here the group $N_{\tilde\xi} (\AAA) A_{M_{\tilde \xi}}^\infty  M_{\tilde\xi}(F)$ is equipped with  a  right-invariant Haar measure. However this measure is not left-invariant: the modular character is given by
  \begin{align*}
    \delta_{P,\tilde\xi}:x\mapsto \exp(\bg \rho_P,H_P(\tilde\xi x {\tilde \xi}^{-1})\bd)
  \end{align*}
  see \cite[VII p.221]{JLR}. The integral in \eqref{eq:Jxila} is understood as a right-$G'(\AAA)$-invariant linear form on the space of $(N_{\tilde\xi} (\AAA) A_{M_{\tilde \xi}}^\infty  M_{\tilde\xi}(F),\delta_{P,\tilde\xi})$-equivariant functions which contains
  \begin{align*}
    \exp(\bg \la,H_P(\tilde\xi \cdot)\bd) \varphi(\tilde\xi\cdot)
  \end{align*}
  for  $\la\in \ago_{M,\CC}^{L,*}$.  The integral \eqref{eq:Jxila} makes sense at least formally.  It is in fact absolutely convergent for $\la$ such that $\bg \la,\al^\vee\bd $ is large enough for any $\al^\vee\in \Delta_P^{Q,\vee}$ and it admits a  meromorphic continuation to $\ago_{M,\CC}^{L,*}$. It does not depend on a specific choice of $\tilde\xi$. For all these properties, we refer the reader to \cite[theorem 23 and (proof of) lemma 32 ]{JLR}.   Using this reference and also  \cite[lemma 8.1 case 2]{LapFRTF},  we deduce that  the analytic continuation of  $J(\xi,\varphi,\la)$ coincides with some integral of the mixed truncation of an Eisenstein series $E^Q(\varphi',\la)$ where $\varphi'$ is obtained by integral over the maximal compact subgroup $K'$  of $G'(\AAA)$. In particular,  $J(\xi,\varphi,\la)$ is holomorphic on  $i\ago_M^{L,*}$ and for $\la\in i\ago_M^{L,*}$ the map $\varphi\mapsto J(\xi,\varphi,\la)$ is continuous on $\Ac_{P,\pi}(G)$.
\end{paragr}

\begin{paragr}[Rankin-Selberg period.] --- Let $\varphi\in \Ac_{P,\pi}(G)$  and let  $\PP(E(\varphi,\la))$ be the regularized Rankin-Selberg period of the Eisenstein series $E(\varphi,\la)$ defined by Ichino-Yamana in \cite{IY}. Because we assume that $(M,\pi)$ is $H$-regular, this period is given by the truncated integral
\begin{align}\label{eq:RSperiod-Eis}
  \PP(E(\varphi,\la))=  \int_{[H]} \Lambda_r^T E(h,\varphi,\la)\,dh
  \end{align}
 for some parameter $T\in \ago_{n+1}^+$ and where  $\Lambda_r^T $ is the Ichino-Yamana truncation parameter (whose definition is recalled in \cite[eq. (3.3.2.1)]{BCZ}. One can show that the the right-hand side of \eqref{eq:RSperiod-Eis} does not depend on $T$ (the proof of this property is the same as the proof of  \cite[proposition 5.1.4.1]{BCZ}). In particular, $ \PP(E(\varphi,\la))$ inherits analytic properties of the Eisenstein series $E(\varphi,\la)$. Thus it is holomorphic on $i\ago_M^{*}$ and the map  $\varphi\mapsto \PP(E(\varphi,\la))$ is continuous.
\end{paragr}

\begin{paragr}[The relative character.] ---   \label{S:RelcharJPpi} For any $\psi$ and $\varphi\in \Ac_{P,\pi}(G)$ the expression 
$$ \PP(E(\psi,\la))\cdot J(\xi, \bar\varphi,-\la)$$
is holomorphic on $i\ago_M^{L,*}$ and  gives a continuous pairing on  $\Ac_{P,\pi}(G)$. By \cite[proposition 2.8.4.1]{BCZ}, for any $f\in \Sc(G(\AAA))$ and any $\la\in i\ago_M^{L,*}$ we can define the relative character   
\begin{align*}
    I_{P,\pi}(\la,f)=\sum_{\varphi\in \bc_{P,\pi}}  \PP(E(I_P(\la,f)\varphi,\la))\cdot J(\xi, \bar\varphi,-\la).
\end{align*}
where the sum is over a $K$-basis $\bc_{P,\pi}$.  We get a continuous map $f \mapsto (\la\mapsto  I_{P,\pi}(\la,f))$ from $\Sc(G(\AAA))$ to  the space of Schwartz functions on $i\ago_M^{L,*}$: this is a slight extension of \cite[proposition 4.1.10.1]{BCZ} which basically relies on bounds of Eisenstein series due to Lapid in \cite[proposition 6.1]{LapFRTF}. Note also that this is also an easy consequence of the inequality given in \cite[proof of proposition 7.2.3.2]{AspecExp}.
 \end{paragr}

 \begin{paragr}[The $(G,H)$-regular contribution to the Jacquet-Rallis trace formula.] --- The contribution to the Jacquet-Rallis trace formula of a cuspidal datum $\chi$ is a distribution denoted by $I_\chi$ and defined in  \cite[theorem 3.2.4.1]{BCZ}. In the case of a Hermitian $(G,H)$-regular cuspidal datum, the next theorem relates this contribution  to  the relative character defined above.

\begin{theoreme} \label{thm:Ichi} Let $\chi\in \Xgo(G)$ be a  $(G,H)$-regular cuspidal datum.
  \begin{enumerate}
  \item If $\chi$ is not Hermitian we have $I_\chi=0$.
  \item Assume that $\chi$ is moreover Hermitian and let $(M_P,\pi)$ and $L$ as in §\ref{S:choiceL}. For all $f\in \Sc(G(\AAA))$, we have
  \begin{align}\label{eq:Ichi}
    I_\chi(f)=2^{-\dim(\ago_L)} \int_{i\ago_M^{L,*}}I_{P,\pi}(\la,f)\,d\la
  \end{align}
  where the integral in the right-hand side is absolutely convergent.
\end{enumerate}
\end{theoreme}
\end{paragr}

\begin{paragr}[Proof of theorem \ref{thm:Ichi}.] ---  Let $\chi\in \Xgo(G)$ be a  $(G,H)$-regular cuspidal datum. Let $f\in \Sc(G(\AAA))$ and let $K_{\chi}$ be  the kernel of the right convolution by $f$ on  $L^2_\chi([G])$. By \cite[proposition 3.3.8.1 and theorem 3.3.9.1]{BCZ}, the contribution $I_\chi(f)$ is the constant term in the asymptotic expansion in $T\in \ago_{n+1}^+$ of the absolutely convergent integral
  \begin{align}\label{eq:intLaTK}
      \int_{[H]} \int_{[G'] }   \Lambda_r^TK_\chi (h,g)\,\eta_{G'}(g)dgdh.
    \end{align}

By \cite[lemma 5.2.2.1]{BCZ} we have  for all $h\in [H]$, 
\begin{align}\label{eq:step1}
  \int_{[G'] }   (\Lambda_r^TK_\chi) (h,g)\eta_{G'}(g)\,dhdg=   \Lambda_r^T\left(\int_{[G'] } K_\chi (\cdot,g)\eta_{G'}(g)\,dg\right)(h).
\end{align}

Then we have the following lemma:

\begin{lemme}
  \label{lem:expansion-generic} Let $\chi\in \Xgo(G)$ be a  $(G,H)$-regular cuspidal datum. For all $x\in G(\AAA)$, the absolutely convergent integral  
  \begin{align}\label{eq:intG' Kchi}
    \int_{[G'] } K_\chi (x,g)\eta_{G'}(g)\,dg
  \end{align}
  vanishes unless  $\chi$ is Hermitian. If $\chi$ is  Hermitian and represented by  $(M_P,\pi)$ we have:
  \begin{align}\label{eq:step2}
   \int_{[G']} K_\chi(x,g)\,\eta_{G'}(g)dg=2^{-\dim(\ago_{L})}   \int_{i\ago_{M}^{L,*}} (\sum_{\varphi\in \bc_{P,\pi}}  E(x,I_P(\la,f)\varphi,\la)\cdot J(\xi, \bar\varphi,-\la) \, d\la.
  \end{align}
  where the Levi subgroup $L$ is that defined  in §\ref{S:choiceL} and $\bc_{P,\pi}$ is  a $K$-basis.
\end{lemme}

\begin{preuve}   For brevity reasons, we shall get the lemma as a simple application of a much more difficult result namely \cite[theorem 7.2.4.1]{AspecExp}.  Note that the absolute convergence of \eqref{eq:intG' Kchi} results from \cite[Lemma 2.10.1.1]{BCZ}. Assume that  $(M_P,\pi)$ represents the cuspidal datum $\chi$.  Because $\chi$ is $G$-regular, for any standard parabolic sugroup $P_1$ of $G$, the intersection of  $\Ac_{P_1}(G)$ with  $L^2_\chi([G]_{P_1,0})$ is reduced to 
    \begin{align*}
  \bigoplus_{w\in W(P,P_1)} \Ac_{P_1,w\pi}(G).
  \end{align*}
  In particular this  space  is trivial unless $P$ and $P_1$ are associated. Using this description, we deduce from  (a slight variant of)  \cite[theorem 7.2.4.1]{AspecExp} that the left-hand side of \eqref{eq:step2} is equal to the absolutely convergent expression
   \begin{align}\label{eq:premiere exp}
  |\pc(M)|^{-1} \sum_{P_1}  \sum_{L_1\in \lc_2(M_1)}  2^{-\dim(\ago_{L_1})} \sum_{w\in W(P,P_1)} \int_{i\ago_{M_1}^{L_1,*}} \Ic_{L_1,w\pi}(x,f,\la)\, d\la
    \end{align}
    where the sum is over  the set of standard parabolic sugroups $P_1$ of $G$, we set $M_1=M_{P_1}$, the number of parabolic subgroups of Levi $M$ is denoted by $|\pc(M)|$, the set $\lc_2(M_1)$ is a subset of the set of Levi subgroups of $G$ containing $M_1$ (see \cite[§ 2.2.3]{AspecExp}) and $\Ic_{L_1,w\pi}(x,f,\la)$ is the relative character that essentially appears in the statement of  \cite[theorem 7.2.4.1]{AspecExp}. The only difference here is that we are working on a product of linear groups and the relative character is built upon Eisenstein series and the intertwining periods  defined in \cite[§ 5.1.4]{AspecExp} are twisted by the character $\eta_{G'}$. Note that the apparent discrepancy of a factor $2$ is due to a different choice of  the integration  domain (see proof of \cite[lemma 7.3.1.1]{AspecExp} for a discussion on this point).

    By the fact that $\chi$ is $G$-regular and by the basic properties of the cuspidal intertwining periods (see \cite[theorem 23]{JLR}) we see that the relative character  $\Ic_{L,w\pi}(x,f,\la)$ vanishes unless $\chi$ is Hermitian. So let's assume that $\chi$ is Hermitian. Let $M_1,L_1$ and $w\in  W(P,P_1)$ be as in \eqref{eq:premiere exp}. The intertwining operator $M(w,\la)$ induces  a unitary isomorphism from $\Ac_{P,\pi}$ onto $\Ac_{P_1,w\pi}$, which sends $K$-bases to $K$-bases. Using the functional equations of Eisenstein series on one hand and of the definition of  intertwining periods given in \cite[§ 5.1.4]{AspecExp} on the other hand, we get the functional equation of relative characters (see    \cite[proposition 8.2]{LapFRTF} for a close statement) :
    \begin{align*}
   \forall \la\in i\ago_{M_1}^{L_1,*}\ \    \Ic_{L_1,w\pi}(x,f,\la)= \Ic_{w^{-1}L_1,\pi}(x,f,w^{-1}\la)
    \end{align*}
    where   $w^{-1}L_1$ is the image of $L_1$ by the action of $w^{-1}$ on $\lc(M)$. Still because $\chi$ is $G$-regular,  $\Ic_{w^{-1}L_1,\pi}(x,f,w^{-1}\la)=0$ unless $ w^{-1}L_1=L$. At this point we get by an obvious change of variables that  \eqref{eq:premiere exp} is equal to
    \begin{align}\label{eq:seconde  exp}
      &2^{-\dim(\ago_{L})}|\pc(M)|^{-1} \sum_{P_1}   |W(P,P_1)| \int_{i\ago_{M}^{L,*}} \Ic_{L,\pi}(x,f,\la)\, d\la\\
      &=2^{-\dim(\ago_{L})} \int_{i\ago_{M}^{L,*}} \Ic_{L,\pi}(x,f,\la)\, d\la.
    \end{align}
    Since we have
    \begin{align*}
       \Ic_{L,\pi}(x,f,\la)=\sum_{\varphi\in \bc_{P,\pi}}  E(x,I_P(\la,f)\varphi,\la)\cdot J(\xi, \bar\varphi,-\la) 
    \end{align*}
the lemma is proved.    
\end{preuve}

It follows from lemma \ref{lem:expansion-generic} that,  for $h\in [H]$, we have
\begin{align}\label{eq:step3}
  \Lambda_r^T\left(\int_{[G'] } K_\chi (\cdot,g)\eta_{G'}(g)\,dg\right)(h)= 2^{-\dim(\ago_{L})}   \int_{i\ago_{M}^{L,*}} \sum_{\varphi\in \bc_{P,\pi}}  (\Lambda_r^TE(h,I_P(\la,f)\varphi,\la)\cdot J(\xi, \bar\varphi,-\la) \, d\la.
\end{align}
To get this, we have to permute  the truncation operator with the sum and the integral. However one observes that for a fixed $h\in [H]$ the truncation operator is a finite sum of constant terms. By Fubini theorem, it suffices to show that any parabolic subgroup $Q$ of $G$, the sum
\begin{align*}
  \sum_{\varphi\in \bc_{P,\pi}}  \int_{N_Q}| E(nh,I_P(\la,f)\varphi,\la)| \, dn \cdot |J(\xi, \bar\varphi,-\la) |
\end{align*}
is a Schwartz function in the variable $\la\in i\ago_M^{L,*}$. But this is an easy consequence of \cite[proposition 7.2.3.2]{AspecExp}.
The last step is to observe that
\begin{align}\label{eq:step4}
\int_{[H]} \int_{i\ago_{M}^{L,*}} \sum_{\varphi\in \bc_{P,\pi}}  \Lambda_r^TE(h,I_P(\la,f)\varphi,\la)\cdot J(\xi, \bar\varphi,-\la) \, d\la dh\\
\nonumber = \int_{i\ago_{M}^{L,*}} \sum_{\varphi\in \bc_{P,\pi}}  \int_{[H]} \Lambda_r^TE(h,I_P(\la,f)\varphi,\la)\, dh\cdot J(\xi, \bar\varphi,-\la) \, d\la
\end{align}
and to use \eqref{eq:RSperiod-Eis} to regognize the relative character $I_{P,\pi}(\la, f)$ in the inner sum. Once again as a consequence of \cite[proposition 7.2.3.2]{AspecExp} and also here the basic properties of truncation operator (see \cite[proposition 3.3.2.1]{BCZ}), we have
\begin{align*}
  \int_{i\ago_{M}^{L,*}} \sum_{\varphi\in \bc_{P,\pi}}  \int_{[H]}| \Lambda_r^TE(h,I_P(\la,f)\varphi,\la)|\, dh\cdot |J(\xi, \bar\varphi,-\la) \, |d\la <\infty.
\end{align*}
We can conclude with Fubini's theorem. 
\end{paragr}

\subsection{The relative character in terms of Whittaker functions}\label{ssec:char-Whittaker}

\begin{paragr}
We keep the notations of the previous subsection.  Let $N=N_n\times N_{n+1}$ and $N_H=N_n$ viewed as a diagonal subgroup of $N$. Let $N'=N\cap G'$.
\end{paragr}

\begin{paragr}
  We fix a nontrivial additive character $\psi':\AAA/F\to\CC^\times$. We deduce a character $\psi:\AAA_E/E\to \CC^\times$ trivial on $\AAA$ by $\psi(z)=\psi'(  \mathrm{Tr}_{E/F}(\tau z))$ where $\tau\in E^\times$ is such that $\tau^c=-\tau$. We define a regular character $\psi_n:[N_n]\to \CC^\times$ by
\begin{align*}
  \psi_n(u)=\psi((-1)^n \sum_{i=1}^{n-1}u_{i,i+1})
\end{align*}
for any $u\in [N_n]$. In the same way we get a character $\psi_{n+1}$ of  $[N_{n+1}]$. Thus we have  a character $\psi_{N}=\psi_n\boxtimes\psi_{n+1}$ of  $[N]$. By construction $\psi_N$ is trivial on the subgroups $N'$ and $N_H$.
\end{paragr}

\begin{paragr} Recall that we have fixed a pair $(M,\pi)$ with $M=M_P$ (see §§\ref{S:GHgeneric} and \ref{S:choiceL}).  Let $\la\in i\ago_{P}^{G,*}$.  Let $\pi_\la$ be the representation $\pi$ twisted by the character $m\mapsto \exp(\bg \la,H_M(m\bd)$ and let  $\Pi_\la=\Ind_{P(\AAA)}^{G(\AAA)}(\pi_\la)$ be the induced representation. The representation $\Pi_\la$ is irreducible, unitary and generic. Let $\wc(\Pi_\la,\psi_N)$ be its Whittaker model with respect to the character $\psi_N$.

For any $\varphi\in \Ac_{P,\pi}(G)$ and $g\in G(\AAA)$ let 
$$W(g,\varphi,\la)=\int_{[N]} E(ug,\varphi,\la) \psi_N^{-1}(u)\,du.
$$
We may and shall identify $g\mapsto W(g,\varphi,\la)$ with an element of  $\wc(\Pi_\la,\psi_N)$.
\end{paragr}

\begin{paragr} Let $W\in \wc(\Pi_\la,\psi_N)$.   By \cite{JPSS} and \cite{Jacarch}, the integral 
$$\displaystyle Z^{\RS}(s,W)=\int_{N_{H}(\AAA)\backslash H(\AAA)} W(h) \lvert \det h\rvert_{\AAA_E}^s dh$$
converges for $\Re(s)\gg 0$ and extends to a meromorphic function on $\CC$ which is holomorphic  at $s=0$. Let $\pc=\pc_n\times\pc_{n+1}$ be  the product of the respective mirabolic subgroups $\pc_n$ and $\pc_{n+1}$ of $G_n$ and $G_{n+1}$. Let $\pc'=\pc\cap G'$. We put $\As_G=\As^{(-1)^{n+1}}\boxtimes\As^{(-1)^n}$.  For any finite set $S\subset V_F$, we set
$$\displaystyle \beta_\eta(W)=(\Delta^{S,*}_{G'})^{-1}L^{S,*}(1,\Pi_\la,\As_G)\int_{N'(F_S)\backslash \pc'(F_S)} W(p_S) \eta_{G'}(p_S) dp_S $$
and
$$\displaystyle \langle W,W\rangle_{\Whitt}=(\Delta^{S,*}_{G})^{-1}L^{S,*}(1,\Pi_\la,\Ad)\int_{N(F_S)\backslash \pc(F_S)} \lvert W(p_S)\rvert^2 dp_S.$$
Here and hereafter, $L^*(1)$  means  the leading coefficient in the Laurent expansion of the meromorphic function $L(s)$ at $s=1$. The above expressions converge and are independent of $S$ as soon as it is chosen sufficiently large according to the level of $W$ (see \cite{Flicker} and \cite{JS}). 
\end{paragr}

\begin{paragr} 
  \begin{proposition}\label{prop:ZRS}
For any $\varphi\in \Ac_{P,\pi}(G)$, we have 
\begin{align*}
   Z^{\RS}(0,W(\varphi,\la))=\mathbf{P}(E(\varphi,\la)).
\end{align*}
  \end{proposition}

  \begin{preuve}
This is a straightforward application of results of Ichino-Yamana and the fact that $(M,\pi)$ is $H$-regular.
    First for any $T\in \ago_{n+1}^+ $ and any $s\in \CC$, the integral
    \begin{align}
      \label{eq:LaTs}
\int_{[H]} \Lambda_r^T E(h,\varphi,\la) |\det(h)|^s\,dh
    \end{align}
    converges and defines a holomorphic function in the variable $s$. Moreover since $(M,\pi)$ is $H$-regular, it does not depend on $T$ (see \cite[proof of proposition 5.1.4.1]{BCZ}). Because of this, \eqref{eq:LaTs}  is the regularized Rankin-Selberg period of $h\mapsto E(h,\varphi,\la) |\det(h)|^s$ defined in \cite{IY}. By \cite[theorem 1.1]{IY}, we deduce that \eqref{eq:LaTs} is equal to $Z^{\RS}(s, W(\varphi,\la))$. It suffices to take $s=0$ to get the result.
  \end{preuve}
\end{paragr}

\begin{paragr}
  \begin{proposition}
 \label{prop:Pet-Whit}
For any $\varphi\in \Ac_{P,\pi}(G)$, we have 
$$ \bg\phi,\phi\bd_{\Pet}=\bg W(\varphi,\la), W(\varphi,\la)\bd_{\Whitt}.$$
  \end{proposition}
  
\begin{preuve}
    This is  the proof of \cite[proposition 8.1.2.1]{BCZ}),  the main assumption there being that $(M,\pi)$ is regular.
  \end{preuve}
\end{paragr}

\begin{paragr}
  \begin{proposition}
    \label{prop:Flicker-Whitt}
For any $\varphi\in \Ac_{P,\pi}(G)$, we have 
\begin{align}\label{eq:Jbeta}
  J(\xi,\varphi,\la)=\beta_\eta(W(\varphi,\la)).
\end{align}
  \end{proposition}

  \begin{preuve}
    Let $Q$ be the standard parabolic subgroup of $G$ whose standard Levi factor is the group $L$ defined in  §\ref{S:choiceL}. Then $Q'=Q\cap G'$  is a standard parabolic subgroup of $G'$ of Levi factor $L'=L\cap G$. We shall follow the notations of §\ref{S:intertwining period}. For $\la$ in some cone in $\ago_{M,\CC}^{L,*}$, the intertwining period  $J(\xi,\varphi,\la)$ is given by the integral \eqref{eq:Jxila}. We may and will choose $\tilde\xi\in L(F)$. Observe then that we have $ A_{M_{\tilde\xi}}^\infty=A_L^\infty$, $M_{\tilde\xi}\subset L'$ and $ N_{\tilde\xi}= N_{Q'}$. Using the Iwasawa decomposition, one gets for a suitable measure on $K'$:
   
    \begin{align*}
      J(\xi,\varphi,\la)&=J^L(\xi,\varphi^{K'},\la)
      \end{align*}
where $$\varphi^{K'}(g)=\exp(-\bg \rho_Q,H_Q(g)\bd)\int_{K'} \varphi(gk)\eta_{G'}(gk)\, dk$$ and for $\phi\in \Ac_{L\cap P,\pi}(L)$ one defines the intertwining period $J^L(\xi,\phi,\mu)$ by analytic continuation of the integral
\begin{align*}
      \int_{ A_{M_{\tilde \xi}}^\infty  M_{\tilde\xi}(F) \back L'(\AAA)}\exp(\bg \mu,H_P(\tilde\xi h)\bd) \phi(\tilde\xi h )\eta_{G'}(h) \,dh
\end{align*}
which is convergent for $\mu$ in some cone in $\ago_{P,\CC}^{Q,*}$.

Let $\pc_L$ the mirabolic subgroup of $L$ (defined as the product of the mirabolic subgroups of its linear factors). Let $\pc_{L'}=\pc_L\cap L'$. Let $\As_{L_n}$ (resp. $\As_{L_{n+1}}$) be the tensor product of $\As^{(-1)^{n+1}}$ (resp. $\As^{(-1)^{n}}$) for each linear factor. Let 
$\As_L=\As_{L_n}\boxtimes\As_{L_{n+1}}$.  Let 
$$W_L(g,\varphi,\la)=\int_{[N\cap L]} E^Q(ug,\varphi,\la) \psi_N^{-1}(u)\,du
$$
where $E^Q$ denoted the Eisenstein series
\begin{align*}
  \sum_{\delta\in P(F)\back Q(F)}  \exp(\bg \la-\rho_Q, H_P(\delta g)\bd)\varphi(\delta g).
\end{align*}
Let $\Pi^Q_\la=\Ind_{P(\AAA)}^{Q(\AAA)}(\pi_\la)$ be the induced representation. Using \cite[proposition 3.2]{Zhang2}, \cite{Flicker} and theorem \ref{thm:J-W} in section \ref{sec:Inter} below, we can compute the intertwining period $J^L(\xi,\varphi,\la)$ in terms of $W_L$. We deduce that  for $S$ large enough $J(\xi,\varphi,\la)$ is equal to 
\begin{align*}
 &(\Delta_{L'}^{S,*})^{-1}  L^{S,*}(1,\Pi^Q_\la,\As_{L}) \int_{K'} \int_{(N'\cap L)(F_S)\back \pc'_L(F_S)}  \delta_Q(q_S k)^{-1}{W}_L(q_S k,\varphi,\la)  \eta_{G'}(q_S k)    \, dq_Sdk\\
&= (\Delta_{G'}^{S,*})^{-1}L^{S,*}(1,\Pi^Q_\la,\As_{L}) \int_{Q'(F_S)\back G'(F_S)} \int_{(N'\cap L)(F_S)\back \pc'_L(F_S)}\delta_Q(q_S g_S)^{-1}  {W}_L(q_S g_S,\varphi,\la) \eta_{G'}(q_S g_S)   \, dq_Sdg_S.
    \end{align*}
    Now we can appeal to \cite[theorem 9.1.7.1]{BCZ} to get that the last line is equal to
    \begin{align*}
  (\Delta_{G'}^{S,*})^{-1}L^{S,*}(1,\Pi^Q_\la,\As_{L}) \int_{N'(F_S)\back\pc'(F_S)}\mathbb{W}_S(p_S,\la,\varphi)\, dp_S,
\end{align*}
where $\mathbb{W}_S(g_S,\la,\varphi)$ stands for  the Jacquet integral given by the value at $s=1$ of the holomorphic continuation of the integral:
\begin{align*}
  \int_{  (w_L^{-1}L w_L\cap N)(F_S) \back N(F_S)  } \delta_Q(w_L u g_S)^{-s}  {W}_L(w_Lu g_S ,\varphi,\la) \eta_{G'}(g_S) \psi_N(u)^{-1}\,du , \Re(s)\gg 1.
\end{align*}
To conclude it suffices to observe that for any $g_S\in G(\AAA_S)$
\begin{align*}
L^{S,*}(1,\Pi^Q_\la,\As_{L})  \mathbb{W}_S(g_S,\la,\varphi)= L^{S,*}(1,\Pi_\la,\As) W(g_S,\la,\varphi).
\end{align*}
Indeed this follows from computations of \cite[section 4]{ShaLfunction} and the fact that the functions $L^{S}(s,\Pi^Q_\la,\As_{L})$ and  $L^{S}(s,\Pi_\la,\As) $ have a pole of same order at $s=1$  since $(M,\pi)$ is regular.

\end{preuve}
\end{paragr}

\begin{paragr} Let $\bc_{P,\pi}$ be a $K$-basis of $\Ac_{P,\pi}(G)$. For any $f\in \Sc(G(\AAA))$, we define
  \begin{align}\label{eq:IPIla}
    I_{\Pi_\la}(f)=\sum_{\varphi\in \bc_{P,\pi}} \frac{Z^{RS}(0, W(I_P(\la,f)\varphi,\la))  \overline{\beta_\eta(W(\varphi,\la))}       }{ \bg W(\varphi,\la),W(\varphi,\la)\bd_{\Pet}}
  \end{align}

  \begin{theoreme}\label{thm:IPi=Ipi}
    \begin{enumerate}
    \item The series \eqref{eq:IPIla} converges, does not depend on the choice of $\bc_{P,\pi}$ and defines a continuous distribution on $\Sc(G(\AAA))$
    \item We have   $I_{\Pi_\la}(f)=I_{P,\pi}(\la,f)$.
    \end{enumerate}
\end{theoreme}

\begin{preuve}
The two assertions follow from the fact that  the relative characters $I_{\Pi_\la}$ and $I_{\pi,\la}$ can be identified term by term by the combination of propositions \ref{prop:ZRS}, \ref{prop:Pet-Whit} and \ref{prop:Flicker-Whitt}.
\end{preuve}
\end{paragr}

%% file: Inter-Whittaker.tex
\section{Intertwining periods and Whittaker functions}\label{sec:Inter}

\subsection{Notations}

\begin{paragr}\label{S:notations-W} We follow the notations of §\ref{S:Gn}. However, in this section, we set $G=G_{2n}$ and $G'=G'_{2n}$. Let $P=MN$ be the  maximal standard parabolic subgroup of $G$ where its  standard Levi factor $M$ is $G_n\times G_n$. Let $\sigma$ be an irreducible cuspidal automorphic representation of $G_n$ with central character trivial on $A_{G_n}^\infty$. Let $\pi=\sigma\boxtimes \sigma^*$: this a cuspidal representation of $M$. 
\end{paragr}

\begin{paragr}  We set:
  \begin{align*}
    _PW_P=\{w\in W \mid M\cap w^{-1}B_{2n}w=M \cap B_{2n}=M\cap wB_{2n}w^{-1}\}
  \end{align*}
and the subset of involutions:
  \begin{align*}
    _PW_{P,2}=\{w\in \, _PW_{P}\mid w^2=1\}.
  \end{align*}

We have the following lemma:

\begin{lemme}(Jacquet-Lapid-Rogawski, see \cite[proposition 20]{JLR}). \label{lem:PGG'}
Any double coset in $P(F)\back G(F)/ G'(F)$ has a representative $\tilde\xi$ such that $\xi=\tilde\xi c(\tilde\xi)^{-1}$ belongs to $\,    _PW_{P,2}$. The map $P(F)\tilde\xi G'(F) \mapsto \xi$ is well-defined and induces a bijection from $P(F)\back G(F)/ G'(F)$ onto $  _PW_{P,2}$.
\end{lemme}

For any $\tilde\xi\in G(F)$ such that $\tilde\xi c(\tilde\xi)^{-1}$ belongs to $\,    _PW_{P,2}$, we set $ P_{\tilde \xi}=G'\cap \tilde\xi^{-1}P \tilde\xi $ and $M_{\tilde \xi}=G'\cap \tilde\xi^{-1}M \tilde\xi$. Note that  $M_{\tilde \xi}$ is a Levi factor of   $P_{\tilde \xi}$.
\end{paragr}

\begin{paragr}
    We fix $\tau\in E$ such that $c(\tau)=-\tau$. Let 
\begin{align*}
  \tilde\xi_0=
  \begin{pmatrix}
    I_n & \tau I_n\\   I_n & -\tau I_n
  \end{pmatrix}
      \text{ and }                  \xi_0=\tilde\xi_0 c(\tilde\xi_0)^{-1}.       
\end{align*}
We have  $\xi_0\in  \,    _PW_{P,2}$ and $  P_{\tilde \xi_0}=M_{\tilde \xi_0}$. 

Let $\varphi\in \Ac_{P,\pi}(G)$. For any $\la\in \ago_{P,\CC}^{G,*}$ we consider the intertwining period (due to Jacquet-Lapid-Rogawski, see \cite{JLR})
 \begin{align}\label{eq:periode entrelacement}
    J(\xi_0,\varphi,\la)=\int_{ A_{M_{\tilde \xi_0}}^\infty  M_{\tilde\xi_0}(F) \back G'(\AAA)} \exp(\bg \la,H_P(\tilde\xi h)\bd) \varphi(\tilde\xi h)\,dh
  \end{align}
Note that we have $A_{M_{\tilde \xi_0}}^\infty=A_{G'}^\infty  $. Let $\al$ be the unique root in $\Delta_P$. The integral above is absolutely convergent if $\Re(\bg \al,\la\bd)\gg 0$. Moreover it admits a meromorphic continuation to $\ago_{P,\CC}^{G,*}$, see \cite[theorem 23]{JLR}. 
\end{paragr}

\subsection{Epstein series and intertwining periods}

\begin{paragr} \label{S:Epstein}The group $G$ acts on the right on the space  of rows of size $2n$ (identified to $E^{2n}$). Let $\pc$ be the stabilizer of $e_{2n}=(0,\ldots,0,1)$. Let $\pc'=\pc\cap G'$.

  Let $\Phi\in \Sc(\AAA^{2n})$. For any $s\in \CC$ such that $\Re(s)>1$ and any $h\in [G']$,  the Epstein series is defined by the following absolutely convergent integral
  \begin{align*}
E(h,\Phi,s)=\int_{A_{G'}^\infty} \sum_{\gamma\in \pc'(F)\back  G'(F)} \Phi(e_{2n}\gamma ah) |\det(ah)|^s da.
  \end{align*}
Here $|\cdot|$ is the product over all places $v$ of $F$ of normalized absolute values of the completions $F_v$. The map $s\mapsto E(\Phi,s)$ extends to a meromorphic function valued in $\tc([G'])$ with simple poles at $s=0$, $1$ of respective residues $\Phi(0)$ and $\widehat{\Phi}(0)$ (cf. \cite[Lemma 4.2]{JS}). 
\end{paragr}

\begin{paragr}
  Let $\tilde\xi\in G(F)$ such that $\xi=\tilde\xi c(\tilde\xi)^{-1}$  belongs to $\,    _PW_{P,2}$. We assume also $\xi\not=1$. Let $\varphi\in \Ac_{P,\pi}(G)$. We define

\begin{align}\label{eq:JPhis}
     J(\xi,\varphi,\la,\Phi,s)=\int_{ A_{G'}^\infty  P_{\tilde\xi}(F) \back G'(\AAA)} \exp(\bg \la,H_P(\tilde\xi h)\bd) \varphi(\tilde\xi h) E(h,\Phi,s)\,dh.
\end{align}

If it is well-defined, this integral does not depend  on $\tilde\xi$  provided that we have $\xi=\tilde\xi c(\tilde\xi)^{-1}$, hence the notation.

  \begin{proposition}\label{prop:JlaPhi}
Assume $\xi\not=1$. 
\begin{enumerate}
\item  There exists $r>0$ such that for each for any $\la$ in the domain 
$$\dc_r=\{\la\in \ago_{P,\CC}^{G,*} \mid \Re(\bg \la,\al^\vee\bd)>r\}$$
and any $s\in \CC\setminus \{0,1\}$  the integral \eqref{eq:JPhis} converges absolutely. It also converges uniformly for  $\Re(\la)$ in a compact subset of $\dc_r$ and $s$ in a compact subset of  $\CC\setminus \{0,1\}$.
\item  The map  $(\la,s)\mapsto J(\xi,\varphi,\la,\Phi,s)$ is holomorphic on $ \dc_r \times\CC\setminus \{0,1\}$.
\item If $\xi\not=\xi_0$, the map has a holomorphic extension to $\dc_r \times\CC$.
\item If $\xi=\xi_0$, the map has simple poles at $s=0,1$ with respective residues:
\begin{align*}
  \Phi(0)  J(\xi_0,\varphi,\la) \text{  and  }  \hat\Phi(0)  J(\xi_0,\varphi,\la).
\end{align*}
\end{enumerate}
  \end{proposition}

  \begin{preuve}
We fix a compact subset $\Om $ of $\CC$.  We can write
\begin{align*}
     J(\xi,\varphi,\la,\Phi,s)=\int_{ [G']_0}   \left(\sum_{ \delta \in P_{\tilde\xi}(F) \back G'(F) }  \exp(\bg \la,H_P(\tilde\xi \delta h)\bd) \varphi(\tilde\xi \delta h) \right) E(h,\Phi,s)\,dh.
\end{align*}
We shall use the following two facts:
\begin{itemize}
\item there exist $C,t>0$  such that for all  $h\in G'(\AAA)^1$   and $s\in \Om$
  \begin{align}\label{eq:maj-epstein}
     |s(s-1)E(h,\Phi,s)| \leq C   \|h\|_{[G']}^t.
  \end{align}
\item there exists $r>0$ such that for any $t>0$ and any $\la\in  \dc_r$ there exists $C>0$ such that  for all  $h\in G'(\AAA)^1$ 
  \begin{align}\label{eq:maj-eisenstein-lacunaire}
     \sum_{ \delta \in P_{\tilde\xi}(F) \back G'(F) }  |\exp(\bg \la,H_P(\tilde\xi \delta h)\bd) \varphi(\tilde\xi \delta h)| \leq C  \|h\|_{[G']}^{-t}.
  \end{align}
\end{itemize}
Moreover in  \eqref{eq:maj-eisenstein-lacunaire} the constant $C$ may be chosen uniformly for $\la$ such that $\Re(\la)$ belongs to a compact subset of $\dc_r$.  This  can be extracted from the proof of \cite[theorem 23]{JLR}, more precisely from the  combination of proposition 24 and the lines following (58) of \cite{JLR}.

We get the holomorphic continuation with  at most simple  poles at $s=0,1$. Up to a factor $\Phi(0)  $ or $\hat\Phi(0)  $, the residue is given by:
\begin{align*}
  \int_{ A_{G'}^\infty  P_{\tilde\xi}(F) \back G'(\AAA)} \exp(\bg \la,H_P(\tilde\xi h)\bd) \varphi(\tilde\xi h) \,dh.
\end{align*}
This integral converges absolutely thanks to the majorization \eqref{eq:maj-eisenstein-lacunaire}. According to the proof of \cite[theorem 23]{JLR}, the integral vanishes unless $\xi=\xi_0$. In this case, the integral is nothing else  $J(\xi_0,\varphi,\la)$.
  \end{preuve}
\end{paragr}

\subsection{Period of a pseudo-Eisenstein series: first computation}

\begin{paragr}
 Let $\om$ be the central character of $\pi$. By restriction it induces  a unitary character of $Z_G'(\AAA)$. Let 
\begin{align*}
  \tilde\Phi(h,\om, s)=\int_{Z_{G'}(\AAA)} \Phi(e_{2n} zh) \om(z) |\det(zh)|^s dz.
\end{align*}
The integral is convergent for $\Re(s)>\frac1n$. Let $P'_1=\pc' Z_{G'}$. This is a parabolic subgroup of $G'$ of type $(2n-1,1)$.  Then we have
\begin{align*}
  \int_{A_{G'}^\infty Z_{G'}(F) \back Z_{G'}(\AAA)} E(zh,\Phi,s)\om(z)\, dz=E(h,\tilde\Phi(\om, s))
\end{align*}
where
\begin{align*}
  E(h,\tilde\Phi(\om, s))= \sum_{\gamma\in P_1'(F)\back  G'(F)} \tilde\Phi(\gamma h,\om, s)
\end{align*}
is a usual Eisenstein series, which is    convergent for  $\Re(s)>1$. By the classical computation of the constant term of an Eisenstein series, there exist  a finite set $I$ and  families $(\varphi_{i,s})_{i\in I}\in \Ac_{P'}(G')^I$ and $(\mu_{i,s})_{i\in I} \in  \ago_{P,\CC}^{G,*}$ for any $s\in \CC$ such that  the following conditions are satisfied for any $i\in I$:
\begin{itemize}
\item the map $s\mapsto \mu_{i,s}$ is affine;
\item the map $s\mapsto \varphi_{i,s}$ is a meromorphic function;
\item we have $\varphi_{i,s}(ah)=\varphi_{i,s}(h)$ for any $h\in G'(\AAA)$ and $a\in A_{M'}^{\infty}$;
\item we have for any $a\in  A_{M'}^{\infty}$, $m\in M'(\AAA)^1$ and $k\in K'$
  \begin{align}\label{eq:cst-term}
 \int_{[N_{P'}]}  E(namk ,\tilde\Phi(\om, s))dn= \sum_{i\in I} \varphi_{i,s}(mk) \exp(\bg \mu_{i,s},H_{P}(a)\bd).
\end{align}
\end{itemize}
\end{paragr}

\begin{paragr}
  We shall simply say that a map $\be:\ago_{P,\CC}^{G,*}\to \CC$ is a Paley-Wiener function if it is given by the Fourier-Laplace transform of a compactly supported smooth function on $\ago_P^{G,*}$. 
\end{paragr}

\begin{paragr}\label{S:les-objets}
  In the following we consider the following objects:
  \begin{itemize}
  \item $\be$ a Paley-Wiener function ;
  \item $f\in \Cc(G(F_\infty))$ a decomposable function;
  \item $\kappa\in \ago_P^*$;
  \item $\varphi\in \Ac_{P,\pi}(G)$.
  \end{itemize}
From these, one defines the pseudo-Eisenstein series
\begin{align*}
\theta(g)=  \sum_{\delta\in P(F)\back G(F)} B(\delta g).
\end{align*}
where $B$ is the function  on $A_G^\infty N(\AAA)M(F)\back G(\AAA)$ given by
  \begin{align*}
    B(g)= \int_{ \kappa+i\ago_P^{G,*}}    \exp(\bg \la,H_P(g)\bd) (I_P(\la,f)\varphi)(g)  \be(\la)  \,d\la.
  \end{align*}

In the following, we shall assume that  $\kappa$ is in the region of convergence of the Eisenstein series $E(\varphi,\la)$. Then  we have:
\begin{align*}
\theta(g)= \int_{\kappa+i\ago_P^{G,*}}E(g, I_P(\la,f)\varphi,\la) \be(\la)\, d\la.
\end{align*}
for any $g\in G(\AAA)$.

We fix $r>2$ that satisfies the conditions of  proposition \ref{prop:JlaPhi}.  In the following, we assume moreover that $\kappa$ is such that $\bg \kappa, \al^\vee \bd >r$.
\end{paragr}

\begin{paragr}\label{S:premier-calcul}

  \begin{proposition}\label{prop:premier-calcul}
Let $s\in \CC \setminus\{0,1\}$.  Assume that $\be$ vanishes at the points $-\mu_{i,s}$ for $i\in I$. Then we have
    \begin{align}\label{eq:theta-J}
      \int_{[G']_0} \theta(h) E(h,\Phi,s)\,dh=   \sum_{\xi\in \,_PW_{P,2}, \xi\not=1  }\int_{\kappa+i\ago_P^{G,*}} J(\xi,I_P(\la,f)\varphi,\la,\Phi,s) \be(\la)\, d\la
    \end{align}
where both sides are absolutely convergent.
  \end{proposition}

  \begin{preuve}
Because  $\theta(h)$ is rapidly decreasing, the left-hand side is  absolutely convergent. Using the majorizations  \eqref{eq:maj-epstein} and \eqref{eq:maj-eisenstein-lacunaire} given in the proof of proposition \ref{prop:JlaPhi}, we see that the right-hand side of \eqref{eq:theta-J} is also absolutely convergent.

By the lemma \ref{lem:PGG'}, the left-hand side of \eqref{eq:theta-J} is given by 
\begin{align}\label{eq:thetaP}
  \sum_{\xi\in  \,_PW_{P,2}} \int_{ A_{G'}^{\infty} P_{\tilde\xi}(F)\back G'(\AAA) } B(\tilde\xi h) E(h,\Phi,s)\,dh
 \end{align}
where $\tilde\xi\in G(F)$ is any element such that  $\xi=\tilde\xi c(\tilde\xi)^{-1}$.

Assume $\xi\not=1$. Using the definition of $B$ and permuting the adelic and the complex integrals we get that 
\begin{align*}
  \int_{ A_{G'}^{\infty} P_{\tilde\xi}(F)\back G'(\AAA) } B(\tilde\xi h) E(h,\Phi,s)\,dh=\int_{\kappa+i\ago_P^{G,*}} J(\xi,I_P(\la,f)\varphi,\la,\Phi,s) \be(\la)\, d\la.
\end{align*}
This permutation is easily justified by the majorization \eqref{eq:maj-eisenstein-lacunaire} and the fact that $\be$ is a Paley-Wiener function. 
We have to compute the term corresponding to $\xi=1$ (for which we take $\tilde\xi=1$) namely
\begin{align*}
  \int_{ A_{G'}^{\infty} P'(F)\back G'(\AAA) } B( h) E(h,\Phi,s)\,dh.
\end{align*}
We will show that this integral vanishes.  Using Iwasawa decomposition, we can write  it  as follows:
\begin{align}\label{eq:thetaP2}
  \int_{K'} \int_{ M'(F)Z_{G'}(\AAA)^1\back M'(\AAA)^1}\int_{ A_{G'}^{\infty}\back A_{M'}^\infty} \exp( \bg -\rho_P,H_P(a)\bd) B(a m k) \int_{[N_{P'}]}E( n amk ,\tilde\Phi(\om,s))  dn dm dk.
 \end{align}
According to the shape of the constant term in \eqref{eq:cst-term}, we are reduced to fix $i\in I$, $m\in  M'(\AAA)^1$, $k\in K'$ and to show the vanishing of 
\begin{align*}
  \int_{ A_{G'}^{\infty}\back A_{M'}^\infty}   \exp(\bg \mu_{i,s} ,H_P(a)\bd)\int_{ \kappa+i\ago_P^{G,*}} (I_P(\la,f)\varphi)(mk)  \be(\la)\exp(\bg \la,H_P(a)\bd)     \,d\la da.
\end{align*}
for any $m\in M'(\AAA)^1$ and $k\in K'$. By Fourier inversion,  this is, up to a constant, 
$$ (I_P(-\mu_{i,s},f)\varphi) (mk)  \be(-\mu_{i,s})$$ 
and we are done.
  \end{preuve}
\end{paragr}

\subsection{Period of a pseudo-Eisenstein series: computation in terms of the Whittaker functional}

\begin{paragr}
  We fix a non-trivial character $\psi':\AAA/F\to \CC^\times$.  We define $\psi:[N_{2n}]\to \CC^\times$ by
  \begin{align*}
    \psi(u)=\psi'(  \mathrm{Tr}_{E/F}(\tau\sum_{i=1}^{2n-1}u_{i,i+1}))
  \end{align*}
for $u=(u_{i,j})\in N_{2n}(\AAA)$. Note that $\psi$ is trivial on $N_{2n}'(\AAA)$ where $N_{2n}'=N_{2n}\cap G'$.
\end{paragr}

\begin{paragr}[A variant of mirabolic subgroups.] ---\label{S:mirabol-variant}
  For all $1\leq i \leq 2n$, we define the following subgroup of $G_{2n}$:
  \begin{align*}
    \pc_i=\{
  \begin{pmatrix}
    g & * \\ 0 & u 
  \end{pmatrix}\mid g\in G_{2n-i}, u\in N_{i}
\} 
  \end{align*}
  Note that $\pc_1$ is the mirabolic subgroup $\pc$ defined in §\ref{S:Epstein} and $\pc_{2n}=B_{2n}$. We denote by $N_{\pc_i}$ the unipotent radical of $\pc_i$. Let $P_i$ be the standard parabolic subgroup of $G$ of type $(2n-i,i)$. We have $\pc_i\subset P_i$ and $P_n$ is the  parabolic subgroup $P$ defined in §\ref{S:notations-W}. We denote by an upper  script $'$ the subgroups obtained by intersection with $G'$ that is $\pc_i'=\pc_i\cap G'$.

For any smooth function $\phi$  on $\pc_i(F)\back G(\AAA)$ we put for $g\in G(\AAA)$
\begin{align*}
  W_i(g,\phi)=\int_{[N_{\pc_i}]} \phi(ug)\psi(u)^{-1}\,du.
\end{align*}
When $i=2n$, we have $N_{\pc_i}=N_{2n}$ and we omit the subscript: we set  $W=W_{2n}$.
\end{paragr}

\begin{paragr}
  For any $\varphi\in \Ac_{P,\pi}(G)$ and $\la\in \ago_{P,\CC}^{G,*}$ such that $\Re(\bg \la,\al^\vee\bd)$ is large  enough, we define
\begin{align*}
  \mathbb{W}(g,\varphi,\la)=\int_{(N_{2n}\cap M)(F) \back N_{2n}(\AAA)} \exp(\bg \la,H_P(\xi_0 ug) \bd) \varphi(\xi_0 u g) \psi(u)^{-1}\, du.
\end{align*}
One has in the convergence region
\begin{align}\label{eq:Whit-Eis-Jacquet}
     \mathbb{W}(g, \varphi,\la)=W(g, E(\varphi,\la))
\end{align}
and so the integral has a meromorphic continuation to  $\ago_{P,\CC}^{G,*}$. It factors through the Fourier coefficient
\begin{align*}
  W_M^\psi(g,\varphi)=\int_{ [N_{2n}\cap M]} \varphi(ug)\psi(u)^{-1}\, du.
\end{align*}
\end{paragr}

\begin{paragr} We view $W_M^\psi(1,\cdot)$ as a Whittaker functional on $\pi$.  Let $\pi=\otimes_{v\in V_F}'\pi_v$ be a decomposition of $\pi$ as a restricted tensor product of representations $\pi_v$ of $M(F_v)$. According to this decomposition, we fix a Whittaker functionnal  $W_{M,v}^\psi$  on $\pi_v$ such that  $W_M^\psi(1,\cdot)=\otimes_{v\in V_F} W_{M,v}^\psi$.   Let $V_{F,\infty}\subset S\subset V_F$ be a finite set such that all objects  $\pi$, $\Phi$, $\psi$, $E/F$ are unramified outside $S$.  We put $W_{M,S}^{\psi}=\otimes_{v\in S} W_{M,v}^\psi$. For any $\la\in i\ago_{P,\CC}^{G,*}$, let $\pi_{v,\la}$ be the representation of $M(F_v)$ given by  $\pi_{v,\la}(m)=\exp(\bg\la, H_M(m)\bd) \pi_{v}(m)$ for $m\in M(F_v)$.  Let $\Pi_{v,\la}=\Ind_{P(F_v)}^{G(F_v)} (\pi_{v,\la})$ and  $\Pi_v=\Pi_{v,\la=0}$.  We put  $\Pi_S=\otimes_{v\in S}\Pi_v$.

 For any $g_S\in G(F_S)$ and $\phi\in \Pi_S$, we define the Jacquet integral by the analytic continuation of
  \begin{align}\label{eq:Jacquet integral}
     \mathbb{W}_S(g_S,\phi,\la)=\int_{(N_{2n}\cap M)(F_S) \back N_{2n}(F_S)} W_{M,S}^{\psi}(\Pi_\la(\xi_0 ug_S)\phi) \, du.
  \end{align}
  
  For any $v\notin S$, let $W_v^\psi(\pi_\la)$ be the $K_v$-invariant  Whittaker function such that $W_v^\psi(1,\pi_\la)=1$.   We shall identify  $\varphi$ with $\varphi_S \otimes \varphi^S$ where $\varphi_S\in \Pi_S $ and  $\varphi^S\in  \otimes_{v\notin S}'\Pi_v$ is $K^S$-invariant. Then we have 
  
\begin{align}\label{eq:decom-non-cst}
  \mathbb{W}(g, I_P(\la,f)\varphi,\la)=  \frac1{L^S(1+\bg\la,\al^\vee\bd, \sigma\times \sigma^c)}\mathbb{W}_S(g_S,I_P(\la,f)\varphi_S,\la)  \prod_{v\notin S} W_v^\psi(g_v,\pi_\la)
\end{align}
for all $g=(g_S, (g_v)_{v\notin S})\in G(\AAA)$.

  \begin{proposition}\label{prop:cv-integ-whitt}
    There is $s_0\in \RR$ such that for any $\la\in \kappa +i\ago_P^{G,*}$ and any $s\in \hc_{>s_0}$ the integral 
    \begin{align*}
      \int_{N_{2n}'(\AAA)\back G'(\AAA)}| \mathbb{W}(h, I_P(\la,f)\varphi,\la) \Phi(e_{2n}h)  |\det(h)|^s |\,dh  
    \end{align*}
is absolutely convergent and uniformly bounded on compact  subsets of $\hc_{>s_0}$.
  \end{proposition}

  \begin{preuve}
   We shall use \eqref{eq:decom-non-cst}. First the factor  ${L^S(1+\bg\la,\al^\vee\bd, \sigma\times \sigma^c)}^{-1}$ is bounded uniformy for $\la$ such $\bg\Re(\la),\al^\vee\bd\geq 2$. We have $\Phi=\Phi_S\otimes\Phi^S$ where $\Phi^S$ is the characteristic function of $(\oc^S)^{2n}$. By (the proof of) \cite[proposition 2.6.1 and lemma 3.3.1]{BPAsai}, there exists $s_0\in \RR$  such that for any $\la\in \kappa +i\ago_P^{G,*}$ and any $s\in \hc_{>s_0}$ the integral 
     \begin{align*}
      \int_{N_{2n}'(F_S)\back G'(F_S)}| \mathbb{W}_S(h_S, I_P(\la,f)\varphi_S,\la) \Phi_S(e_{2n}h_S)  |\det(h_S)|^s |\,dh_S
    \end{align*}
is convergent and uniformly bounded on compact  subsets of $\hc_{>s_0}$. So we are left with 
\begin{align*}
        \prod_{v\notin S} \int_{N_{2n}'(F_v)\back G'(F_v)}   |W_v^\psi(h_v,\pi_\la)| \mathbf{1}_{\oc_v^{2n}}(e_{2n}h_v)  |\det(h_v)|_v^s |\,dh_v.  
    \end{align*}
Then we can use the Iwasawa decomposition and the bound of \cite[proposition 2.4.1]{JPSS-GL3-I} since we may assume that $S$ is large enough so that for $v\notin S$  the cardinality of the residue field of $F_v$ is bigger than $n$. The details are left to the reader (see also \cite[proof of lemma 4.5]{IY}).
  \end{preuve}
\end{paragr}

\begin{paragr} We consider the situation of §\ref{S:les-objets}.

 \begin{proposition}\label{prop:second-calcul}
 There is $s_0\in \RR$ such that for any $s\in \hc_{>s_0}$ and any Paley-Wiener function $\be$ that vanishes at the points $\pm s\al/2$, we have
    \begin{align}\label{eq:theta-W}
      \int_{[G']_0} \theta(h) E(h,\Phi,s)\,dh=   \int_{\kappa+i\ago_P^{G,*}} \left(\int_{N_{2n}'(\AAA)\back G'(\AAA)}\mathbb{W}(h, I_P(\la,f)\varphi,\la) \Phi(e_{2n}h)|\det(h)|^s \,dh   \right) \be(\la)\, d\la,
    \end{align}
where the two sides are given by absolutely convergent integrals.
  \end{proposition}

Before giving the proof (which is to find in  §§\ref{S:proof-2nd-calcul}-\ref{S:vanish-n-1} below), we shall give a corollary.

\begin{corollaire}\label{cor:intert-zeta}
  There is $s_0\in \RR$ such that for any $\la\in \kappa+i\ago_P^{G,*}$ and $s\in \hc_{>s_0}$  we have
  \begin{align*}
     \sum_{\xi\in \,_PW_{P,2}, \xi\not=1  }J(\xi,\varphi,\la,\Phi,s)= \int_{N_{2n}'(\AAA)\back G'(\AAA)}\mathbb{W}(h, \varphi,\la) \Phi(e_{2n}h)|\det(h)|^s \,dh
  \end{align*}
where the right-hand side is given by an absolutely and uniformly convergent integral  on compact subsets of  $\hc_{>s_0}$.
\end{corollaire}

\begin{preuve}
  The combination of propositions    \ref{prop:premier-calcul} and \ref{prop:second-calcul} implies  (see  \cite[lemma 9.1.2]{LR} for a simple argument) 
 \begin{align*}
     \sum_{\xi\in \,_PW_{P,2}, \xi\not=1  }J(\xi,I_P(\la,f)\varphi,\la,\Phi,s)= \int_{N_{2n}'(\AAA)\back G'(\AAA)}\mathbb{W}(h, I_P(\la,f)\varphi,\la) \Phi(e_{2n}h)|\det(h)|^s \,dh
  \end{align*}
for any $\la\in \kappa+i\ago_P^{G,*}$. Since we can find $f$ such that $I_P(\la,f)\varphi=\varphi$ we get the result.
\end{preuve}

  \end{paragr}

\begin{paragr}[Proof of proposition \ref{prop:second-calcul}.] ---\label{S:proof-2nd-calcul} We first  recall the well-known computation of the constant term of a pseudo  Eisenstein series.

    \begin{lemme}\label{lem:vanish-cst-term}
    For $1\leq i \leq 2n$, the constant term of $\theta$ along $P_i$, defined by
    \begin{align*}
  \forall g\in G(\AAA)\ \ \    \theta_{P_i}(g)=\int_{[N_{P_i}]} \theta(ng)\, dn,
    \end{align*}
vanishes unless $i\in \{n,2n\}$. Moreover we have
\begin{align*}
   \theta_{P}(g)=\sum_{w\in W(M)} \int_{\kappa+i\ago_P^{G,*}} \exp(w\la, H_P(g)\bd  (M(w,\la)I_P(\la,f)\varphi)(g) \be(\la)\, d\la.
\end{align*}
  \end{lemme}

Unfolding the Epstein series, we get
  \begin{align*}
    \int_{[G']_0} \theta(h) E(h,\Phi,s)\, dh&=\int_{\pc'_1(F)\back G'(\AAA)} \theta(h) \Phi(e_{2n}h)|\det(h)|^s \,dh\\
&= \int_{\pc'_1(F)N'_{\pc_1}(\AAA)\back G'(\AAA)}  \left(\int_{[N'_{\pc_1}]} \theta(nh) \, dn\right)\Phi(e_{2n}h)|\det(h)|^s \,dh.
  \end{align*}
  The Fourier expansion of the map $n\in [N_{\pc_1}]\mapsto \theta(nh)$ gives:
  \begin{align*}
    \theta(h)=\theta_{P_1}(h)+ \sum_{\gamma\in  \pc_{2}(F)\back \pc_{1}(F)}  W_1(\gamma h, \theta)
  \end{align*}
from which we deduce
\begin{align*}
  \int_{[N'_{\pc_1}]} \theta(nh) \, dn=\theta_{P_1}(h)+ \sum_{\gamma\in  \pc_{2}'(F)\back \pc_{1}'(F)}   W_1(\gamma h, \theta).
\end{align*}
If $n>1$ then $\theta_{P_1}=0$ (see lemma \ref{lem:vanish-cst-term}). In particular, we get
\begin{align*}
    &\int_{[G']_0} \theta(h) E(h,\Phi,s)\, dh\\
  &=\int_{\pc'_2(F) N'_{\pc_1}(\AAA) \back G'(\AAA)} W_1(h, \theta)\Phi(e_{2n}h)|\det(h)|^s \,dh\\
&=\int_{\pc'_2(F)N'_{\pc_2}(\AAA)\back G'(\AAA)}  \left(\int_{[ N'_{\pc_1}\back  N'_{\pc_2}]} W_1(nh,\theta) \, dn\right)\Phi(e_{2n}h)|\det(h)|^s \,dh.
\end{align*}
Next, using  the Fourier expansion of $n\in  [N_{\pc_2}\cap M_{P_1}]\mapsto W_1(nh,\theta) $ we get:
\begin{align*}
\int_{[ N'_{\pc_1}\back  N'_{\pc_2}]} W_1(nh,\theta) \, dn=W_1(h,\theta_{P_2}) + \sum_{\gamma\in  \pc_{3}'(F)\back \pc_{2}'(F)}   W_2(\gamma h, \theta).
\end{align*}
If $n>2$ we have $\theta_{P_2}=0$ (see lemma \ref{lem:vanish-cst-term}). By recursion we get that the left-hand side of \eqref{eq:theta-W} is the sum of 
\begin{align}\label{eq:n-1}
\int_{\pc'_n(F) N'_{\pc_n}(\AAA) \back G'(\AAA)} W_{n-1}(h, \theta_{P_n})\Phi(e_{2n}h)|\det(h)|^s \,dh
\end{align}
and
\begin{align}\label{eq:n}
  \int_{\pc'_{n+1}(F) N'_{\pc_n}(\AAA) \back G'(\AAA)} W_n(h,\theta) \Phi(e_{2n}h)|\det(h)|^s \,dh  .
\end{align}
The manipulation is justified as in \cite[corollary 4.3 and bottom of p. 697]{IY}. The next step is to compute both expressions. Let's start with the second.
\end{paragr}

\begin{paragr}[Computation of \eqref{eq:n}.] --- We can continue the process. Since $\theta_{P_k}=0$ for $k>n$, $W=W_{2n}$ and $\pc'_{2n}(F) N'_{\pc_{2n}}(\AAA)=N'_{2n}(\AAA)$, we get that \eqref{eq:n} is equal to 
\begin{align*}
  \int_{N_{2n}(\AAA)\back G'(\AAA)} W(h,\theta) \Phi(e_{2n}h)|\det(h)|^s \,dh  .
\end{align*}
Using  \eqref{eq:Whit-Eis-Jacquet}, we also have
\begin{align*}
   W_{2n}(h,\theta)=  \int_{\kappa+i\ago_P^{G,*}}  \mathbb{W}(h, I_P(\la,f)\varphi,\la)    \be(\la)\, d\la.
\end{align*}
To get the right-hand side of \eqref{eq:theta-W},we just need to permute the adelic integral and the integral over $\la$. This  is justified by proposition \ref{prop:cv-integ-whitt}. To conclude, it suffices to show that  \eqref{eq:n-1} vanishes: this is done in the next §.
\end{paragr}

\begin{paragr}[Vanishing of \eqref{eq:n-1}.] --- \label{S:vanish-n-1}Recall that $P_n=P$. Using the Iwasawa decomposition, we get that the expression \eqref{eq:n-1} is equal to
  \begin{align}\label{eq:after-Iwa}
    \int_{A_{M'}^{\infty}}   \exp(-\bg \rho_P,H_P(a)\bd)  \int_{     M'(\AAA)\cap (\pc'_n(F) N'_{\pc_n}(\AAA) )\back  M'(\AAA)} \int_{K'}W_{n-1}(amk, \theta_{P})\Phi(e_{2n}ak)|\det(a)|^s \,dh.
  \end{align}
By lemma \ref{lem:vanish-cst-term}, we see that the expression $W_{n-1}(amk, \theta_{P})$ is the sum over $w\in W(M)$ of

\begin{align*}
 \int_{\kappa+i\ago_P^{G,*}} \exp(\bg w\la,H_P(a)\bd)  \int_{[N_{\pc_{n-1}}]}(M(w,\la)I_P(\la,f)\varphi  )(uamk)\psi(u)^{-1}\, du \be(\la) d\la.
\end{align*}

Writing $u=u_P u_M$ with $u_M\in (N_{\pc_{n-1}}\cap M)$ and $u_P\in (N_{\pc_{n-1}}\cap N_P)$ we see that
\begin{align*}
  (M(w,\la)I_P(\la,f)\varphi  )(uamk)= \exp(\bg \rho_P,H_P(a)\bd)  (M(w,\la)I_P(\la,f)\varphi  )(u_M mk).
\end{align*}
Let $A^{\infty}_{\pc_n'}$ be the stabilizer of $e_{2n}$ in $A^{\infty}_{M'}$. We have $|\det(a)|^s=\exp(s\bg \al, H_P(a)\bd/2)$. The contribution in \eqref{eq:after-Iwa} corresponding to $w\in W(M)$  factors through the integral:
\begin{align}\label{eq:preFourier}
   \int_{A_{\pc'_n}^{\infty}}   \exp(s\bg \al, H_P(a)\bd/2)   \int_{\kappa+i\ago_P^{G,*}} (M(w,\la)I_P(\la,f)\varphi  )(umk)  \exp(\bg w\la,H_P(a)\bd)\be(\la) \, d\la da.
\end{align}
for some Haar measure on $A_{\pc'}^{\infty}$. If $w=1$, the expression \eqref{eq:preFourier} is simply
\begin{align*}
 \int_{A_{\pc'_n}^{\infty}}   \exp(s\bg \al, H_P(a)\bd/2)   \int_{\kappa+i\ago_P^{G,*}}  (I_P(\la,f)\varphi  )(umk)  \exp(\bg \la,H_P(a)\bd)\be(\la) \, d\la da.
\end{align*}
 By Fourier inversion, it  is, up to a constant, $(I_P(-s\al/2,f)\varphi  )(umk)\be(-s\al/2)$ and thus vanishes.

If $w\not=1$, then $w\la=-\la$ and the expression \eqref{eq:preFourier} can be written as
\begin{align*}
   \int_{A_{\pc'_n}^{\infty}}   \exp(s\bg \al, H_P(a)\bd/2)   \int_{\kappa+i\ago_P^{G,*}} (M(w,\la)I_P(\la, f)\varphi  )(umk)  \exp(- \bg \la,H_P(a)\bd)\be(\la) \, d\la da.
\end{align*}
Assume that $\Re(s)$ is large enough so that $\Re(s)  \bg \al,\al^\vee\bd > 2 \bg \kappa,\al^\vee\bd$. One can check that there exists $c$ such that the inner integral vanishes unless $\bg \al, H_P(a)\bd\leq c$. Thus one can restricts the outer integral to this ``half-line''. Then we can permute the two integrals. By Cauchy formula, we get that it is, up to a constant, $(M(w,s \al/2) I_P(s\al/2, f) \varphi  )(umk)\be( s\al/2)$ and thus it also vanishes. This concludes the proof of proposition \ref{prop:second-calcul}.
\end{paragr}

\subsection{Final result}

\begin{paragr} We keep the notations of previous sections. Let $\la\in \ago_{P,\CC}^{G,*}$. For $s\in \CC$ let $L^S(s,\Pi_\la,\As)$ be the Asai $L$-function ``outside $S$''. We have
\begin{align*}
  {L^S(s,\Pi_\la,\As) }= L^S(s+\bg\la,\al^\vee\bd,\sigma,\As)L^S(s-\bg\la,\al^\vee\bd,\sigma^*,\As) { L^S(s, \sigma\times \sigma^\vee)   }.
\end{align*}
The factor  $L^S(s, \sigma\times \sigma^\vee)  $ has a simple pole at $s=1$ and we set:
\begin{align*}
  L^{S,*}(1, \sigma\times \sigma^\vee)  =\lim_{s\to 1}(s-1) L^S(s, \sigma\times \sigma^\vee)  .
\end{align*}
Then we get an analytic function in $\la$ by setting:
\begin{align*}
  {L^{S,*}(1,\Pi_\la,\As) }= L^S(1+\bg\la,\al^\vee\bd,\sigma,\As)L^S(1-\bg\la,\al^\vee\bd,\sigma^*,\As) { L^{S,*}(1, \sigma\times \sigma^\vee)   }.
\end{align*}

\begin{theoreme}\label{thm:J-W}
We have the following equality of meromorphic  functions on  $i\ago_{P,\CC}^{G,*}$:
\begin{align*}
  J(\xi_0,\varphi,\la) &= (\Delta_{G'}^{S,*})^{-1}  L^{S,*}(1,\Pi_\la,\As) \int_{N_{2n}'(F_S)\back \pc'(F_S)}\mathbb{W}(p_S,\varphi,\la)\, dp_S.
\end{align*}
where the integral is absolutely convergent, the left-hand side and the integrand in the right-hand are respectively defined in \eqref{eq:periode entrelacement} and \eqref{eq:Whit-Eis-Jacquet}.
Moreover both sides of the equality are holomorphic at  $\la\in i\ago_{P}^{G,*}$ in the following cases:
\begin{itemize}
\item $\sigma$ is not $G_n'$-distinguished;
\item $\sigma$ is $G_n'$-distinguished  and $\la\not=0$.
\end{itemize}
  \end{theoreme}
\end{paragr}

\begin{paragr}[Proof of theorem \ref{thm:J-W}.] ---
  By \eqref{eq:decom-non-cst}, we have for all $g_S\in G(F_S)$
  
\begin{align*}
  \mathbb{W}(g_S, \varphi,\la)=  \frac1{L^S(1+\bg\la,\al^\vee\bd, \sigma\times \sigma^c)}\mathbb{W}_S(g_S,\varphi_S,\la).
\end{align*}
Let $\Phi_S$ a test function in  the Schwartz space $\Sc(F_S^{2n})$. Let's consider the following integrals:
\begin{align}\label{eq:integ-S-W}
      \int_{N_{2n}'(F_S)\back \pc'(F_S)} \mathbb{W}_S(h_S, \varphi_S,\la)   \,dh_S
\end{align}
and
 \begin{align}\label{eq:integ-S-W-dets}
      \int_{N_{2n}'(F_S)\back G'(F_S)} \mathbb{W}_S(h_S, \varphi_S,\la) \Phi_S(e_{2n}h_S)  |\det(h_S)|^s  \,dh_S.
    \end{align}

\begin{lemme}\label{lem:cv-integW}
  \begin{enumerate}
  \item  There exists $\eta>0$ (resp. and $\eps>0$) such that the integral \eqref{eq:integ-S-W}, resp. \eqref{eq:integ-S-W-dets}, is absolutely convergent and holomorphic  on the subset of $\la\in \ago_{P}^{G,*}$ such that $|\bg \Re (\la),\al^\vee\bd|<\eta$, resp. and $s\in \hc_{1-\eps}$.
  \item The integral \eqref{eq:integ-S-W-dets} admits a meromorphic continuation to $\CC\times \ago_{P,\CC}^{G,*}$ denoted by $Z_S(s,\la,\Phi_S)$.
  \item For any  $c> 0$ there exists $s_1\in \RR$ such that for $s\in \hc_{s_1}$ and  $|\bg \Re (\la),\al^\vee\bd|<c$ the integral \eqref{eq:integ-S-W-dets} is absolutely convergent and coincides with  $Z_S(s,\la,\Phi_S)$.
  \end{enumerate}
  \end{lemme}

  \begin{preuve}
    All the results are slight variations on  \cite[lemma 3.3.1, lemma 3.3.2 and section 3.10]{BPAsai}. To get the convergence of \eqref{eq:integ-S-W} or the precise lower bound $1-\eps$, we need to observe that $\varphi_S$ belongs to the induced representation of the $S$-component of an irreducible automorphic cuspidal representation of $M(\AAA)$. As such, it is an irreducible, generic and unitary representation. We can now use the classification of local irreductible generic and unitary representations as fully induced from essentially discrete series with exponents in $]-1/2;1/2[$ (see e.g. \cite{Tadic-unit-local}, \cite{ZelevII} and \cite[Theorem 8.2.]{Ba-Renard}).
  \end{preuve}

Let $\Phi=\Phi_S\otimes \Phi^S$ where $\Phi_S$ is the Schwartz space $\Sc(F_S^{2n})$ and $\Phi^S$ is the characteristic function of $(\oc^{S})^{2n}$. Let $c> \bg \kappa,\al^\vee\bd$.  By lemma \ref{lem:cv-integW} assertion 3 there exists $s_1\in \RR$ such that for $s\in \hc_{s_1}$ and  $|\bg \Re (\la),\al^\vee\bd|<c$ the integral \eqref{eq:integ-S-W-dets} is absolutely convergent. We may and shall assume that $s_1$ is large enough so that proposition  \ref{prop:cv-integ-whitt} holds for $s_1$. Then by  the factorization \eqref{eq:decom-non-cst}, some local computations \cite{Flicker} and lemma \ref{lem:cv-integW} assertion 2  we have
  \begin{align*}
    \int_{N_{2n}'(\AAA)\back G'(\AAA)}\mathbb{W}(h, \varphi,\la) \Phi(e_{2n}h)|\det(h)|^s \,dh=   (\Delta_{G'}^{S,*})^{-1} \frac{L^S(s,\Pi_\la,\As) }{ L^S(1+\bg\la,\al^\vee\bd, \sigma\times \sigma^c)}Z_S(s,\la,\Phi_S)\\
  \end{align*}
for any $\la\in\kappa+i\ago_P^{G,*}$ and $s\in \hc_{>s_1}$. By analytic continuation of the equality of corollary  \ref{cor:intert-zeta}, we have:
\begin{align*}
     \sum_{\xi\in \,_PW_{P,2}, \xi\not=1  }J(\xi,\varphi,\la,\Phi,s)=    (\Delta_{G'}^{S,*})^{-1} \frac{L^S(s,\Pi_\la,\As) }{ L^S(1+\bg\la,\al^\vee\bd, \sigma\times \sigma^c)} Z_S(s,\la,\Phi_S) 
  \end{align*}
for all $s\in \CC$ and $\la\in \dc_r$ given by proposition \ref{prop:JlaPhi}. We can compute the residue of the left-hand side at $s=1$ following proposition \ref{prop:JlaPhi}. We get:

\begin{align}\label{eq;proportionalite}
     \hat\Phi(0)  J(\xi_0,\varphi,\la) =   (\Delta_{G'}^{S,*})^{-1}  Z_S(1,\la,\Phi_S) \frac{L^{S,*}(1,\Pi_\la,\As) }{ L^S(1+\bg\la,\al^\vee\bd, \sigma\times \sigma^c)}
  \end{align}
  Both sides are analytic in $\la$.  Thus the equality for   $\la\in i\ago_{P}^{G,*}$. But by lemma \ref{lem:cv-integW} assertion 1 and  assertion 2,  one has
\begin{align*}
    Z_S(1,\la,\Phi_S)=  \int_{N_{2n}'(F_S)\back G'(F_S)} \mathbb{W}_S(h_S, \varphi_S,\la) \Phi_S(e_{2n}h_S)  |\det(h_S)| \,dh_S.
    \end{align*}  
Let $\bar{N}_1$ be the unipotent radical of the opposite of the parabolic subgroup $P_1$ defined in §\ref{S:mirabol-variant}. The standard Levi factor of $P_1$ decomposes as $G_{2n-1}\times G_1$. Thus we have $N_{2n}' \back  P_1'\simeq N_{2n}' \back \pc'\times G_1'$. By a usual decomposition of measures we get
\begin{align*}
     Z_S(1,\la,\Phi_S) &= \int_{ \bar{N}_1(F_S) } \int_{G_1'(F_S)}\left(\int_{N_{2n}'(F_S)\back \pc'(F_S)} \mathbb{W}_S(h  t n , \varphi_S,\la) \,dh\right)  \Phi_S(e_{2n}tn )|t|^{2n} \, dt dn.
\end{align*}
However we have also
\begin{align*}
    \hat\Phi_S(0) &=\int_{F_S^{2n}} \Phi_S(X) \,dX\\
&= \int_{ \bar{N}_1(F_S) } \int_{G_1'(F_S)} \Phi_S(e_{2n}tn )|t|^{2n} \, dt dn.
\end{align*}
Since \eqref{eq;proportionalite} holds for any Schwartz function $\Phi_S$, we get that $J(\xi_0,\varphi,\la) $ is equal to 
\begin{align*}
  & \frac{L^{S,*}(1,\Pi_\la,\As) }{L^S(1+\bg\la,\al^\vee\bd, \sigma\times \sigma^c)}\int_{N_{2n}'(F_S)\back \pc'(F_S)}\mathbb{W}_S(g_S,\varphi_S,\la)\\
  &= L^{S,*}(1,\Pi_\la,\As) \int_{N_{2n}'(F_S)\back \pc'(F_S)}\mathbb{W}(p_S,\varphi,\la)\, dp_S
 \end{align*}
where we have used the factorization  \eqref{eq:decom-non-cst}. In the first expression above, the  integral is holomorphic on $i\ago_P^{G,*}$, see lemma \ref{lem:cv-integW}. Using the factorization of $L^S(s, \sigma\times \sigma^c)$ in terms of Asai $L$-functions $L(s,\sigma,As^{\pm})$ we get:
\begin{align*}
  \frac{L^{S,*}(1,\Pi_\la,\As) }{L^S(1+\bg\la,\al^\vee\bd, \sigma\times \sigma^c)}=   L^{S,*}(1, \sigma\times \sigma^\vee)  \frac{L^S(1-\bg\la,\al^\vee\bd,\sigma^*,\As)}{L^S(1+\bg\la,\al^\vee\bd,\sigma,\As^-)}.
\end{align*}
On $i\ago_P^{G,*}$,  the $L$-function $L^S(1+\bg\la,\al^\vee\bd,\sigma,\As^-)$ does not vanish by  \cite[theorem 5.1]{ShaLfunction} and  $L^S(1-\bg\la,\al^\vee\bd,\sigma^*,\As)$ is holomorphic unless $\la=0$ and  $\sigma^*$ (thus $\sigma$) is $G_n'$-distinguished (see \cite{Flicker}). On the other hand, if $\sigma$  is not $G_n'$-distinguished then $J(\xi_0,\varphi,\la) $ is known to be holomorphic on $i\ago_{P}^{G,*}$, see \cite[lemma 8.1]{LapFRTF}. Otherwise,  $J(\xi_0,\varphi,\la) $ is holomorphic on $i\ago_{P}^{G,*}\setminus\{0\}$ but it may have a simple pole  at $\la=0$.

  \end{paragr}

%% file: FRspectral.tex
\section{The $(G,H)$-regular contribution to the Jacquet-Rallis trace formula: alternative proof}\label{chap: Alt proof}

\subsection{Statement}

\begin{paragr}
The goal of this section is to provide an alternative proof of the following combination of Theorem \ref{thm:Ichi} and Theorem \ref{thm:IPi=Ipi}.

\begin{theoreme}\label{theo alt GHreg contributions}
Let $\chi\in \Xgo(G)$ be a $(G,H)$-regular cuspidal datum and let $f\in \Sc(G(\AAA))$. Then:
\begin{enumerate}
\item If $\chi$ is not Hermitian, we have $I_\chi(f)=0$.
		
\item If $\chi$ is Hermitian, we have
\begin{equation}
\displaystyle I_\chi(f)=2^{-\dim(\mathfrak{a}_L)}\int_{i\mathfrak{a}_M^{L,*}} I_{\Pi_\lambda}(f) d\lambda
\end{equation}
where we recall that $(M,\pi)$ is a pair representing $\chi$, $L\supset M$ is the Levi subgroup defined by \eqref{eq:leLevi}, $\Pi_\lambda$ stands for the induced representation $\Ind_{P(\AAA)}^{G(\AAA)}(\pi_\lambda)$, for $P$ a chosen parabolic subgroup with Levi factor $M$, and $I_{\Pi_\lambda}$ is the relative character defined by \eqref{eq:IPIla}.
\end{enumerate}
\end{theoreme}
More precisely, the proof will be very similar to that given in \cite[Section 8]{BCZ} and is based on two ingredients of independent interests. The first one is that the Rankin-Selberg period (over $H$) admits a continuous extension to the space $\mathcal{T}_\chi([G])$ of functions of uniform moderate growth supported on a $H$-regular cuspidal datum $\chi$ and that this extension can moreover be described in terms of the analytic continuation of Zeta integrals of Rankin Selberg type. This was already established in \cite[Section 7]{BCZ}. The second ingredient is an explicit spectral decomposition of the Flicker-Rallis period (over $G'$) restricted to $\mathcal{S}_\chi([G'])$ when the cuspidal datum $\chi$ is ($G$-)regular. This was already done in \cite[Section 6]{BCZ} under the stronger assumption that $\chi$ is $*$-regular. The aim of the Subsection \ref{S: Spectral FR} is to state and prove the extension of this result to the regular case. Once established, we will be able to give a proof of Theorem \ref{theo alt GHreg contributions} in Subsection \ref{S: altproof} in much the same lines as \cite[\S 8.2]{BCZ}. 
\end{paragr}

\subsection{Spectral decomposition of the Flicker-Rallis period for regular cuspidal data}\label{S: Spectral FR}

\begin{paragr}[Zeta integrals]
Let $n\geqslant 0$ be an integer. We will freely use the notation introduced in Section \ref{sec:GH-contrib-JR}. For every $f\in \tc([G_n])$, we denote by
$$\displaystyle W_f(g)=\int_{[N_n]} f(ug) \psi_n(u)^{-1} du,\;\;\; g\in G_n(\AAA),$$
its Whittaker function and, for $\phi\in \Sc(\AAA^n)$, we set
$$\displaystyle Z_\psi^{\FR}(s,f,\phi):=\int_{N_n'(\AAA)\backslash G_n'(\AAA)} W_f(h) \phi(e_nh)\lvert \det h\rvert^s dh.$$
This expression is absolutely convergent for $\Re(s)$ sufficiently large. More precisely, for every $N>0$, there exists $c>0$ such that $Z_\psi^{\FR}(s,f,\phi)$ converges for $s\in \cH_{>c}$ (see \cite[Theorem 6.2.5.1]{BCZ}).
\end{paragr}

\begin{paragr}[Flicker-Rallis period]
Recall that for every $f\in \cC([G_n])$, the period integral
$$\displaystyle P_{G'_n}(f)=\int_{[G'_n]} f(h)dh$$
is convergent \cite[Theorem 6.2.6.1]{BCZ}.
\end{paragr}

\begin{paragr}[Hermitian cuspidal data]
Let $\chi\in \Xgo(G_n)$ be a cuspidal datum. Then, we can find a pair $(M,\pi)$ representing $\chi$ with
$$\displaystyle M=\prod_{i\in I} G_{n_i}^{\times d_i} \times \prod_{j\in J} G_{n_j}^{\times d_j} \times \prod_{k\in K} G_{n_k}^{\times d_k}$$
and
$$\displaystyle \pi=\bigotimes_{i\in I} \pi_i^{\boxtimes d_i} \boxtimes \bigotimes_{j\in J} \pi_j^{\boxtimes d_j} \boxtimes \bigotimes_{k\in K} \pi_k^{\boxtimes d_k}$$
for some disjoint finite sets $I$, $J$, $K$, families of positive integers $(n_l)_{l\in I\cup J\cup K}$, $(d_l)_{l\in I\cup J\cup K}$ and a family of distinct cuspidal automorphic representations $(\pi_l)_{l\in I\cup J\cup K}$ satisfying:
\begin{itemize}
	\item For every $i\in I$, $\pi_i\not\simeq \pi_i^*$;
	
	\item For every $j\in J$, $\pi_j\simeq \pi_j^*$ and $L(s,\pi_j,\As^{(-1)^{n+1}})$ has no pole at $s=1$;
	
	\item For every $k\in K$, $\pi_k\simeq \pi_k^*$ and $L(s,\pi_k,\As^{(-1)^{n+1}})$ has a pole at $s=1$.
\end{itemize}

Fixing data as above (which are unique up to reordering), we recall that $\chi$ is said {\em Hermitian} (see \S \ref{S:GHgeneric}) if the following condition is satisfied:
\begin{num}
\item For every $i\in I$, there exists $i^*\in I$ such that $\pi_{i^*}\simeq \pi_i^*$ and for every $j\in J$, $d_j$ is even.
\end{num}
\end{paragr}

\begin{paragr}\label{S LeviL}
Assume furthermore that $\chi$ is {\em regular} in the sense of \S \ref{S:GHgeneric} or \cite[\S 2.9.7]{BCZ} and fix a pair $(M,\pi)$ representing $\chi$ together with data $I$, $J$, $K$, $(n_l)_{l\in I\cup J\cup K}$, $(d_l)_{l\in I\cup J\cup K}$, $(\pi_l)_{l\in I\cup J\cup K}$ as in the previous paragraph. By the regularity assumption, we have $d_l=1$ for every $l\in I\cup J\cup K$. Moreover, $\chi$ is Hermitian if and only if $J=\emptyset$ and there exists an involution $i\mapsto i^*$ of $I$ without fixed point such that $\pi_{i^*}=\pi_i^*$ for every $i\in I$. If this is the case, we choose a subset $I'\subset I$ such that $I$ is the disjoint union of $I'$ and $(I')^*=\{i^*\mid i\in I' \}$ and we define a Levi subgroup $L\supset M$ by
$$\displaystyle L:=\prod_{i\in I'} G_{n_i+n_{i^*}}\times \prod_{k\in K} G_{n_k}.$$
\end{paragr}

\begin{paragr}\label{S betan}
Let $P$ be a parabolic subgroup with Levi factor $M$. For every $\lambda\in \ago_{M,\CC}^*$ and $f\in \cC([G_n])$, we set
$$\displaystyle \Pi_\lambda:=\Ind_{P(\AAA)}^{G_n(\AAA)}(\pi_\lambda)$$
and
$$\displaystyle W_{f,\Pi_\lambda}:=W_{f_{\Pi_\lambda}}$$
where $f_{\Pi_\lambda}\in \tc([G_n])$ is defined as in \cite[Eq. (2.9.8.14)]{BCZ}.

Assuming that $\chi$ is Hermitian, we define for every $\lambda\in i\ago_{M}^{L,*}$ a linear form $\beta_n$ on $\wc(\Pi_\lambda,\psi_n)$ by setting
$$\displaystyle \beta_n(W)=(\Delta^{S,*}_{G_n'})^{-1} L^{S,*}(1,\Pi,\As) \int_{N_n'(F_S)\backslash \pc_n'(F_S)} W(p_S)dp_S$$
for every $W\in \wc(\Pi_\lambda,\psi_n)$ where $S$ is a sufficiently large finite set of places of $F$ (depending on $W$). That the above integral is convergent follows from \cite[Proposition 2.6.1, Lemma 3.3.1]{BPAsai}, \cite{JS} and moreover the product stabilizes for $S$ sufficiently large by the unramified computation of \cite[Proposition 3]{Flicker}.
\end{paragr}

\begin{paragr}
\begin{theoreme}\label{theo FR generic}
Let $\chi\in \Xgo^{\reg}(G_n)$ for which we adopt the notation introduced in the previous three paragraphs. Then, for every $f\in \cC_\chi([G_n])$ and $\phi\in \Sc(\AAA^n)$ we have:
\begin{enumerate}
	\item If $\chi$ is Hermitian, the function $\lambda\in i\mathfrak{a}_M^{L,*}\mapsto \beta_n(W_{f,\Pi_\lambda})$ is Schwartz and the resulting map
	\begin{equation}\label{eq0 theo FR generic}
	\displaystyle \cC_\chi([G_n])\to \Sc(i\mathfrak{a}_M^{L,*}),\;\; f\mapsto \left( \lambda\in i\mathfrak{a}_M^{L,*}\mapsto \beta_n(W_{f,\Pi_\lambda})\right)
	\end{equation}
	is continuous.
	\item The function $s\mapsto (s-1)Z_\psi^{\FR}(s,{}^0 f,\phi)$ admits an analytic continuation to $\cH_{>1}$ with a limit at $s=1$ and setting $Z_\psi^{\FR,*}(1,{}^0f,\phi):=\lim\limits_{s\to 1^+} (s-1)Z_\psi^{\FR}(s,{}^0f,\phi)$, we have
	\begin{equation}\label{eq1 theo FR generic}
	\displaystyle Z_\psi^{\FR,*}(1,{}^0f,\phi)=\left\{\begin{array}{ll}
	2^{1-\dim(\mathfrak{a}_L)}\widehat{\phi}(0) \int_{i\ago_{M}^{L,*}}\beta_n(W_{f,\Pi_\lambda}) d\lambda \mbox{ if } \chi \mbox{ is Hermitian}, \\ \\
	0 \mbox{ otherwise.}
	\end{array}\right.
	\end{equation}
	
	\item The equality
	\begin{equation}\label{eq FR period Zeta integral}
	\displaystyle \widehat{\phi}(0)P_{G'_n}(f)=\frac{1}{2}Z_\psi^{\FR,*}(1,{}^0f,\phi).
	\end{equation}
\end{enumerate}
 
\end{theoreme}
\end{paragr}

\begin{paragr}
We note the following immediate corollary of Theorem \ref{theo FR generic}.

\begin{corollaire}\label{cor FR generic}
Let $\chi\in \Xgo^{\reg}(G_n)$ for which we adopt the notation introduced in the paragraphs \ref{S LeviL} and \ref{S betan}. Then, for every $f\in \cC_\chi([G_n])$ we have
\begin{equation}
\displaystyle P_{G'_n}(f)=\left\{\begin{array}{ll}
	2^{-\dim(\mathfrak{a}_L)} \int_{i\ago_{M}^{L,*}}\beta_n(W_{f,\Pi_\lambda}) d\lambda \mbox{ if } \chi \mbox{ is Hermitian}, \\ \\
	0 \mbox{ otherwise.}
\end{array}\right.
\end{equation}
\end{corollaire}
\end{paragr}

\begin{paragr}[Proof of Theorem \ref{theo FR generic} 1. and 2.]
Here we prove parts 1. and 2. of Theorem \ref{theo FR generic}. We will explain how to deduce the last part in the next subsection. The proof is actually along the same lines as that of \cite[Theorem 6.2.5.1]{BCZ}. Therefore, we will be brief on the parts that are similar.

Set
$$\displaystyle \Ac:=(i\RR)^{I\cup J\cup K}$$
and let $\Ac_0$ be the subspace of vectors $\underline{x}=(x_{\ell})_{\ell\in I\cup J\cup K}\in \Ac$ such that $\sum_{\ell\in I\cup J \cup K} x_\ell=0$. We equip $\Ac$ with the product of Lebesgue measures and $\Ac_0$ with the unique measure inducing on $\Ac/\Ac_0\simeq i\RR$ the Lebesgue measure via the map $\underline{x}\mapsto \sum_{\ell\in I\cup J \cup K} x_\ell$. We identify $\Ac$ with $i\ago_M^*$ by sending $\underline{x}\in \Ac$ to the unramified character
$$\displaystyle g=(g_{\ell})_{\ell\in I\cup J\cup K}\in M(\AAA)=\prod_{\ell\in I\cup J\cup K} G_{n_\ell}(\AAA)\mapsto \prod_{\ell\in I\cup J\cup K} \lvert \det g_\ell\rvert_{\AAA_E}^{x_\ell/n_{\ell}}.$$
We note that this isomorphism sends $\Ac_0$ onto $i\mathfrak{a}_M^{G_r,*}$. Furthermore, by our choice of measure on $i\mathfrak{a}_M^*$ (see \S \ref{S:Haar} as well as \cite[Eq. (2.3.1.1)]{BCZ}), this isomorphism sends the Haar measure just described on $\Ac$ to $(2\pi)^{r}P$ times the Haar measure on $i\mathfrak{a}_M^*$ where we have set
$$\displaystyle r=\dim(A_M)=\lvert I\cup J\cup K\rvert,\;\;\; P=\prod_{\ell\in I\cup J\cup K} n_\ell.$$
Set $f_{\underline{x}}=f_{\Pi_{\underline{x}}}$ for every $\underline{x}\in \Ac$. Then, by a computation completely similar to that leading to \cite[Eq. (6.3.0.3)]{BCZ} (using one of the main results of \cite{LapHC}), we have
\begin{equation*}
\displaystyle Z_\psi^{\FR}(s,{}^0f,\phi)=\frac{n}{P}(2\pi)^{-r+1}\int_{\Ac_0} Z_\psi^{\FR}(s,f_{\underline{x}},\phi) d\underline{x}
\end{equation*} 
for $\Re(s)\gg 1$.

Let $S_0$ be a finite set of places of $F$ including the Archimedean ones and outside of which $\pi$ is unramified and let $S_{0,f}\subset S_0$ be the subset of finite places. We fix, for every $l\in I\cup J\cup K$ and $v\in S_{0,f}$, polynomials $Q_l(T), Q_{l,v}(T)\in \CC[T]$ with all their roots in $\cH_{]0,1[}$ and $\cH_{]q_v^{-1},1[}$ respectively such that the products $s\mapsto Q_l(s) L_\infty(s,\pi_l,\As)$ and $s\mapsto Q_{l,v}(q_v^{-s})L_v(s,\pi_l,\As)$ have no pole in $\cH_{]0,1[}$. Let $I_0$ be the set of $i\in I$ for which there exists $i^*\in I$ (necessarily unique by regularity of $\chi$) satisfying $\pi_{i^*}=\pi_i^*$ and fix a subset $I'\subset I_0$ such that $I_0$ is the disjoint union of $I'$ and $(I')^*$. Set
\[\begin{aligned}
\displaystyle P(s,\underline{x})= & \prod_{i\in I'} (s+\frac{x_i+x_{i^*}}{n_i})(s+\frac{x_i+x_{i^*}}{n_i}-1)\prod_{k\in K} (s+\frac{2x_k}{n_k})(s-1+\frac{2x_k}{n_k})\\
& \times \prod_{l\in I\cup J\cup K} Q_l(s+\frac{2x_l}{n_l}) \prod_{\substack{l\in I\cup J \cup K\\ v\in S_{0,f}}} Q_{l,v}(q_v^{-s-\frac{2x_l}{n_l}})
\end{aligned}\]
Then, we can prove exactly as in \cite[\S 6.3]{BCZ} that the two functions
\begin{equation*}
(s,\underline{x})\in \CC\times \Ac_0\mapsto P(s+\frac{1}{2},\underline{x})Z_\psi^{\FR}(s+\frac{1}{2},f_{\underline{x}},\phi)
\end{equation*}
and
\begin{equation*}
(s,\underline{x})\in \CC\times \Ac_0\mapsto P(\frac{1}{2}-s,\underline{x})Z_{\psi^{-1}}^{\FR}(s+\frac{1}{2},\widetilde{f}_{\underline{x}},\widehat{\phi})
\end{equation*}
where $\widetilde{f}_{\underline{x}}(g)=f_{\underline{x}}({}^t g^{-1})$ satisfy the conditions of \cite[Corollary A.0.11.1]{BCZ}. Therefore, by the conclusion of this corollary, the map
$$\displaystyle s\mapsto F_s:=\left(\underline{x}\mapsto \prod_{i\in I'} (s+\frac{x_i+x_{i^*}}{n_i}-1)\prod_{k\in K} (s+\frac{2x_k}{n_k}-1) Z_\psi^{\FR}(s,f_{\underline{x}},\phi) \right)$$
induces a holomorphic map $\cH_{1-\epsilon}\to \Sc(\Ac_0)$ for some $\epsilon>0$. This already implies that
$$\displaystyle s\mapsto Z_\psi^{\FR}(s,{}^0f,\phi)=\frac{n}{P}(2\pi)^{-r+1}\int_{\Ac_0} \prod_{i\in I'} (s+\frac{x_i+x_{i^*}}{n_i}-1)^{-1}\prod_{k\in K} (s+\frac{2x_k}{n_k}-1)^{-1} F_s(\underline{x}) d\underline{x}$$
has an analytic continuation to $\cH_{>1}$. Moreover, if $\chi$ is not Hermitian the linear forms $\underline{x}\in \Ac_0\mapsto x_i+x_{i^*}$, $i\in I'$, and $\underline{x}\in \Ac_0\mapsto x_k$, $k\in K$, are linearly independent and thus by \cite[Lemma 3.11]{BPPlanch}, $Z_\psi^{\FR}(s,{}^0f,\phi)$ has a limit as $s\to 1$.

We assume from now on that $\chi$ is Hermitian so that $I=I'\cup (I')^*$ and $J=\emptyset$. Set $\Ac^L=(i\RR)^{I'}$, $\Ac'=(i\RR)^{I'\cup K}$ and let $\Ac'_0$ be the subspace of $\underline{x}\in \Ac'$ such that $\sum_{i\in I'} x_i+\sum_{k\in K} x_k=0$. We have a short exact sequence
$$\displaystyle 0\to \Ac^L\to \Ac_0\to \Ac'_0\to 0$$
where the first map sends $\underline{x}\in \Ac^L$ to the vector $\underline{y}\in \Ac_0$ with coordinates $y_i=x_i$ for $i\in I'$, $y_i=-x_{i^*}$ for $i\in (I')^*$ and $y_k=0$ for $k\in K$ whereas the second map is given by
$$\displaystyle \underline{x}\in \Ac_0\mapsto \left((x_i+x_{i^*})_{i\in I'}, (x_k)_{k\in K} \right).$$
We equip $\Ac^L$ and $\Ac'$ with the products of Lebesgue measures and $\Ac'_0$ as before with the unique measure inducing on $\Ac'/\Ac'_0\simeq i\RR$ the Lebesgue measure. Then, it is easy to see that the above exact sequence is compatible with the different Haar measures. In particular, we have
$$\displaystyle Z_\psi^{\FR}(s,{}^0f,\phi)=\frac{n}{P}(2\pi)^{-r+1}\int_{\Ac'_0} \prod_{i\in I'} (s+\frac{y_i}{n_i}-1)^{-1}\prod_{k\in K} (s+\frac{2y_k}{n_k}-1)^{-1} \int_{\Ac^L}  F_s(\underline{x}+\underline{y}) d\underline{x} d\underline{y}.$$
From this and \cite[Proposition 3.12]{BPPlanch}, we obtain
\begin{align}\label{eq1 proof of theo FR generic}
\displaystyle \lim\limits_{s\to 1^+} (s-1)Z_\psi^{\FR}(s,{}^0f,\phi) & =\frac{n}{P}(2\pi)^{-r+1} \frac{\prod_{i\in I'} n_i \prod_{k\in K} \frac{n_k}{2}}{\sum_{i\in I'} n_i+\sum_{k\in K} \frac{n_k}{2}} (2\pi)^{\lvert I'\cup K\rvert-1}\int_{\Ac^L}  F_1(\underline{x}) d\underline{x} \\
\nonumber & =\frac{(2\pi)^{-\lvert I'\rvert}}{P'} 2^{1-\lvert K\rvert}\int_{\Ac^L}  F_1(\underline{x}) d\underline{x}
\end{align}
where we have set $P'=\prod_{i\in I'} n_i$. Furthermore, from the unramified computation of \cite[Proposition 3]{Flicker} and \cite[Lemma 2.16.3]{BPPlanch} we have
\begin{align}\label{eq2 proof of theo FR generic}
\displaystyle F_1(\underline{x})=\lim\limits_{s\to 1^+} (s-1)^{\lvert I'\cup K\rvert}Z_\psi^{FR}(s,f_{\underline{x}},\phi)=\widehat{\phi}(0)\beta_n(W_{f_{\underline{x}}})
\end{align}
for $\underline{x}\in \Ac^L$. On the other hand, the isomorphism $\Ac\simeq i\mathfrak{a}_M^*$ sends $\Ac^L$ onto $i\mathfrak{a}_M^{L,*}$ and sends the measure on $\Ac^L$ to $\frac{(2\pi)^{\lvert I'\rvert} P'}{2^{\lvert I'\rvert}}$ times the measure on $i\mathfrak{a}_M^{L,*}$. As, by \cite[Corollary A.0.11.1]{BCZ}, the function $F_1$ is Schwartz, so is the function $\lambda\in i\mathfrak{a}_M^{L,*}\mapsto \beta_n(W_{f,\Pi_\lambda})$. Moreover, that the map \eqref{eq0 theo FR generic} is continuous follows from \eqref{eq2 proof of theo FR generic} together with \cite[Theorem 6.2.5.1 1., Eq. (A.0.4.3)]{BCZ} and the closed graph theorem. Finally, combining \eqref{eq1 proof of theo FR generic} with \eqref{eq2 proof of theo FR generic} and the above comparison of measures readily gives the identity \eqref{eq1 theo FR generic} and ends the proof of Theorem \ref{theo FR generic} 2.
\end{paragr}

\begin{paragr}[Proof of Theorem \ref{theo FR generic} 3.]

The proof of \cite[Theorem 6.2.6.1]{BCZ} applies verbatim noting that in {\em loc. cit.} the condition that $\chi$ is $*$-regular is only used at the end of the proof to show that the family of bilinear forms denoted by
$$\displaystyle s\mapsto P_{G'_r}\widehat{\otimes} \mathcal{Z}^{FR}_{n-r}(s)$$
on $\mathcal{C}_\chi([M_r])\times \Sc(\AAA^{n-r})$ extends holomorphically to $\cH_{s>1}$. However, thanks to the proof of part 2. of Theorem \ref{theo FR generic} given in the previous subsection, we know that this property continues to hold for cuspidal data that are only assumed to be regular.

\end{paragr}

\subsection{Proof of Theorem \ref{theo alt GHreg contributions}}\label{S: altproof}
Since $\chi$ is $H$-regular, by \cite[Theorem 7.1.3.1]{BCZ}, the period integral
$$\displaystyle P_H: \phi\in \Sc_\chi([G])\mapsto \int_{[H]} \phi(h)dh$$
extends continuously to a linear form on $\tc_\chi([G])$ that we shall denote by $P_H^*$. Then, by the very same argument as for \cite[Eq. (8.2.3.5)]{BCZ} we have
\begin{equation}\label{eq1 theo alt GHreg contributions}
\displaystyle I_\chi(f)=P_H^*\left(\int_{[G']} K_{f,\chi}(.,g') \eta_{G'}(g') dg' \right)
\end{equation}
where the inner integral is taken in $\tc_N([G])$ for $N$ large enough. Furthermore, applying Corollary \ref{cor FR generic} instead of \cite[Corollary 6.2.7.1]{BCZ}, the discussion of \cite[\S 8.2.4]{BCZ} shows that this inner integral is identically zero if $\chi$ is not Hermitian whereas otherwise it leads to the following expansion
\begin{equation}\label{eq2 theo alt GHreg contributions}
\displaystyle \int_{[G']} K_{f,\chi}(.,g') \eta_{G'}(g') dg'=2^{-\dim(\mathfrak{a}_L)}\int_{i\mathfrak{a}_M^{L,*}} \sum_{\phi\in \mathcal{B}_{P,\pi}} E(.,I_P(\lambda,f)\phi,\lambda) \overline{\beta_\eta(W(\phi,\lambda))} d\lambda
\end{equation}
for $\mathcal{B}_{P,\pi}$ a $K$-basis of $I_{P(\mathbb{A})}^{G(\mathbb{A})}(\pi)$. We assume from now on that $\chi$ is Hermitian. We claim:
\begin{num}
\item\label{eq3 theo alt GHreg contributions} There exists $N>0$ such that for every $\lambda\in i\mathfrak{a}_M^{L,*}$ the series
\begin{equation}\label{eq4 theo alt GHreg contributions}
\displaystyle \sum_{\phi\in \mathcal{B}_{P,\pi}} E(.,I_P(\lambda,f)\phi,\lambda) \overline{\beta_\eta(W(\phi,\lambda))}
\end{equation}
converges absolutely in $\tc_N([G])$ and furthermore the integral over $i\mathfrak{a}_M^{L,*}$ in the right hand side of \eqref{eq2 theo alt GHreg contributions} converges absolutely in $\tc_N([G])$.
\end{num}
Indeed, that the series \eqref{eq4 theo alt GHreg contributions} converges absolutely in $\tc_N([G])$ for large enough $N$ (independent of $\lambda$) follows from \cite[Theorem 2.9.8.1]{BCZ}. Note that for every $\lambda\in i\mathfrak{a}_M^{L,*}$ and $g\in [G]$ we have
$$\displaystyle \sum_{\phi\in \mathcal{B}_{P,\pi}} E(g,I_P(\lambda,f)\phi,\lambda) \overline{\beta_\eta(W(\phi,\lambda))}=\beta_\eta(K_f(g,.)_{\Pi_\lambda})$$
where $K_f(g,.)_{\Pi_\lambda}$ denotes the projection of the function $K_{f}(g,.)\in \Sc([G])$ to $\Pi_\lambda$ as defined in \cite[Eq. (2.9.8.14)]{BCZ}. Moreover, by \cite[Lemma 2.10.1.1]{BCZ}, for every continuous semi-norm $\mu$ on $\Sc([G])$ we can find $N_\mu>0$ such that for every $f'\in \Sc(G(\AAA))$ we have
$$\displaystyle \mu(K_{f'}(g,.))\ll_{\mu,f'} \lVert g\rVert_G^{N_\mu}, \mbox{ for } g\in [G].$$
In particular, combining this estimates with the first part of Theorem \ref{theo FR generic} we deduce that for each $d>0$ we can find $N_d>0$ such that, for $f'\in \Sc(G(\AAA))$, we have
$$\displaystyle \left\lvert \sum_{\phi\in \mathcal{B}_{P,\pi}} E(g,I_P(\lambda,f')\phi,\lambda) \overline{\beta_\eta(W(\phi,\lambda))} \right\rvert=\lvert \beta_\eta(K_{f'}(g,.)_{\Pi_\lambda})\rvert \ll_{d,f'} \lVert g\rVert_G^{N_d} (1+\lVert \lambda\rVert)^{-d}$$
for $(\lambda,g)\in i\mathfrak{a}_M^{L,*}\times [G]$. Choosing $d$ such that $\lambda\mapsto (1+\lVert \lambda\rVert)^{-d}$ is integrable on $i\mathfrak{a}_M^{L,*}$ and applying the above inequality to $f'=R(X)f$ for $X\in \mathcal{U}(\ggo)$ gives the second part of \eqref{eq3 theo alt GHreg contributions}.

We now need to recall how the extension $P_H^*$ is defined in \cite[Theorem 7.1.3.1]{BCZ}: for $\Phi\in \tc_\chi([G])$ the integral
$$\displaystyle Z^{RS}_\psi(s,\Phi)=\int_{N_H(\AAA)\backslash H(\AAA)}  W_\Phi(h) \lvert \det h\rvert_{\AAA_E}^s dh,$$
where $W_\Phi(h):=\int_{[N]}\Phi(uh) \psi_N(u)^{-1}du$, a priori only defined when $\Re(s)$ is sufficiently large, admits an analytic continuation to $\CC$ and we have
$$\displaystyle P_H^*(\Phi)=Z^{RS}_\psi(0,\Phi).$$
Therefore, by continuity of $P_H^*$, \eqref{eq2 theo alt GHreg contributions} and \eqref{eq3 theo alt GHreg contributions}, we obtain
\[\begin{aligned}
\displaystyle P_H^*\left(\int_{[G']} K_{f,\chi}(.,g') \eta_{G'}(g') dg' \right) & =2^{-\dim(\mathfrak{a}_L)}\int_{i\mathfrak{a}_M^{L,*}} \sum_{\phi\in \mathcal{B}_{P,\pi}} P_H^*(E(.,I_P(\lambda,f)\phi,\lambda)) \overline{\beta_\eta(W(\phi,\lambda))} d\lambda \\
 & =2^{-\dim(\mathfrak{a}_L)}\int_{i\mathfrak{a}_M^{L,*}} \sum_{\phi\in \mathcal{B}_{P,\pi}} Z^{RS}_\psi(0,W(I_P(\lambda,f)\phi,\lambda)) \overline{\beta_\eta(W(\phi,\lambda))} d\lambda \\
 & =2^{-\dim(\mathfrak{a}_L)}\int_{i\mathfrak{a}_M^{L,*}} I_{\Pi_\lambda}(f) d\lambda.
\end{aligned}\]
Together with \eqref{eq1 theo alt GHreg contributions} this ends the proof of the theorem.

%% file: Preuve.tex
\section{Proof of Gan-Gross-Prasad conjecture : Eisenstein case}\label{sec:Preuve Eis}

\subsection{Global comparison  of relative characters}\label{ssec:Globalcomparison}

\begin{paragr}[Notations.] --- We follow the notations of sections \ref{sec:spectral-unitary} and \ref{sec:GH-contrib-JR}. In particular, we fix an integer $n\geq 1$. We consider the group $G=G_n\times G_{n+1}$ and for $h\in \hc_n$ the groups $U_h'\subset U_h$. Since $n$ is fixed we drop it from the notation: $\hc=\hc_n$.  The various Haar measures are those considered  in §\ref{S:Haar}.
  \end{paragr}

\begin{paragr}\label{S:finitesetS0}  Let $V_{F,\infty}\subset S_0\subset V_F$ be a finite set of places containing all the places that are ramified in $E$.  For every $v\in V_F$, we set $E_v=E\otimes_F F_v$ and when $v\notin V_{F,\infty}$ we denote by $\oc_{E_v}\subset E_v$ its ring of integers. Let $\hc^\circ\subset \hc$ be the (finite) subset of Hermitian spaces of rank $n$ over $E$ that admits a selfdual $\oc_{E_v}$-lattice for every $v\notin S_0$. For each $h\in \hc^\circ$, the group $U_h$ is then  defined over $\oc_F^{S_0}$ and we fix a choice of such a model.  For $v\notin S_0$, we define the open compact subgroups $K_{h,v}=U_h(\oc_v)$ and $K_v=G(\oc_v)$ respectively of $U_h(F_v)$ and $G(F_v)$. We set
        \begin{align*}
          K_h^\circ=\prod_{v\notin S_0} K_{h,v} \text{ and } K^\circ=\prod_{v\notin S_0} K_v.
        \end{align*}
We choose also for each $v\in S_0$ some maximal compact subgroup $K_{h,v}\subset U_h(F_v)$  and $K_v\subset G(F_v)$ (see §\ref{S:unitary-choices}).
	
	Let $v\notin S_0$. We denote by $\Sc^\circ(U_h(F_v))$, resp.  $\Sc^{\circ}(G(F_v))$, the corresponding spherical Hecke algebra. We have the base change homomorphism 
	\begin{align*}
	BC_{h,v}: \Sc^{\circ}(G(F_v))\to \Sc^{\circ}(U_h(F_v)).
	\end{align*}
	We denote by $\Sc^{\circ}(U_h(\AAA^{S_0}))$, resp.  $\Sc^{\circ}(G(\AAA^{S_0}))$, the restricted tensor product of $\Sc^\circ(U_h(F_v))$, resp. $\Sc^{\circ}(G(F_v))$, for $v\notin S_0$. We set  $BC_h^{S_0}=\otimes_{v\notin S_0}   BC_{h,v}$.
	
	We also denote by $ \Sc^\circ (G(\AAA))\subset  \Sc (G(\AAA))$ and $\Sc^\circ(U_h(\AAA))\subset \Sc(U_h(\AAA))$, for $h\in \hc^\circ$, the subspaces of functions that are respectively bi-$K^\circ$-invariant and bi-$K_h^\circ$-invariant. 
	\end{paragr}

\begin{paragr}[Transfer.] --- Let $h\in  \hc^\circ$. We shall say that $f_{S_0}\in  \Sc(G(F_{S_0}))$  and  $f_{S_0}^h \in \Sc(U_h(F_{S_0}))$ are transfers  if the functions $f_{S_0}$ and $f_{S_0}^h$ have matching regular orbital integrals in the sense of  \cite[Definition 4.4]{BPLZZ}. 
\end{paragr}

\begin{paragr}[Cuspidal datum $\chi_0$ and $H$-regular Hermitian Arthur parameter $\Pi$.] --- \label{S:lesdonnees} Let $P$ be a standard parabolic subgroup of $G$ and $\pi$ be a cuspidal automorphic representation of $M=M_P$. Let $\chi_0\in \Xgo(G)$ be the class of the pair $(M_P,\pi)$. We assume henceforth that $\chi_0$ is a Hermitian $(G,H)$-regular cuspidal datum in the sense of \S \ref{S:GHgeneric}. We assume also that $(M,\pi)$ satisfies the conditions of §\ref{S:choiceL} and that $L$ is the standard Levi subgroup containing $M$ defined in this §. Set $\Pi=\Ind_P^G(\pi)$: this is a $H$-regular Hermitian Arthur parameter in the sense of §\ref{S:AparamG}. We assume  that $S_0$ is large enough so that $\Pi$ admits $K^\circ$-fixed vectors. Let $\Pi_0$ be the discrete component of $\Pi$ (see §\ref{S:Aparam}).

The group $S_\Pi=S_{\Pi_0}$ is defined in §§ \ref{S:Aparam} and \ref{S-intro:product}. Its order can be computed as follows:
\begin{align}
  \label{eq:card-SPi}
|S_\Pi|=2^{\dim(\ago_L)-\dim(\ago_M^L)}.
\end{align}

For any $\la\in i\ago_M^{L,*}$, we shall use the  distribution $I_{P,\pi}(\la,\cdot)$ on $\Sc(G(\AAA))$  defined in  §\ref{S:RelcharJPpi}.
\end{paragr}

\begin{paragr} Let $S'_0$ be the union of $S_0\setminus V_{F,\infty}$ and the set of all finite places of $F$ that are inert in $E$.   Let $h\in \hc^\circ$ and let $\Xgo_\Pi^h\subset \Xgo(U_h)$  be the set of cuspidal data represented by  pairs $(M_h,\sigma)$ such that $M_h$ is the standard Levi factor of a standard parabolic sugroup $P_h$ of $U_h$ and $\sigma$ is a cuspidal automorphic representation of $M_h(\AAA)$ such that 
	\begin{itemize}
        \item $\Pi$ is a weak base change of $(P_h,\sigma)$ (in the sense of §§\ref{S:weakBC} and \ref{S-intro:product});
        \item $\sigma$ is $\prod_{v\notin S_0} (M_h(F_v)\cap K_{h,v})$-unramified;
        \item with the identifications  $M_h=G^\sharp\times  U$ where $G^\sharp$ is a product of linear groups and  $U=U(h_{n_0})\times U(h_{n_0'})$ for some $h_{n'_0} \in \hc_{n_0}$ and  $h_{n_0} \in \hc_{n_0'}$ and   $\sigma=\sigma^\sharp\boxtimes \sigma_0$ accordingly (where $\sigma_0$ is a cuspidal automorphic representation of $U(\AAA)$), for all $v\notin S_0'\cup V_{F,\infty}$,    the representation $\Pi_{0,v}$  is  the split base change of the representation $\sigma_{0,v}$.
        \end{itemize}

Since $\Pi$ is a $H$-regular Hermitian Arthur parameter,  the class of  $(M_h,\sigma)$ is $(U,U')$-regular in the sense of §\ref{S:UU'reg}. Moreover  there is  a natural  isomorphism 
\begin{align}
  \label{eq:aMh-aML}
bc=bc_{(M_h,\sigma)}:\ago_M^{L,*}\to \ago_{M_h}^*.
\end{align}
More precisely, we wrote $M_h=G^\sharp\times  U$ where $G^\sharp$ is a product of linear groups $G_{m_i}$ for some integers $m_i$ and $1\leq i \leq s$ so that $\ago_{M_h}^*=\oplus_{i=1}^s \ago_{G_{m_i}}^*$. But the product $\prod_{i=1}^s  (G_{m_i}\times G_{m_i})$ is also a factor in the decomposition of $M$ into a product of  linear groups. So the space $\ago_{M_h}^*\oplus \ago_{M_h}^*=\oplus_{i=1}^s (\ago_{G_{m_i}}^*\oplus \ago_{G_{m_i}}^*)$ is a subspace of $\ago_M^{*}$. The antidiagonal map $x\mapsto (x,-x)$ then identifies $\ago_{M_h}^*$ with a subspace of $\ago_{M_h}^*\oplus \ago_{M_h}^*$ and thus of $\ago_M^{*}$ and this subspace is precisely $\ago_M^{L,*}$.

\begin{remarques}\label{rq:bc-notation}  The base change map $bc$ depends implicitly on the choice of the pair $(M_h,\sigma)$. Indeed,  if $(M_{1,h},\sigma_1)$ is a pair equivalent to $(M_h,\sigma)$ and if $w\in W(M_h,M_{1,h})$ is such that $\sigma_1=w\cdot \sigma$ then $bc_{(M_{1,h},\sigma_1)}=w \circ bc_{(M_h,\sigma)}$ (where we view $w$ as a map $\ago_{M_h}^*\to \ago_{M_{1,h}}^*$). In order to not burden the notation we shall omit the subscript $(M_h,\sigma)$ from the notation  $bc$ in \eqref{eq:aMh-aML}.

We warn also the reader that the map \eqref{eq:aMh-aML} does not preserve the various choices of Haar measures. In fact the pullback of the measure on $\ago_{M_h}^*$ is $2^{\dim(\ago_{M}^L)}$ times the  measure on $\ago_M^{L,*}$.
\end{remarques}

For any $\la\in i\ago_{M}^{L,*}$, we set 
\begin{align}\label{eq:JPih}
  J^{h}_{\Pi}(\la,f)=\sum_{(M_h,\sigma) }J^h_{P_h,\sigma}(bc(\la),f)
\end{align}
where the sum is over a set of representatives $(M_h,\sigma)$ of classes in $\Xgo_\Pi^h$ and  $J^{h}_{P_h,\sigma}(bc(\la),f)=J^{U_h}_{P_h,\sigma}(bc(\la),f)$ is the distribution introduced in  §\ref{S:rel-charU}. By proposition \ref{prop:JU-eqF} and remarks \ref{rq:bc-notation}, the distribution does not depend on the choice of the representative.

The sum above is \emph{a priori} absolutely convergent and the convergence is uniform on compact subset of $i\ago_M^{L,*}$ (as follows from from a reinforcement of \cite[proposition 2.8.4.1]{BCZ} based on results of Müller, see \cite[corollary 0.3]{Muller-traceII}). In particular, the expression $J^{h}_{\Pi}(\la,f)$ is holomorphic for $\la\in i\ago_{M}^{L,*}$.
\end{paragr}

\begin{paragr}[A global relative characters identity.] ---

  	\begin{theoreme}\label{thm:comparison}
		Let $f\in \Sc^\circ(G(\AAA))$ and $f^h\in \Sc^\circ(U_h(\AAA))$ for every $h\in \hc^\circ$. Assume that the following properties are satisfied for every $h\in \hc^\circ$:
		\begin{enumerate}
			\item $f=(\Delta_H^{S_0,*} \Delta_{G'}^{S_0,*} ) f_{S_0}\otimes f^{S_0}$ with $f_{S_0}\in \Sc(G(F_{S_0}))$ and $f^{S_0}\in  \Sc^{\circ}(G(\AAA^{S_0}))$.
			\item $f^h= (\Delta_{U_h'}^{S_0})^2f^h_{S_0}\otimes f^{h,S_0}$ with $f_{S_0}^h \in \Sc(U_h(F_{S_0}))$ and $f^{h,S_0} \in  \Sc^{\circ}(U_h(\AAA^{S_0}))$. 
			\item The functions $f_{S_0}$ and $f_{S_0}^h$ are transfers.
			\item $f^{h,S_0}=BC^{S_0}_h(f^{S_0})$
			\item The function $f^{S_0}$ is a product of a smooth compactly supported function on the restricted product $\prod_{v\notin S_0'\cup V_{F,\infty}}' G(F_v)$ by the characteristic function of $\prod_{v\in S_0'\setminus S_0} G(\oc_v)$.
		\end{enumerate}
		
		Then for any $\la\in i\ago_M^{L,*}$, we have:

		\begin{align}\label{eq:comparison}
		\sum_{h\in\hc^\circ} J_{\Pi}^h(\la, f^h)=|S_\Pi|^{-1}  I_{P,\pi}(\la,f).
		\end{align}
	\end{theoreme}
	
	\begin{remarque} If $\Pi$ is a discrete Arthur parameter (that is $L=M$), the statement  reduces to \cite[proposition 10.1.6.1]{BCZ}. As we observed in \cite{BCZ}, if the assumptions hold for the set $S_0$, they also hold for any large enough finite set containing $S_0$.
	\end{remarque}
	
\end{paragr}

\begin{paragr}[Proof of theorem \ref{thm:comparison}.]  --- As in \cite{BPLZZ} and \cite{BCZ}, our proof is based on the global comparison of Jacquet-Rallis trace formulas and the use of multipliers to isolate some spectral contributions. However the spectral contributions we consider here are  continuous in nature and we need further considerations.

	In \cite[Theorem 3.2.4.1]{BCZ}, we defined a distribution $I$ on $\Sc(G(\AAA))$: this is the ``Jacquet-Rallis trace formula'' for $G$. For unitary groups, we have an analogous distribution  $J^h=J^{U_h}$ on  $\Sc(U_h(\AAA))$ for each $h\in \hc$ (see theorem \ref{thm:jfDefU}).  By \cite[théorème 1.6.1.1]{CZ}, we have for functions $f$ and $f^h$ as in the statement
		\begin{align}\label{eq:thm:CZ}
		I_{}(f)=\sum_{h\in\hc^\circ}  J_{}^h(f^h).
		\end{align}

In the following, for each finite place outside $S_0$,  we fix open compact subgroups $K_v'\subset K_v$ and $K_{h,v}'\subset K_{h,v}$ of finite index. We set
\begin{align*}
  K^\infty_1=\prod_{v\in S_0\setminus V_{F,\infty} } K_v'   \prod_{v\notin S_0} K_v\text{ and } K^\infty_{h,1}=\prod_{v\in S_0\setminus V_{F,\infty} } K_{h,v}' \prod_{v\notin S_0}K_{h,v}.
\end{align*}

We denote by $\Sc(G(\AAA),K_1^\infty)\subset  \Sc^\circ(G(\AAA))$, resp.  $\Sc(U_h(\AAA),K_{h,1}^\infty)\subset \Sc^\circ(U_h(\AAA))$, the subalgebras of bi-$K_1^\infty$ (resp. bi-$K_{h,1}^\infty$) invariant functions. Since we can shrink $ K^\infty_1$ and $K^\infty_{h,1}$ if necessary, it suffices to prove the theorem for functions in these subalgebras.

	We denote by $\mc^{S_0'}(G(\AAA))$, resp.  $\mc^{S_0'}(U_h(\AAA))$, the algebra of $S_0'$-multipliers defined in \cite[definition 3.5]{BPLZZ} relatively to the subgroup $\prod_{v\notin S_0'} K_v$, resp. $\prod_{v\notin S_0'} K_{h,v}$. Any multiplier $\mu\in \mc^{S_0'}(G(\AAA))$, resp. $\mu\in \mc^{S_0'}(U_h(\AAA))$, gives rise to a linear operator $\mu \ast$ of the algebra  $\Sc(G(\AAA),K_1^\infty)$, resp. $\Sc(U_h(\AAA),K_{h,1}^\infty)$. For every admissible irreducible representation $\rho$ of $G(\AAA)$, resp. of $U_h(\AAA)$, there exists a complex number $\mu(\rho)$ such that $\rho(\mu \ast f)=\mu(\rho) \rho(f)$ for all $f\in \Sc^\circ(G(\AAA),K_1^\infty)$, resp. $f\in \Sc^\circ(U_h(\AAA),K_{h,1}^\infty)$. Note that $\mu(\rho)$ depends only on the infinitesimal character of the archimedean component of $\rho$ and on the components outside $S_0'$. Moreover if  $Q$ is a standard parabolic  subgroup of $G$, resp. $U_h$, and if $\rho$ is an admissible irreducible representation of the Levi component $M_Q(\AAA)$, then we have for any $\la\in \ago_{Q,\CC}^*$ and $R_\la=\Ind_{Q}^G(\rho_\la)$, resp. $\Ind_{Q}^{U_h}(\rho_\la)$,
        \begin{align*}
          R_\la (\mu \ast  f)=\mu( R_\la)  R_\la( f)
        \end{align*}
where $\mu( R_\la) \in \CC$ and the map $\la\mapsto \mu( R_\la) $ is holomorphic.

By \cite{GRS}, there exists at least one $h^\flat\in\hc^\circ$ such that the set $\Xgo_\Pi^{h^\flat}$ is nonempty. We fix such a form $h^\flat$ and a pair $(M_{h^\flat},\sigma^\flat)$. The following lemmas are  based on  \cite[theorem 3.17]{BPLZZ} and a strong multiplicity one theorem of Ramakrishnan.

\begin{lemme}\label{lem:multU}
  Let $h\in \hc^\circ$ and $\la\in i\ago_{M}^{L,*}$ in general  position. Then there exists a multiplier $\mu_h\in \mc^{S_0'}(U_h(\AAA))$ such that
  \begin{enumerate}
  \item For all $f^h\in \Sc(U_h(\AAA),K_{h,1}^\infty)$, the right convolution by $\mu_h*f^h$ sends $L^2([U_h])$ into
    \begin{align*}
     \widehat{ \oplus}_{\chi\in \Xgo^h_\Pi} L^2_\chi([U_h]).
    \end{align*}
\item For $(M_h,\sigma)\in \Xgo_{\Pi}^h$,    we have $\mu_h(\Ind_{P_h}^{U_h}(\sigma_{bc(\la)}))=1$.
  \end{enumerate}
\end{lemme}

\begin{preuve} Let $\la\in i\ago_{M}^{L,*}$ in general  position. Let  $h\in \hc^\circ$ and let $\Xgo^h_\la$ be the set of cuspidal data  of $U_h$ represented by pairs $(M_Q,\rho)$ with $Q\subset U_h$ a standard parabolic subgroup for which there exists $\la'\in \ago_{Q,\CC}^*$ such that
\begin{itemize}
\item[(A)] the Archimedean infinitesimal characters of  $\Ind_{P_{h^\flat}}^{U_{h^\flat}}(\sigma_{bc(\la)}^\flat)$ and $\Ind_{Q}^{U_h}(\rho_{\la'})$ are the same;
\item[(B)] for all  finite places $v$ outside $S_0'$, the irreducible representation $\Ind_{P_{h^\flat}(F_v)}^{U_{h^\flat}(F_v)}(\sigma_{bc(\la),v}^\flat)$ is a constituent of $\Ind_{Q(F_v)}^{U_h(F_v)}(\rho_{\la',v})$.
\end{itemize}

By   \cite[theorem 3.19]{BPLZZ}, there exists a $S_0'$-multiplier  $\mu_h$ such that
 \begin{itemize}
  \item[(C)] For all $f^h\in \Sc(U_h(\AAA),K_{h,1}^\infty)$, the right convolution by $\mu_h*f^h$ sends $L^2([U_h])$ into
    \begin{align*}
     \widehat{ \oplus}_{\chi\in \Xgo^h_\la} L^2_\chi([U_h]).
    \end{align*}
\item[(D)] $\mu_h(\Ind_{P_{h^\flat}}^{U_{h^\flat}}(\sigma_{bc(\la)}^\flat))=1$.
\end{itemize}

Let $\la'\in \ago_{Q,\CC}^*$  and $(M_Q,\rho)$ be a representative of an element $\Xgo^h_\la$ such that conditions (A) and (B) hold. We can write the Levi factor $M_Q$ of $Q$ as a product $G^\sharp\times U_0$ where $G^\sharp$ is a product of linear groups and $U_0$ is a product of two unitary  groups. Accordingly we write $\rho=\rho^\sharp\times \rho_0$. Let  $Q_E=\Res_{E/F}(Q\times_F E)$ and let $\la_E'\in \ago_{Q_E}^*$ be the base change of $\la'$. Observe that $\Pi_v$ is generic for all places $v$ and  is the split base of the representation $\Ind_{P_{h^\flat}(F_v)}^{U_{h^\flat}(F_v)}(\sigma_{bc(\la),v}^\flat)$  for finite $v\notin S_0'$. So the representation   $\Ind_{P_{h^\flat}(F_v)}^{U_{h^\flat}(F_v)}(\sigma_{bc(\la),v}^\flat)$ is also generic and thus $\Ind_{Q(F_v)}^{U_h(F_v)}(\rho_{\la',v})$ and $\rho_{\la',v}$ are also generic for finite $v\notin S_0'$ . By \cite[Theorem 4.14 (1)]{BPLZZ}, there exists an isobaric automorphic representation $\rho_E$ of $M_{Q_E}$ such that the split base change of $\rho$  is  $\rho_E$ at almost all split places. By the strong multiplicity one theorem of Ramakrishnan (see \cite{Ram2}), we deduce that  the isobaric automorphic representations of $G(\AAA)$ associated to $(P,\pi\otimes {\la})$ and $(Q_E,\rho_{E}\otimes \la'_E)$  are isomorphic.  Using \cite[theorem 4.4]{JS2} and the fact that $\la$ is in general position, we  conclude that, up to a change of representative, we have $P \subset Q_E$, the inclusion given by base change $\ago_Q^*\subset\ago_{Q_E}^* $ induces an isomorphism of $\ago_Q^*$ onto $\ago_{M}^{L,*}$ which identifies $\la'$ to $bc(\la)$,  the representation $\Pi_0$ is the weak base change of $\rho_0$ and $\Ind_{Q_E}^G(\rho^\sharp\boxtimes \Pi_0 \boxtimes \rho^{\sharp,*})=\Pi$. Using condition (B), we deduce that $\Pi_{0,v}$ is also the split base change of $\rho_{0,v}$ for all finite places $v$ outside $S_0'$. We get that  $\Xgo^h_\la\subset  \Xgo^h_\Pi$ and so (C) implies assertion 1.

Still we have to check assertion 2. Let $(M_{Q},\rho)$ be a representative of an element in $ \Xgo^h_\Pi$. We claim that $\Ind_{Q}^{U_h}(\rho)$ and $\Ind_{P_{h^\flat}}^{U_{h^\flat}}(\sigma_{bc(\la)}^\flat)$ have the same Archimedean infinitesimal character and have the same local component for  all finite places $v$ outside $S_0'$. The latter condition follows directly from the definition of the set $ \Xgo^h_\Pi$ and the fact that the split base is injective. The former condition follows from \cite[Theorem 4.14 (4)]{BPLZZ} (applied to $\rho$) and the fact that the base change map is injective at the level of archimedean infinitesimal characters. Since  condition 2 depends only on the components of $\Ind_{Q}^{U_h}(\rho)$ on $V_F\setminus S_0'$, we see that (D) implies assertion 2.
\end{preuve}

\begin{lemme}\label{lem:multG}
  Let $\la\in i\ago_{M}^{L,*}$. Then there exists a multiplier $\mu\in \mc^{S_0'}(G(\AAA))$ such that
  \begin{enumerate}
  \item For all $f\in \Sc(G(\AAA),K_{h,1}^\infty)$, the right convolution by $\mu*f$ sends $L^2(G(F)A_G(\AAA)\back G(\AAA) )$ into
  $L^2_{\chi_0}([G])$.
\item We have    $\mu(\Pi_\la)=1$.
  \end{enumerate}
\end{lemme}

\begin{preuve}
 The proof is similar to (but simpler than) that of lemma \ref{lem:multU} and is based on \cite[Theorem 1.3]{BPLZZ}).
\end{preuve}

Let $\la_0$ be an element of $i\ago_{M}^{L,*}$ in general position. Let $\mu\in \mc^{S_0'}(G(\AAA))$ and $\mu_h\in \mc^{S_0'}(U_h(\AAA))$  for $h\in \hc^\circ$ satisfying assertions of lemmas \ref{lem:multU}  and \ref{lem:multG} for $\la=\la_0$.   We may and we shall also assume that $\mu^h$ is the ``base change'' of $\mu$ (see  \cite[Lemma 4.12]{BPLZZ}). Let $f\in \Sc(G(\AAA),K_1^\infty)$ and for $h\in \hc^\circ$ let $f^h\in  \Sc(U_h(\AAA),K_{h,1}^\infty)$ that satisfy the hypotheses of theorem \ref{thm:comparison}. Then $\mu\ast f$ and $\mu^h\ast f^h$, for $h\in \hc^\circ$ still satisfy the hypotheses (see \cite[Proposition 4.8]{BPLZZ}). Hence, by \eqref{eq:thm:CZ}, we have:
\begin{align}\label{thm:CZ2}
		I_{}(\mu\ast f)=\sum_{h\in\hc^\circ}  J_{}^h(\mu^h\ast f^h).
		\end{align}
It follows from lemma \ref{lem:multG} and \cite[Proposition 3.3.3.1 and Theorem 3.3.9.1]{BCZ}   that we have 	$I_{}(\mu\ast f)=I_{\chi_0}(\mu\ast f)$. In the same way, lemma \ref{lem:multU}  and corollary  \ref{cor:LaTu-cstterm} show that the right-hand side of \eqref{thm:CZ2} reduces to

\begin{align*}
  \sum_{h\in\hc^\circ}  \sum_{\chi\in \Xgo_\Pi^h  } J_{\chi}^h(\mu_h*f^h).
\end{align*}
By theorems \ref{thm:JchiU} and \ref{thm:Ichi} and by an elementary change of variables we get:
\begin{align}
 \nonumber  2^{-\dim(\ago_L)}\int_{i\ago_M^{L,*}} I_{P,\pi}(\la,\mu*f)\, d\la&=\sum_{h\in\hc^\circ} \sum_{(M_h,\sigma)\in  \Xgo_\Pi^h} \int_{i\ago_{M_h}^*} J_{P_h,\sigma}^h(bc(\la),\mu_h*f^h)\, d\la\\
 \label{eq:integ-la}&=2^{-\dim(\ago_M^L)}\sum_{h\in\hc^\circ } \int_{i\ago_{M}^{L,*}} J_{\Pi}^h(\la,\mu_h*f^h)\, d\la.
\end{align}

Let $v_1,v_2$ two finite places outside  $S_0'$     with distinct residual characteristics.  Let  $S_1\subset V_F\setminus S_0'$  be a finite set of finite places. We assume that $S_1$ contains $v_1$ and $v_2$. Let $\Ac_{S_1}$ be the spherical algebra $\otimes_{v\in S_1}\Sc^\circ(G(F_{v}))$. Let $g\in \Ac_{S_1}$. For all $h\in \hc^\circ$, let $g^h=(\otimes_{v\in S_1} BC_{h,v})(g)$ be its base change to $\otimes_{v\in S_1}\Sc^\circ(U_h(F_{v}))$. The assumptions  still hold for  the convolutions $f*g$ and $f^h*g^h$. Note that we have for any $\la\in i\ago_{M}^{L,*}$,
\begin{align*}
  J_{\Pi}^h(\la,\mu_h*(f^h*g^h))= J_{\Pi}^h(\la,\mu_h*f^h) \hat g(\Pi_{\la,S_1}) 
\end{align*}
where $\hat g$ is the Satake transform and $\hat g(\Pi_{\la,S_1})$ is the scalar by which $g$ acts on  $\Pi_{\la,S_1}=\otimes_{v\in S_1} \Pi_{\la,v} $. For all $\la\in i\ago_M^{L,*}$ we  set 
\begin{align*}
  h(\la)=|S_\Pi|^{-1} I_{P,\pi}(\la,\mu*f)  - \sum_{h\in\hc^\circ}J_{\Pi}^h(\la,\mu_h*f^h).
\end{align*}
The equality \eqref{eq:integ-la} implies that for all  $g\in \Ac_{S_1}$ we have
\begin{align*}
    \int_{i\ago_M^{L,*}} \hat g(\Pi_{\la,S_1}) h(\la) \, d\la =0.
\end{align*}
 Let $\mathcal{T}_{S_1}:\la\in i\ago_M^{L,*} \mapsto \Pi_{\la,S_1}$. The map  $\mathcal{T}_{\{v_1,v_2\}}$ has finite fibers of  uniformly bounded cardinality. In particular, the map  $\mathcal{T}_{S_1}$ has the same property with the same bound. By Stone-Weierstrass theorem, the set $\{\hat g \mid g\in \Ac_{S_1}\}$ is dense in the set of continuous fonctions on the unramified unitary dual of $\prod_{v\in S_1} G(F_v)$.The push-forward of the measure $h(\la)\, d\la $ by the map   $\mathcal{T}_{S_1}$   is  thus zero. In this way, we get that for almost all $\la \in i\ago_M^{L,*}$ 
\begin{align}\label{eq:vanish-sum}
  \sum_{\la_1\in \tc^{-1}_{S_1}(\Pi_{\la,S_1})  } h(\la_1) =0
\end{align}
for all finite sets $S_1$ as above.  By continuity of $h$ this equality holds for all $\la \in i\ago_M^{L,*}$ and in particular for $\la=\la_0$.

By the strong multiplicity one theorem of Ramakrishnan (see \cite{Ram2}), we may enlarge  $S_1$ so that for any $\la\in \tc^{-1}_{S_1}(\Pi_{\la_0,S_1})$ we have $\Pi_\la=\Pi_{\la_0}$. By \cite[theorem 4.4]{JS2} and the fact that $(M,\pi)$ is regular, we have even $\pi_\la=\pi_{\la_0}$ and thus $\la=\la_0$. So  \eqref{eq:vanish-sum} reduces to $h(\la_0)=0$. Because $\la_0$ is in general position, we may use condition 2 of  lemmas \ref{lem:multU} and \ref{lem:multG}: we get
\begin{align*}
  h(\la_0)=|S_\Pi|^{-1} I_{P,\pi}(\la_0,f)  - \sum_{h\in\hc^\circ}J_{\Pi}^h(\la_0,f^h).
\end{align*}
So we get \eqref{eq:comparison} for $\la$ in general position and so for  all $\la\in i\ago_M^{L,*}$ since both members of \eqref{eq:comparison}   are  continuous (and even analytic). 
\end{paragr}

\subsection{Proof of Theorem \ref{thm:GGP}}

\begin{paragr} Once we have theorem \ref{thm:comparison}, the proof of theorem \ref{thm:GGP} is very similar to the proof of \cite[theorem 1.1.5.1]{BCZ}. For the reader's convenience, we recall some steps.  We use notations of section \ref{ssec:Globalcomparison}. Let $\Pi=\Ind_P^G(\pi)$ be a $H$-regular Hermitian Arthur parameter of $G$ and let $\la\in i\ago_{M}^{L,*}$.  The relative character  $I_{\Pi_\la}$ defined in \eqref{eq:IPIla} is built upon two linear forms namely $Z^{RS}(0)$ and $\beta_\eta$.  The linear form $\beta_\eta$  is not identically zero (as a consequence of \cite{GelKa}, \cite[Proposition 5]{JacDist} and \cite{Kemar}) and $Z^{RS}(0)$  is nonzero if and only if $L(\frac12,\Pi_\la)\not=0$ (as follows  from the work Jacquet, Piatetski-Shapiro and Shalika \cite{JPSS}, \cite{JacIntegral}). Since by theorem \ref{thm:IPi=Ipi}, we have $I_{\Pi_\la}=I_{P,\pi}(\la)$, we deduce that $I_{P,\pi}(\la)$ is non zero  if and only if $L(\frac12,\Pi_\la)\not=0$.

  Let's consider $h\in \hc$, a parabolic subgroup  $P_h=M_{h}N_{P_h}$ of $U_h$ and a cuspidal \emph{subrepresentation} $\sigma$ of $M_h$. Let $\Ac_{P_h,\sigma_h}(U_h)\subset \Ac_{P_h}(U_h)$  be the space of forms $\varphi\in \Ac_{P}(G)$ such that
$$m\in [M_P]\mapsto \exp(-\langle \rho_P, H_P(m)\rangle) \varphi(mg)$$
belongs to the space of $\sigma$   for every $g\in G(\AAA)$.   By a variation on §\ref{S:rel-charU}, we define for $\mu\in\ago_{M_h}^*$ the relative character
\begin{align*}
  J_{P_h,\sigma}^{U_h}(\mu, f)=\sum_{\varphi\in \bc_{P_h,\sigma}}   \pc_{U'_h}(I_{P_h}(\mu,f)\varphi,\mu)\overline{   \pc_{U'_h}(\varphi,\mu)  }
  \end{align*}
where $\bc_{P_h,\sigma}$ is a $K$-basis of $\Ac_{P_h,\sigma_h}(U_h)$ that is it is the union over of $\tau\in \hat K_h$ of orthonormal bases $\bc_{P_h,\sigma,\tau}$  for the Petersson inner product of  the finite dimensional  subspaces $\Ac_{P_h,\sigma}(U_h,\tau)$ of functions in $\Ac_{P_h,\sigma}(U_h)$   which transform under $K_h$ according to $\tau$.

  Let's assume that  $\Pi$  is the weak base change of $(P_h, \sigma)$. Then we have the map $bc:a_M^{L,*}\to \ago_{M_h}^*$, see \eqref{eq:aMh-aML}. It is clear that the distribution  $J_{P_h,\sigma}^h(bc(\la))$ is non-zero if and only if the period integral  $\pc_{U'_h}(\cdot,bc(\la))$ induces a non-zero linear form on the space of $\sigma$. Then theorem \ref{thm:GGP} reduces  to the equivalence between the two assertions:
	\begin{enumerate}
		\item[(A)] The distribution $I_{P,\pi}(\la)$  is non-zero.
		\item[(B)] There exist $h\in \hc$, a parabolic subgroup  $P_h=M_{h}N_{P_h}$ of $U_h$ and a cuspidal subrepresentation $\sigma$ of $M_h$   such that $\Pi$  is the weak base change of $(P_h, \sigma)$ and $J_{P_h,\sigma}^h(bc(\la),f) \not=0$.
			\end{enumerate}
	
                      \end{paragr}
                      
\begin{paragr}[Proof of $(A)\Rightarrow (B)$.] --- Let  $S_0$ be the finite set of  §\ref{S:finitesetS0} such that $I_\Pi$ is not identically zero on $\Sc^\circ(G(\AAA))$. Then by  results of \cite{Xue} towards the archimedean transfert  and the existence of $p$-adic transfer of \cite{Z1}, we see that there exist functions $f$ and $f^h$ for $h\in \hc^\circ$ satisfying the hypotheses of theorem  \ref{thm:comparison} and such that $I_{P,\pi}(\la,f)\not=0$. Then $(B)$ follows from \eqref{eq:comparison}.	
\end{paragr}

\begin{paragr}[Proof of $(B) \Rightarrow (A)$.] ---   We may choose the set $S_0$ so that there exist $h_0\in \hc^\circ$, a parabolic subgroup  $P_{h_0}=M_{h_0}N_{P_{h_0}}\subset U_{h_0}$, a cuspidal subrepresentation $\sigma_0$ of $M_{h_0}(\AAA)$ and   a function  $\xi \in \Sc^\circ(U_{h_0}(\AAA))$ such that
  \begin{itemize}
  \item the class of $(M_{h_0},\sigma_0)$ belongs to $\Xgo_\Pi^{h_0}$;
  \item $J_{P_{h_0},\sigma_0}^{h_0}(bc_{}(\la),\xi) \not=0$.
  \end{itemize}

  We set $f_0^{h_0}=\xi*\xi^\vee$ where $\xi^\vee(g)=\overline{\xi(g^{-1})}$. Then (see \cite[p. 993]{Z1}) we have  $J_{P,\sigma}^{h_0}(bc(\la),f^{h_0}_0)\geqslant 0$ for all pairs $(P,\sigma)$ whose class belongs to $\Xgo_\Pi^{h_0}$.  Moreover $J_{P_{h_0},\sigma_0}^{h_0}(bc_{h_0}(\la),f^{h_0}_0)>0$.   For any $h\in \hc^\circ$ such that $h\not=h_0$ we set $f^h_0=0$. Up to enlarging $S_0$, we may and shall assume that the family $(f^h_0)_{h\in \hc^\circ}$ satisfies conditions 2  and 5  of theorem \ref{thm:comparison}. The left hand side of \eqref{eq:comparison} for the family $(f^h_0)_{h\in \hc^\circ}$ is nonzero. By the existence of transfer in \cite{Z1} and the results towards archimedean transfer in \cite{Xue},  we can find test functions $f\in \Sc^\circ(G(\AAA))$ and $f^h\in \Sc^\circ(U_h(\AAA))$, for $h\in \hc^\circ$, satisfying all the conditions of theorem \ref{thm:comparison} and such that the left hand side of \eqref{eq:comparison} is nonzero.  Assertion (A) is then clear.
  \end{paragr}

\subsection{Proof of Theorem \ref{thm:II}}

\begin{paragr}\label{S:reformulation II}
	Let $h\in \hc$. Let  $P=M_PN_P$ be a parabolic subgroup of $U_h$ and $\sigma$ be a cuspidal automorphic subrepresentation of $M_P(\AAA)$ which is tempered everywhere. Then the group $\Res_{E/F}(P\times_F E)$ obtained by extension to $E$ and restriction of scalars to $F$ can be identified to a parabolic subgroup $Q=M_QN_Q$ of $G$. Then by \cite{Mok}, \cite{KMSW}, $\sigma$ admits a strong base-change $\pi$ to $M_Q$ namely for every place $v$ of $F$, the local base-change of $\sigma_{v}$ (defined in \cite{Mok} and \cite{KMSW}) coincides with $\pi_v$. Let $\Pi=\Ind_Q^G(\pi)$. It follows that $\pi$ and $\Pi$ are also tempered everywhere. We assume that $\Pi$ is a $H$-regular hermitian Arthur parameter. As in §\ref{S:weakBC}, we shall not distinguish in the notation the spaces $\ago_\Pi^*$ and $\ago_P^*$. Let $\la\in i\ago_\Pi^*$.
	
	We choose a finite set of places $S_0$ as in \S \ref{S:finitesetS0} such that $h\in \hc^\circ$ and $\sigma$ as well as the additive character $\psi'$ used to normalize local Haar measures in §\ref{S:Haar} are unramified outside of $S_0$. 
	
We have a decomposition $\sigma=\otimes_{v\in V_F}\sigma_v$. Let $\Sigma_{\la,v}$ be the full induced representation $\Ind_{P(F_v)}^{U_h(F_v)}(\sigma_v\otimes \la)$. Let $v\in V_F$.  We define a distribution $J_{\Sigma_{\la,v}}$ on $\Sc(U_{h}(F_v))$ by
        \begin{align*}
J_{\Sigma_{\la,v}}(f^{h}_v)=\int_{U'_{h}(F_v)} \mathrm{Trace}(\Sigma_{\la,v}(h_v)\Sigma_{\la,v}(f^{h}_v)) dh_v,\;\; f^{h}_v\in \Sc(U_{h}(F_v)),
        \end{align*}
        where
        \begin{align*}
          \displaystyle \Sigma_{\la,v}(f^{h}_v)=\int_{U_{h}(F_v)} f^{h}_v(g_v) \Sigma_{\la,v}(g_v)dg_v.
        \end{align*}
By \cite{NHar},  since the representations $\Sigma_{\la,v}$ are all tempered, the expression defining $J_{\Sigma_{\la,v}}$ is absolutely convergent and for every $v\notin S_0$ we have
	$$\displaystyle J_{\Sigma_{\la,v}}(\mathbf{1}_{U_{h}(\oc_v)})=\Delta_{U'_{h},v}^{-2} \frac{L(\frac12,\Pi_{\la,v})}{L(1,\Pi_{\la,v},\As')}.$$

	If there exists a place $v\in S_0$ such that $\sigma_{v}$ does not support any nonzero continuous $U'_{h}(F_v)$-invariant functional, both sides of \eqref{intro-eq:II-E} in  theorem \ref{thm:II} are clearly automatically zero. So we shall assume that for every $v\in S_0$, the local representation $\sigma_{v}$ supports a nonzero continuous $U'_{h}(F_v)$-invariant functional. Then the semi-local  distribution $\prod_{v\in S_0} J_{\Sigma_{\la, v}}(f_v^{h})$ does not vanish identically by \cite[theorem 8.2.1]{BP3}. According to our choice of local measure, theorem \ref{thm:II} is then equivalent to the following assertion: for all factorizable test function $f^{h}\in \Sc(U_{h_0}(\AAA))$ of the form $f^{h}=(\Delta_{U'_{h}}^{S_0})^2\prod_{v\in S_0} f^{h}_v\times \prod_{v\notin S_0} \mathbf{1}_{U_{h}(\oc_v)}$, we have
	\begin{align}\label{eq:reformulation II}
	\displaystyle J^{h}_{P,\sigma}(\la,f^{h})=\lvert S_\Pi\rvert^{-1} \frac{L^{S_0}(\frac12,\Pi_\la)}{L^{S_0}(1,\Pi_\la,\As')}\prod_{v\in S_0} J_{\Sigma_{\la, v}}(f_v^{h}).
	\end{align}
\end{paragr}

\begin{paragr} By theorem \ref{thm:IPi=Ipi} and the definition of $I_{\Pi_\la}$ there,  for every factorizable test function $f\in \Sc(G(\AAA))$ of the form $f=\Delta_H^{S_0,*}\Delta_{G'}^{S_0,*}\prod_{v\in S_0} f_v \times \prod_{v\notin S_0} \mathbf{1}_{G(\oc_v)}$, we have
	\begin{align}\label{eq:factorization RS}
	I_{Q,\pi}(\la,f)=\frac{L^{S_0}(\frac12,\Pi_\la) }{ L^{S_0,*}(1,\Pi_\la,\As')  }     \prod_{v\in S_0} I_{\Pi_{\la,v}}(f_v).
	\end{align}
where for for every place $v\in S_0$, we introduce the local relative character $I_{\Pi_{\la,v}}$ defined by
        \begin{align*}
I_{\Pi_{\la,v}}(f_v)=\sum_{W_v} \frac{\al_v(\Pi_v(f_v)W_v) \overline{\beta_{\eta,v}(W_v)}}{\langle W_v,W_v\rangle_{\Whitt, v}},\;\; f_v\in \Sc(G(F_v)).
        \end{align*}
Here the sum runs over a $K_v$-basis of the Whittaker model $\wc(\Pi_{\la,v},\psi_{N,v})$  and $\lambda_v$, $\beta_{\eta,v}$, $\langle .,.\rangle_{\Whitt,v}$ are  given by
\begin{align*}
  &   \al_v(W_v)=\int_{N_H(F_v)\backslash H(F_v)} W_v(h_v)dh_v \\
  & \beta_{\eta,v}(W_v)=\int_{N'(F_v)\backslash \pc'(F_v)} W_v(p_v) \eta_{G',v}(p_v) dp_v, \\
  & \langle W_v,W_v\rangle_{\Whitt,v}=\int_{N(F_v)\backslash \pc(F_v)} \lvert W_v(p_v)\rvert^2 dp_v.
\end{align*}
The above expressions, especially $\al_v(W_v)$, are all absolutely convergent due to the fact that $\Pi_{\la,v}$ is tempered (see \cite[Proposition 8.4]{JPSS}). The above definition implicitely depends on the choice of an additive character $\psi$ of $\AAA_E/E$ trivial on $\AAA$ which,  up to enlarging $S_0$, we may assume to be  unramified outside of $S_0$. 
\end{paragr}

\begin{paragr} 	Let $f^{h}$ be a test function as in \eqref{eq:reformulation II}. Since both sides of \eqref{eq:reformulation II} are continuous functionals in $f^{h}_v$ for $v\in S_0$, we may assume that the function $f^{h}_v$ admits a transfer $f_v\in \Sc(G(F_v))$ for every $v\in S_0$ using results of \cite{Xue}  and \cite{Z1}. Moreover, by the results of those references, we may also assume that for every $h'\in \hc^\circ$ with $h'\neq h$, the zero function on $U_{h'}(F_{S_0})$ is a transfer of $f_{S_0}=\prod_{v\in S_0} f_v$. We set $f=\Delta_H^{S_0,*}\Delta_{G'}^{S_0,*}f_{S_0} \times \prod_{v\notin S_0} \mathbf{1}_{G(\oc_v)}$. Then, setting $f^{h'}=0$ for every $h'\in \hc^\circ \setminus \{h \}$, the functions $f$ and $(f^{h'})_{h'\in \hc^\circ}$ satisfy the assumptions of theorem \ref{thm:comparison}. Therefore, we have
	\begin{align}\label{eq:comparison II}
	 J_{\Pi}^{h}(\la,f^{h})=|S_\Pi|^{-1}I_{Q,\pi}(\la,f).
	\end{align}
\end{paragr}

\begin{paragr} By  the local Gan-Gross-Prasad conjecture \cite{BP3}, the classification of cuspidal automorphic representations of $U_{h}$ in terms of local $L$-packets \cite{Mok}, \cite{KMSW},  all the terms in the definition of $J_{\Pi}^{h}(\la,f^{h})$, see  \eqref{eq:JPih}, vanish   except possibly $J_{P,\sigma}^{h}(\la, f^{h})$. So \eqref{eq:comparison II} reduces to
  \begin{align}\label{eq:comparison III}
    J_{P,\sigma}^{h}(\la,f^{h})=|S_\Pi|^{-1}I_{Q,\pi}(\la,f).
  \end{align}
By \cite[Theorem 5.4.1]{BPPlanch} and since $\Pi_v$ is the local base-change of $\sigma_{v}$, there are explicit constants $\kappa_v\in \CC^\times$ for $v\in S_0$ satisfying $\prod_{v\in S_0} \kappa_v=1$ and such that
	\begin{align}\label{eq:comparison local}
 I_{\Pi_{\la,v}}(f_v)=\kappa_v J_{\Sigma_{\la,v}}(f^{h}_v)
	\end{align}
	for every $v\in S_0$. Now  \eqref{eq:reformulation II}  results from the combination of \eqref{eq:comparison III}, the factorization \eqref{eq:factorization RS} and the local comparison \eqref{eq:comparison local}.
\end{paragr}

%% file: Besselcomputation.tex
	\section{Application to the Gan-Gross-Prasad and Ichino-Ikeda conjectures for Bessel periods}\label{sect: application}
	
	General notation (probablement d\'ej\`a introduites dans le reste de l'article): $E/F$ is a quadratic extension of number fields. For every place $v$ of $F$ we set $E_v=F_v\otimes_F E$ and we let $\cO_{E_v}$ be the ``ring of integers'' in the quadratic \'etale extension $E_v$ of $F_v$. We denote by $\AAA=\prod_v' F_v$ is the ring of adeles of $F$, $\AAA_E=\prod_v' E_v$ that of $E$.  Also for any linear algebraic group $G$ defined over $F$ we set $[G]=G(F)\backslash G(\AAA)$ and we fix a {\it height} function $\lVert .\rVert: G(\AAA)\to \RR_{\geqslant 1}$ as in \cite[\S 2.4.1]{BCZ} and we set
	$$\displaystyle \lVert x\rVert_G=\inf_{\gamma\in G(F)} \lVert \gamma x\rVert,\; \mbox{ for } x\in G(\AAA),$$
	$$\displaystyle \sigma(x)=1+\log \lVert x\rVert,\; \mbox{ for } x\in G(\AAA).$$
	
	If $V$ is a finite dimensional vector space over $F$, we also fix a height $\lVert .\rVert_{V_{\AAA}}: V_{\AAA}=V\otimes_F \AAA\to \RR_{\geqslant 1}$ as in \cite[\S 2.4.2]{BCZ} and for every place $v$ of $F$, we denote by $\lVert .\rVert_{V_v}$ the restriction of $\lVert .\rVert_{V_{\AAA}}$ to $V_v:=V\otimes_F F_v$.
	
	\subsection{Groups}

        \begin{paragr}[Notations.] --- They are as in sections \ref{sec:prelim}, \ref{sec:spectral-unitary} and \ref{sec:GH-contrib-JR}. In particular,  $E/F$ is a quadratic extension of number fields. For every place $v$ of $F$ we set $E_v=F_v\otimes_F E$ and we let $\cO_{E_v}$ be the ``ring of integers'' in the quadratic \'etale extension $E_v$ of $F_v$. We denote by $\AAA_E=\prod_v' E_v$ the ring of adeles of $E$. We fix a non-trivial character $\psi: \AAA_E/E\to \CC^\times$ that is trivial on $\AAA$. We also set, see §\ref{S:hauteurs},
          \begin{align*}
            \sigma(x)=1+\log \lVert x\rVert,\; \mbox{ for } x\in G(\AAA).
          \end{align*}
	If $V$ is a finite dimensional vector space over $F$, we also fix a height $\lVert .\rVert_{V_{\AAA}}: V_{\AAA}=V\otimes_F \AAA\to \RR_{\geqslant 1}$ as in \cite[\S 2.4.2]{BCZ} and for every place $v$ of $F$, we denote by $\lVert .\rVert_{V_v}$ the restriction of $\lVert .\rVert_{V_{\AAA}}$ to $V_v:=V\otimes_F F_v$.
          
        \end{paragr}
	
        \begin{paragr}[Linear groups.] --- For every $k\geqslant 0$ we set $G_k=\Res_{E/F}\GL_{k,E}$ equipped with the pair $(B_k,T_k)$ as in §\ref{S:Gn}. Let $N_k\subset B_k$ be the unipotent radical of $ B_k$. We define two generic characters $\psi_k,\psi_{-k}:[N_k]\to \CC^\times$ by
	\begin{equation*}
		\displaystyle \psi_{k}(u)=\psi(\sum_{i=1}^{k-1}u_{i,i+1}).
	\end{equation*}
	\begin{equation*}
		\displaystyle \psi_{-k}(u)=\psi(-\sum_{i=1}^{k-1}u_{i,i+1}),\;\; u\in [N_k].
	\end{equation*}
	We also let $P_k\subset G_k$ be the {\it mirabolic} subgroup consisting of matrices with last row $(0,\ldots,0,1)$ and $K_k=\prod_{v}K_{k,v}$ be the standard maximal compact subgroup of $G_k(\AAA)$.
      \end{paragr}

        \begin{paragr}[Unitary groups.] ---	
	We define the Hermitian form $h_s\in \cH_2$ by $h_s=h_0\oplus -h_0$. We also fix a basis $(x,y)$ of $E^2$ consisting of isotropic vectors for $h_s$ (i.e. $h_s(x,x)=h_s(y,y)=0$) such that $h_s(x,y)=1$ and we set $X=Ex$, $Y=Ey$.
	
	In this chapter we fix two positive integers $n\geqslant m$ of the same parity. Thus $n=m+2r$ for some $r\geqslant 0$. Let $h_m\in \cH_m$. We set $h_{m+1}=h_m\oplus h_0$, $h_n=h_s^{\oplus r}\oplus h_m$, $h_{n+1}=h_n\oplus h_0=h_s^{\oplus r}\oplus h_{m+1}$ and
	$$\displaystyle U_m=U(h_m),\; U_{m+1}=U(h_{m+1}),\; U_n=U(h_n),\; U_{n+1}=U(h_{n+1}).$$
	Note that we have natural inclusions $U_m\hookrightarrow U_{m+1}\hookrightarrow U_n\hookrightarrow U_{n+1}$. We also define the following products
	\begin{equation}\label{product unitary groups}
		\displaystyle \mathcal{G}=U_m\times U_{n+1},\;\; U=U_n\times U_{n+1},\;\; \widetilde{U}=U_{m}\times U_{m+1}.
	\end{equation}
	
	For every $0\leqslant i\leqslant r$, we let
	$$\displaystyle x_i=(\underbrace{0,\ldots,0}_{2i-2},x,\underbrace{0,\ldots,0}_{1+m+2r-2i}),\; y_i=(\underbrace{0,\ldots,0}_{2i-2},y,\underbrace{0,\ldots,0}_{1+m+2r-2i}),\; v_0=(\underbrace{0,\ldots,0}_{2r+m},1)\in E^{n+1}$$
	and we set $X_i=\langle x_1,\ldots,x_i\rangle$. Then, $0=X_0\subsetneq X_1\ldots \subsetneq X_r$ is a flag of $h_{n+1}$-isotropic subspaces and we let $P\subset U_{n+1}$ be the parabolic subgroup stabilizing it. Let $N$ be the unipotent radical of $P$. We define a character $\psi_N:[N]\to \CC^\times$ by
	$$\displaystyle \psi_N(u)=\psi\left(\sum_{i=1}^{r-1} h_{n+1}(ux_{i+1},y_i)+h_{n+1}(uv_0,y_r) \right),\;\;\; u\in [N].$$
	The subgroup $U_m\subset U_{n+1}$ normalizes $N$ and $\psi_N$ is invariant by $U_m(\AAA)$-conjugation. We define the following three subgroups of the products \eqref{product unitary groups}
	\begin{equation}\label{Bessel subgroups}
		\displaystyle \mathcal{B}=U_m\ltimes N,\;\; U'=U_n,\;\;, \widetilde{U}'=U_{m}
	\end{equation}
	where the embedding $\cB\subset \cG$ is the product of inclusion $\cB\subset U_{n+1}$ and the natural projection $\cB\to U_m$ whereas the embeddings $U'\subset U$, $\widetilde{U}'\subset \widetilde{U}$ are the diagonal ones. We also let $\psi_{\mathcal{B}}: [\mathcal{B}]\to \CC^\times$ be the character that coincides with $\psi_N$ on $[N]$ and is trivial on $[U_m]$.
	
	We will also need some auxiliary parabolic subgroups. First, we let $P_n=P\cap U_n$ be the stabilizer of the flag $X_0\subsetneq \ldots\subsetneq X_r$ in $U_n$ and we will also write $P'$ for $P_n$ when we consider it as a subgroup of $U'$. We also define $Q_n\subset U_n$ to be the parabolic subgroup stabilizing the $h_n$-isotropic subspace $X_r$, we write $V_n$ for its unipotent radical and we let $L_n\subset Q_n$ be the Levi factor stabilizing $Y_r$. Note that, as $X_r$ comes with a basis (namely $(x_1,\ldots,x_r)$), we have natural identifications
	\begin{equation}\label{decomp Ln}
		\displaystyle L_n=R_{E/F}GL(X_r)\times U(h_m)=G_r\times U_m.
	\end{equation}
	We will also use the natural identification $a_{Q_n}^*\simeq a_{G_r}^*\simeq \CC$ sending $s\in \CC$ to the unramified character $(g_r,h)\in L_n(\AAA)\mapsto \lvert \det(g_r)\rvert_{\AAA_E}^s$.
	
	We also set $Q=Q_n\times U_{n+1}$ (a maximal parabolic subgroup of $U$) and
	\begin{equation*}
		\displaystyle L=L_n\times U_{n+1}=G_r\times \cG
	\end{equation*}
	(a Levi factor of $Q$). We define the following subgroup of $L$
	\begin{equation*}
		\displaystyle \cB^L=N_r\times \cB
	\end{equation*}
	as well as the character $\psi^L_{\cB}=\psi_{-r}\boxtimes \psi_{\cB}$ of $[\cB^L]$. We also set
	\begin{equation*}
		\displaystyle \cB'=U'\cap \cB=(N_r\times U_m)\ltimes V_n.
	\end{equation*}
	
	We fix a finite set $S$ of places of $F$ containing the Archimedean places as well a the places dividing $2$ and such that the character $\psi$ and the Hermitian form $h_m$ are both unramified outside $S$ i.e. there exists a lattice $\Lambda_m\subset E^m$ such that $\Lambda_{m,v}$ is self-dual with respect to $h_m$ for every $v\notin S$. Then, the same holds for the lattices $\Lambda_{m+1}=\Lambda_{m}\oplus \cO_E$, $\Lambda_n=\cO_E^{\oplus 2r}\oplus \Lambda_m$ and $\Lambda_{n+1}=\Lambda_n\oplus \cO_E$ with respect to the Hermitian forms $h_{m+1}$, $h_n$ and $h_{n+1}$ respectively.
	
	For $\ell\in \{m,m+1,n,n+1 \}$ we fix a maximal compact subgroup $K^U_{\ell}=\prod_v K^U_{\ell,v}\subset U_{\ell}(\AAA)$ such that for every $v\notin S$, $K^U_{\ell,v}$ is the stabilizer of the lattice $\Lambda_{\ell}\otimes_{\cO_E} \cO_{E_v}$.
      \end{paragr}

      \subsection{Measures}
	
      \begin{paragr}
        	For every linear algebraic group $\mathbb{G}$ defined over $F$, we have equipped $\mathbb{G}(\AAA)$ with its (left) Tamagawa measure $dg$, see §\ref{S:Haar}. Also, for each $\mathbb{G}\in \{U_m,U_n,U_{n+1},U_{n+2}, G_r, V_n \}$, we fix a factorization $dg=\prod_v dg_v$ into local Haar measures giving $\mathbb{G}(\cO_v)$ measure one for almost all $v$. In the case where $\mathbb{G}=N_k$ it will be convenient to fix more precisely the local measures as follows: for every place $v$ of $F$, let $d_{\psi_v}x_v$ be the Haar measure on $E_v$ that is self-dual with respect to $\psi_v$ then we equip $N_k(F_v)$ with the product measure
	$$\displaystyle du_v=\prod_{1\leqslant i<j\leqslant k} d_{\psi_v} u_{v,i,j}.$$
	It is well-known, and easy to check, that $du=\prod_v du_v$ is indeed the global Tamagawa measure on $N_k(\AAA)$ i.e. it gives $[N_k]$ volume one.

	Similarly, for every $v$ there is another natural measure on $G_k(F_v)=\GL_k(E_v)$ defined by
	$$\displaystyle d_{\psi_v}g_v=\frac{\prod_{1\leqslant i,j\leqslant k} d_{\psi_v} g_{v,i,j}}{\lvert \det(g_v)\rvert_{E_v}^k}.$$
	We will denote by $\nu(G_{k,v})\in \RR_{>0}$ the quotient $dg_v(d_{\psi_v}g_v)^{-1}$ between the Haar measure we have fixed on $G_r(F_v)$ and the above one. Set
	$$\displaystyle \Delta^*_{G_k}=\zeta^*_E(1)\prod_{i=2}^k \zeta_{E}(i),$$
	where $\zeta_E(s)$ denotes the completed Dedekind Zeta function of $E$ and $\zeta_E^*(1)$ its residue at $s=1$. Similarly for every place $v$ (resp. every finite set $T$ of places) we set
	$$\displaystyle \Delta_{G_k,v}=\prod_{i=1}^k \zeta_{E_v}(i) \;\;\; (\mbox{resp. } \Delta^{T,*}_{G_k}=\zeta^{T,*}_E(1)\prod_{i=2}^k \zeta^T_{E}(i))$$
	where $\zeta_{E_v}$ denotes the local Eulerian factor \footnote{More precisely, this is really a product of two such factors when $v$ splits in $E$.}  of $\zeta_E$ at $v$(resp. $\zeta^T_E(s)$ denotes the corresponding partial Dedekind zeta function and $\zeta^{T,*}_E(1)$ is its residue at $s=1$). Then, for every non-Archimedean $v$ where $\psi_v$ is unramified, we have
	\begin{equation*}
		\displaystyle \vol(K_{k,v},d_{\psi_v}g_v)=\Delta_{G_k,v}^{-1}
	\end{equation*}
	and the global Tamagawa measure on $G_k(\AAA)$ is, by definition,
	\begin{equation*}
		\displaystyle dg=(\Delta^*_{G_k})^{-1}\prod_v \Delta_{G_k,v} d_{\psi_v}g_v.
	\end{equation*}
	In particular, it follows that if $T$ is a sufficiently large finite set of places, we have $\nu(G_{k,v})=\Delta_{G_k,v}$ for $v\notin T$ and
	\begin{equation}\label{eq product nu}
		\displaystyle \prod_{v\in T} \nu(G_{k,v})=(\Delta^{T,*}_{G_k})^{-1}.
	\end{equation}
      \end{paragr}	

      \begin{paragr}	Finally, we record the following Fourier inversion formula: for every $f\in C_c^\infty(P_{k+1}(F_v))$ setting
	$$\displaystyle W_f(g_1,g_2)=\int_{N_{k+1}(F_v)} f(g_1^{-1}u_vg_2) \psi_{k+1}(u_v) du_v,\;\;\; g_1,g_2\in P_{k+1}(F_v),$$
	we have
	\begin{equation*}
		\displaystyle f(p)=\int_{N_k(F_v)\backslash G_k(F_v)} W_f(\gamma,\gamma p) d_{\psi_v}\gamma
	\end{equation*}
	for every $p\in P_{k+1}(F_v)$ where $d_{\psi_v}\gamma$ denotes the quotient of the Haar measure $d_{\psi_v}g_v$ on $G_k(F_v)$ by the Haar measure $du_v=d{\psi_v}u_v$ on $N_k(F_v)$. In particular, replacing $d_{\psi_v}\gamma$ by the quotient measure $d\gamma$ of $dg_v$ by $du_v$ we obtain the following renormalized inversion formula
	\begin{equation}\label{eq1 Fourier inversion}
		\displaystyle f(p)=\nu(G_{k,v})^{-1}\int_{N_k(F_v)\backslash G_k(F_v)} W_f(\gamma,\gamma p) d\gamma.
	\end{equation}
      \end{paragr}

	\subsection{Global periods}
	
        \begin{paragr}
          	We define the {\it Whittaker period} on $G_r$ as the linear form
	$$\displaystyle \cP_{N_r,\psi_{-r}}: \cA([G_r])\ni \phi\mapsto \int_{[N_r]} \phi(u) \psi_{-r}(u) du.$$
	(Note the minus sign.)
      \end{paragr}

      \begin{paragr}      Under extra assumptions ensuring absolute convergence (e.g. cuspidality), we will consider the following {\it global Bessel period} for $\phi\in \cA([\cG])$ (resp. $\phi\in \cA([U])$):
	$$\displaystyle \mathcal{P}_{\mathcal{B},\psi_{\cB}}(\phi)=\int_{[\cB]} \phi(s) \psi_{\cB}(s) ds.$$
	$$\displaystyle (\mbox{resp. } \mathcal{P}_{U'}(\phi)=\int_{[U']} \phi(h)dh.)$$
	For example, it is readily seen that the period $\cP_{U'}(\phi)$ is absolutely convergent for $\phi=\phi_n\otimes \phi_{n+1}$ where $\phi_n\in \cA([U_n])$ and $\phi_{n+1}\in \cA_{cusp}([U_{n+1}])$.
      \end{paragr}

      \begin{paragr}      Again provided it is convergent, for $\phi\in \cA([L])$, we define the {\it mixed Whittaker-Bessel period}
	$$\displaystyle \cP_{\cB^L,\psi^L_{\cB}}(\phi)= \int_{[\cB^L]} \phi(s) \psi^L_{\cB}(s) ds.$$
	For example, this period is absolutely convergent when $\phi=\phi_r\otimes \phi'$ where $\phi_r\in \cA([G_r])$ and $\phi'\in \cA_{cusp}([\cG])$ in which case we have
	$$\displaystyle \cP_{\cB^L,\psi^L_{\cB}}(\phi)=\cP_{N_r,\psi_{-r}}(\phi_r)\cP_{\cB,\psi_{\cB}}(\phi').$$
      \end{paragr}

      \subsection{Local Bessel periods}\label{Sect local periods}
	
        \begin{paragr}	Let $v$ be a place of $F$. For every connected reductive group $\mathbb{G}$ defined over $F_v$, we let $\cC^w(\mathbb{G}(F_v))$ be the space of {\it tempered functions} on $\mathbb{G}(F_v)$ as defined in \cite[\S 1.5]{BP3} (where it is called the {\it weak Harish-Chandra Schwartz space}). Since we will need it, let us recall quickly its definition. Let $\Xi^{\mathbb{G}}$ Harish-Chandra special spherical function on $\mathbb{G}(F_v)$; which strictly speaking depends on the choice of a maximal compact subgroup $\mathbb{K}_v\subset \mathbb{G}(F_v)$ (such a choice has already been done for all the groups we will have to consider). We let $\cC^w_d(\mathbb{G}(F_v))$ be the space of functions $f:\mathbb{G}(F_v)\to \CC$ such that:
	\begin{itemize}
		\item If $v$ is non-Archimedean, $f$ is biinvariant by a compact-open subgroup and we have
		$$\displaystyle \lvert f(g)\rvert\ll \Xi^{\mathbb{G}}(g)\sigma(g)^d, \mbox{ for } g\in \mathbb{G}(F_v);$$
		\item If $v$ is Archimedean, $f$ is $C^\infty$ and for every $X,Y\in \mathcal{U}(\Lie(\mathbb{G}(F_v)))$,
		$$\displaystyle \lvert (L(X)R(Y)f)(g)\rvert\ll_{X,Y} \Xi^{\mathbb{G}}(g)\sigma(g)^d, \mbox{ for } g\in \mathbb{G}(F_v).$$
	\end{itemize}
	When $v$ is Archimedean, $\cC^w_d(\mathbb{G}(F_v))$ is naturally a Fr\'echet space whereas if $v$ is non-Archimedean, $\cC^w_d(\mathbb{G}(F_v))$ is a strict LF space. By definition the space of tempered functions is $\cC^w(\mathbb{G}(F_v))=\bigcup_{d>0} \cC^w_d(\mathbb{G}(F_v))$. It is equipped with the direct limit locally convex topology and it contains $C_c^\infty(\mathbb{G}(F_v))$ as a dense subspace. (However, note that $C_c^\infty(\mathbb{G}(F_v))$ is not dense in $\cC^w_d(\mathbb{G}(F_v))$ for any $d>0$.)
      \end{paragr}

      \begin{paragr}By definition the local period $\cP_{U',v}$ is the linear form
	$$\displaystyle \cP_{U',v}:\cC^w(U(F_v))\ni f\mapsto \int_{U'(F_v)} f(h)dh$$
	the integral being absolutely convergent by \cite[Lemma 6.5.1(i)]{BP3}.
	
	On the other hand, it is shown in \cite[Proposition 7.1.1]{BP3} that the linear form
	$$\displaystyle C_c^\infty(\cG(F_v))\ni f_v\mapsto \int_{\cB(F_v)} f_v(s) \psi_{\cB,v}(s) ds$$
	extends by continuity to $\cC^w(\cG(F_v))$. We denote this unique continuous extension by $\mathcal{P}_{\cB,\psi_{\cB},v}$ and call it the {\it local Bessel period}.
	
	A similar argument shows that the linear forms
	$$\displaystyle C_c^\infty(G_r(F_v))\ni f_v\mapsto \int_{N_r(F_v)} f_v(u) \psi_{-r,v}(u) du,$$
	$$\displaystyle C_c^\infty(L(F_v))\ni f_v\mapsto \int_{\cB^L(F_v)} f_v(s) \psi^L_{\cB,v}(s) ds$$
	extend by continuity to $\cC^w(G_r(F_v))$ and $\cC^w(L(F_v))$ respectively. We denote these unique continuous extensions by $\mathcal{P}_{N_r,\psi_{-r},v}$, $\mathcal{P}_{\cB^L,\psi^L_{\cB},v}$ and call them the local {\it Whittaker period} and {\it Whittaker-Bessel period} respectively.
      \end{paragr}	

      \begin{paragr}
              Let $(\mathbb{G},\mathbb{H})\in \{(U,U'),(\cG,(\cB,\psi_{\cB})),(L,(\cB^L,\psi^L_{\cB})), (G_r,(N_r,\psi_{-r})) \}$. Then, for any tempered irreducible representation $\sigma_v$ of $\mathbb{G}(F_v)$ equipped with an invariant inner product $(\cdot,\cdot)_v$ and vectors $\phi_v,\phi'_v\in \sigma_v$, the matrix coefficient
	$$\displaystyle f_{\phi_v,\phi'_v}: g\in \mathbb{G}(F_v)\mapsto (\sigma_v(g) \phi_v,\phi'_v)_v$$
	belongs to $\cC^w(\mathbb{G}(F_v))$. We set
	\begin{equation*}
		\displaystyle \mathcal{P}_{\mathbb{H},v}(\phi_v,\phi'_v):=\mathcal{P}_{\mathbb{H},v}(f_{\phi_v,\phi'_v}).
	\end{equation*}
      \end{paragr}
      
      \subsection{Relation between global periods}

      \begin{paragr}
        
	\begin{proposition}\label{prop relation global periods}
		Let $\phi_n\in \cA_{Q_n}(U_n)$ and $\phi_{n+1}\in \cA_{cusp}(U_{n+1})$. Then, there exists $c>0$ such that for $s\in \cH_{>c}$ we have the identity
		\begin{equation*}
			\displaystyle \cP_{U'}(E_{Q_n}^{U_n}(\phi_n,s)\otimes \phi_{n+1})=\int_{\cB'(\AAA)\backslash U'(\AAA)} \cP_{\cB^L,\psi^L_{\cB}}\left(R(h)(\phi_{n,s}\otimes \phi_{n+1})\right) dh
		\end{equation*}
		where the right expression is absolutely convergent.
	\end{proposition}
	
	\begin{preuve}
		For $\Re(s)$ sufficiently large, we have
		$$\displaystyle E_{Q_n}^{U_n}(h,\phi_n,s)=\sum_{\gamma\in Q_n(F)\backslash U_n(F)} \phi_{n,s}(\gamma h),\;\;\; h\in [U_n],$$
		so that (the resulting expression being absolutely convergent by cuspidality of $\phi_{n+1}$)
		\begin{align}\label{Unfolding PU' period}
			\displaystyle \cP_{U'}(E_{Q_n}^{U_n}(\phi_n,s)\otimes \phi_{n+1}) & =\int_{[U_n]} \sum_{\gamma\in Q_n(F)\backslash U_n(F)} \phi_{n,s}(\gamma h) \phi_{n+1}(h) dh \\
			\nonumber & =\int_{Q_n(F)\backslash U_n(\AAA)} \phi_{n,s}(h) \phi_{n+1}(h) dh  \\
			\nonumber & =\int_{L_n(F)V_n(\AAA)\backslash U_n(\AAA)} \phi_{n,s}(h) \phi_{n+1,V_n}(h) dh
		\end{align}
		where
		$$\displaystyle \phi_{n+1,V_n}(h):=\int_{[V_n]} \phi_{n+1}(vh)dv.$$
		Note that the subgroup $G_r\ltimes N$ of $U_{n+1}$ contains $V_n$ as a normal subgroup and that the quotient $G_r\ltimes N/V_n$ can be identified with the mirabolic subgroup $P_{r+1}$ of $G_{r+1}$ via restriction to the subspace $\langle x_1,\ldots,x_r,x_{r+1}\rangle$, where we have set $x_{r+1}=v_0$. Moreover, for any $h\in U_{n+1}(\AAA)$, the function $p\in [P_{r+1}]\mapsto \phi_{n+1,V_n}(ph)$ is cuspidal. It follows that we have a Fourier expansion
		\begin{equation}\label{Fourier expansion phin+1Vn}
			\displaystyle \phi_{n+1,V_n}(h)=\sum_{\delta\in N_r(F)\backslash G_r(F)} \phi_{n+1,N,\psi_N}(\delta h)
		\end{equation}
		where
		$$\displaystyle \phi_{n+1,N,\psi_N}(h):=\int_{[N_{r+1}]} \phi_{n+1,V_n}(uh) \psi_{r+1}(u)du=\int_{[N]} \phi_{n+1}(uh) \psi_N(u) du.$$
		Replacing $\phi_{n+1,V_n}$ by its Fourier expansion \eqref{Fourier expansion phin+1Vn} in \eqref{Unfolding PU' period}, we formally obtain (remembering the decomposition \eqref{decomp Ln} for $L_n$)
		\[\begin{aligned}
			\displaystyle \cP_{U'}(E_{Q_n}^{U_n}(\phi_n,s)\otimes \phi_{n+1}) & =\int_{N_r(F)U_m(F)V_n(\AAA)\backslash U_n(\AAA)} \phi_{n,s}(h) \phi_{n+1,N,\psi_N}(h) dh \\
			& =\int_{\cB'(\AAA)\backslash U'(\AAA)} \cP_{\cB^L,\psi^L_{\cB}}\left(R(h)(\phi_{n,s}\otimes \phi_{n+1})\right) dh.
		\end{aligned}\]
		To justify this formal computation, it remains to check that the integral
		$$\displaystyle \int_{N_r(F)U_m(F)V_n(\AAA)\backslash U_n(\AAA)} \left\lvert\phi_{n,s}(h) \phi_{n+1,N,\psi_N}(h) \right\rvert dh$$
		converges for $\Re(s)$ sufficiently large. This can be rewritten as
		$$\displaystyle \int_{P'(\AAA)\backslash U'(\AAA)}\int_{[U_m]\times T_r(\AAA)\times [N_r]} \left\lvert\phi_{n,s}(uagh) \phi_{n+1,N,\psi_N}(agh) \right\rvert \delta_{Q_n}(a)^{-1}\delta_{B_r}(a)^{-1} dudadgdh.$$
		Thus, as $P'(\AAA)\backslash U'(\AAA)$ is compact, it suffices to check the convergence of
		\[\begin{aligned}
			\displaystyle & \int_{[U_m]\times T_r(\AAA)\times [N_r]} \left\lvert\phi_{n,s}(uag) \phi_{n+1,N,\psi_N}(ag) \right\rvert \delta_{Q_n}(a)^{-1}\delta_{B_r}(a)^{-1} dudadg= \\
			& \int_{[U_m]\times T_r(\AAA)\times [N_r]} \left\lvert\phi_{n}(uag) \phi_{n+1,N,\psi_N}(ag) \right\rvert \delta_{Q_n}(a)^{-1}\delta_{B_r}(a)^{-1} \lvert \det a\rvert_{\AAA_E}^s dudadg.
		\end{aligned}\]
		Let us embed $T_r(\AAA)$ into $\AAA_E^r$ in the natural way. Then, by \cite[Lemma 2.6.1.1]{BCZ} and since $\phi_{n+1}$ is cuspidal, for every $R>0$ we have
		\begin{equation}\label{ineq phin+1NpsiN}
			\displaystyle \left\lvert \phi_{n+1,N,\psi_N}(ag) \right\rvert\ll_R \lVert a\rVert_{\AAA_E^r}^{-R} \lVert g\rVert_{U_m}^{-R},\;\; (a,g)\in T_r(\AAA)\times [U_m].
		\end{equation}
		On the other hand, since $\phi_n$ is of moderate growth, we can find $D>0$ such that
		\begin{equation}\label{ineq phin}
			\displaystyle \left\lvert\phi_{n}(uag)\right\rvert\ll \lVert a\rVert^D \lVert g\rVert_{U_m}^{D},\;\; (u,a,g)\in [N_r]\times T_r(\AAA)\times [U_m].
		\end{equation}
		Combining \eqref{ineq phin+1NpsiN} with \eqref{ineq phin}, we are eventually reduced to the following readily checked property: for every $s$ large enough there exists $R>0$ such that the integral
		$$\displaystyle \int_{T_r(\AAA)} \lVert a\rVert^D \lVert a\rVert_{\AAA_E^r}^{-R} \delta_{Q_n}(a)^{-1} \delta_{B_r}(a)^{-1} \lvert \det a\rvert_{\AAA_E}^s da$$
		converges.
	\end{preuve}
      \end{paragr}

      \subsection{Relations between local periods}
	
      \begin{paragr}
        	Let $v$ be a place of $F$ and let $\tau$, $\sigma_{m}$, $\sigma_{n+1}$ be irreducible representations of $G_r(F_v)$, $U_m(F_v)$ and $U_{n+1}(F_v)$ respectively. We set $\sigma=\sigma_{m}\boxtimes \sigma_{n+1}$ (an irreducible representation of $\cG(F_v)=U_m(F_v)\times U_{n+1}(F_v)$). Let
	$$\displaystyle \cL \in \Hom_{\cB^L(F_v)}(\tau\boxtimes \sigma, \psi^L_{\cB,v}).$$
	By multiplicity one results \cite{AGRS}, \cite[Corollary 15.3]{GGP}, \cite{JSZ} $\cL$ factors as $\cL=\cL^W\otimes \cL^B$ where $\cL^W\in \Hom_{N_r(F_v)}(\tau,\psi_{-r,v})$ and $\cL^B\in \Hom_{\cB(F_v)}(\sigma,\psi_{\cB,v})$.
	
	For every $s\in \CC$, we set $\tau_s=\tau \lvert \det \rvert_{E_v}^s$ and we denote by $I_{Q_n(F_v)}^{U_n(F_v)}(\tau_s \boxtimes \sigma_{m})$ the normalized parabolic induction of $\tau_s \boxtimes \sigma_{m}$ to $U_n(F_v)$.
	
	\begin{proposition}\label{proposition induced local intertwining}
		There exists $c\in \RR$ such that the functional
		$$\displaystyle \mathcal{L}_s^{U'}: \phi_{n,s}\otimes \phi_{n+1}\in I_{Q_n(F_v)}^{U_n(F_v)}(\tau_s \boxtimes \sigma_{m})\otimes \sigma_{n+1}\mapsto \int_{\cB'(F_v)\backslash U'(F_v)} \cL(\phi_{n,s}(h)\otimes \sigma_{n+1}(h)\phi_{n+1}) dh$$
		converges absolutely for $s\in \cH_{>c}$. If moreover both $\sigma$ and $\tau$ are tempered we may take $c=-1/2$. Furthermore, for $s$ with sufficiently large real part the following assertions are equivalent:
		\begin{enumerate}
			\item There exist $\phi_{\tau}\in \tau$, $\phi_m\in \sigma_{m}$ and $\phi_{n+1}\in \sigma_{n+1}$ such that $\mathcal{L}(\phi_{\tau}\otimes \phi_m\otimes \phi_{n+1})\neq 0$;
			
			\item There exist $\phi_{n,s}\in I_{Q_n(F_v)}^{U_n(F_v)}(\tau_s\boxtimes \sigma_{m})$ and $\phi_{n+1}\in \sigma_{n+1}$ such that
			$\mathcal{L}_s^{U'}(\phi_{n,s}\otimes \phi_{n+1})\neq 0$.
		\end{enumerate}
		
	\end{proposition}
	
	\begin{preuve}
		Since $P'=T_r\cB'$ is a parabolic subgroup of $U'$, it suffices to check the convergence of the integral
		\[\begin{aligned}
			\displaystyle & \int_{T_r(F_v)} \cL(\phi_{n,s}(a)\otimes \sigma_{n+1,v}(a)\phi_{n+1}) \delta_{P'}(a)^{-1}da= \\
			& \int_{T_r(F_v)} \cL((\tau(a)\otimes \sigma(a))\phi_{n,s}(1)\otimes \phi_{n+1}) \lvert \det a\rvert_{E_v}^s\delta_{Q_n}(a)^{1/2}\delta_{P'}(a)^{-1}da
		\end{aligned}\]
		for $\Re(s)\gg 0$.
		
		\begin{lemme}
			There exists $D>0$ such that for every $\phi\in \tau\boxtimes \sigma$ and $R>0$, we have
			\begin{equation}
				\displaystyle \left\lvert \mathcal{L}((\tau(a)\otimes\sigma(a))\phi)\right\rvert\ll_R  \lVert a\rVert_{E_v^r}^{-R} \lVert a\rVert^D,\;\; \mbox{ for } a \in T_r(F_v).
			\end{equation}
			If moreover $\sigma$ and $\tau$ are tempered, for every $\phi\in \tau\boxtimes \sigma$ and $R>0$, we have
			\begin{equation}
				\displaystyle \left\lvert \mathcal{L}((\tau(a)\otimes\sigma(a))\phi)\right\rvert\ll_R  \lVert a\rVert_{E_v^r}^{-R} \Xi^{G_{r,v}}(a)\Xi^{U_{n+1,v}}(a),\;\; \mbox{ for } a \in T_r(F_v).
			\end{equation}
		\end{lemme}
		
		\begin{preuve}
			It suffices to prove the lemma when $\phi$ is a pure tensor i.e. it is of the form $\phi=\phi_\tau\otimes \phi_m\otimes \phi_{n+1}$ where $\phi_\tau\in \tau$, $\phi_m\in \sigma_{m}$ and $\phi_{n+1}\in \sigma_{n+1}$. Indeed, in the case where $v$ is non-Archimedean every vector in $\tau\boxtimes \sigma$ is a sum of pure tensors whereas if $v$ is Archimedean the claimed inequalities would automatically extend from the algebraic tensor product to the completed tensor product by the Banach-Steinhaus theorem (see e.g. \cite[Theorem 34.1]{Treves}).
			
			Using the equality $\cL=\cL^W\otimes \cL^B$, we are therefore reduced to show the existence of $D>0$ such that for every $R>0$  (resp. that when $\tau$, $\sigma$ are tempered for every $R>0$) we have
			\begin{equation}\label{ineq LW}
				\displaystyle  \left\lvert \mathcal{L}^W(\tau(a)\phi_\tau)\right\rvert\ll \lVert a\rVert^D \;\; (\mbox{resp. }  \left\lvert \mathcal{L}^W(\tau(a)\phi_\tau)\right\rvert\ll \Xi^{G_{r,v}}(a)\; )
			\end{equation}
			and
			\begin{equation}\label{ineq LB}
				\displaystyle  \left\lvert \mathcal{L}^B(\phi_m\otimes \sigma_{n+1}(a)\phi_{n+1})\right\rvert\ll_R  \lVert a\rVert_{E_v^r}^{-R} \lVert a\rVert^D
			\end{equation}
			$$\displaystyle (\mbox{resp. } \left\lvert \mathcal{L}^B(\phi_m\otimes \sigma_{n+1}(a)\phi_{n+1})\right\rvert\ll_R  \lVert a\rVert_{E_v^r}^{-R} \Xi^{U_{n+1,v}}(a)\; )$$
			for $a\in T_r(F_v)$.
			
			The estimates \eqref{ineq LW} and \eqref{ineq LB} can be established along the same lines as Lemma B.2.1 and Lemma 7.3.1 (i) of \cite{BP3} respectively. More precisely, in the tempered case \eqref{ineq LW} is a direct application of Lemma B.2.1 of loc.cit. The same inequality for general representations (and for a suitable $D$) is a consequence of the continuity of $\cL^W$ when $v$ is Archimedean whereas, for $v$ non-Archimedean, by the same argument as in \cite[Lemma B.2.1]{BP3} we can bound the function $a\in T_r(F_v)\mapsto \left\lvert \mathcal{L}^W(\tau(a)\phi_\tau)\right\rvert$ by a matrix coefficient of $\tau_v$ which is in turn essentially bounded by $\lVert a\rVert^D$ for some $D>0$.
			
			As for \eqref{ineq LB}, we first note that
			\begin{equation}\label{eq LB1}
				\displaystyle \psi_{\cB,v}(e^{\Ad(a) X})\cL^B(\phi_m\otimes \sigma_{n+1}(a)\phi_{n+1})=\cL^B(\phi_m\otimes \sigma_{n+1}(ae^X)\phi_{n+1})
			\end{equation}
			for every $(a,X)\in T_r(F_v)\times \Lie(\cB(F_v))$. Furthermore, when $v$ is Archimedean, differentiating the above identities yields
			\begin{equation}\label{eq LB2}
				\displaystyle d\psi_{\cB,v}(\Ad(a) X)\cL^B(\phi_m\otimes \sigma_{n+1}(a)\phi_{n+1})=\cL^B(\phi_m\otimes \sigma_{n+1}(a)\sigma_{n+1}(X)\phi_{n+1})
			\end{equation}
			where $d\psi_{\cB,v}: Lie(\cB(F_v))\to \CC$ denotes the differential of $\psi_{\cB,v}$ at $1$. From \eqref{eq LB2}, we deduce that for every linear form $\lambda: E_v^r\to F_v$ there exists $X_\lambda\in Lie(\cB(F_v))$ such that
			$$\displaystyle \lambda(a) \cL^B(\phi_m\otimes \sigma_{n+1}(a)\phi_{n+1})=\cL^B(\phi_m\otimes \sigma_{n+1}(a)\sigma_{n+1}(X_\lambda)\phi_{n+1})$$
			for every $a\in T_r(F_v)$. Similarly, when $v$ is non-Archimedean, we deduce from \eqref{eq LB1} and the smoothness of $\phi_{n+1}$ the existence of $C>0$ such that $\cL^B(\phi_m\otimes \sigma_{n+1}(a)\phi_{n+1})=0$ unless $\lVert a\rVert_{E_v^r}\leqslant C$. Combining these two facts, we are now reduced to show the existence of $D>0$ such that (resp. that for $\sigma$ tempered we have)
			\begin{equation*}
				\displaystyle  \left\lvert \mathcal{L}^B(\phi_m\otimes \sigma_{n+1}(a)\phi_{n+1})\right\rvert\ll \lVert a\rVert^D \;\; (\mbox{resp. } \left\lvert \mathcal{L}^B(\phi_m\otimes \sigma_{n+1}(a)\phi_{n+1})\right\rvert\ll \Xi^{U_{n+1,v}}(a))
			\end{equation*}
			for $a\in T_r(F_v)$. The case where $\sigma$ is tempered is a direct application of \cite[Lemma 7.3.1(i)]{BP3}  of whereas the case of a general representation follows from continuity of $\cL^B$ in the Archimedean case or an argument similar to that of loc. cit. to show that $a\mapsto \mathcal{L}^B(\phi_m\otimes \sigma_{n+1}(a)\phi_{n+1})$ is bounded by a matrix coefficient of $\sigma$ in the non-Archimedean case.
		\end{preuve}
		
		By the lemma, and since there exists $d>0$ such that (see \cite[Lemme II.1.1]{WaldPlanch})
		$$\displaystyle \Xi^{G_{r,v}}(a)\Xi^{U_{n+1,v}}(a)\ll \delta_{B_r}(a)^{1/2}\delta_P(a)^{1/2}\sigma(a)^d,\;\; \mbox{ for } a\in T_r(F_v),$$
		the convergence part of the proposition reduces to the two following readily checked facts:
		\begin{itemize}
			\item Let $D>0$. Then, we can find $c>0$ such that for every $s>c$ and for suitable $R>0$ the integral
			$$\displaystyle \int_{T_r(F_v)} \lVert a \rVert^D \lVert a\rVert_{E_v^r}^{-R} \lvert \det a\rvert_{E_v}^s \delta_{Q_n}(a)^{1/2} \delta_{P'}(a)^{-1} da$$
			converges;
			
			\item For every $s>-1/2$ we can find $R>0$ such that the integral
			\[\begin{aligned}
				\displaystyle & \int_{T_r(F_v)} \lVert a\rVert_{E_v^r}^{-R} \lvert \det a\rvert_{E_v}^s \delta_{B_r}(a)^{1/2}\delta_P(a)^{1/2}\delta_{Q_n}(a)^{1/2} \delta_{P'}(a)^{-1} \sigma(a)^dda= \\
				& \int_{T_r(F_v)} \lVert a\rVert_{E_v^r}^{-R} \lvert \det a\rvert_{E_v}^{s+1/2} \sigma(a)^dda
			\end{aligned}\]
			converges.
		\end{itemize}
		
		The implication $2.\Rightarrow 1.$ is clear from the definition of $\cL_s^{U'}$. Let us show the converse. Thus we assume that $\cL$ is not identically zero and we aim to prove that the same holds for $\cL_s^{U'}$ for $\Re(s)\gg 0$. First, by the equality
		$$\displaystyle \cL_s^{U'}(\phi_{n,s}\otimes \phi_{n+1})=\int_{Q_n(F_v)\backslash U_n(F_v)} \int_{N_r(F_v)\backslash G_r(F_v)} \cL((\tau(g)\otimes \sigma(gh))\phi_{n,s}(h)\otimes \phi_{n+1}) \lvert \det(g)\rvert_{E_v}^s\delta_{Q_n}(g)^{-1/2}dgdh$$
		and since the space of $I_{Q_n(F_v)}^{U_n(F_v)}(\tau_{s}\boxtimes \sigma_{m})$ is stable by multiplication by functions in $C^\infty(Q_n(F_v)\backslash U_n(F_v))$, it suffices to show the existence $\phi_\tau\in \tau$, $\phi_m\in \sigma_{m}$ and $\phi_{n+1}\in \sigma_{n+1}$ such that
		$$\displaystyle\int_{N_r(F_v)\backslash G_r(F_v)} \cL(\tau(g)\phi_\tau\otimes \phi_m\otimes  \sigma_{n+1}(g)\phi_{n+1}) \lvert \det(g)\rvert_{E_v}^s\delta_{Q_n}(g)^{-1/2}dg\neq 0.$$
		Let $f\in \cS(V_{n+1}(F_v))$, $g\in G_r(F_v)$, $\phi_m\in \sigma_m$ and $\phi_{n+1}\in \sigma_{n+1}$. From the equivariance property of $\cL^B$ we deduce
		$$\displaystyle \cL^B(\phi_m\otimes \sigma_{n+1}(g)\sigma_{n+1}(f)\phi_{n+1})=\hat{f}(g^*y_r)\cL^B(\phi_m\otimes \sigma_{n+1}(g)\phi_{n+1})$$
		where for $y\in Y_r(F_v)$ we have set
		$$\displaystyle \hat{f}(y)=\int_{V_{n+1}(F_v)} f(u) \psi_v(h_{n+1}(uv_0,y))du.$$
		By the theory of Fourier transform, $f\mapsto \hat{f}$ induces a surjective map $\cS(V_{n+1}(F_v))\to  \cS(Y_r(F_v))$ whereas the map $g\mapsto g^* y_r$ induces an embedding $\cS(P_r(F_v)\backslash G_r(F_v))\hookrightarrow \cS(Y_r(F_v))$ (where we recall thet $P_r$ denotes the mirabolic subgroup of $G_r$). Consequently, it suffices to prove the existence of $\phi_\tau\in \tau$, $\phi_m\in \sigma_{m}$, $\phi_{n+1}\in \sigma_{n+1}$  and $f\in \cS(P_r(F_v)\backslash G_r(F_v))$ such that
		\[\begin{aligned}
			\displaystyle & \int_{P_r(F_v)\backslash G_r(F_v)} f(g) \int_{N_r(F_v)\backslash P_r(F_v)} \cL^W(\tau(pg)\phi_\tau) \\
			& \cL^B(\phi_m\otimes  \sigma_{n+1}(pg)\phi_{n+1}) \lvert \det(pg)\rvert_{E_v}^{s} \delta_{Q_n}(pg)^{-1/2} \lvert \det(p)\rvert_{E_v}^{-1} dpdg\neq 0
		\end{aligned}\]
		or, equivalently, the existence of $\phi_\tau\in \tau$, $\phi_m\in \sigma_{m}$, $\phi_{n+1}\in \sigma_{n+1}$ such that
		$$\displaystyle \int_{N_{r-1}(F_v)\backslash G_{r-1}(F_v)} \cL^W(\tau(h)\phi_\tau)\cL^B(\phi_m\otimes  \sigma_{n+1}(h)\phi_{n+1}) \lvert \det(h)\rvert_{E_v}^{s-1} \delta_{Q_n}(h)^{-1/2} dh\neq 0.$$
		By \cite[Theorem 6]{GKaz} and \cite[Theorem 1]{Kemar}, for every $f\in C_c^\infty(N_{r-1}(F_v)\backslash G_{r-1}(F_v))$, we can find $\phi_\tau\in \tau$ such that $\cL^W(\tau(h)\phi_\tau)=f(h)$ for $h\in G_{r-1}(F_v)$ and from this the claim follows readily from the non-vanishing of $\cL^B$.
	\end{preuve}
	
      \end{paragr}

        \begin{paragr}

	\begin{proposition}\label{prop relation local periods}
		Assume that $\sigma$ and $\tau$ are tempered. Then, for every $\phi_n,\phi'_n\in I_{Q_n(F_v)}^{U_n(F_v)}(\tau\boxtimes \sigma_{m})$ and $\phi_{n+1},\phi_{n+1}'\in \sigma_{n+1}$, we have
		\begin{align}
			\displaystyle & \mathcal{P}_{U',v}(\phi_n\otimes \phi_{n+1},\phi'_n\otimes \phi'_{n+1})= \\
			\nonumber & \nu(G_{r,v})^{-1}\int_{(\cB'(F_v)\backslash U'(F_v))^2} \cP_{\cB^L,\psi^L_{\cB},v}(\phi_n(h_1)\otimes \sigma_{n+1}(h_1)\phi_{n+1},\phi'_n(h_2)\otimes \sigma_{n+1}(h_2)\phi'_{n+1}) dh_1dh_2.
		\end{align}
	\end{proposition}
	
	\begin{preuve}
		By definition of the local period $\mathcal{P}_{U',v}$ and of the invariant inner product on $I_{Q_n(F_v)}^{U_n(F_v)}(\tau\boxtimes \sigma_{m})$, we have
		\begin{align}\label{eq1 unfold local periods}
			\displaystyle & \mathcal{P}_{U'_v}(\phi_n\otimes \phi_{n+1},\phi'_n\otimes \phi'_{n+1})= \\
			\nonumber & \int_{U'(F_v)}\int_{Q_n(F_v)\backslash U_n(F_v)} (\phi_n(h_2h_1),\phi'_n(h_2)) dh_2 (\sigma_{n+1}(h_1)\phi_{n+1},\phi'_{n+1})dh_1.
		\end{align}
		The above double integral is absolutely convergent. Indeed, from \cite[Theorem 2]{CoHaHo} and \cite[Lemme II.1.6]{WaldPlanch} we have
		$$\displaystyle \int_{Q_n(F_v)\backslash U_n(F_v)} \left\lvert (\phi_n(h_2h_1),\phi'_n(h_2))\right\rvert dh_2 \ll \Xi^{U_{n,v}}(h_1),\;\; h_1\in U_n(F_v),$$
		and $(\sigma_{n+1}(h_1)\phi_{n+1},\phi'_{n+1})\ll \Xi^{U_{n+1,v}}(h_1)$ whereas $\Xi^{U_v}=\Xi^{U_{n,v}}\Xi^{U_{n+1,v}}$ is integrable on $U'(F_v)$ (\cite[Lemme 6.5.1(i)]{BP3}).
		
		To finish the proof, we need the following lemma.
		
		\begin{lemme}\label{lemma local Fourier inversion}
			For every $f\in \cC^w(G_r(F_v)\times \cG(F_v))$ we have the identity
			\begin{equation}\label{eq lemma local Fourier inversion}
				\displaystyle \int_{Q_n(F_v)} f(q) \delta_{Q_n}(q)^{1/2}d_Lq=\nu(G_{r,v})^{-1}\int_{(N_r(F_v)\backslash G_r(F_v))^2} \cP_{\cB^L,\psi^L_{\cB},v}(L(g_1)R(g_2)f) \delta_{Q_n}(g_1g_2)^{-1/2}dg_1dg_2
			\end{equation}
			where both sides are absolutely convergent.
		\end{lemme}
		
		Indeed, assuming the lemma for the moment, from \eqref{eq1 unfold local periods} we obtain
		\[\begin{aligned}
			\displaystyle & \mathcal{P}_{U'_v}(\phi_n\otimes \phi_{n+1},\phi'_n\otimes \phi'_{n+1}) \\
			& =\int_{(Q_n(F_v)\backslash U_n(F_v))^2} \int_{Q_n(F_v)} (\phi_n(qh_1),\phi'_n(h_2)) (\sigma_{n+1}(qh_1)\phi_{n+1},\sigma_{n+1}(h_2)\phi'_{n+1})d_Lq dh_1dh_2 \\
			& =\int_{(Q_n(F_v)\backslash U_n(F_v))^2} \int_{(N_r(F_v)\backslash G_r(F_v))^2} \\
			& \cP_{\cB^L,\psi^L_{\cB},v}(\phi_n(g_1h_1)\otimes \sigma_{n+1}(g_1h_1)\phi_{n+1} ,\phi'_n(g_2h_2)\otimes \sigma_{n+1}(g_2h_2)\phi'_{n+1}) \delta_{Q_n}(g_1g_2)^{-1}dg_1dg_2 dh_1dh_2 \\
			& =\int_{(\cB'(F_v)\backslash U'(F_v))^2} \cP_{\cB^L,\psi^L_{\cB},v}(\phi_n(h_1)\otimes \sigma_{n+1}(h_1)\phi_{n+1},\phi'_n(h_2)\otimes \sigma_{n+1}(h_2)\phi'_{n+1}) dh_1dh_2
		\end{aligned}\]
		where in the second equality we have applied the lemma thanks to the fact that the function
		$$(g_r,g_m,g_{n+1})\mapsto \left((\tau(g_r)\otimes \sigma_m(g_m))\phi_n(h_1),\phi'_n(h_2)\right)\left(\sigma_{n+1}(g_{n+1})\sigma_{n+1}(h_1)\phi_{n+1},\sigma_{n+1}(h_2)\phi'_{n+1}\right)$$
		belongs to $\cC^w(G_r(F_v)\times \cG(F_v))$ and we are done.
		
		\begin{preuve}
			(of Lemma \ref{lemma local Fourier inversion}). First, we check that both sides are convergent and define continuous functionals on $\cC^w(G_r(F_v)\times \cG(F_v))$.
			
			For the left hand side, we can use the identity (\cite[Lemme II.1.6]{WaldPlanch})
			$$\displaystyle \Xi^{U_v}(g)=\int_{K^{U'}_{v}} \Xi^{G_{r,v}\times \cG_v}(l(kg)) \delta_{Q}(l(kg))^{1/2}dk,\;\; \mbox{ for } g\in U(F_v),$$
			where, for $g\in U(F_v)$, $l(g)$ denotes any element in $L(F_v)=G_r(F_v)\times \cG(F_v)$ such that $l(g)^{-1}g\in V_n(F_v)K^U_{v}$, which in turn implies
			\[\begin{aligned}
				\displaystyle \int_{Q_n(F_v)} \Xi^{G_{r,v}\times \cG_v}(q) \sigma(q)^d \delta_{Q_n}(q)^{1/2}d_Lq & =\int_{U'(F_v)} \Xi^{G_{r,v}\times \cG_v}(l(h)) \sigma(l(h))^d \delta_{Q}(l(h))^{1/2} dh \\
				& \ll \int_{U'(F_v)} \Xi^{U_v}(h) \sigma(h)^d dh<\infty
			\end{aligned}\]
			for every $d>0$.
			
			For the right hand side, as in the proof of Proposition \ref{proposition induced local intertwining}, it suffices to show for every $d>0$ the existence of $d'>0$ as well as a continuous semi-norm $\nu_d$ on $\cC^w_d(G_r(F_v)\times \cG(F_v))$ such that
			\begin{equation}
				\displaystyle \left\lvert \cP_{\cB_U,\psi_{\cB_U},v}(L(a_1)R(a_2)f)\right\rvert\ll \Xi^{G_{r,v}\times U_{n+1,v}}(a_1)\Xi^{G_{r,v}\times U_{n+1,v}}(a_2) \sigma(a_1)^{d'}\sigma(a_2)^{d'} \nu_d(f)
			\end{equation}
			for $f\in \cC^w_d(G_r(F_v)\times \cG(F_v))$ and $a_1,a_2\in T_r(F_v)$. Such an inequality can be proved along the same lines as \cite[Lemma 7.3.1(ii)]{BP3}.
			
			Now that we know that both sides of \eqref{eq lemma local Fourier inversion} define continuous functionals on $\cC^w(G_r(F_v)\times \cG(F_v))$, by density it suffices to check the claimed identity for $f=f_1\otimes f_2\in C_c^\infty(G_r(F_v))\otimes C_c^\infty(\cG(F_v))$. The subgroup $G_r\cB$ of $U_{n+1}$ contains $U_m\ltimes V_n$ as a normal subgroup and we may define a function $f_3\in C_c^\infty(G_r(F_v)\cB(F_v)/U_m(F_v)V_n(F_v))$ by
			$$\displaystyle f_3(p)=\delta_{Q_n}(p)^{1/2}\int_{U_m(F_v)\times V_n(F_v)} f_2(phu) dudh,\;\; p\in G_r(F_v)\cB(F_v).$$
			Moreover, restriction to the susbpace $\langle x_1,\ldots,x_{r+1}\rangle$, where we have again set $x_{r+1}=v_0$, induces an identification $G_r\cB/U_mV_n\simeq P_{r+1}$. Seeing $f_3$ as a test function on $P_{r+1}(F_v)$ in this way the identity of the lemma becomes
			\begin{align}\label{eq2 Fourier inversion}
				\displaystyle \int_{G_r(F_v)} f_1(g)f_3(g) dg=\nu(G_{k,v})^{-1}\int_{(N_r(F_v)\backslash G_r(F_v))^2} W_{f_1}(g_1,g_2) W_{f_3}(g_1,g_2) dg_2dg_1
			\end{align}
			where we have set
			$$\displaystyle W_{f_1}(g_1,g_2)=\int_{N_r(F_v)} f_1(g_1^{-1}u g_2) \psi_r(u)^{-1}du,\;\; W_{f_3}(g_1,g_2)=\int_{N_{r+1}(F_v)} f_3(g_1^{-1}ug_2)\psi_{r+1}(u) du.$$
			But \eqref{eq2 Fourier inversion} is now a direct consequence of the Fourier inversion formula \eqref{eq1 Fourier inversion}.

		\end{preuve}

              \end{preuve}

      \end{paragr}

      \subsection{Unramified computation}\label{ssect: unr computation}
	
      \begin{paragr}

	We continue with the setting of the previous subsection assuming moreover that $v\notin S$ and the representations $\tau$, $\sigma_{m}$, $\sigma_{n+1}$ are unramified in the sense that $\tau^{K_{r,v}}\neq 0$, $\sigma_m^{K^U_{m,v}}\neq 0$ and $\sigma_{n+1}^{K^U_{n+1,v}}\neq 0$. Note that this implies $I_{Q_n(F_v)}^{U_n(F_v)}(\tau_s\boxtimes \sigma_{m})^{K^U_{n,v}}\neq 0$ i.e. $I_{Q_n(F_v)}^{U_n(F_v)}(\tau_s\boxtimes \sigma_{m})$ is also unramified.
      \end{paragr}

      \begin{paragr}
        
      \begin{proposition}\label{prop unr computation}
		For $\Re(s)$ sufficiently large, $\phi^\circ_{n,s}\in I_{Q_n(F_v)}^{U_n(F_v)}(\tau_s\boxtimes \sigma_{m})^{K^U_{n,v}}$ and $\phi^\circ_{n+1}\in \sigma_{n+1}^{K^U_{n+1,v}}$, we have
		\begin{align}
		\displaystyle & \mathcal{L}_s^{U'}(\phi^\circ_{n,s}\otimes \phi^\circ_{n+1})= \\
		\nonumber & \frac{\vol(K^U_{n,v})}{\vol(K^U_{n,v}\cap \cB'(F_v))}\frac{L(\frac{1}{2}+s,\tau\times \sigma_{n+1})}{L(1+s,\tau\times \sigma_{m})L(1+2s,\tau,\As^{(-1)^m})}\mathcal{L}(\phi^\circ_{n,s}(1)\otimes \phi^\circ_{n+1}).
		\end{align}
	\end{proposition}
	
	\begin{preuve}
		Let $\phi^\circ_\tau\in \tau^{K_{r,v}}$ and $\phi_{m}^\circ\in \sigma_m^{K^U_{m,v}}$ be such that $\phi_n^\circ(1)= \phi^\circ_\tau\otimes \phi_{m}^\circ$. Recall the factorization $\cL=\cL^W\otimes \cL^B$. Up to scaling, we may assume, without loss of generality, that
		$$\displaystyle \cL^W(\phi_\tau^\circ)=1 \mbox{ and } \cL^B(\phi_m^\circ\otimes \phi_{n+1}^\circ)=1.$$
		
		Let $P'=P\cap U'\subset U'$ be the parabolic subgroup stabilizing the flag
		$$\displaystyle 0=X_0\subset X_1\subset\ldots\subset X_r.$$
		Then, we have $P'=\cB'\rtimes T_r$ and from the Iwasawa decomposition $U'(F_v)=P'(F_v)K_{n,v}^{U}$, we obtain
		\begin{align}\label{eq1 unr computation}
		\displaystyle & \cL_s^{U'}(\phi_{n,s}^\circ\otimes \phi_{n+1}^\circ) =\frac{vol(K_{n,v}^U)}{vol(K_{n,v}^U\cap \cB'(F_v))} \sum_{t\in \Lambda_r} \cL(\phi_{n,s}^\circ(t)\otimes \sigma_{n+1}(t)\phi_{n+1}^\circ) \delta_{P'}(t)^{-1} \\
		\nonumber & =\frac{vol(K_{n,v}^U)}{vol(K_{n,v}^U\cap \cB'(F_v))} \sum_{t\in \Lambda_r} \cL^W(\tau(t)\phi_{\tau}^\circ) \cL^B(\phi_{m}^\circ\otimes \sigma_{n+1}(t)\phi_{n+1}^\circ) \lvert \det t\rvert_{E_v}^s \delta_{Q_n}(t)^{1/2}\delta_{P'}(t)^{-1}
		\end{align}
		where we have set $\Lambda_r=T_r(F_v)/T_r(\cO_v)$ which we will identify with the cocharacter lattice $X_*(T_{r,v})$ via the map $\lambda\mapsto \lambda(\varpi_F)$. Let $\Lambda_r^+\subset \Lambda_r$ be the cone of dominant cocharacters with respect to $B_r$ and $\Lambda_r^{++}\subset \Lambda_r^+$ the subcone of cocharacters that are moreover dominant with respect to $B_{r+1}$ through the embedding $g\in G_r\mapsto \begin{pmatrix} g & \\ & 1 \end{pmatrix}\in G_{r+1}$. To lighten a bit the computations we will assume from now on that the local measures at $v$ have been chosen such that $vol(K_{n,v}^U)=vol(K_{n,v}^U\cap \cB'(F_v))=1$.
		
		Before proceeding, it will be convenient to introduce some notation pertaining to complex dual groups:
		\begin{itemize}
			\item For $\ell\in\{m,m+1,n+1\}$, we identify the dual group $\widehat{U}_\ell$ of $U_\ell$ (resp. $\widehat{G}_r$ of $G_r$) with $\GL_\ell(\CC)$ (resp. $\GL_r(\CC)\times \GL_r(\CC)$) equipped with its standard pinning. Its $L$-group can be written
			$$\displaystyle {}^L U_\ell=\GL_\ell(\CC)\rtimes \Gal(E/F) \;\; (\mbox{resp. } {}^L G_r=(\GL_r(\CC)\times \GL_r(\CC))\rtimes \Gal(E/F))$$
			where the Galois action is given by $c(g)=g^\star$ (resp. $c(g_1,g_2)=(g_2,g_1)$) where the involution $g\in\GL_\ell(\CC)\mapsto g^\star$ is
			$$\displaystyle g^\star=J_\ell {}^t g^{-1} J_\ell^{-1},\;\;\; J_\ell=\begin{pmatrix} & & 1 \\ & \iddots & \\ (-1)^{\ell-1} & & \end{pmatrix}.$$
			We will also denote by $S\mapsto S^\star$ the automorphism of ${}^L G_r$ which is the identity on $\Gal(E/F)$ and given by $(g_1,g_2)\mapsto (g_2^\star, g_1^\star)$ on $\widehat{G}_r$.
			
			\item We will write $(\widehat{T}_r,\widehat{B}_r)$ for the standard Borel pair in $\widehat{G}_r$ and $(\widehat{T}^U_\ell,\widehat{B}_\ell^U)$ for that in $\widehat{U}_\ell$. The corresponding semi-direct products with $\Gal(E/F)$ will be denoted
			$$\displaystyle {}^L T_r=\widehat{T}_r\rtimes \Gal(E/F),\; {}^L B_r=\widehat{B}_r\rtimes \Gal(E/F),\; {}^L T^U_\ell=\widehat{T}^U_\ell\rtimes \Gal(E/F),\; {}^L B^U_\ell=\widehat{B}^U_\ell\rtimes \Gal(E/F).$$
			
			\item The dual groups of $\cG$, $\widetilde{U}$ are
			$$\displaystyle \widehat{\cG}=\widehat{U}_m\times \widehat{U}_{n+1},\; \widehat{\widetilde{U}}=\widehat{U}_m\times \widehat{U}_{m+1}$$
			and, writing $\times_\Gamma$ for the fibered products over $\Gal(E/F)$, their $L$-groups are
			$$\displaystyle {}^L \cG={}^L U_m\times_\Gamma {}^L U_{n+1},\; {}^L \widetilde{U}={}^L U_m\times_\Gamma {}^L U_{m+1}$$
			respectively. We let $\widehat{B}=\widehat{B}^U_m\times \widehat{B}^U_{n+1}$, $\widehat{T}=\widehat{T}^U_m\times \widehat{T}^U_{n+1}$ (resp. $\widehat{\widetilde{B}}=\widehat{B}^U_m\times \widehat{B}^U_{m+1}$, $\widehat{\widetilde{T}}=\widehat{T}^U_m\times \widehat{T}^U_{m+1}$) be the standard Borel and maximal torus in $\widehat{\cG}$ (resp. in $\widehat{\widetilde{U}}$). We also write
			$$\displaystyle {}^L T=\widehat{T}\rtimes \Gal(E/F),\; {}^L \widetilde{T}=\widehat{\widetilde{T}}\rtimes \Gal(E/F).$$
			
			\item The parabolic subgroup $Q_{n+1}\subset U_{n+1}$ stabilizing the isotropic subspace $X_r$ corresponds to a standard parabolic subgroup $\widehat{Q}_{n+1}$ of $\widehat{U}_{n+1}$ with standard Levi
			$$\displaystyle \widehat{L}_{n+1}=\begin{pmatrix} \GL_r(\CC) & & \\ & \GL_{m+1}(\CC) & \\ & & \GL_r(\CC) \end{pmatrix}.$$
			The corresponding $L$-group ${}^L L_{n+1}=\widehat{L}_{n+1}\rtimes \Gal(E/F)$ is isomorphic to ${}^L G_r\times_\Gamma {}^L U_{m+1}$ via the map which is the identity on $\Gal(E/F)$ and
			$$\displaystyle \begin{pmatrix} g^{(1)}_r & & \\ & g_{m+1} & \\ & & g_r^{(2)}\end{pmatrix}\mapsto ((g_r^{(1)},g_r^{(2)\star}),g_{m+1})$$
			on $\widehat{L}_{n+1}$. For $S\in {}^L L_{n+1}$, we denote by $(S^{(r)}, S^{(m+1)})\in {}^L G_r\times_\Gamma {}^L U_{m+1}$ its image by this isomorphism.
			
			\item It will also be convenient to use the parabolic subgroup $\cQ=U_m\times Q_{n+1}$ of $\cG$. We set $\widehat{\cQ}=\widehat{U}_m\times \widehat{Q}_{n+1}$, $\widehat{\cL}=\widehat{U}_m\times \widehat{L}_{n+1}$ and ${}^L \cL=\widehat{\cL}\rtimes \Gal(E/F)$. There are two natural decompositions
			$$\displaystyle {}^L \cL={}^L U_m\times_{\Gamma}{}^L L_{n+1} \mbox{ and } {}^L \cL={}^L \widetilde{U}\times_\Gamma {}^L G_r$$
			and for $S\in {}^L \cL$, we will denote by $(S_m,S_{n+1})$ and $(\widetilde{S}, S_{n+1}^{(r)})$ the corresponding respective decompositions of $S$.
			
			\item For every complex Lie group ${}^L \mathbb{G}$ with a subgroup ${}^L \mathbb{Q}$ and respective identity components $\widehat{\mathbb{G}}$, $\widehat{\mathbb{Q}}$, we set
			$$\displaystyle D_{\widehat{\mathbb{G}}/\widehat{\mathbb{Q}}}(S)=\det(1-Ad(S)\mid Lie(\widehat{\mathbb{G}})/Lie(\widehat{\mathbb{Q}})), \mbox{ for } S\in {}^L \mathbb{Q}.$$
			
			\item For $\mathbb{G}\in \{G_r,U_\ell, \cG, \cL, L_{n+1}, T_r, T^U_\ell, T  \}$, we denote by ${}^L \mathbb{G}_v$ the $L$-group of $\mathbb{G}_v=\mathbb{G}\times_F F_v$ that is:
			$$\displaystyle {}^L \mathbb{G}_{v} =\left\{\begin{array}{ll}
			{}^L \mathbb{G} & \mbox{ if } v \mbox{ is inert in } E; \\
			
			\widehat{\mathbb{G}} & \mbox{ if } v \mbox{ splits in } E.
		\end{array} \right.$$
		Also, for $\mathbb{G}\in \{G_r,U_\ell, \cG, \cL, L_{n+1}\}$, we write $W(\mathbb{G}_v)$ for the Weyl group $Norm_{\widehat{\mathbb{G}}}({}^L \mathbb{T}_v)/\widehat{\mathbb{T}}$ where $\widehat{\mathbb{T}}\subset \widehat{\mathbb{G}}$ is the standard maximal torus.
		
		\item The choice of the Borel pair $(T_r,B_r)$ allows to identify $\Lambda_r$ with the group of characters of ${}^L T_{r,v}$ that are trivial on $\Gal(E/F)$ and we will denote by $\Lambda_r\ni t\mapsto \chi_t$ this identification. For $t\in \Lambda_r^+$, we write $ch_t$ for the character of the irreducible representation of ${}^L G_r$ with highest weight $\chi_t$ (see Appendix \ref{appendix Weyl character formula}).

		\item For $k,\ell\in \NN$, we define the representation
		$$\displaystyle {}^L(G_k\times G_\ell)={}^LG_k\times_\Gamma {}^L G_\ell\to \GL(\CC^k\otimes \CC^\ell\oplus \CC^k\otimes \CC^\ell)$$
		$$\displaystyle (S_k,S_\ell)\mapsto S_k\itimes S_\ell$$
		which sends $((g_k^{(1)},g_k^{(2)}),(g^{(1)}_\ell,g_\ell^{(2)}))\in \widehat{G_k}\times \widehat{G_\ell}$ to $g_k^{(1)}\otimes g_\ell^{(1)}\oplus g_k^{(2)}\otimes g_\ell^{(2)}$ and $c\in \Gal(E/F)$ to the operator $\iota$ that swaps the two copies of $\CC^k\otimes \CC^\ell$.
		
		\item Composing this representation with the embeddings ${}^L U_i\to {}^L G_i$, $g\in \widehat{U}_i\mapsto (g,g^\star)$, ($i=k,\ell$) we obtain representations
		$$\displaystyle {}^L(U_k\times G_\ell)\to \GL(\CC^k\otimes \CC^\ell\oplus \CC^k\otimes \CC^\ell)$$
		and
		$$\displaystyle {}^L(U_k\times U_\ell)\to \GL(\CC^k\otimes \CC^\ell\oplus \CC^k\otimes \CC^\ell)$$
		that for simplicity we will also denote by the symbol $\itimes$.
		
		\item In particular, we have two representations
		$$\displaystyle \cR: {}^L \cG={}^L(U_m\times U_{n+1})\to \GL(\CC^m\otimes \CC^{n+1}\oplus \CC^m\otimes \CC^{n+1}),\; (S_m,S_{n+1})\mapsto S_m \itimes S_{n+1}$$
		$$\displaystyle \widetilde{\cR}: {}^L \widetilde{U}={}^L(U_m\times U_{m+1})\to \GL(\CC^m\otimes \CC^{m+1}\oplus \CC^m\otimes \CC^{m+1}),\; (S_m,S_{m+1})\mapsto S_m \itimes S_{m+1}$$
		that we will denote by $\cR$ and $\widetilde{\cR}$ respectively. The subspace
		$$\displaystyle V_-:=\left\langle (e_i\otimes e_j,0) \mid \substack{1\leqslant i\leqslant m \\ 1\leqslant j\leqslant n+1 \\ i+j> m+r+1}\right\rangle\oplus \left\langle (0,e_i\otimes e_j) \mid \substack{1\leqslant i\leqslant m \\ 1\leqslant j\leqslant n+1 \\ i+j> m+r+1}\right\rangle $$
		$$\displaystyle \left(\mbox{resp. } \widetilde{V}_-:=\left\langle (e_i\otimes e_j,0) \mid \substack{1\leqslant i\leqslant m \\ 1\leqslant j\leqslant m+1 \\ i+j> m+1}\right\rangle\oplus \left\langle (0,e_i\otimes e_j) \mid \substack{1\leqslant i\leqslant m \\ 1\leqslant j\leqslant m+1 \\ i+j> m+1}\right\rangle \right),$$
		where we denote by $(e_i)_{i=1}^k$ the standard basis of $\CC^k$ for any $k$, is stable by ${}^L T$ (resp. by ${}^L \widetilde{T}$) and we set $\cR_-(S):=\cR(S)\mid_{V_-}$ (resp. $\widetilde{\cR}_-(\widetilde{S}):=\widetilde{\cR}(\widetilde{S})_{\mid \widetilde{V}_-}$) for $S\in {}^L T$ (resp. $\widetilde{S}\in {}^L \widetilde{T}$).
		
		\item We will also denote by $\As_m$ the representation
		$$\displaystyle \As_m: {}^L G_r\to \GL(\CC^r\otimes \CC^r)$$
		given by $\As_m(g^{(1)},g^{(2)})=g^{(1)}\otimes g^{(2)}$ for $(g^{(1)},g^{(2)})\in\widehat{G}_r$ and $\As_m(c)=(-1)^m s$ where $s: \CC^r\otimes \CC^r\to \CC^r\otimes \CC^r$ is defined by $s(u\otimes v)=v\otimes u$ for $u,v\in \CC^r$.
		
		\item The Satake parameters of the unramified representations $\tau$, $\sigma_m$, $\sigma_{n+1}$ will be denoted by $S_\tau$, $S_m$ and $S_{n+1}$ respectively. These are semisimple conjugacy classes in ${}^L G_r$, ${}^L U_m$ and ${}^L U_{n+1}$ and to simplify some arguments we will choose representatives of them in ${}^L T_r$, ${}^L T_m^U$ and ${}^L T_{n+1}^U$ respectively. Thus, denoting by $Frob_v\in \Gal(E/F)$ the Frobenius at $v$, we have
		$$\displaystyle S_\tau\in \widehat{T}_r Frob_v,\;\; S_m\in \widehat{T}^U_m Frob_v \; \mbox{ and }\; S_{n+1}\in \widehat{T}^U_{n+1}Frob_v.$$
		We will also write $S=(S_m,S_{n+1})\in {}^L T$ for the Satake parameter of $\sigma=\sigma_m\boxtimes \sigma_{n+1}$.
	\end{itemize}
	
	By Shintani and Casselman-Shalika's formula \cite{Shint} \cite{CasSha}, we have
	\begin{equation}\label{eq2 unr computation}
	\displaystyle \cL^W(\tau(t)\phi_{\tau}^\circ)=\left\{\begin{array}{ll}
	\delta_{B_r}(t)^{1/2} ch_t(S_\tau) & \mbox{ if } t\in \Lambda_r^+, \\
	0 & \mbox{ otherwise.}
	\end{array} \right.
	\end{equation}
	On the other hand, according to the formulas given in \cite[Theorem 11.4]{Khou}, \cite[Proposition 6.4]{Liu} and \cite{LZhangWS}, when $S\in {}^L T$ is regular we have:
	\begin{align}\label{eq3 unr computation}
	\displaystyle & \cL^B(\phi_{m}^\circ\otimes \sigma_{n+1}(t)\phi_{n+1}^\circ)= \\
	\nonumber  & \left\{\begin{array}{ll}
	\displaystyle (\Delta_{m,v}^U)^{-1} \sum_{w\in W(\cG_v)} \frac{\det(1-q^{-1/2} \cR_-(wS))}{D_{\widehat{\cG}/\widehat{B}}(wS)} \chi_{t}((wS_{n+1})^{(r)}) \delta_P(t)^{1/2} & \mbox{ if } t\in \Lambda_r^{++} \\ \\
	0 & \mbox{ otherwise.}
	\end{array}
	\right.
	\end{align}
	Moreover, the above sum over $W(\cG_v)$ extends to a regular function on ${}^L T$ and the formula is still valid when we interpret the right-hand side in term of this extension. We will now prove the formula of the proposition assuming that $S$ is regular but, as it is an identity between rational functions, the extension to the non-regular case will follow.
	
	\begin{lemme}
		For $t\in \Lambda_r^{++}$, we have
		\begin{align}\label{eq5 unr computation}
		\displaystyle & \cL^B(\phi_{m}^\circ\otimes \sigma_{n+1}(t)\phi_{n+1}^\circ)= \\
		\nonumber & \sum_{w\in W(U_{n+1,v})} \frac{\det(1-q^{-1/2} S_m\itimes (wS_{n+1})^{(r)\star})}{D_{\widehat{U}_{n+1}/\widehat{B}^U_{n+1}}(wS_{n+1})} ch_t((wS_{n+1})^{(r)}) \delta_P(t)^{1/2}.
		\end{align}
	\end{lemme}
	
	\begin{preuve}
		First, we note that the function
		$$\displaystyle S_{n+1}\in \widehat{T}_{n+1}^U Frob_v\mapsto \det(1-q^{-1/2} S_m\itimes S_{n+1}^{(r)\star})ch_t(S_{n+1}^{(r)})$$
		is invariant under the action of $W(L_{n+1,v})$. Combining this with the identity (see Corollary \ref{corollary appendix})
		$$\displaystyle \sum_{w\in W(L_{n+1,v})} D_{\widehat{U}_{n+1}/\widehat{B}^U_{n+1}}(wS_{n+1})^{-1}=D_{\widehat{U}_{n+1}/\widehat{Q}_{n+1}}(S_{n+1})^{-1}$$
		we see that the right hand side of \eqref{eq5 unr computation} can be rewritten as
		$$\displaystyle \sum_{w\in W(L_{n+1,v})\backslash W(U_{n+1,v})} \frac{\det(1-q^{-1/2} S_m\itimes (wS_{n+1})^{(r)\star})}{D_{\widehat{U}_{n+1}/\widehat{Q}_{n+1}}(wS_{n+1})} ch_t((wS_{n+1})^{(r)}) \delta_P(t)^{1/2}.$$
		The natural projection $W(\cG_v)\to W(U_{n+1,v})$ induces a bijection
		$$\displaystyle W(\cL_v)\backslash W(\cG_v)\simeq W(L_{n+1,v})\backslash W(U_{n+1,v})$$
		and thus by \eqref{eq3 unr computation} we just need to establish, for every regular $S=(S_m,S_{n+1})\in {}^L T$ and $t\in \Lambda_r$, the identity
		\begin{equation}\label{eq3bis unr computation}
		\displaystyle (\Delta_{m,v}^U)^{-1} \sum_{w\in W(\cL_v)} \frac{\det(1-q^{-1/2} \cR_-(wS))}{D_{\widehat{\cG}/\widehat{B}}(wS)} \chi_{t}((wS_{n+1})^{(r)})=\frac{\det(1-q^{-1/2} S_m\itimes S_{n+1}^{(r)\star})}{D_{\widehat{U}_{n+1}/\widehat{Q}_{n+1}}(S_{n+1})} ch_t(S_{n+1}^{(r)}).
		\end{equation}
		We have decompositions
		\[\begin{aligned}
		\displaystyle D_{\widehat{\cG}/\widehat{B}}(S) =D_{\widehat{\cG}/\widehat{\cQ}}(S) D_{\widehat{\cL}/\widehat{B}_{\cL}}(S)=D_{\widehat{U}_{n+1}/\widehat{Q}_{n+1}}(S_{n+1}) D_{\widehat{\widetilde{U}}/\widehat{\widetilde{B}}}(\widetilde{S}) D_{\widehat{G}_r/\widehat{B}_r}(S_{n+1}^{(r)}),
		\end{aligned}\]
		\[\begin{aligned}
		\displaystyle \cR_-(S)=S_m\itimes S_{n+1}^{(r)\star} \oplus \widetilde{\cR}_-(\widetilde{S}).
		\end{aligned}\]
		and
		$$\displaystyle W(\cL_v)=W(\widetilde{U}_v)\times W(G_{r,v}).$$
		This leads to the following expression for the left hand side of \eqref{eq3bis unr computation}:
		\[\begin{aligned}
		\displaystyle & \frac{\det(1-q^{-1/2} S_m\itimes S_{n+1}^{(r)\star})}{D_{\widehat{U}_{n+1}/\widehat{Q}_{n+1}}(S_{n+1})} \times (\Delta_{m,v}^U)^{-1} \sum_{\widetilde{w}\in W(\widetilde{U}_v)} \frac{\det(1-q^{-1/2} \widetilde{\cR}_-(\widetilde{w}\widetilde{S}))}{D_{\widehat{\widetilde{U}}/\widehat{\widetilde{B}}}(\widetilde{w}\widetilde{S})} \\
		& \times \sum_{w_r\in W(G_{r,v})} \frac{\chi_{t}(w_rS_{n+1}^{(r)})}{D_{\widehat{G_r}/\widehat{B_r}}(w_rS_{n+1}^{(r)})}.
		\end{aligned}\]
		Furthermore, Weyl's character formula implies (see Proposition \ref{prop Weyl character formula})
		$$\displaystyle \sum_{w_r\in W(G_{r,v})} \frac{\chi_{t}(w_rS_{n+1}^{(r)})}{D_{\widehat{G_r}/\widehat{B_r}}(w_rS_{n+1}^{(r)})}=ch_t(S_{n+1}^{(r)})$$
		while, according to \cite[Proposition 6.4]{Liu}, we have
		$$\displaystyle (\Delta_{m,v}^U)^{-1} \sum_{\widetilde{w}\in W(\widetilde{U}_v)} \frac{\det(1-q^{-1/2} \widetilde{\cR}_-(\widetilde{w}\widetilde{S}))}{D_{\widehat{\widetilde{U}}/\widehat{\widetilde{B}}}(\widetilde{w}\widetilde{S})}=1.$$
		This shows the formula \eqref{eq3bis unr computation} and ends the proof of the lemma.
	\end{preuve}
	
	The next lemma is a consequence of the Cauchy identity \cite[Theorem 43.3]{BumpLieGroups}.
	
	\begin{lemme}
		Let $S_1,S_2\in \widehat{T_{r}}Frob_v$. Then, for $\Re(s)$ sufficiently large we have
		$$\displaystyle \sum_{t\in \Lambda_r^{++}} \lvert \det t\rvert_{E_v}^{1/2+s}ch_t(S_{1}) ch_t(S_2) =\det(1-q^{-1/2} S_{1,s}\itimes S_{2})^{-1}$$
		where we have set $S_{1,s}:=q^{-s}S_1$ 
	\end{lemme}
	
	Combining the two above lemmas with \eqref{eq1 unr computation}, \eqref{eq2 unr computation} as well as the identity
	$$\displaystyle \delta_P(t)^{1/2}\delta_{B_r}(t)^{1/2} \delta_{Q_n}(t)^{1/2} \delta_{P'}(t)^{-1}=\lvert \det t\rvert_{E_v}^{1/2},\;\; t\in \Lambda_r,$$
	we obtain
	\begin{align}\label{eq4 unr computation}
     &\cL_s^{U'}(\phi_{n,s}^\circ\otimes \phi_{n+1}^\circ) \\
          \nonumber & =\sum_{w\in W(U_{n+1,v})} \frac{\det(1-q^{-1/2} S_m\itimes (wS_{n+1})^{(r)\star})}{D_{\widehat{U}_{n+1}/\widehat{B}^U_{n+1}}(wS_{n+1})} \det(1-q^{-1/2} S_{\tau,s}\itimes (wS_{n+1})^{(r)})^{-1}.
	\end{align}
	For every $w\in W(U_{n+1,v})$ we have the identity
	$$\displaystyle S_{\tau,s}\itimes wS_{n+1}=S_{\tau,s}\itimes (wS_{n+1})^{(r)}\oplus S_{\tau,s}\itimes (wS_{n+1})^{(m+1)}\otimes S_{\tau,s}\itimes (wS_{n+1})^{(r)\star}.$$
	It follows that
        \begin{align*}
           & L(\frac{1}{2},\tau_s\times \sigma_{n+1})^{-1}=\det(1-q^{-1/2} S_{\tau,s}\itimes S_{n+1}) \\
	& =\det(1-q^{-1/2} S_{\tau,s}\itimes (wS_{n+1})^{(r)})\det(1-q^{-1/2} S_{\tau,s}\itimes (wS_{n+1})^{(m+1)})\det(1-q^{-1/2} S_{\tau,s}\itimes (wS_{n+1})^{(r)\star}).
	\end{align*}
      
	Thus, \eqref{eq4 unr computation} can be rewritten as
	\begin{align*}
	\displaystyle & \cL_s^{U'}(\phi_{n,s}^\circ\otimes \phi_{n+1}^\circ) = L(\frac{1}{2},\tau_s\times \sigma_{n+1}) \sum_{w\in W(U_{n+1,v})}\\
	\nonumber &  \frac{\det(1-q^{-1/2} S_m\itimes (wS_{n+1})^{(r)\star})}{D_{\widehat{U}_{n+1}/\widehat{B}^U_{n+1}}(wS_{n+1})} \det(1-q^{-1/2} S_{\tau,s}\itimes (wS_{n+1})^{(r)\star})\det(1-q^{-1/2} S_{\tau,s}\itimes (wS_{n+1})^{(m+1)}).
	\end{align*}
	Since we have
	$$\displaystyle L(1,\tau_s\times \sigma_m)=\det(1-q^{-1}S_{\tau,s}\itimes S_m)^{-1},\;\; L(1,\tau_s,\As^{(-1)^m})=\det(1-q^{-1}\As_m(S_{\tau,s}))^{-1}$$
	we see that the proposition is now reduced to the equality
	\begin{align}\label{eq6 unr computation}
	\displaystyle & \det(1-q^{-1}S_{\tau,s}\itimes S_m)\det(1-q^{-1}\As_m(S_{\tau,s}))=\sum_{w\in W(U_{n+1,v})}  \\
	\nonumber & \frac{\det(1-q^{-1/2} S_m\itimes (wS_{n+1})^{(r)\star})}{D_{\widehat{U}_{n+1}/\widehat{B}^U_{n+1}}(wS_{n+1})} \det(1-q^{-1/2} S_{\tau,s}\itimes (wS_{n+1})^{(r)\star})\det(1-q^{-1/2} S_{\tau,s}\itimes (wS_{n+1})^{(m+1)}).
	\end{align}
	
	To prove the above identity, we first show that the right hand side of \eqref{eq6 unr computation} considered as a function of the regular element $S_{n+1}\in \widehat{T}_{n+1}^U Frob_v$ is constant. For this, we first note that the function
        \begin{align*}
        &  S_{n+1}\in \widehat{T}_{n+1}^U Frob_v\mapsto \\
          &\det(1-q^{-1/2} S_m\itimes (wS_{n+1})^{(r)\star})\det(1-q^{-1/2} S_{\tau,s}\itimes (wS_{n+1})^{(r)\star})\det(1-q^{-1/2} S_{\tau,s}\itimes (wS_{n+1})^{(m+1)})
        \end{align*}
	is (the restriction of) a linear combination of characters of ${}^L T_{n+1}^U$. Thus, by Proposition \ref{prop Weyl character formula}, it suffices to check that for every character $\chi$ appearing in this linear combination the sum $\rho+\chi\mid_{\widehat{T}_{n+1}^U}$, where $\rho\in X^*(\widehat{T}_{n+1}^U)$ stands for the half-sum of positive roots (with respect to $\widehat{B}_{n+1}^U$), is either singular or conjugate, under the full Weyl group of $(\widehat{U}_{n+1},\widehat{T}_{n+1}^U)$, to $\rho$. Let $\chi$ be such a character. Using the natural isomorphism $X^*(\widehat{T}_{n+1}^U)\simeq \mathbb{Z}^{n+1}$, we have $\chi\mid_{\widehat{T}_{n+1}^U}=(\lambda_1,\ldots,\lambda_{n+1})$ where the $\lambda_i$'s are integers satisfying
        \begin{align*}
          -m-r\leqslant \lambda_i\leqslant 0 &\text{ for } 0\leqslant i\leqslant r,\;\\
          -r\leqslant \lambda_i\leqslant r &\text{ for } r+1\leqslant i\leqslant r+m+1,\; \\
          0\leqslant \lambda_i\leqslant m+r &\text{ for } r+m+2\leqslant i\leqslant n+1.
        \end{align*}
	Furthermore, we have $\rho=(\frac{n}{2},\frac{n-2}{2},\ldots, -\frac{n}{2})$ and from the above inequality we see that the coordinates of $\rho+\chi\mid_{\widehat{T}_{n+1}^U}$ are all integers or half-integers between $-\frac{n}{2}$ and $\frac{n}{2}$. The claim about $\chi\mid_{\widehat{T}_{n+1}^U}+\rho$ readily follows and we therefore deduce that the right-hand side of \eqref{eq6 unr computation} is indeed independent of $S_{n+1}$.
	
	Let $A,B\in \GL_r(\CC)$ be such that $S_{\tau,s}=(A,B)Frob_v$. Then, plugging
	$$\displaystyle S_{n+1}=\begin{pmatrix}q^{1/2}A^\star & & \\ & A_{m+1} & \\ & & q^{-1/2} B \end{pmatrix} Frob_v\in \widehat{T}_{n+1}^U Frob_v$$
	in the right-hand side of \eqref{eq6 unr computation}, where the matrix $A_{m+1}\in \GL_{m+1}(\CC)$ is chosen such that $S_{n+1}$ is regular, the term indexed by $w\in W(U_{n+1,v})$ is nonzero only when $w\in W(L_{n+1,v})$. It follows that this sum equals
        \begin{align*}
        &   \sum_{w\in W(L_{n+1,v})} \frac{\det(1-q^{-1/2} S_m\itimes (wS_{n+1})^{(r)\star})}{D_{\widehat{U}_{n+1}/\widehat{B}^U_{n+1}}(wS_{n+1})} \det(1-q^{-1/2} S_{\tau,s}\itimes (wS_{n+1})^{(r)\star})\times\\
          &\det(1-q^{-1/2} S_{\tau,s}\itimes (wS_{n+1})^{(m+1)}) \\
        &=\det(1-q^{-1/2} S_m\itimes S_{n+1}^{(r)\star})\det(1-q^{-1/2} S_{\tau,s}\itimes S_{n+1}^{(r)\star})\det(1-q^{-1/2} S_{\tau,s}\itimes S_{n+1}^{(m+1)}) \times\\
        &\sum_{w\in W(L_{n+1,v})} D_{\widehat{U}_{n+1}/\widehat{B}^U_{n+1}}(wS_{n+1})^{-1} \\
	& =\det(1-q^{-1/2} S_m\itimes S_{n+1}^{(r)\star})\det(1-q^{-1/2} S_{\tau,s}\itimes S_{n+1}^{(r)\star})\times\\
          &\det(1-q^{-1/2} S_{\tau,s}\itimes S_{n+1}^{(m+1)}) D_{\widehat{U}_{n+1}/\widehat{Q}_{n+1}}(S_{n+1})^{-1}
	\end{align*}
      
	where the last equality follows from Corollary \ref{corollary appendix}. By direct computation, we have
        \begin{align*}
          \det(1-q^{-1/2} S_m\itimes S_{n+1}^{(r)\star})=\det(1-q^{-1}S_m\itimes S_{\tau,s}),\; \\
          \det(1-q^{-1/2} S_{\tau,s}\itimes S_{n+1}^{(r)\star})=\det(1-q^{-1} S_{\tau,s}\itimes S_{\tau,s}),
        \end{align*}
	$$\displaystyle \det(1-q^{-1/2} S_{\tau,s}\itimes S_{n+1}^{(m+1)})=\det(1-q^{-1/2} S_{\tau,s}\itimes S_{m+1}),$$
	$$\displaystyle D_{\widehat{U}_{n+1}/\widehat{Q}_{n+1}}(S_{n+1})=\det(1-q^{-1}As_{m+1}(S_{\tau,s}))\det(1-q^{-1/2} S_{\tau,s}\itimes S_{m+1})$$
	where we have set $S_{m+1}=A_{m+1}Frob_v\in \widehat{T}_{m+1}^U Frob_v$.
	
	
\end{preuve}
\end{paragr}	

\subsection{Reduction to the corank one case}\label{ssect: reduction}
	
\begin{paragr}  	In this subsection, we let:
	\begin{itemize}
		\item $\sigma_m$ and $\sigma_{n+1}$ be cuspidal automorphic representations of $U_m(\AAA)$ and $U_{n+1}(\AAA)$  respectively;
		
		\item $\tau$ be an irreducible automorphic representation of $G_r(\AAA)$ that is induced from a unitary cuspidal representation meaning that there exist a parabolic subgroup $R=MN_R\subset G_r$ as well as a unitary cuspidal automorphic representation $\kappa$ of $M(\AAA)$ such that
		$$\displaystyle \tau=\{E_R^{G_r}(\phi,0)\mid \phi\in I_{R(\AAA)}^{G_r(\AAA)}(\kappa) \}.$$
	\end{itemize}
	 Moreover, we henceforth identify $\cA_{Q_n,\tau\boxtimes \sigma_m}(U_n)$ with the parabolic induction $\sigma_n:=I_{Q_n(\AAA)}^{U_n(\AAA)}(\tau\boxtimes \sigma_m)$ so that for $\phi_n \in \sigma_n$ and any $s\in \CC$ we have $\phi_{n,s}\in \sigma_{n,s}:=I_{Q_n(\AAA)}^{U_n(\AAA)}(\tau_s\boxtimes \sigma_m)$. Then, by Proposition \ref{prop relation global periods} and Proposition \ref{prop unr computation}, for $\phi_n\in \sigma_n$ and $\phi_{n+1}\in \sigma_{n+1}$ we have
	\begin{align}\label{eq1 relation global periods}
		\displaystyle \cP_{U'}(E_{Q_n}^{U_n}(\phi_n,s)\otimes \phi_{n+1})= & \frac{\vol(K_n^{U,T})}{\vol(K_n^{U,T}\cap \cB'(\AAA^T))}\frac{L^T(\frac{1}{2}+s,\tau\times \sigma_{n+1})}{L^T(1+s,\tau\times \sigma_m) L^T(1+2s,\tau,\As^{(-1)^m})} \\
		\nonumber & \times \int_{\cB'(F_T)\backslash U'(F_T)} \cP_{\cB^L,\psi^L_{\cB}}(\phi_{n,s}(h)\otimes \sigma_{n+1}(h)\phi_{n+1})dh
	\end{align}
	for $\Re(s)\gg 0$ and where $T\supset S$ is any sufficiently large finite set of places such that $\phi_n$ is $K_n^{U,T}$-invariant and $\phi_{n+1}$ is $K_{n+1}^{U,T}$-invariant.
      \end{paragr}	

      \begin{paragr}

      \begin{proposition}\label{prop:reduc-GGP-cork1}
		The following are equivalent:
		\begin{enumerate}
			\item There exist $\phi_m\in \sigma_m$ and $\phi_{n+1}\in \sigma_{n+1}$ such that $\cP_{\cB,\psi_{\cB}}(\phi_m\otimes \phi_{n+1})\neq 0$;
			
			\item There exist $\phi_n\in \sigma_n$, $\phi_{n+1}\in \sigma_{n+1}$ and $s\in \CC$ such that $E_{Q_n}^{U_n}(\phi_n,.)$ has no pole at $s$ and $\cP_{U'}(E_{Q_n}^{U_n}(\phi_n,s)\otimes \phi_{n+1})\neq 0$.
		\end{enumerate}
	\end{proposition}
	
	\begin{preuve}
		First, we remark that since there exists $\phi_{\tau}\in \tau$ whose Whittaker period $\int_{[N_r]} \phi_{\tau}(u)\psi_{-r}(u)du$ is nonzero, assertion 1. is equivalent to the non-vanishing of $\cP_{\cB^L,\psi^L_{\cB}}$ on $\tau\boxtimes \sigma_m\boxtimes \sigma_{n+1}$. Then, the equivalence $2.\Leftrightarrow 1.$ follows from \eqref{eq1 relation global periods} and the last part of Proposition \ref{proposition induced local intertwining}. 
	\end{preuve}
      \end{paragr}

      \begin{paragr}
              Choose isomorphisms $\tau\simeq \bigotimes_v' \tau_v$, $\sigma_m\simeq \bigotimes_v' \sigma_{m,v}$ an $\sigma_{n+1}\simeq\bigotimes'_v \sigma_{n+1,v}$. This induces an isomorphism $\sigma_n\simeq \bigotimes_v' \sigma_{n,v}$ where $\sigma_{n,v}:= I_{Q_n(F_v)}^{U_n(F_v)}(\tau_v\boxtimes \sigma_{m,v})$. We assume henceforth that for every place $v$ the representations $\tau_v$, $\sigma_{m,v}$ and $\sigma_{n+1,v}$ are all tempered. Let $\phi_\tau\in \tau$, $\phi_m\in \sigma_m$, $\phi_n\in \sigma_n$ and $\phi_{n+1}\in\sigma_{n+1}$ be factorizable vectors i.e.
	$$\displaystyle \phi_\tau=\bigotimes_v' \phi_{\tau,v},\; \phi_m=\bigotimes'_v \phi_{m,v},\; \phi_n=\bigotimes_v' \phi_{n,v},\; \phi_{n+1}=\bigotimes_v' \phi_{n+1,v}$$
	where $\phi_{\tau,v}\in \tau_v$, $\phi_{m,v}\in \sigma_{m,v}$, $\phi_{n,v}\in \sigma_{n,v}$ and $\phi_{n+1,v}\in \sigma_{n+1,v}$. 
	
	We equip $\sigma_m$, $\sigma_{n+1}$, $\tau$ and $\sigma_n$ with invariant inner products as follows:
\begin{itemize}
	\item We endow $\sigma_m$, $\sigma_{n+1}$ with the Petersson inner products $\langle .,.\rangle_{Pet}$ i.e. the $L^2$ inner products with respect to the Tamagawa measures on $[U_m]$ and $[U_{n+1}]$ respectively.
	
	\item On $\tau$ we put the inner product defined by
	$$\displaystyle \langle E_R^{G_r}(\phi,0), E_R^{G_r}(\phi',0)\rangle_\tau=\int_{R(\AAA)\backslash G_r(\AAA)} \int_{[M]^1} \phi(mg) \overline{\phi'(mg)} dm dg$$
	for $\phi,\phi'\in I_{R(\AAA)}^{G_r(\AAA)}(\kappa)$.
	
	\item $\sigma_n$ is equipped with the inner product induced from that on $\tau\boxtimes \sigma_m$, here denoted $\langle .,.\rangle_{\tau\boxtimes \sigma_m}$, that is:
	$$\displaystyle \langle \phi,\phi'\rangle_{\sigma_n}=\int_{Q_n(\AAA)\backslash U_n(\AAA)} \langle \phi(g),\phi'(g)\rangle_{\tau\boxtimes \sigma_m} dg,\mbox{ for } \phi,\phi'\in \sigma_n.$$
\end{itemize}
	We also fix factorizations of these inner products on $\sigma_{m}$, $\sigma_{n+1}$, $\tau$ and $\sigma_n$ into local invariant inner products. Following Section \ref{Sect local periods}, this allows to define, for every place $v$ of $F$, local periods $\cP_{U',v}$, $\cP_{\cB,\psi_{\cB},v}$, $\cP_{N_r,\psi_{-r},v}$ and $\cP_{\cB^L,\psi^L_{\cB},v}$ on $\sigma_{n,v}\boxtimes \sigma_{n+1,v}$, $\sigma_{m,v}\boxtimes \sigma_{n,v}$, $\tau_v$ and $\tau_v\boxtimes \sigma_{m,v}\boxtimes \sigma_{n+1,v}$ respectively. Furthermore, for almost all $v$ we have
	\[\begin{aligned}
		\displaystyle \cP_{U',v}(\phi_{n,v}\otimes \phi_{n+1,v},\phi_{n,v}\otimes \phi_{n+1,v})=\Delta_{U_{n+1},v} \frac{L(\frac{1}{2},\sigma_{n,v}\times \sigma_{n+1,v})}{L(1,\sigma_{n,v},\Ad)L(1,\sigma_{n+1,v},\Ad)}
	\end{aligned}\]
	\begin{equation*}
		\displaystyle \cP_{\cB,\psi_{\cB},v}(\phi_{m,v}\otimes \phi_{n+1,v},\phi_{m,v}\otimes \phi_{n+1,v})=\Delta_{U_{n+1},v} \frac{L(\frac{1}{2},\sigma_{m,v}\times \sigma_{n+1,v})}{L(1,\sigma_{m,v},\Ad)L(1,\sigma_{n+1,v},\Ad)}
	\end{equation*}
	\begin{equation*}
		\displaystyle \cP_{N_r,\psi_{-r},v}(\phi_{\tau,v},\phi_{\tau,v})=\frac{\Delta_{G_r,v}}{L(1,\tau_{v},\Ad)}
	\end{equation*}
	where we have set $\Delta_{U_{n+1},v}=\prod_{i=1}^{n+1} L(i,\eta_{E/F,v}^i)$ and $\Delta_{G_r,v}=\prod_{i=1}^r \zeta_{E_v}(i)$. Note that the last two equalities above imply that
	\[\begin{aligned}
		\displaystyle & \cP_{\cB^L,\psi^L_{\cB},v}(\phi_{\tau,v}\otimes \phi_{m,v}\otimes \phi_{n+1,v},\phi_{\tau,v}\otimes \phi_{m,v}\otimes \phi_{n+1,v})= \\
		& \Delta_{U_{n+1},v} \Delta_{G_r,v} \frac{L(\frac{1}{2},\sigma_{m,v}\times \sigma_{n+1,v})}{L(1,\sigma_{m,v},\Ad)L(1,\sigma_{n+1,v},\Ad)L(1,\tau_{v},\Ad)}
	\end{aligned}\]
	for almost all $v$. Given all these identities, it makes sense to define
	\[\begin{aligned}
		\displaystyle & \prod_v' \cP_{U',v}(\phi_{n,v}\otimes \phi_{n+1,v},\phi_{n,v}\otimes \phi_{n+1,v}):= \\
		& \Delta^T_{U_{n+1}} \frac{L^T(\frac{1}{2},\sigma_{n}\times \sigma_{n+1})}{L^{T,*}(1,\sigma_{n},\Ad)L^T(1,\sigma_{n+1},\Ad)} \prod_{v\in T} \cP_{U',v}(\phi_{n,v}\otimes \phi_{n+1,v},\phi_{n,v}\otimes \phi_{n+1,v}),
	\end{aligned}\]
	\[\begin{aligned}
		\displaystyle & \prod_v'\cP_{\cB,\psi_{\cB},v}(\phi_{m,v}\otimes \phi_{n+1,v},\phi_{m,v}\otimes \phi_{n+1,v}):= \\
		& \Delta^T_{U_{n+1}} \frac{L^T(\frac{1}{2},\sigma_{m}\times \sigma_{n+1})}{L^T(1,\sigma_{m},\Ad)L^T(1,\sigma_{n+1},\Ad)}\prod_{v\in T} \cP_{\cB,\psi_{\cB},v}(\phi_{m,v}\otimes \phi_{n+1,v},\phi_{m,v}\otimes \phi_{n+1,v}),
	\end{aligned}\]
	\begin{equation*}
		\displaystyle \prod_v'\cP_{N_r,\psi_{-r},v}(\phi_{\tau,v},\phi_{\tau,v}):=\frac{\Delta^{T,*}_{G_r}}{L^{T,*}(1,\tau,\Ad)}\prod_{v\in T} \cP_{N_r,\psi_{-r},v}(\phi_{\tau,v},\phi_{\tau,v})
	\end{equation*}
	for any sufficiently large finite set $T$ of places where we have set
	$$\displaystyle \Delta^T_{U_{n+1}}=\prod_{i=1}^{n+1} L^T(i,\eta_{E/F}^i),\;\;\; \Delta^{T,*}_{G_r}=\zeta_E^{T,*}(1)\prod_{i=2}^r \zeta_E^T(i)$$
	and $L^{T,*}(1,\sigma_{n},\Ad)$, $L^{T,*}(1,\tau,\Ad)$ stand for the regularized values
	$$\displaystyle L^{T,*}(1,\sigma_{n},\Ad):=\left((s-1)^a L^T(s,\sigma_n,\Ad) \right)_{s=1},\;\;\; L^{T,*}(1,\tau,\Ad):=\left((s-1)^a L^T(s,\tau,\Ad) \right)_{s=1}$$
	with $a=\dim(A_M)$
	respectively.
	We define similarly
	$$\displaystyle \prod_v' \cP_{\cB^L,\psi^L_{\cB},v}(\phi_{\tau,v}\otimes \phi_{m,v}\otimes \phi_{n+1,v},\phi_{\tau,v}\otimes \phi_{m,v}\otimes \phi_{n+1,v}).$$
	
	Of course, the previous discussion applies verbatim when we replace $\tau$ and $\sigma_n$ by $\tau_s$ and $\sigma_{n,s}$ for every $s\in i\RR$.
	
	\begin{proposition}
		Let $c\in \CC$ and assume that for every factorizable vectors $\phi_n\in \sigma_n$, $\phi_{n+1}\in \sigma_{n+1}$ and every $s\in i\RR$ we have
		$$\displaystyle \left\lvert \cP_{U'}(E_{Q_n}^{U_n}(\phi_n,s)\otimes \phi_{n+1})\right\rvert^2=c \prod_v' \cP_{U',v}(\phi_{n,s,v}\otimes \phi_{n+1,v},\phi_{n,s,v}\otimes \phi_{n+1,v}).$$
		Then, for every factorizable vectors $\phi_m\in \sigma_m$ and $\phi_{n+1}\in \sigma_{n+1}$ we have
		$$\displaystyle \left\lvert \cP_{\cB,\psi_{\cB}}(\phi_m\otimes \phi_{n+1})\right\rvert^2=c \prod_v' \cP_{\cB,\psi_{\cB},v}(\phi_{m,v}\otimes \phi_{n+1,v},\phi_{m,v}\otimes \phi_{n+1,v}).$$
	\end{proposition}
	
	\begin{preuve}
		Let $\phi_n\in \sigma_n$ and $\phi_{n+1}\in \sigma_{n+1}$ be factorizable vectors and let $T$ be a finite sett of places of $F$ that we will assume throughout to be sufficiently large. By \eqref{eq1 relation global periods}, we have
		\begin{align}\label{eq1 preuve reduction II}
			\displaystyle & \left\lvert \cP_{U'}(E_{Q_n}^{U_n}(\phi_n,s)\otimes \phi_{n+1})\right\rvert^2= \left\lvert \frac{L^T(\frac{1}{2}+s,\tau\times \sigma_{n+1})}{L^T(1+s,\tau\times \sigma_m) L^T(1+2s,\tau,\As^{(-1)^m})}\right\rvert^2 \\
			\nonumber &  \int_{(\cB'(F_T)\backslash U'(F_T))^2} \cP_{\cB^L,\psi_{\cB}^L}(\phi_{n,s}(h_1)\otimes \sigma_{n+1}(h_1)\phi_{n+1}) \overline{\cP_{\cB^L,\psi_{\cB}^L}(\phi_{n,s}(h_2)\otimes \sigma_{n+1}(h_2)\phi_{n+1})} dh_2dh_1
		\end{align}
		for $s\in i\RR$. On the other hand, from the hypothesis, Proposition \ref{prop relation local periods} and \eqref{eq product nu} we obtain
		\begin{align}\label{eq2 preuve reduction II}
			\displaystyle & \left\lvert \cP_{U'}(E_{Q_n}^{U_n}(\phi_n,s)\otimes \phi_{n+1})\right\rvert^2=c \Delta^{T,*}_{G_r}\Delta_{U_{n+1}}^T \frac{L^T(\frac{1}{2},\sigma_{n,s}\times \sigma_{n+1})}{L^{T,*}(1,\sigma_{n,s},\Ad)L^T(1,\sigma_{n+1},\Ad)} \\
			\nonumber &  \prod_{v\in T} \int_{(\cB'(F_v)\backslash U'(F_v))^2} \cP_{\cB^L,\psi^L_{\cB},v}(\phi_{n,s,v}(h_1)\otimes \sigma_{n+1,v}(h_2)\phi_{n+1,v},\phi_{n,s,v}(h_1)\otimes \sigma_{n+1,v}(h_2)\phi_{n+1,v} ) dh_2dh_1
		\end{align}
		for $s\in i\RR$. From \eqref{eq1 preuve reduction II}, \eqref{eq2 preuve reduction II}, the last part of Proposition \ref{proposition induced local intertwining} as well as the identity of partial $L$-functions for $s\in i\RR$
		$$\displaystyle \frac{L^T(\frac{1}{2},\sigma_{n,s}\times \sigma_{n+1})}{L^{T,*}(1,\sigma_{n,s},\Ad)}=\frac{L^T(\frac{1}{2},\sigma_{m}\times \sigma_{n+1})}{L^{T,*}(1,\tau,\Ad)L^{T}(1,\sigma_{m},\Ad)} \left\lvert \frac{L^T(\frac{1}{2}+s,\tau\times \sigma_{n+1})}{L^T(1+s,\tau\times \sigma_m) L^T(1+2s,\tau,\As^{(-1)^m})}\right\rvert^2,$$
		we deduce the equality
		\begin{align}\label{eq3 preuve reduction II}
			\displaystyle \left\lvert \cP_{\cB^L,\psi_{\cB}^L}(\phi\otimes \phi_{n+1})\right\rvert^2= & c\Delta^{T,*}_{G_r}\Delta_{U_{n+1}}^T \frac{L^T(\frac{1}{2},\sigma_{m}\times \sigma_{n+1})}{L^{T,*}(1,\tau,\Ad)L^{T}(1,\sigma_{m},\Ad)L^T(1,\sigma_{n+1},\Ad)}  \\
			\nonumber & \times \prod_{v\in T} \cP_{\cB^L,\psi^L_{\cB},v}(\phi_v\otimes \phi_{n+1,v},\phi_{v}\otimes \phi_{n+1,v} ) \\
			\nonumber & =c\prod_v' \cP_{\cB^L,\psi^L_{\cB},v}(\phi_v\otimes \phi_{n+1,v},\phi_{v}\otimes \phi_{n+1,v} )
		\end{align}
		for every factorizable vector $\phi=\bigotimes_v' \phi_v\in \tau\boxtimes \sigma_m$. Furthermore, for $\phi_\tau\in \tau$ and $\phi_m\in \sigma_m$ two factorizable vectors, we have
		$$\displaystyle \cP_{\cB^L,\psi_{\cB}^L}(\phi_\tau\otimes \phi_m\otimes \phi_{n+1})=\cP_{N_r,\psi_{-r}}(\phi_\tau)\cP_{\cB,\psi_{\cB}}(\phi_m\otimes \phi_{n+1})$$
		and
		$$\displaystyle \cP_{\cB^L,\psi^L_{\cB},v}(\phi_{\tau,v}\otimes \phi_{m,v}\otimes \phi_{n+1,v},\phi_{\tau,v}\otimes \phi_{m,v}\otimes \phi_{n+1,v} )=\cP_{N_r,\psi_{-r},v}(\phi_{\tau,v},\phi_{\tau,v}) \cP_{\cB,\psi_{\cB},v}(\phi_{m,v}\otimes \phi_{n+1,v},\phi_{m,v}\otimes \phi_{n+1,v} )$$
		for every place $v$ as well as the identity (cf. \cite[Eq. (11.3)]{FLO})
		\begin{equation}\label{eq4 preuve reduction II}
			\displaystyle \left\lvert \cP_{N_r,\psi_{-r}}(\phi_\tau)\right\rvert^2=\prod_v' \cP_{N_r,\psi_{-r},v}(\phi_{\tau,v},\phi_{\tau,v}).
		\end{equation}
		The proposition can now be deduced by dividing the identity \eqref{eq4 preuve reduction II} by \eqref{eq3 preuve reduction II}. Note that such a deduction is valid as the Whittaker period $\cP_{N_r,\psi_{-r}}$ is known to be nonzero on $\tau$.
	\end{preuve}
      \end{paragr}
      
\appendix

\section{Weyl character formula for non-connected groups}\label{appendix Weyl character formula}

\begin{paragr}
  
Let $\widehat{G}$ be a connected complex reductive group and $\Gamma$ be a finite group acting on $\widehat{G}$ by holomorphic automorphisms preserving a given Borel pair $(\widehat{B},\widehat{T})$. Set ${}^L G=\widehat{G}\rtimes \Gamma$, ${}^L T=\widehat{T}\rtimes \Gamma$ and ${}^L B=\widehat{B}\rtimes \Gamma$. Let $\widehat{B}\subset \widehat{G}$ be a Borel subgroup and ${}^L B\subset {}^L G$ be the normalizer of $\widehat{B}$. We say that an element $S\in {}^L T$ is {\em regular} if the neutral component of its centralizer is contained in $\widehat{T}$.

Let $X^*(\widehat{T})$ be the group of algebraic characters of $\widehat{T}$ and $X^*(\widehat{T})^+\subset X^*(\widehat{T})$ be teh subset of dominant elements. Note that the subgroup of $\Gamma$-invariant characters $X^*(\widehat{T})^\Gamma$ can be identified with the set of algebraic characters ${}^L T\to \CC^\times$ that are trivial on $\Gamma$. For $\chi\in X^*(\widehat{T})^+\cap X^*(\widehat{T})^\Gamma$, any irreducible representation $(\pi_\chi,V_\chi)$ of $\widehat{G}$ with highest weight $\chi$ can be extended in a unique way to a representation of ${}^L G$ such that $\Gamma$ acts trivially on the line of highest weight vectors. We shall denote by $ch_\chi$ the character of that representation of ${}^L G$ (i.e. $ch_\chi(S)=\Tr \pi_\chi(S)$ for $S\in {}^L G$).
\end{paragr}

\begin{paragr}
  
Let $\widehat{W}=\Norm_{\widehat{G}}(\widehat{T})/\widehat{T}$ be the Weyl group of $\widehat{T}$ in $\widehat{G}$. Then, $\Gamma$ acts in a natural way on $\widehat{W}$ and the subgroup of fixed points $\widehat{W}^\Gamma$ can be identified with $W=\Norm_{\widehat{G}}({}^L T)/\widehat{T}$. Let $\rho\in \frac{1}{2}X^*(\widehat{T})$ the half-sum of the positive roots with respect to $\widehat{B}$. The dot action of $\widehat{W}$ on $X^*(\widehat{T})$ is defined by
$$\displaystyle w\cdot \chi=w(\chi+\rho)-\rho,\;\; (w,\chi)\in \widehat{W}\times X^*(\widehat{T}).$$
For $\chi\in X^*(\widehat{T})$, we recall that the following alternative. Either:
\begin{itemize}
	\item $\chi+\rho$ is singular, i.e. there exists a coroot $\alpha^\vee$ such that $\langle \alpha^\vee, \chi+\rho\rangle=0$;
	
	\item or there exists a unique $w_\chi\in \widehat{W}$ such that $w_\chi\cdot \chi\in X^*(\widehat{T})$.
\end{itemize}
Note that, by unicity, for $\chi\in X^*(\widehat{T})^\Gamma$ such that $\chi+\rho$ is non singular, we have $w_\chi\in W$.

For every standard parabolic subgroup $\widehat{B}\subset \widehat{Q}\subset \widehat{G}$ that is $\Gamma$-stable, we set
$$\displaystyle D_{\widehat{G}/\widehat{Q}}(S)= \det(1-Ad(S))\mid Lie(\widehat{G})/Lie(\widehat{Q})$$
for $S\in {}^LT$.

\begin{proposition}\label{prop Weyl character formula}
	Let $\chi\in X^*(\widehat{T})^\Gamma$ and $F\in \Gamma$.
	\begin{enumerate}
		\item If $\chi+\rho$ is singular, we have
		$$\displaystyle \sum_{w\in W} \frac{\chi(wS)}{D_{\widehat{G}/\widehat{B}}(wS)}=0$$
		for every regular $S\in {}^L T$.
		
		\item If $\chi+\rho$ is non singular, setting $\chi^+=w_\chi\cdot \chi$, there exists a root of unity $\epsilon_\chi\in \CC^\times$ such that
		$$\displaystyle \sum_{w\in W} \frac{\chi(wS)}{D_{\widehat{G}/\widehat{B}}(wS)}=\epsilon_\chi ch_\chi(S)$$
		for every regular $S\in \widehat{T} F$. Moreover, if $\chi$ is dominant we have $\epsilon_\chi=1$.

	\end{enumerate}
\end{proposition}

\begin{preuve}
	Let $X=\widehat{G}/\widehat{B}$ be the flag variety of $\widehat{G}$. The action of $\widehat{G}$ on $X$ naturally extends to ${}^L G$ e.g. because we can also write $X={}^L G/{}^L B$. Let $\mathcal{L}_\chi$ be the ${}^L G$-equivariant line bundle on $X$ such that the action of ${}^LB$ on the fiber above $1\in X$ is given by $\chi$. Then, by the Borel-Weil-Bott theorem we have:
	\begin{itemize}
		\item If $\chi+\rho$ is singular, $H^i(X,\mathcal{L}_\chi)=0$ for $i\geqslant 0$,
		
		\item otherwise (i.e. if $\chi+\rho$ is non singular), there exists a unique $i\geqslant 0$ such that $H^i(X,\mathcal{L}_\chi)\neq 0$ and moreover $H^i(X,\mathcal{L}_\chi)\simeq V_\chi$ as $\widehat{G}$-modules.		
	\end{itemize}
	Let $S\in \widehat{T} F$ be regular. We apply Atiyah-Bott fixed point theorem \cite[Theorem 4.12]{AtiBott} to the action of $S$ on the pair $(X,\mathcal{L}_\chi)$. First note that the set of fixed points of $S$ in $X$ is precisely the image of the natural embedding $W\subset X$. Moreover, for $w\in W$ the action of $S$ on the fiber $(\mathcal{L}_\chi)_w$ (resp. on the tangent space $T_wX$ identified in a natural way with $Lie(\widehat{G})/Lie(\widehat{B})$) is the multiplication by $\chi(w^{-1}S)$ (resp. the adjoint operator $Ad(w^{-1}S)$). Given this as well as the above description of the cohomology groups of $\mathcal{L}_\chi$, the Atiyah-Bott fixed point theorem implies directly the proposition when $\chi+\rho$ is singular whereas in the non-singular case it gives
	$$\displaystyle \sum_{w\in W} \frac{\chi(wS)}{D_{\widehat{G}/\widehat{B}}(wS)}=\Tr(S\mid H^i(X,\mathcal{L}_\chi)).$$
	However, as $H^i(X,\mathcal{L}_\chi)\simeq V_\chi$ as $\widehat{G}$-modules and $V_\chi$ is irreducible, we see that $H^i(X,\mathcal{L}_\chi)$ is isomorphic as a ${}^L G$-representation to a twist of $V_\chi$ by a character of $\Gamma$. Denoting by $\epsilon_\chi$ the value of this character on $F$ we obtain the second formula of the proposition.
\end{preuve}

\begin{corollaire}\label{corollary appendix}
	Let $\widehat{Q}\subset \widehat{G}$ be a $\Gamma$-stable standard parabolic subgroup. Let $\widehat{L}\subset \widehat{Q}$ be the unique Levi component containing $\widehat{T}$ and set $W^L=\Norm_{{}^L L}(\widehat{T})/\widehat{T}\subset W$ where we have set ${}^L L=\widehat{L}\rtimes \Gamma$. Then, we have
	$$\displaystyle \sum_{w\in W^L} D_{\widehat{G}/\widehat{B}}(wS)^{-1}=D_{\widehat{G}/\widehat{Q}}(S)^{-1}$$
	for every regular $S\in {}^L T$.
\end{corollaire}

\begin{preuve}
	Note that
	$$\displaystyle D_{\widehat{G}/\widehat{B}}(S)=D_{\widehat{G}/\widehat{Q}}(S)D_{\widehat{L}/\widehat{B}_L}(S)$$
	for $S\in {}^L T$, where we have set $\widehat{B}_L=\widehat{B}\cap \widehat{L}$. The corollary now follows readily from the previous proposition applied to the trivial character and ${}^L L$ instead of ${}^L G$.
\end{preuve}
\end{paragr}